\documentclass[11pt]{report}

\usepackage[a4paper,left=3cm,right=3cm,top=2cm,bottom=2cm]{geometry}
\usepackage{setspace}
\usepackage{url}

\usepackage{pgf}
\usepackage{mathtools}
\usepackage{graphicx}
\usepackage{amsthm}
\usepackage{amsmath}
\usepackage{amsfonts}
\usepackage{MnSymbol}
\usepackage[utf8]{inputenc}
\usepackage[english]{babel}
\usepackage{cite}
\usepackage[toc,page]{appendix}
\usepackage{youngtab}
\usepackage{ytableau}
\usepackage{comment}
\usepackage{floatrow}
\usepackage{relsize}
\usepackage{float}
\usepackage{dsfont}

\usepackage{adjustbox}

\usepackage{enumitem}

\usepackage{scrextend}
\addtokomafont{labelinglabel}{\sffamily}

\usepackage[hidelinks]{hyperref}

\usepackage{tikz}
\usetikzlibrary{arrows}

\usepackage{cellspace}
\newtheorem{thm}{Theorem}
\numberwithin{thm}{section}

\newtheorem{corollary}[thm]{Corollary}
\newtheorem{lemma}[thm]{Lemma}

\newtheorem{conj}[thm]{Conjecture}

\theoremstyle{definition}
\newtheorem{defn}[thm]{Definition}
\newtheorem{example}[thm]{Example}
\newtheorem{rem}[thm]{Remark}

\newcommand{\E}{\mathrm{E}}
\newcommand{\Var}{\mathrm{Var}}
\newcommand{\sh}{\Phi}
\newcommand{\rtrsh}{\textnormal{sh}}

\newcommand{\pro}{\kappa}

\newcommand{\iso}{\pi}

\newcommand{\best}{T_\lambda^{\rightarrow}}
\newcommand{\worst}{T_\lambda^{\downarrow}}
\newcommand{\SEQ}{T_\lambda^{\rightarrow}}
\newcommand{\OSEQ}{T_\lambda^{\downarrow}}
\newcommand{\YT}{\textnormal{YT}}
\newcommand{\SYT}{\textnormal{SYT}}
\newcommand{\eig}{\textnormal{eig}}
\newcommand{\eval}{\textnormal{eval}}

\newcommand{\righthand}{R^{i}}
\newcommand{\lefthand}{L^{i}}
\newcommand{\NA}{N_{\alpha}(n)}
\newcommand{\NW}{N_{w}(n)}
\newcommand{\D}{\textnormal{Diag}}

\newcommand{\ind}{\mathds{1}}

\newcommand{\dt}{d_{\textnormal{\tiny TV}}}
\newcommand{\ds}{d_{\textnormal{\tiny Sep}}}
\newcommand{\mt}{t_{\textnormal{\tiny TV}}}
\newcommand{\ms}{t_{\textnormal{\tiny Sep}}}

\newcommand{\Tr}{\textnormal{Tr}}
\newcommand{\support}{\textnormal{Supp}}
\newcommand{\I}{\textnormal{Id}}

\newcommand{\RT}{\textnormal{RT}_{n}}
\newcommand{\ART}{\textnormal{ART}}
\newcommand{\ttr}{\textnormal{TTR}_{n}}
\newcommand{\rtt}{\textnormal{RTT}_{n}}
\newcommand{\rtr}{\textnormal{RTR}_{n}}

\newcommand{\artr}{\textnormal{ARTR}_{n}}
\newcommand{\artrp}{\textnormal{ARTR}_{n+1}}
\newcommand{\artrn}{\textnormal{ARTR}}

\newcommand{\RL}{\textnormal{OST}}
\newcommand{\LR}{\textnormal{OST}}
\newcommand{\ALR}{\textnormal{AOST}}
\newcommand{\BLR}{\textnormal{OST}}
\newcommand{\ABLR}{\textnormal{AOST}}

\newcommand{\brt}{\overline{\textnormal{RT}}_{n}}
\newcommand{\abrt}{\overline{\textnormal{ART}}_{n}}
\newcommand{\abrtp}{\overline{\textnormal{ART}}_{n+1}}
\newcommand{\abrtn}{\overline{\textnormal{ART}}}

\newcommand{\bblr}{\overline{\textnormal{OST}}}
\newcommand{\abblr}{\overline{\textnormal{AOST}}_{n,w}}
\newcommand{\abblrp}{\overline{\textnormal{AOST}}_{n+1,w}}

\newcommand{\ad}{\Phi}

\newcommand{\rep}{\textnormal{Irr}}

\newcommand{\Triv}{\textnormal{Triv}}

\newcommand{\gpid}{e}

\newcommand{\sign}{\textnormal{sgn}}

\newcommand{\Ker}{\textnormal{Ker}}
\newcommand{\im}{\textnormal{Im}}

\newcommand{\bl}{\overline{\lambda}}
\newcommand{\po}{\lambda^{1}}
\newcommand{\pt}{\lambda^{2}}
\newcommand{\bm}{\overline{\mu}}
\newcommand{\bt}{\overline{T}}
\newcommand{\tl}{T_{\overline{\lambda}}}
\newcommand{\tm}{T_{\overline{\mu}}}

\newcommand{\res}{\textnormal{Res}}
\newcommand{\induce}{\textnormal{Ind}}

\begin{document}

\pagenumbering{roman}	
	
\begin{titlepage}
	\centering
	{\huge\bfseries Random Walks on the Symmetric Group: \\
		\LARGE Cutoff for One-sided Transposition Shuffles\par}
	\vspace{5cm}
	{\Large\itshape Oliver Thomas Matheau-Raven\par}
	\vspace{5cm}
	{\scshape\Large PhD\par}
	\vspace{2cm}
	{\scshape\LARGE University of York \par}
	\vspace{0.1cm}
	{\scshape\LARGE Mathematics \par}
	\vspace{5cm}
	{\scshape\LARGE August 2020. \par}	
	\vfill
\end{titlepage}

\newpage

\addcontentsline{toc}{section}{Abstract}

\begin{abstract}
	In this thesis we introduce a new type of card shuffle called the \emph{one-sided 
		transposition shuffle}. At each step a card is chosen uniformly from the pack and then transposed with another card chosen uniformly from \emph{below} it. This defines a random walk on the symmetric group generated by a distribution which is non-constant on the conjugacy class of transpositions. Nevertheless, we provide an explicit formula for all eigenvalues of the shuffle by demonstrating a useful correspondence between eigenvalues and standard Young tableaux. This allows us to prove the existence of a total-variation cutoff for the one-sided transposition shuffle at time $n\log n$. We also study weighted 
	generalisations of the one-sided transposition shuffle called \emph{biased one-sided transposition shuffles}. We compute the full spectrum for every biased one-sided transposition shuffle, and  prove the existence of a total variation cutoff for certain choices of weighted distribution. In particular, we recover the eigenvalues and well known mixing time of the classical random transposition shuffle. We study the hyperoctahedral group as an extension of the symmetric group, and formulate the one-sided transposition shuffle and random transposition shuffle as random walks on this new group. We determine the spectrum of each hyperoctahedral shuffle by developing a correspondence between their eigenvalues and standard Young bi-tableaux. We prove that the one-sided transposition shuffle on the hyperoctahedral group exhibits a cutoff at $n\log n$, the same time as its symmetric group counterpart. We conjecture that this results extends to the biased one-sided transposition shuffles and the random transposition shuffle on the hyperoctahedral group.  
\end{abstract}

\newpage
\tableofcontents

\newpage
\addcontentsline{toc}{section}{List of Tables}
\listoftables

\newpage
\addcontentsline{toc}{section}{List of Figures}
\listoffigures

\newpage
\section*{Introduction}
\addcontentsline{toc}{section}{Introduction}

\paragraph{}
Consider a stacked deck of $n$ distinct cards, whose positions are labelled by elements of the set $[n]:=\{1,\dots,n\}$ from bottom to top. Any shuffle 
which involves choosing two positions and switching the cards found there (if the two positions coincide then no cards are 
moved)
is 	called a \emph{transposition shuffle}, and may be viewed as a random walk on the 
symmetric group $S_n$. We refer to the two positions as being picked by our right and left hands, associating with them random variables, $\righthand$ and $\lefthand$ respectively, which take values in $[n]$. Diaconis and Shahshahani \cite{diaconis1981generating} were the 
first to study random walks on groups  using representation theory; they famously showed that the \emph{random transposition shuffle}, 
in which the two positions are chosen independently and uniformly on $[n]$, 
takes $(n/2)\log n$ steps to randomise the order of the deck. The time taken to randomise the order of the deck is known as the \emph{mixing time} of the shuffle. 

We may expect random walks to converge to their stationary distribution smoothly as time increases but there are many classes of random walks for which convergence happens only once a certain threshold time has been met. This behaviour was first observed in the random transposition shuffle by Diaconis and Shahshahani and has been given the name \emph{the cutoff phenomenon}. 
Suppose the mixing time of a succession of random walks on groups $\{G_{n}\}_{n \in \mathbb{N}}$ may be tightly bounded above and below by a time of the form $t_{n} \pm c \,w_{n}$, where $c$ is of constant order and $w_{n} = o(t_{n})$, then the sequence of random walks is said to exhibit a \emph{cutoff} at time $t_{n}$.
The mixing time of the random transposition shuffle may be shown to be upper and lower bounded around at time $(n/2)\log n \pm c\, n$, thus it shows a cutoff at time $(n/2)\log n$. 
The cutoff phenomenon is prevalent throughout the study of random walks on groups with many examples shown to exhibit a cutoff. However,  there is currently no known sufficient criterion for the existence of a cutoff, and proving a cutoff for any particular random walk is generally difficult.

\paragraph{}
Since the analysis of the random transposition shuffle a variety of algebraic and probabilistic techniques have been 
employed to study the mixing time of different transposition shuffles. Notable 
examples include the \emph{top-to-random transposition} shuffle 
\cite{saloff2004random}, the \emph{adjacent transposition} shuffle 
\cite{lacoin2016mixing}, and the \emph{random transposition  shuffle on a cycle} where our positions are constrained to be within some fixed distance \cite{berestycki2012}. All of these shuffles have the property 
that, at each step, the transposition to be applied is chosen uniformly 
from a subset of transpositions which generate the entire group $S_n$.

An interesting class of transposition shuffles which are not chosen uniformly from a set are the \emph{semi-random transposition shuffles}  (see e.g. \cite{Mossel2004, Pymar2011, pinsky2012cyclic}).
Any semi-random transposition shuffle is driven by the following procedure: pick one position with our right hand uniformly at random, and then independently choose a second position with our left hand following some stochastic (possibly time-inhomogeneous) process on the set $[n]$, then transpose the cards in the chosen positions. This class includes the random transposition shuffle by letting the stochastic process for the left hand be repeated uniform draws from $[n]$.  Mossel, Peres, and 
Sinclair \cite{Mossel2004} were able to establish a universal upper bound of $O(n \log n)$ on the mixing time of any semi-random transposition shuffle. The key idea behind this proof relies on the positions chosen being independent. The card picked by our stochastic process is uniformly transposed with a card in the deck, therefore, once our stochastic process has picked every possible card
our deck will be at a uniformly chosen permutation of $S_{n}$. This reduces the analysis of the mixing time to a coupon collector's problem, and thus we find an upper bound of order $n\log n$.

\paragraph{}
In this thesis we introduce a new class of shuffles called 
the \emph{one-sided transposition shuffles}: these have the defining property 
that at step $i$ the right hand's position ($\righthand$) is chosen 
according to some arbitrary distribution on $[n]$, and given the value of 
$\righthand$ the distribution of the left hand's position ($\lefthand$) is 
supported on the set $\{1,\dots,\righthand\}$.  When our right hand is chosen uniformly from $[n]$ the one-sided transposition shuffle follows a similar description to that of a semi-random transposition shuffle however with the key difference that our choice for $L^{i}$ now depends on our first uniform choice. This change of dependence greatly influences the behaviour of the one-sided transposition shuffles. There does not exist a universal upper bound on the mixing time of the one-sided transposition shuffles without imposing further constraints on the distribution of the left hand, since a arbitrarily slow shuffle can be created by increasing the probability that the two hands choose the same position. For the majority of Chapter 5 we focus on the case when the right and left hands are both chosen uniformly from their possible ranges: we call this the \emph{(unbiased) one-sided transposition shuffle}. Although the support for the one-sided transposition shuffle is the entire conjugacy class of transpositions, our probability distribution on 
this set is in general far from uniform. (E.g. when right and left hands are both uniformly chosen from their permitted ranges, the probabilities attached to  different transpositions range from $1/n^2$ to $1/2n$)

The main results of this thesis are as follows: we recover the eigenvalues for the one-sided transposition shuffle using the technique of lifting eigenvectors which was recently constructed by Dieker and Saliola \cite{dieker2018spectral}. Analysis of the spectrum for the one-sided transposition shuffle allows us to find tight bounds on the mixing time of the shuffle and prove that it exhibits a cutoff at time $n\log n$. Afterwards we modify our analysis to study one-sided transposition shuffles where the right hand  $\righthand$ is chosen via a weighted distribution: we call these the \emph{biased one-sided transposition shuffles}.  We extend the technique of lifting eigenvectors to compute the full spectrum of all biased one-sided transposition shuffles. Furthermore, for particular choices of weighted distribution we are able to  prove the existence of a cutoff for the biased one-sided transposition shuffles. We go on to study random walks on the hyperoctahedral group $B_{n}$ as an extension of the symmetric group where each card now has a distinguishable back and front. We create transposition shuffles on $B_{n}$ from existing shuffles on $S_{n}$ by allowing ourselves a chance of turning cards over. We show that the method of lifting eigenvectors may be modified to the setting of the hyperoctohedral group, and compute the spectrum of the random transposition shuffle and one-sided transposition shuffle on $B_{n}$. We end by showing that the unbiased one-sided transposition shuffle on the hyperoctahedral group exhibits a cutoff at time $n\log n$, the same time as its symmetric group counterpart.

\paragraph{}
The structure of this thesis is as follows:
Chapter \ref{chpt1:chpt} gives an introduction to the topic of random walks on groups, including all the necessary probabilistic and algebraic background. We restrict our attention to random walks which converge to a unique stationary distribution and introduce total variation distance to  measure the convergence of random walks to their stationary distributions. We define the mixing time of a random walk, and what it means for a random walk to exhibit a cutoff in total variation distance. We recall the definitions of representations, modules and characters for a group $G$, and state classical results including Schur's Lemma and Maschke's Theorem. We also explore the properties of discrete Fourier transforms on a group $G$. Finally we recall the \emph{upper bound lemma} which allows us to analyse the mixing time of a random walk on a group $G$ using the  irreducible representations of $G$. 

In Chapter \ref{chpt4:chpt} we specialise the results of Chapter \ref{chpt1:chpt} to the symmetric group $S_{n}$. We build up a picture of the structure of the symmetric group including a construction of its simple modules. To demonstrate a variety of techniques used to analyse random walks on groups we present three longer examples: the random transposition shuffle\cite{diaconis1981generating}, the top-to-random shuffle \cite{aldous1986shuffling}, and the random-to-random shuffle\cite{dieker2018spectral}.  	
The most recent of these is the random-to-random shuffle where Dieker and Saliola first used the technique of lifting eigenvectors. We explore their work showing that the eigenvalues of the random-to-random shuffle may be computed by transforming eigenvalues of the random-to-random shuffle on $n$ cards to those of the shuffle on $n+1$ cards.

The main body of Chapter \ref{chpt5:chpt} is comprised of a paper jointly written with Michael E. Bate and Stephen B. Connor \cite{bate2019cutoff}. In this chapter we explore the one-sided transposition shuffle which stems from a modification of the classical random transposition shuffle. Our main result is proving that the one-sided transposition shuffle exhibits a cutoff in total variation distance at time $ n \log n$. To prove this result we analyse the spectrum of the shuffle which we compute by inductively lifting the eigenvectors for the shuffle. We present the details of our lifting in full as it differs in interesting ways from that of the random-to-random shuffle, in particular we are able to recover all eigenvectors and eigenvalues for the one-sided transposition shuffle. We generalise our lifting to recover eigenvalues for biased one-sided transposition shuffles and for particular choices of weighted distribution show them to exhibit a cutoff in total variation distance. Finally we show the one-sided transposition shuffle to exhibit a cutoff in separation distance at time $n\log n$ via the use of a strong stationary time.

Chapter \ref{chpt6:chpt} establishes the hyperoctahedral group $B_{n}$ as an extension of the symmetric group. We recall facts about the hyperoctahedral group, relating them to the definitions and results of the symmetric group in Chapter \ref{chpt4:chpt}. We examine the module structure of $B_{n}$ and construct its permutation and simple modules which are now indexed by bi-partitions of $n$. We extend the random transposition shuffle and one-sided transposition shuffle to random walks on the group $B_{n}$ by adding in a chance to flip cards over during our shuffling procedure. We generalise the technique of lifting eigenvectors to the hyperoctahedral group allowing us to recover the entire spectrum of the random transposition shuffle and one-sided transposition shuffle. We conjecture that the two shuffles on the hyperoctohedral group exhibit a cutoff in total variation distance at the same time as their symmetric group counterparts.

\newpage
\section*{Acknowledgements}
\addcontentsline{toc}{section}{Acknowledgements}

First and Foremost, I would like to thank Michael Bate and Stephen Connor for introducing me to this wonderful area of mathematics. Thanks to your attentive supervision and patience I have grown as a mathematician and as a communicator in general. This thesis would not have been possible without your help.

Thank you to my examiners, Harry Geranios and Evita Nestoridi, for making my viva an enjoyable experience and for your comments on this thesis. Thank you to EPRSC and the department of mathematics for their the financial support which made my PhD possible.
I also want to thank all the staff at York for their support, particularly Vicky Gould and Alet Roux for being on my TAP.  

I want to give a big thank you to the rats of Millfield Road: Allan, Chris, Christoph, and Scott, for being the best housemates, board game players and quarantine companions I could ask for. Thank you to the whole mathematics PhD community at York, including those that have left and those that have just joined, for all the Wednesday talks, pub trips and tea breaks, which provided needed respite from mathematics. Thank you to the archery society for fostering a hopefully life long hobby in me. I also extend my gratitude to everyone who has been a part of my PhD journey over the last three years.

Finally, I want to thank my family for all the love and support they have always given me.

\newpage
\section*{Declaration}
\addcontentsline{toc}{section}{Declarations}

I declare that this thesis is a presentation of original work and I am the sole author. This work has been carried out under the supervisor of Dr. Michael E. Bate, and Dr. Stephen B. Connor, and has not previously been presented for an award at this, or any other, University. All sources are acknowledged as References.
Chapters 1 and 2, form a review of background material needed for this thesis, and the appropriate literature is referenced before each chapter. Chapter \ref{chpt5:chpt} is based on a paper written in collaboration with  Michael E. Bate and  Stephen B. Connor, which has been accepted for publication in the Annals of Applied Probability \cite{bate2019cutoff}. Chapter \ref{chpt6:chpt} consists of original work by the author.

\newpage \vspace*{8cm}
\thispagestyle{empty}

\emph{\Large This thesis is dedicated to the memory of little gran.}

\chapter{Preliminaries}
\label{chpt1:chpt}

\pagenumbering{arabic}

\section{Random Walks on Finite Groups}

There are many different ways we may formulate a random walk on a group $G$. We could designate one step transition probabilities for every two states $g\to h$, or even have the probabilities depend on our past like a self-avoiding random walk. However, in the first case we have not used any of the symmetry a group has, and in the second case our walk is not Markovian. 
We want to exploit the structure of groups to make our random walk inherit sensible properties and we also want the walk to form a Markov chain. In this thesis we focus on one  particular Markovian description of a random walk where the walk is driven by a single probability distribution $P:G\to [0,1]$. All groups we consider are finite and have identity element $\gpid$, unless otherwise stated. A detailed account of the theory of Markov chains presented in this section can be found in \emph{Markov Chains and Mixing Times} by Levin, Peres, and Wilmer \cite{Levin2017}. We begin this section with the definition of a random walk on a group before moving on to recalling important definition and facts about Markov chains. The section ends with the introduction of the cutoff phenomenon for random walks on groups. 

\subsection{Markov Chains from Random Walks}

\begin{defn}
	Let $G$ be a finite group, suppose $P:G\to [0,1]$, and $\mu:G \to [0,1] $, are probability distributions on $G$. Define each individual step of our random walk as sequence of i.i.d. (independent and identically distributed) random 
	variables $\{\xi^{t}\}_{t\geq 1}$ on $G$ distributed according to $P$. Create a sequence of random variables 
	$\{X^{t}\}_{t\geq 	0}$, by setting $X^{0} \sim_{d} \mu$ and inductively defining
	\[X^{t+1} = \xi^{t+1} X^{t} \]
	for all $t>0$. Then $\{X^{t}\}_{t\geq 0}$ defines a Markov chain with one step transition probabilities $P(g,h) :=P(hg^{-1})$. We call this the \emph{random walk} on $G$ \emph{driven} by $P$.
\end{defn}

Given a random walk on $G$, its driving probability $P$ gives rise to the transition matrix of our random walk defined using one step transition probabilities $P(g,h)$. We distinguish between these two interpretations of the symbol $P$ by the number of arguments each one takes.  Note that from the transition matrix we may recover the driving probability by setting $P(g) = P(\gpid,g)$.

\begin{defn}
	\label{chpt2:def:conv}
	Let $P, Q$ be probabilities on a finite group $G$, define the \emph{convolution} 
	of $P$ with $Q$, denoted 
	$P \star Q: G \to [0,1]$ as
	\[(P \star Q)(g) :=\sum_{h\in G}P(gh^{-1})Q(h).\]
	This is the probability of a random walk ending at element $g$ if we start at $\gpid$ and take the first step according to $Q$ and the second step to $P$.
	We denote $P \star P = P^{\star 2}$, and the $t^{\textnormal{th}}$ convolution of $P$ with itself as $P^{\star t}$.
\end{defn}
\begin{lemma}
	Let $P^{t}$ be the $t^{\textnormal{th}}$ power of our transition matrix 
	$P$. Then $P^{t} (g,h)= P^{\star t}(hg^{-1})$. 
	
	Hence, we drop the 	$\star$ in all subsequent work and just write $P^{t}$.
\end{lemma}
\begin{proof}
	We proceed by induction, by definition the statement holds for $t=1$. Now 
	\begin{eqnarray*}
		P^{t+1}(g,h) =  \sum_{a\in G} 
		P(g,a)
		P^{t}(a,h) &= &\sum_{a\in G} P(ag^{-1}) P^{\star t}(ha^{-1})\\
		& = & \sum_{b\in G} P(b) P^{\star t}(h g^{-1} b^{-1}) = 
		P^{\star (t+1)}(hg^{-1}). \qedhere
	\end{eqnarray*}
\end{proof}

Lemma 3 tells us that convolution and matrix multiplication amount to the same thing. We use $P^{t}$ to stand for both the $t$-step 
transition matrix and the $t$-fold convolution of our probability $P$, with the 
assumption that if a starting state is not specified we assume it to be the 
identity.

Our starting state $X^{0}$ is picked via a 
probability distribution $\mu$, i.e. $\mathbb{P}(X^{0}=g) = \mu(g)$. Given the 
starting distribution $\mu$ the probability of being in state $g$ at time $t$ 
is $(P^{t} \mu)(g)$. In practice we often fix the starting 
distribution to be a single element of our state space, usually the identity of our finite group. We note, however, that the choice of 
starting distribution can have a large impact on the behaviour of a random walk, as the following example shows.

\begin{example}
	The simple random walk on $\mathbb{Z}_{n}$, is the random walk generated by probability $P:\mathbb{Z}_{n} \to [0,1]$ with $P(1) = P(-1) =1/2$.
	
	Consider $P$ to be the simple random walk on $\mathbb{Z}_{6}$.
	Suppose we start at state $0$ then after $3$ steps we have $P^{3}(0,0) = 
	0$ and $P^{3}(0,1) = 6/16$. This is because we always take one step at 
	each time meaning at odd times we are on odd elements on $\mathbb{Z}_{6}$ 
	and at even times we are on even elements.
	Instead now suppose we start at element $1$, then $P^{3}(1,0) = 6/16$, 
	and $P^{3}(1,0) = 0$ for the same reason as above.
	Combining the two starting distributions let $\mu$ now have $\mu(0) = 
	\mu(1) 
	=1/2$ then $(P^{3}\mu)(0)  = \frac{1}{2}(P^{3}(0,0) + 
	P^{3}(1,0)) = 3/16 =  (P^{3} \mu)(1).$
\end{example} 

\paragraph{}
Understanding how the distribution $P^{t}$ evolves in time is key to the study of Markov chains. In theory we could 
always compute $P^{t}$ given the probability $P$ but this becomes 
impractical as we consider larger state spaces and  behaviour at times with $t$ large. Instead, we often look to bound these probabilities to understand their behaviour without exact computation. 
Development of tools to help understand how $P^{t}\mu$ behaves are an active area of research. Often we can use 
characteristics of our random walk: the starting distribution $\mu$, the 
transition matrix $P$, and the group $G$, in order to form bounds on the 
probability distribution $P^{t} \mu$. In particular the eigenvalues 
of the transition matrix $P$ play a big part in its behaviour. Furthermore, random 
walks on groups allow the use of algebraic tools 
using the representations of $G$, which we explore in section \ref{chpt3:chpt}. Next we recall the definition of stationary 
distributions for a Markov chain.

\begin{defn}
	Let $\{X^{t}\}_{t\geq0}$ be a Markov chain with 
	transition matrix $P$. A \emph{stationary distribution} for the Markov chain is a probability distribution $\pi$ such that $ P\pi  = \pi$. Note that for any stationary distribution $\pi$ we have $P^{t}\pi  = \pi$ for all $t \geq 1$.
\end{defn}

Generic Markov chains do not necessarily have stationary distributions nor do they have to be unique.
\begin{example}
	Let $G = \mathbb{Z}$, consider a Markov chain started 
	at $0$ and driven by probability $P(1) =1 $. This random walk has no stationary distribution because for 
	any distribution $\mu$ on $\mathbb{Z}$ we have $( P \mu)(i) = \mu(i+1)$.
	
	Let $G= \mathbb{Z}_{2}^{2}$, consider a random walk with probability 
	$P((0,0)) = P((0,1)) = 1/2$. Suppose we have $\mu_{0},\mu_{1}$ defined by 
	$\mu_{0}((i,j)) = 1/2$ if $i=0$, and $\mu_{1}((i,j))= 
	1/2$ if $i=1$. Then both $\mu_{0},\mu_{1}$ are stationary distributions for 
	$P$. In fact let $a \in [0,1]$ and $\mu_{a}$ be defined by
	\[\mu_{a}((i,j)) = \begin{cases}
	\frac{1-a}{2} & \textnormal{ if } i=0\\
	\frac{a}{2} & \textnormal{ if } i=1 
	\end{cases}.\]
	Then $\mu_{a}$  is a stationary distribution for $P$, so we have infinitely many stationary distributions for this random walk.

\end{example}

\begin{lemma}
	\label{chpt2:lem:uniform}
	Let $\pi$ denote the uniform distribution on a finite group $G$, that is for all $g \in G$,  
	$\pi(g) = 1/|G|$. Then $\pi$ is a stationary distribution for any random walk on $G$.
\end{lemma}

\begin{proof}
	To prove this all we need to show is that $P \pi = P \star \pi = \pi$ for any probability $P$ defined on $G$. 
	Following from the definition
	\[(P \star \pi)(g) = \sum_{h \in G} P(gh^{-1}) \pi(h) 
	=\frac{1}{|G|}\sum_{h\in G} P(gh^{-1}) = \pi(g) .\]
\end{proof}

Under mild assumptions we can prove that a Markov chain has a unique stationary 
distribution which we denote by $\pi$. Below we list some of the 
properties a Markov chain may exhibit. 

\begin{defn} 
	\label{chpt2:prop:properties}
	Let $\{X^{t}\}$ be a Markov chain on a space $\mathcal{X}$ with transition 
	matrix $P$. The Markov chain may have the following properties.

	\setlist[description]{font=\normalfont \underline}
	\begin{description}
		\item[Irreducible:] $\{X^{t}\}$ is called 
		\emph{irreducible} if for all $x,y$, there exists $t>0$ such that 
		$P^{t}(x,y) >0$.
		\item[Aperiodic:] Let $r(x) = \{t\geq 1 : P^{t}_{x,x}>0\}$ these are 
		the \emph{return times} of the state $x$, the \emph{period} of $x$ 
		is $\gcd( r(x))$. 
		The Markov chain is called \emph{aperiodic} if all states have period 
		$1$.
		\item[Transitive:] A Markov chain is called \emph{transitive} if for all pairs $(x_{1},x_{2})\in \mathcal{X} \times \mathcal{X}$ 
		there exists a bijection $\phi$, such that $\phi(x_{1})=x_{2}$ and 
		$\phi$ preserves all one step transition probabilities, i.e. for all 
		$y,z \in \mathcal{X}$ we have $P(y,z) = P(\phi(y),\phi(z))$. 
		\item[Reversible:]A Markov chain is called \emph{reversible} if there exists a probability 
		distribution 
		$\pi$ on $\mathcal{X}$ such that
		\begin{eqnarray}
		\label{chpt2:eqn:detail}
		\pi(x)P(x,y)  =  \pi(y)P(y,x) 
		\end{eqnarray}
		for all $x,y \in \mathcal{X}$. In this case $\pi$ is 
		a stationary distribution. The set of equations defined by 	\eqref{chpt2:eqn:detail} are called 
		the \emph{detailed balance equations}.
	\end{description}
\end{defn}

All random walks we study in detail in this thesis will be irreducible, 
aperiodic, and transitive. We shall see one by one why these conditions are 
necessary for the study of mixing times of random walks on groups.
First we prove that all random walks on finite groups are transitive.

\begin{lemma}[Section 2.6.2 \cite{Levin2017}]
	\label{chpt2:lem:groupstrans}
	Random walks on finite groups define transitive Markov chains.
\end{lemma}

\begin{proof}
	Let $\{X^{t}\}_{t\geq 0}$ be a random walk on a finite group with driving 
	probability $P$. Take a pair $(x_{1},x_{2}) \in G\times G$, and 
	define 	a bijection $\phi_{x_{1},x_{2}}:G \to G$ by $\phi_{x_{1},x_{2}}(g) =gx_{1}^{-1}x_{2}$. Then 
	$\phi_{x_{1},x_{2}}(x_{1}) 
	= x_{2}$, and for any pair $g,h \in G$ we get
	\[P(g,h) = P(hg^{-1}) =P((hx_{1}^{-1}x_{2})(x_{2}^{-1}x_{1}g^{-1})) 
	=P(\phi_{x_{1},x_{2}}(g),\phi_{x_{1},x_{2}}(h))	\qedhere .\]
\end{proof}

Irreducibility implies that the Markov chain has a unique stationary distribution \cite[Corollary 1.17]{Levin2017}. Therefore, by Lemma \ref{chpt2:lem:uniform} if a 
random walk on a finite group $G$ is irreducible its unique stationary distribution must in fact be the uniform distribution.

\begin{lemma}
	Let $P$ be the transition matrix of an aperiodic irreducible Markov chain, there 
	exists a unique stationary distribution $\pi$ for $P$. 
\end{lemma}

\begin{corollary}
	Let $P$ be a probability on a finite group $G$. Define the \emph{support} of $P$ as the set $\support(P) =\{g : 
	P(x)>0\}$. Suppose $\support(P)$ 
	is a generating set for $G$, then the 
	random walk on $G$ driven by $P$ is irreducible, and hence has the uniform distribution as its unique stationary distribution.
\end{corollary}

Suppose a random walk $P$ on a group $G$ is not 
irreducible, then if we 
restrict our state space to the subgroup $ \langle \support(P) \rangle =H 
\subsetneq G$ we recover an irreducible random walk on group $H$.
When choosing random walks on groups we like to fix the support of a 
probability to be a set of generators $S \subseteq G$ and vary our probability to see what different behaviour can arise.  One important case we focus on is transposition shuffles of the symmetric group.
We shall see that the one-sided transposition shuffle has very different behaviour to the random transposition shuffle despite being generated by the same conjugacy class.
If our support $S$ is a union of conjugacy classes of our group $G$ then we gain access to extra algebraic tools to help our analysis, we explore these in section \ref{chpt3:chpt}. From here onwards we assume all the random walks on 
groups are irreducible with stationary distribution $\pi: G\to G$ being uniform.

Once we know a unique stationary distribution exists our next question is 
whether our probability $P^{t}$ will ever reach this equilibrium? Convergence 
of $P^{t}$ to a stationary distribution generally depends on the starting 
distribution of our Markov chain.

\begin{example}
	\label{chpt2:ex:period}
	Consider the simple random walk on $\mathbb{Z}_{n}$ started at $0$, this 
	Markov chain is irreducible therefore has stationary distribution 
	$\pi_{n}(i) = 1/n$. Suppose $n>2$ is even, then the walk has a period of 
	$2$: at all even times 	$X^{2t}$ we must be at an even integer, and at all 
	odd times $X^{2t+1}$ we are at an odd integer. Therefore, our probability $P^{t}$ will always be distinguishable from $\pi$. However, if we start at 
	the distribution $\mu(0) = \mu(1) = 0.5$, then we have $ (P^{t} \mu )(i) 
	\to 1/n$ for all $i$ as $t\to \infty$.
	
	Consider instead $n>2$ being odd, the simple random walk on $\mathbb{Z}_{n}$ has period $1$. Now for any single starting point $i$ we have 
	$P^{t}(i,j) \to 1/n$ for all $j$ as $t \to \infty$.
\end{example}

The key difference between the two random walks in Example 
\ref{chpt2:ex:period} is their period.  In Section 1.2.2. we establish that aperiodic random walks  always converge to their stationary distributions.
To get rid of periodicity concerns for a random walk, we introduce the notion of a \emph{lazy} random walk.

\begin{lemma}
	\label{chpt2:lem:lazy}
	Let $\{X^{t}\}$ be an periodic Markov chain driven by $P$. Create a new 
	Markov chain $\{Y^{t}\}$ with $Y^{0} =X^{0}$, and new transition 
	probabilities $Q$ created in the following way: flip a fair coin, if heads 
	do nothing, otherwise proceed according to $P$. Therefore, $Q(\gpid) = 
	\frac{1}{2} + \frac{1}{2}P(\gpid)$ and $Q(g) = \frac{1}{2}P(g)$ for 
	$g\neq e$. The Markov chain $\{Y^{t}\}$ is called the \emph{lazy} version 
	of $\{X^{t}\}$, and is aperiodic.
\end{lemma}

\begin{proof}
	The lazy walk has $Q(\gpid) >1/2 $, and so for any $g\in G$ we find $r(g) = \{t\geq 1 : Q^{t}(g,g)>0\} 
	= \mathbb{N}$. Hence, all states have period $1$.
\end{proof}

\begin{example}
	The lazy simple random walk on $\mathbb{Z}_{n}$ has probability $P(0) =1/2$ 
	and $P(-1)= P(1) =1/4$. This walk is now aperiodic for all $n\geq  2$.
\end{example}

The proof of Lemma~\ref{chpt2:lem:lazy} exploits the fact that so long as we 
have $P(\gpid)>0$ our random walk is aperiodic. All random walks we study we detail will be aperiodic because they will have a non-zero probability of remaining still. In the next 
section we a define pair of measures on the space of probability distributions, and define a notion of convergence with respect to these measures.

\subsection{Convergence to a Stationary Distribution}

Consider a Markov chain $\{X^{t}\}$ with transition matrix $P$ and stationary 
distribution $\pi$. We 
know that if our chain starts at $X^{0} \sim_{d} \pi$ then we 
remain at the distribution $\pi$ after every step. We show that for 
any starting distribution an irreducible, aperiodic Markov chain always converges to its unique stationary distributions. To prove this we 
first establish a measure on our space of probability distributions, and then use this to define convergence. To this end 
we introduce two notions of distance on the space of probability 
distributions over a finite group $G$.

\subsubsection{Total Variation Distance}

The first metric we introduce is ubiquitous within the study of Markov chains, it 
is called total variation distance.

\begin{defn}
	\label{chpt2:def:totalvariation}
	Let $\mu,\nu$ be two probabilities on a finite group $G$. Define the \emph{total 
		variation distance} between $\mu, \nu$ as follows:
	\begin{eqnarray}
	\lVert \mu - \nu \rVert_{\textnormal{\tiny TV}} = \sup_{A \subseteq G} | \mu(A) -\nu(A) |. \label{chpt2:eqn:TV1}
	\end{eqnarray}
\end{defn}

The total variation distance between any two probabilities always lies in the 
range $[0,1]$. The formulation of total variation distance provided in 
Definition \ref{chpt2:def:totalvariation} is frequently too cumbersome for use. 
We may reformulate definition \eqref{chpt2:eqn:TV1} to  involve summing over the group $G$ instead of its subsets.

\begin{lemma}[Section 4.1 \cite{Levin2017}]
	\label{chpt2:lem:totalvariation}
	Let $\mu,\nu$ be two probabilities on $G$, then we have
	\[\lVert \mu - \nu \rVert_{\textnormal{\tiny TV}} = \frac{1}{2}\sum_{g\in 
		G} | \mu(g) - \nu(g)| \qedhere .\]
\end{lemma}

\begin{proof}
	Let $A= \{g : \mu(g) \geq \nu(g)\}$, and so $A^{c} = \{g : \nu(g) > 
	\mu(g)\}$. First notice that the supremum in Definition 
	\ref{chpt2:def:totalvariation} is reached by set $A$, and secondly that  
	$\mu(A) - \nu(A) = \nu(A^{c}) - \mu(A^{c})$. Putting these facts together 
	we find,
	\[\lVert \mu - \nu \rVert_{\textnormal{\tiny TV}} = \mu(A)- \nu(A) = 
	\frac{1}{2}\left(\mu(A)- \nu(A) + \nu(A^{c}) - \mu (A^{c})\right) = 
	\frac{1}{2}\sum_{g\in G} | \mu(g) - \nu(g) |. \qedhere \]
\end{proof}

Following from our definitions it is now easy to see that total variation distance forms a metric on the space of probability distributions of $G$ (the triangle inequality follows from \ref{chpt2:lem:totalvariation}).

We are interested in measuring the convergence of a Markov chain $P^{t}$ to its 
stationary distribution. To understand this we need to analyse the total 
variation distance $\lVert P^{t}(g,\cdot) - \pi \rVert_{\textnormal{\tiny TV}}$ 
as a function of $t$. The following lemma demonstrates that $\lVert P^{t}(g,\cdot) - \pi 
\rVert_{\textnormal{\tiny TV}}$ is a non-increasing function; this follows the 
logic that performing another step of our random walk should never bring 
us further away from our equilibrium.

\begin{lemma}
	\label{chpt2:lem:TVdecreasing}
	Let $P$ be a transition matrix and $\mu,\nu $ be probabilities on 
	$G$. Then
	\[\lVert  P\mu -  P \nu\rVert_{\textnormal{\tiny TV}} \leq \lVert \mu - \nu 
	\rVert_{\textnormal{\tiny TV}} .\]
\end{lemma}

\begin{proof}
	Going from left to right we have:
	\begin{eqnarray*}
		\lVert  P \mu -  P \nu \rVert_{\textnormal{\tiny TV}} & = & 
		\frac{1}{2}\sum_{g\in 
			G} |  P\mu(g) -  P\nu(g)|\\
		& = &  \frac{1}{2}\sum_{g\in G} | \sum_{h \in G} 
		P(gh^{-1})(\mu(h) -\nu(h))|\\
		& \leq &  \frac{1}{2}\sum_{g\in G} \sum_{h \in G} 
		|P(gh^{-1}) | | \mu(h) -\nu(h) | =  \lVert \mu - \nu\rVert_{\textnormal{\tiny TV}}.  \qedhere
	\end{eqnarray*}
\end{proof}

\begin{corollary}
	Let $P$ be a transition matrix for a random walk on a finite group $G$ with stationary distribution $\pi$, then for all $t\in\mathbb{N}$ and $g\in G$ we have:
	\[\lVert P^{t+1}(g,\cdot) - \pi \rVert_{\textnormal{\tiny TV}} \leq \lVert 
	P^{t}(g,\cdot) - 
	\pi \rVert_{\textnormal{\tiny TV}}\]
\end{corollary}

\begin{proof}
	This is a consequence of Lemma \ref{chpt2:lem:TVdecreasing} with $\mu= 
	P^{t}(g,\cdot),\, \nu =\pi$. 
\end{proof}

Our definition of total variation distance so far has depended on our starting 
state $g$. When bounding the speed at which a random walk convergences to its stationary distribution we look at the time it takes from its worst possible starting state. For irreducible random walks on finite groups every starting state gives the same total variation distance.

\begin{lemma}
	\label{chpt2:lem:tvtran}
	Let $P$ define an irreducible random walk on $G$. Then 
	for any $g,h\in  G$ we have 
	\[\lVert P^{t}(g,\cdot) - \pi \rVert_{\textnormal{\tiny TV}} = \lVert P^{t}(h,\cdot) 
	- \pi \rVert_{\textnormal{\tiny TV}}.\]
\end{lemma}

\begin{proof}
	Let $P$ be the transition probability for our Markov chain, and $\phi$ be our probability preserving bijection with  $\phi(g) = h$, from the definition of transitivity. It 
	follows that from the 
	proof of Lemma \ref{chpt2:lem:groupstrans} that $P^{t}(a,b) 
	= P^{t}(\phi(a),\phi(b))$ for all $a,b \in G$. Therefore,
	\[\lVert P^{t}(g,\cdot) - \pi \rVert_{\textnormal{\tiny TV}} = 
	\frac{1}{2}\sum_{a\in G}|P^{t}(g,a) - \pi(a) | = \frac{1}{2}\sum_{a\in 
		G}|P^{t}(\phi(g),\phi(a)) - \pi(\phi(a)) | =\lVert P^{t}(h,\cdot) - \pi 
	\rVert_{\textnormal{\tiny TV}}. \qedhere\]
\end{proof}

Following Lemma \ref{chpt2:lem:tvtran} we start all our random walks on 
groups at the identity of the group. We introduce the function $\dt(t)$ as 
condensed notation for the (worst) total variation distance of a Markov chain.
\begin{defn}
	\label{chpt2:def:dt}
	Let $P$ be a probability which defines an irreducible random walk on $G$ 
	and $t\geq 0$. Define the function $\dt(t)$  as follows
	\[\dt(t) := \max_{g\in G} \lVert P^{t}(g,\cdot) - 
	\pi \rVert_{\textnormal{\tiny TV}} = \lVert P^{t} - \pi 
	\rVert_{\textnormal{\tiny TV}}.\]
\end{defn}

\subsubsection{Separation Distance}

We now provide a second notion of distance on probabilities as an alternative to total variation distance.
Separation distance is a commonly used measure of distance between a Markov 
chain driven by $P$ and its stationary distribution. Unlike total variation distance to 
define separation distance we require the stationary distribution to exist and 
be unique.

\begin{defn}[Section 6.4 \cite{Levin2017}]
	Let $P$ be a transition matrix for a Markov chain with unique stationary distribution $\pi$. Define \emph{separation distance} as follows:
	\begin{eqnarray}
	\lVert P^{t}(g,\cdot) -\pi \rVert_{\textnormal{Sep}} = \max_{h \in G} 
	\left( 1 -\frac{P^{t}(g,h)}{\pi(h)}\right).
	\end{eqnarray}
\end{defn}

Separation distance takes values in $[0,1]$ and tells us the maximum ratio of the probability being of in single state $h$ against the uniform distribution.

Similarly to 
total variation distance for a random walk on a group $G$ separation distance does not depend on the starting state $g$, that is  for any $g,h \in G$, we have $\lVert P^{t}(g,\cdot) - \pi \rVert_{\textnormal{Sep}} = \lVert P^{t}(h,\cdot) - \pi \rVert_{\textnormal{Sep}}$.

\begin{defn}
	Let $P$ define by an irreducible random walk on a group $G$ and 
	$t\geq 0$. Define the function $\ds(t)$  as follows
	\[\ds(t) := \max_{g\in G}\lVert P^{t}(g,\cdot)- \pi 
	\rVert_{\textnormal{Sep}} = \lVert P^{t} - \pi 
	\rVert_{\textnormal{Sep}}.\]
\end{defn}

The separation distance of an random walk on a group is always decreasing in time, $\ds(t+1) \leq \ds(t)$, for any $t\geq 0$ \cite[Section 6.4]{Levin2017}.
From our 
definitions it is unclear whether separation distance is statistically different 
from total variation distance. 

\begin{example}
	\label{chpt2:ex:tvsep}
	Consider the simple random walk on $\mathbb{Z}_{5}$ started at point $0$. 
	The table below gives the separation and total variation distance for 
	$0\leq t \leq 6$:
	\begin{table}[H]
		\begin{tabular}{c|c|c|c|c|c|c|c}
			$t$ & $0$ & $1$ & $2$ & $3$ & $4$ & 5 & 6\\ \hline
			$\ds(t)$ & 1 & 1 & 1 & 1 & 
			0.6875 & 0.6875 & 0.453125\\ \hline
			$\dt(t)$ & 0.8 & 0.6 
			& 0.4 & 0.35 & 0.275 & 0.225 & 0.18125
		\end{tabular}
		\caption[The values of separation and total variation distance for the simple random walk on $\mathbb{Z}_{5}$]{The values of separation and total variation distance for the simple random walk on $\mathbb{Z}_{5}$.}
		\label{chpt2:table:tvsep}
	\end{table} 
\end{example}

In Table \ref{chpt2:table:tvsep} we can see that total variation distance is 
always smaller than separation distance. This is not specific to our example 
and holds for all Markov chains for which separation distance is well defined.

\begin{lemma}
	\label{chpt2:lem:sepupperbound}
	Let $P$ define an irreducible random walk on a finite group $G$ with stationary 
	distribution $\pi$. Separation distance forms an upper bound on total 
	variation distance, that is for all $t\geq 0$,
	\begin{eqnarray}
	\dt(t) \leq \ds(t) .
	\end{eqnarray}
\end{lemma}

\begin{proof}
	Following from the definition of total variation distance we find:
	\begin{eqnarray*}
		\lVert P^{t} -\pi \rVert_{\textnormal{\tiny TV}}	= \sum_{\substack{h 
				\in G \\ P(h) < \pi(h)}} 
		\pi(h) - P^{t}(h)  = \sum_{\substack{h \in G \\ 
				P(h) < \pi(h)}} \pi(h)
		\left(1 - \frac{ P^{t}(h)}{\pi(h)} \right) \leq \max_{h \in  
			G} \left(1	-\frac{ P^{t}(h)}{\pi(h)} \right) \qedhere
	\end{eqnarray*}
\end{proof}

\subsubsection{Bounds on Total Variation and Separation Distance}

With our new definitions in hand we may state one of the cornerstone theorems 
in Markov chain theory, the convergence of irreducible aperiodic Markov chains to their stationary distributions. A proof of the following result may be found in \cite[Theorem 4.9]{Levin2017}.

\begin{thm}
	\label{chpt2:thm:convergence}
	Suppose that $P$ is an irreducible, aperiodic random walk on a finite group $G$, with 
	stationary distribution $\pi$. Then there exists some constants $c \in 
	[0,1)$ and 
	$A>0$ such that
	\begin{eqnarray}
	\lVert P^{t} - \pi 
	\rVert_{\textnormal{\tiny TV}} 
	\leq A c^{t}
	\end{eqnarray}
	Therefore, $\dt(t) \to 0 $ as 
	$t\to\infty$. Similarly $\ds(t) \to 0 $ as $t \to 
	\infty$.
\end{thm}

This theorem tells us that any irreducible aperiodic random walk 
necessarily converges to its stationary distribution and 
moreover we can bound the rate of this convergence. The time at which our 
random walk is close to uniform is called the \emph{mixing time} of the Markov 
chain.

\begin{defn}[Section 4.5 \cite{Levin2017}]
	\label{chpt2:def:mixingtime}
	Let $\{X^{t}\}$ be a Markov chain with transition matrix $P$. Define the \emph{$\varepsilon$-total variation mixing time}
	$\mt(\varepsilon)$ of our Markov chain as:
	\begin{eqnarray}
	\mt(\varepsilon) := \min \{t : \dt(t) \leq \varepsilon\}
	\end{eqnarray}
	Now define the \emph{total variation mixing time} $\mt$ of a Markov chain to be 
	$\mt:=\mt(1/4)$ (the choice of $1/4$ is semi-arbitrary, see \cite[Section 4.5]{Levin2017}).
	
	We may define the \emph{separation mixing time} $\ms(\varepsilon)$ and $\ms$, 
	in an analogous way replacing $\dt$ with $\ds$.
\end{defn}

From now on all random walks on groups we study will be irreducible, and aperiodic, meaning that they always have convergence to their uniform distribution. It is then of interest to ask what this convergence looks like for a given random walk.
Theorem \ref{chpt2:thm:convergence} allows us to bound the mixing time of a random walk from 
above. Without any more assumptions about our Markov chain there is not much 
more we can say about the rate of the convergence.

The techniques and tools used to 
bound mixing times have long been developed by probabilists. However, it is often the case that ad hoc techniques need to be developed specific to the random walk being studied. We 
introduce the key probabilistic and algebraic tools in this chapter.  In Chapter \ref{chpt4:chpt} we shall see how we have to use a 
combination of tools to prove precise bounds on the mixing time for the 
random transposition shuffle and the top-to-random shuffle. To begin with we state a 
classical upper bound on total variation distance for reversible random walks. 

\begin{thm}[Classical $\ell^{2}$ bound - see Lemma 12.16 \cite{Levin2017}]
	\label{chpt2:thm:classicL2}
	Let $P$ be the transition matrix for a reversible, transitive, irreducible, 
	aperiodic, 	Markov chain on a finite group $G$, 
	with stationary distribution $\pi$. Then we may label the eigenvalues so that $1=\beta_{1} > \beta_{2} \geq \dots \geq 
	\beta_{|G|} >-1$. Furthermore, we have the following bound,
	\begin{eqnarray}
	\label{chpt2:eqn:classicalL2}
	4\lVert P^{t} - \pi \rVert^{2}_{\textnormal{TV}} \leq \sum_{i\neq 1} 
	\beta_{i}^{2t} .
	\end{eqnarray}	
\end{thm}

For irreducible random walks on groups the 
detailed balanced equations \eqref{chpt2:eqn:detail} reduce to $P(g,h) = 
P(h,g)$, therefore a random walk on a group is reversible if and only if $P(g^{-1}) =P(g)$. We make use of Theorem \ref{chpt2:thm:classicL2} in the analysis of the random-to-random shuffle in Chapter \ref{chpt4:chpt}, and the one-sided transposition shuffle in Chapter \ref{chpt5:chpt}. Another popular method to establish an upper bound on separation distance (and therefore total variation distance) is the use of strong stationary times. To 
define a strong stationary time, we first need to define the notion of a 
stopping time.

\begin{defn}
	Let $\tau$ be a random variable taking values in $\mathbb{N}^{0}$. We call 
	$\tau$ 
	a \emph{stopping} time for a Markov chain $\{X^{t}\}$, if we can decide 
	the event $\{\tau = t\}$ with the knowledge of our Markov chain up to time $t$, that is the states $\{X^{0}, \ldots, X^{t}\}$.
\end{defn}

\begin{example}
	Let $G = S_{5}$ be the symmetric group on $5$ elements, and define a random walk on  $G$ given by the uniform distribution $\pi$.  We may view $S_{5}$ as the permutations of a deck of $5$ cards. We start at the identity permutation, and at every step of our random walk choose a permutation to apply uniformly at random. 
	Let $\tau$ be the first time card $5$ is moved. Then 
	$\tau$ is a stopping time for our random walk because if we observe our 
	walk up to time $t$, we can tell whether card $5$ has been moved and if so when was the first time this happened in our random walk.
	Suppose instead that $\tau$ is the last time card $5$ is moved up until a 
	fixed future time $T$. To even decide the event $\{\tau=1\}$ we would need to know the entire future of our Markov chain up to time $T$, and so $\tau$ is no longer a stopping time.
\end{example}

Strong stationary times are a special kind of stopping times. We impose the extra conditions that once a strong stationary time is met the distribution of $X^{t}$ must be stationary and independent of $\tau$.

\begin{defn}
	\label{chpt2:def:sst}
	Let $\tau$ be a stopping time for an irreducible, aperiodic Markov chain 
	$\{X^{t}\}$ which is a random walk on a group. We say $\tau$ is a 
	\emph{strong stationary time} if

	\begin{eqnarray}
	\mathbb{P}(X^{t} = g \, | \, \tau =t) = \pi(g).
	\end{eqnarray}
\end{defn}

From the definition we swiftly prove how strong stationary 
times may 
be used to bound separation distance and thus total variation distance for random walks.
\begin{lemma}
	\label{chpt2:lem:sst}
	Let $\tau$ be a strong stationary time for an irreducible, aperiodic random 
	walk on $G$ with driving probability $P$. The following holds for all $t\geq 1$:
	\[\lVert P^{t} -\pi \rVert_{\textnormal{\tiny TV}} \leq \lVert P^{t} -\pi 
	\rVert_{\textnormal{Sep}} \leq \mathbb{P}(\tau>t).\]
\end{lemma}
\begin{proof}
	The first inequality comes directly from Lemma 
	\ref{chpt2:lem:sepupperbound}. For the second inequality we start at the 
	definition of separation distance,
	\[1 -\frac{P^{t}(g)}{\pi(g)} \leq 1 -\frac{\mathbb{P}(X^{t}=g, \tau \leq 
		t)}{\pi(g)} = 1 -\frac{\pi(g)\mathbb{P}(\tau \leq 
		t)}{\pi(g)} = \mathbb{P}(\tau >t) \qedhere .\]
\end{proof}

Strong stationary times are useful because they reduce the analysis of 
separation distance based on the probability $P^{t}$, to the analysis of a 
single random variable $\tau$. In practice 
strong stationary times can be tricky to find but we will see 
two examples of them in this thesis: a classical argument for the top-to-random shuffle in Chapter \ref{chpt4:chpt}, and an original argument for the one-sided transposition shuffle in Chapter \ref{chpt5:chpt}.

All the methods we have mentioned so far give upper bounds on the mixing time 
of a random walk. When analysing the rate of convergence of a Markov chain is 
it also useful to bound the mixing time from below. To formulate a lower bound we frequently make use of the following simple lemma.

\begin{lemma}
	\label{chpt2:lem:lowerbound}
	Let $P$ be a random walk on a finite group $G$ with stationary distribution $\pi$. 
	Suppose $A \subseteq G$ then
	\[|P^{t}(A) - \pi(A) | \leq \lVert P^{t} - \pi \rVert_{\textnormal{\tiny 
			TV}}.\]
\end{lemma}
\begin{proof}
	The inequality follows from Definition \ref{chpt2:def:totalvariation}.
\end{proof}

This bound may seem na\"{i}ve but throughout this thesis  we shall see that it is a very adaptable tool. The key idea behind using Lemma 
\ref{chpt2:lem:lowerbound} for a lower bound on the  mixing time of a random walk $P$ is to find a set $A$ 
which tells apart $P^{t}$ from $\pi$. Usually we look for $A$ that has small 
probability under the stationary distribution, and that has high probability 
under the probability $P^{t}$ for time $t< \mt$. The set 
$A$ is often chosen based on information about our group $G$. For random walks 
on the symmetric group a commonly used idea is to take $A$ as a set of permutations with certain fixed points with 
the intention of reducing the calculation of $P^{t}(A)$ to an estimation of how 
long it takes for each fixed point to be eliminated. This allows us to use a modified coupon collectors argument to estimate our probability. We use this idea to formulate lower bounds for both the random transposition shuffle and the one-sided transposition shuffle.

\subsection{The Cutoff Phenomenon}

\begin{figure}[H]
	\includegraphics[scale=0.17]{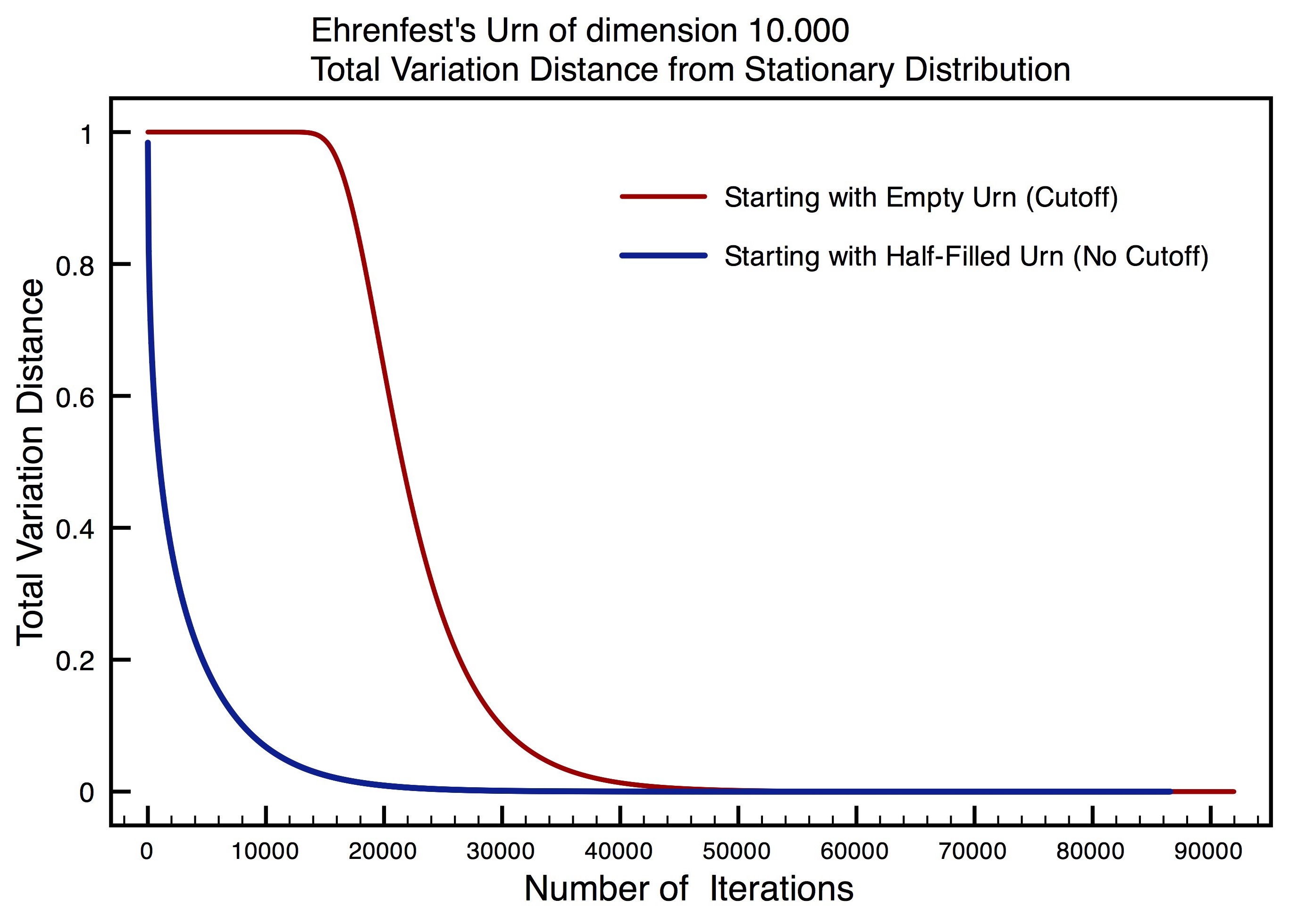}
	\caption[Total variation distance of Ehrenfest's urn from two different starting distributions]{Ehrenfest's urn is a random process starting with $n$ balls distributed between two urns, and at each step  a ball is chosen uniformly at random and moved to the opposite urn. The lines shows the total variation distance of the process from its uniform distribution from different starting states. The blue process starts with both urns containing $1/2$ the total amount of balls, the red process starts with one urn full and the other empty. The figure is taken from \cite{cutoffPresentation}.}
	\label{chpt2:fig:cutoff}
\end{figure}

In practice, for many Markov chains, the total variation distance to uniform does not decrease steadily as suggested by Theorem \ref{chpt2:thm:convergence}. Rather, a phenomenon is often seen where the distance decreases sharply past a critical time as show in Figure \ref{chpt2:fig:cutoff}: The blue line shows an exponential convergence from the start, whereas the red line only shows convergence to $0$ after a certain time threshold has been met.

Furthermore, there exist many families of Markov 
chains $\{X_{n}^{t}\}$ with transition probabilities $\{P_{n}\}$ that have the 
transition from $1$ to $0 $ becoming sharper as $n\to\infty$, 
thus when such a phase transition is 
present we may think about 
the mixing time of a random walk being the time of this phase transition.
We now switch our view to now 
consider a family of Markov chains $\{X^{t}_{n}\}$ with transition matrices 
$\{P_{n}\}$ and stationary distributions $\{\pi_{n}\}$ where $n \in \mathbb{N}$ 
indexes the family. The phase transition which has been observed in many 
families of random walks has been given the name \emph{the cutoff phenomenon}. 

\begin{defn}[Section 18.1 \cite{Levin2017}]
	\label{chpt2:def:cutoff}
	A sequence of random walks $\{X_{n}^{t}\}$ on finite groups $\{G_{n}\}$ with driving probabilities 	$\{P_{n}\}$ and 
	stationary distributions $\{\pi_{n}\}$ exhibits a (total 
	variation/separation) \emph{cutoff} at time 
	$t_{n}$ with a \emph{window} of size $w_{n}$ if $w_{n}=o(t_{n})$ 
	and the 
	following limits hold:
	\begin{eqnarray}
	\lim_{c\rightarrow \infty}\limsup_{n\rightarrow \infty} \lVert P_{n}^{t_{n} 
		+cw_{n}} -\pi_{n} \rVert_{\textnormal{TV/Sep}}&=&0 
	\label{chpt2:eqn:upper} \\
	\lim_{c\rightarrow \infty}\liminf_{n\rightarrow \infty} \lVert P_{n}^{t_{n} 
		-cw_{n}} -\pi_{n} \rVert_{\textnormal{TV/Sep}} &=&1 
	\label{chpt2:eqn:lower}.
	\end{eqnarray}
	The limits \eqref{chpt2:eqn:lower} and \eqref{chpt2:eqn:upper}, respectively define upper and lower bounds on the mixing time of the random walk .
\end{defn}

The heuristic behind this definition is the phase transition for the family of 
chains happens at time $t_{n}$ and this transition becomes sharper as 
$n\to\infty$. If a family shows a cutoff at time $t_{n}$ then we must have that 
$\mt \sim t_{n}$.  An equivalent definition (see \cite[Lemma 18.1]{Levin2017}) of a 
cutoff is the function $\dt(\cdot)$ converges to a step 
function around the mixing time of the random walk. Precisely, a family of random walks $P_{n}$ exhibits a total variation cutoff if and only if the following limit holds:
\[\lim_{n\to\infty} \dt(c\, \mt) = \begin{cases} 
1 & \textnormal{ if } c>1 \\
0 & \textnormal{ if } c<1 \end{cases}.\] 
Once we have established a cutoff for a 
family of Markov chains we have complete asymptotic information about the families convergence to their stationary distributions. Definition \ref{chpt2:def:cutoff} is given in 
terms of both total variation and separation distance. We may hope that a cutoff in one would imply a cutoff in the other, however this is not the case \cite{hermon2016total}.
This does not mean that a cutoff in total variation cannot help us prove 
a cutoff in separation or vice versa. By Lemma 
\ref{chpt2:lem:sepupperbound} if we establish an upper bound on separation 
distance we get an upper bound on total variation distance, similarly a lower bound on 
total variation distance gives a lower bound on separation distance. For the one-sided 
transposition shuffle we first prove a cutoff in total variation distance then a 
cutoff in separation distance with the lower bound following from our first proof. 

To 
establish a cutoff it is not enough to just find the order of $t_{n}$, 
although this is a useful first step. Computing the exact constant factor of the correct order for a cutoff is often where the difficultly in the analysis of mixing times stems from. Frequently finding the  correct cutoff time 
for a random walk involves 
an initial ansatz, and then refinement to make sure both the upper and lower 
bounds hold. We prove the existence of a cutoff for several random walks 
on the symmetric group in Chapter \ref{chpt4:chpt}, for the moment we state some simple cutoff results.

\begin{example}
	\label{chpt2:ex:lazysimple}
	Consider the family of lazy simple random walks on $\mathbb{Z}_{n}$, which are driven by probabilities
	\[P_{n}(i)= \begin{cases}
	\frac{1}{2} & \text{ if } i =0 \\
	\frac{1}{4} & \text{ if } i = \pm 1\\
	0 & \text{ otherwise }
	\end{cases}.\]
	This random walk is irreducible, aperiodic and transitive. We may ask 
	whether $\dt(t)$ shows a cutoff as $n\to\infty$. The answer is no, this 
	family of random walks fails to show a cutoff \cite[Section 
	3.3]{saloff2004random}. The  lazy simple random walk on a circle steadily 
	spreads out from the starting point $0$ meaning there is no phase transition and thus no cutoff. 
\end{example}

\begin{example}
	\label{chpt2:ex:hypercube}
	Consider the hypercube $\mathbb{Z}_{2}^{n}$, define the identity of the group as $e_{0}$, the zero vector. For $i\in\{1,\dots n\}$, define $e_{i}$ as
	the element which has $1$ in position $i$, and $0$ in every other position. The nearest neighbour walk on the hypercube is driven by probability $P_{n}$ defined by
	\[P_{n}(g)= \begin{cases}
	\frac{1}{n+1} & \text{ if }g = e_{i} \text{ for } i\in\{0,1,\ldots,n\}  
	\\0 & \text{ otherwise }\end{cases}.\]
	The probabilities $\{P_{n}\}$ define a family of irreducible, aperiodic, and transitive random walks, 
	which show a cutoff in total variation distance at time 
	$t_{n} = \frac{n+1}{4}\log n$, with a window of $w_{n} = \frac{n+1}{4}$. 
	Therefore we may say the mixing time of this random walk is exactly $\mt = 
	\frac{n+1}{4} \log n$. Details of this random walk may be found at \cite[Chapter 3C]{Diaconis1988}.

\end{example}

\paragraph{}
Since its formulation the cutoff phenomenon has been a great area of interest 
in the study of random walks. Establishing a cutoff for any particular example  
requires detailed knowledge and understanding of the inner workings of the Markov chain. The upper bound \eqref{chpt2:eqn:upper} and the lower bound \eqref{chpt2:eqn:lower} 
have to be proved separately and require the use of a wide variety of 
techniques. We would instead prefer to find a simple criterion for a random walk exhibiting a cutoff; however,  no 
such result currently exists and in fact there are very few global results 
on total variation cutoff of random walks. Contrary to this much is known about cutoffs in $\ell_{p}$-distance with $1<p<\infty$ \cite{chen2008cutoff}, which further highlights why the lack of results for $p=1$ (total variation) is so surprising.  Why any family of Markov chains shows a total variation cutoff is a mystery with many explanations proposed. One of the most popular theories is that the existence of a cutoff is linked to the 
multiplicity and value of the second biggest eigenvalue.

\begin{conj}[Peres's Conjecture \cite{peresconj}]
	Let $\{P_{n}\}$ be a family of transition matrices 
	Markov chains, with 
	second biggest (in absolute value) eigenvalue $\beta_{n}$. The family shows a 
	total variation cutoff at 
	time $t_{n}$ if and only if
	\begin{eqnarray}
	\label{chpt2:eqn:peres}
	t_{n}(1 - \beta_{n})  \rightarrow \infty \textnormal{ as } n 
	\to \infty .
	\end{eqnarray}
\end{conj}
If we measure cutoff instead in $\ell_{p}$-distance for $1<p<\infty$, Chen and Saloff-Coste \cite{chen2008cutoff} have managed to prove a criterion for a cutoff in the spirit of Peres' conjecture. 
If a family of Markov chains shows a total variation cutoff, condition \eqref{chpt2:eqn:peres} is known to be satisfied.  
On the other hand Aldous managed to find a counterexample for Peres' conjecture using random walks on graphs, see \cite[Section 4.2]{chen2006cutoff} for details. However, if we specialise to random walks on finite groups then there is currently no known counterexample to Peres's conjecture. 
Random walks on finite groups are an active area of research, the algebraic setting allowing for specialised techniques in order to prove cutoff results. In the next section we cover the 
algebraic tools available to find the mixing time of a random walk on $G$, these require background about the representations, modules and characters of groups.

\section{Representations, Modules, and Applications to Random Walks}
\label{chpt3:chpt}

Group representation theory aims to  study groups using linear algebra. The goal is to reduce the abstract elements of groups into the well understood elements of linear algebra. To do this we represent the elements of our group as matrices in such a way that the structure of the group is respected. For the applications in this thesis it is enough to develop the representation theory of groups over the complex numbers, so all vectors spaces, matrices, etc., are over the field $\mathbb{C}$ unless otherwise stated. 
This assumption helps to simplify the 
representation theory presented and allows the use of several theorems which do not hold over general fields. A detailed account of all the representation theory presented in this section can be found in \emph{Representations and Characters of Groups} by James and Liebeck\cite{james2001representations}. We begin this section by introducing the representations and characters of a finite group $G$.  We focus on the irreducible representations of $G$ from which any representation may be built. The following section is devoted to the properties of the discrete Fourier transform, which is used to analyse probabilities on $G$ using the group's irreducible characters. Our last section covers the upper bound lemma, which establishes a close link between random walks on $G$ and the representation theory of $G$.

\subsection{Representations, Characters and Modules}

\subsubsection{Representations and Modules}

\begin{defn}
	Let $G$ be a finite group. A \emph{representation} of $G$ is a vector space $V$ 
	together with a \emph{group homomorphism} from $G$ to the general linear group $GL(V)$, denoted $\rho:G \to GL(V)$,
	which has the property 
	$\rho(gh) = \rho(g)\rho(h)$ for all $g,h \in G$. We say the pair $(\rho,V)$ is a representation. The dimension of $V$ is denoted $d_{\rho}$ and is called 
	the \emph{dimension} of the representation $(\rho,V)$.
\end{defn}

The vector space $V$ our matrices $\rho(\cdot)$ act on is critical to our knowledge of the representation. Thus, we write $(\rho,V)$ when we want to emphasise the space $V$. If the space $V$ is clear from 
context we will just talk about the group homomorphism $\rho$ being a representation by itself. 
We now introduce the notion of $G$ modules, by equipping our group with a vector space structure.

\begin{defn}
	Let $G$ be a finite group. The \emph{group algebra} of $G$ is the vector space $\mathbb{C}[G] = \langle g \, | \, g \in G 
	\rangle$, spanned by elements of $G$ with coefficients in $\mathbb{C}$. This vector space has a natural multiplication	given by multiplying complex coefficients and the elements of $G$ separately, e.g. for $g,h,f \in G$ we have $(4g) 
	(ih  - \frac{1}{2}f) = 4i(gh) - 2(gf)$.
\end{defn}

\begin{defn}
	Let $G$ be a finite group. A \emph{$\mathbb{C}[G]$-module}  is a vector space $V$ with a 
	bilinear multiplication $\mu:\mathbb{C}[G]\times V \to V$ such that: $\mu(1,v) = v$ for all 
	$v\in V$, and for all $g,h \in G$, $\mu(h,\mu(g,v)) = \mu(hg,v)$. To simplify notation we frequently write the multiplication $\mu(g,v)$ as just $g \, v$, and say $\mathbb{C}[G]$ acts on $V$. 
\end{defn}

We can think of a $\mathbb{C}[G]-$module as a vector space $V$ which our group $G$ acts on. For any group $G$, and vector space $V$, we may always define a trivial action by setting, $g \, v = v$ for all $g \in G, v\in V$. The \emph{trivial module} of a group $G$ is formed from a one-dimensional vector space with the trivial action. The \emph{trivial representation} of a group $G$ is a one-dimensional vector space with group homomorphism $\rho(g) = \I $ for all $g \in G$. The trivial module and trivial representation are two different ways to view the same algebraic object. The next result established that the representations of $G$ (over field $\mathbb{C}$) and $\mathbb{C}[G]-$modules are in one-to-one correspondence.

\begin{thm}[Theorem 4.4 \cite{james2001representations}]
	\label{chpt3:thm:1to1}
	Let $G$ be a finite group. The representations of $G$ and  $\mathbb{C}[G]$-modules are in one-to-one 
	correspondence. If $(\rho,V)$ is a representation for $G$ then $V$ is a $\mathbb{C}[G]$-module with multiplication defined by $\mu(g,v) = \rho(g)v$ for all $g \in G$. On the 
	other hand if $V$ is a $\mathbb{C}[G]$-module with multiplication $\mu$, define $\rho(g) = 
	\mu(g,\cdot)$, then $(\rho,V)$ is a representation for $G$.

\end{thm}

\begin{proof}
	The constructions in the theorem may be verified by checking the 
	definitions of modules and representations. To prove this is a one-to-one 
	correspondence take a representation $(\rho,V)$ and form the module $V$ 
	with action $\mu$ via 
	the process described in the theorem. Now we turn the module $V$
	into a representation $(\varphi,V)$ with a new mapping $\varphi(g) = 
	\mu(g,\cdot) = \rho(g)$, 
	therefore we recover our original representation $(\rho,V)$. We may perform 
	a similar check starting with a module $V$ and transforming it to a 
	representation and back to a module recovering $V$. Thus, 
	representations and modules are in one-to-one correspondence. 
\end{proof}

From now on out we may talk about representations and modules interchangeably, and every result we state is valid for either view point via the correspondence described in Theorem \ref{chpt3:thm:1to1}. We now give a brief example of some natural representations.

\begin{example}
	Let $S_{3}$ denote the symmetric group on $3$ elements.  	Consider the space $\mathbb{C}^{3}$, we may form an action of $S_{3}$ on $\mathbb{C}^{3}$, by permuting the basis elements. Realising this action in terms of matrices gives us a group homomorphism $\rho$ detailed in Table~\ref{chpt3:table:repexample}, therefore the pair $(\rho, \mathbb{C}^{3})$ defines a representation of $S_{3}$.	For another module of $S_{3}$ take the 2-dimensional subspace $W = \{v =(v_{1},v_{2},v_{3})^{T} \in \mathbb{C}^{3} \, | \, v_{1}+v_{2}+v_{3} = 0\}$, spanned by vectors $w_{1} = (1, -1,0)^{T}$, and $w_{2} = (0, 1,-1)^{T}$. Again consider an action of $S_{3}$ on $W$ by permuting the coordinates of $\mathbb{C}^{3}$, e.g., $(1\,2) \, w_{1} = -w_{1}$ and $(1 \, 2) \,w_{2} = w_{1} + w_{2}$. Realising this action as matrices on $W = \langle w_{1}, w_{2}\rangle$ we find the group homomorphism $\varphi$ shown in Table \ref{chpt3:table:repexample},  hence $(\varphi,W)$ is a representation of $S_{3}$.

	\begin{table}[H]
		\makebox[\linewidth]{
			\renewcommand{\arraystretch}{1.25}
			\begin{tabular}{c|c|c|c|c|c|c}
				$S_{3}$ & $\gpid$ & $(1\,2)$ & $(1 \, 3)$ & $(2 \,3)$ & $(1\,2\,3)$ & $(1\,3\,2)$\\\hline
				$\rho(\cdot)$ & $\left({\begin{array}{ccc}
					1 & 0 & 0  \\
					0 & 1 & 0\\
					0 & 0 & 1
					\end{array}} \right)$	&$\left({\begin{array}{ccc}
					0 & 1 & 0  \\
					1 & 0 & 0\\
					0 & 0 & 1
					\end{array} }\right)$ &$\left({\begin{array}{ccc}
					0 & 0 & 1  \\
					0 & 1 & 0\\
					1 & 0 & 0
					\end{array} }\right)$ &$\left( {\begin{array}{ccc}
					1 & 0 & 0  \\
					0 & 0 & 1\\
					0 & 1 & 0
					\end{array} }\right)$ & $\left({\begin{array}{ccc}
					0 & 1 & 0  \\
					0 & 0 & 1\\
					1 & 0 & 0
					\end{array} }\right)$& $\left({\begin{array}{ccc}
					0 & 0 & 1  \\
					1 & 0 & 0\\
					0 & 1 & 0
					\end{array} }\right)$ \\\hline
				$\varphi(\cdot)$ & $\left({\begin{array}{cc}
					1 & 0 \\
					0 & 1
					\end{array}} \right)$	&$\left({\begin{array}{cc}
					-1 & 0 \\
					1 & 1
					\end{array}}\right)$ &$\left({\begin{array}{cc}
					0 & -1 \\
					-1 & 0
					\end{array}}\right)$ &$\left( {\begin{array}{cc}
					1 & 1 \\
					0 & -1
					\end{array}}\right)$ & $\left({\begin{array}{cc}
					0 & 1 \\
					-1 & -1
					\end{array}}\right)$& $\left({\begin{array}{cc}
					-1 & -1 \\
					1 & 0
					\end{array}}\right)$ 
			\end{tabular}
		}
		\caption[The action of two group homomorphisms of the symmetric group $S_{3}$]{The table shows the group homomorphisms $\rho, \varphi$ on the vectors spaces $\mathbb{C}^{3}$ and $W$ respectively. Note that the elements $\rho(\cdot)$ are the 3-dimensional permutation matrices.

		}
		\label{chpt3:table:repexample}
	\end{table}

\end{example}

\begin{defn}
	Let $V$ be a $\mathbb{C}[G]$-module. We call a subspace $W \subseteq V$ \emph{$G$-stable} if $\{g \,w \, | \, g \in \mathbb{C}[G], \, w \in W\} = W$. A \emph{submodule} of $V$ is a stable subspace $W \subseteq V$. A module $V$ is called \emph{simple} if it has no non-trivial $(W\neq V,\{0\})$ submodules. Equivalently a representation $(\rho,V)$ is called \emph{irreducible} if it has no non-trivial subrepresenations.
\end{defn}

Take any one-dimensional vector space $V$, we know that  $\{0\}$ and $V$ 
are the only subspaces of $V$. 
Therefore, any one dimension module is necessarily simple. 
The simple modules of a group $G$ are the building blocks of its 
representation theory, in the sense that every module may be decomposed into its simple submodules, see Theorem \ref{chpt3:thm:decomposition}. We now introducing mappings 
between modules which preserve their structure, these are called homomorphisms.

\begin{defn}
	Let $V,W$ be $\mathbb{C}[G]$-modules. A \emph{homomorphism of modules} is a linear map $\psi:V\to W$ such that $\psi(g \, v) = g \, \psi 
	(v)$. An \emph{isomorphism of modules} is a bijective homomorphism $\psi:V \to W$, in this case we say $V$ is isomorphic to $W$ as $\mathbb{C}[G]$-modules, denoted $V \cong W$.	
	Define the 
	\emph{kernel} and \emph{image} 
	of a homomorphism $\psi:V\to W$ as follows: $\Ker(\psi) = \{v : \psi(v) = 0\}$, $\im(\psi) = \{w : \psi(v) =w \text{ for some $v$}\}$.
\end{defn}

\begin{lemma}
	\label{chpt3:lem:kerim}
	Let $V,W$ be $\mathbb{C}[G]$-modules, and $\psi:V\to W$ a 
	homomorphism. Then $\Ker(\psi)$ 
	is a submodule of $V$, and $ 
	\im(\psi)$ is 	a submodule of $W$.  If 
	$\Ker(\psi) = \{0\}$ then $\psi$ is an injective map, and if $\im(\psi) = W$ then $\psi$ is a surjective map. 
\end{lemma}

\begin{proof}
	Take $v \in \Ker(\psi)$. For any element $g \in \mathbb{C}[G]$ we have $ gv \in \Ker(\psi)$ because $\psi(g \,v) = g \psi(v) = 0$. Hence, the kernel of $\psi$ is stable under the action of $\mathbb{C}[G]$ and so is a $\mathbb{C}[G]$-module contained in $V$. Similarly the image of $\psi$ is stable under the action of $\mathbb{C}[G]$ and so is a $\mathbb{C}[G]$-module contained in $W$. The last statements are standard facts of linear maps.
\end{proof}

We are now in a position to state Schur's 
Lemma which tells us that homomorphisms between simple modules are trivial.

\begin{lemma}[Schur's Lemma]
	\label{chpt3:lem:schur}
	Let $V,W$ be simple $\mathbb{C}[G]$-modules. If $\psi:V
	\to W$ is a  homomorphism then either $\psi$ is the zero morphism or an 
	isomorphism of representations. 	Furthermore, the only homomorphisms $\psi:V 
	\to V$ are scalar multiplies of the identity map, i.e., $\psi(v) =  \beta v$ for some $\beta \in \mathbb{C}$.
\end{lemma}

\begin{proof}
	By $V$ being a simple module we have two choices for the 
	kernel of $\psi$: $V$ or $\{0\}$. If $\Ker(\psi)= V$ then $\psi$ is clearly the 
	zero morphism. If $\Ker(\psi) = \{0\}$ then our map is injective, moreover 
	$\im(\psi) \neq \{0\}$ and hence must have 
	$\im(\psi) = W$ because $W$ is a simple. Therefore the map $\psi$ is a bijective homomorphism.
	
	Now assume $\psi:V\to V$ is a homomorphism. Since we are working over the field $\mathbb{C}$ the map $\psi$ must 
	have an eigenvalue $\beta$ with eigenvector $v$. Form a new homomorphism 
	$\psi^{\prime} = \psi - \beta \,\I$, with $\I$ being the identity homomorphism. We can see that $\psi^{\prime}(v) = 0$ 
	thus $\Ker(\psi^{\prime}) = V$, and so our new map $\psi^{\prime}$ must be the zero morphism. Hence, we have $\psi = \beta \, \I$.
\end{proof}

Isomorphism forms an equivalence relation on the set of $\mathbb{C}[G]$-modules. For any group $G$ we would like to form a complete collection of its simple modules (equivalently irreducible representations) up to isomorphism. 

\begin{lemma}
	Let $G$ be a finite group. Then there are only finitely many simple $\mathbb{C}[G]$-modules up to isomorphism (see Definition \ref{chpt3:def:regularmod}).
\end{lemma}
\begin{defn}
	
	Let $G$ be a finite group.	Define the \emph{complete collection of simple $\mathbb{C}[G]$-modules} as a set, denoted $\rep(G)$, such that if $V$ is an simple $\mathbb{C}[G]$-module then $V \cong M$ for exactly one $M \in \rep(G)$.
	
\end{defn}

We now look to decompose any reducible module into a direct sum of smaller submodules.

\begin{defn}
	Let $V_{1},V_{2}$ be $\mathbb{C}[G]$-modules. The \emph{direct sum} of vector spaces $V_{1}\oplus V_{2}$ is a $\mathbb{C}[G]$-module under the action $g \,(v_{1},v_{2}) = 
	(g\, v_{1},g \, v_{2})$ for $v_{1}\in V_{1}, \, v_{2} \in V_{2}$ and $g\in G$. Moreover, we have $\dim(V) = \dim(V_{1}) + \dim(V_{2})$.	Conversely, suppose $W$ is a $\mathbb{C}[G]$-module, with submodules $W_{1},W_{2}$, such that $ W_{1} \oplus W_{2} = W$, and $W_{1} \cap W_{2} = \{0\}$. Then our module $W$ may be decomposed into a direct sum $W = W_{1} \oplus W_{2}$. 
\end{defn}

We now prove that every stable subspace of a module $V$ has a stable complement.
Thus, any non-simple module may be split into a direct sum of smaller modules. This allows us to prove the powerful Maschke's Theorem.

\begin{lemma}[See Theorem 8.1 \cite{james2001representations}] 
	\label{chpt3:lem:spacedecomp}
	Let $V$ be a $\mathbb{C}[G]$-module. Suppose $V$ contains a stable 
	subspace $W$. Then $V$ contains a second stable subspace $W^{0}$ such that 
	$W \cap W^{0} = \{0\}$ and $W \oplus W^{0} = V$. The vector space $W^{0}$ 
	is called the \emph{complement} of $W$.
\end{lemma}

\begin{proof}

	Suppose $v_{1},\dots,v_{n}$ form a basis of $V$ with  $v_{1},\ldots v_{l}$ simultaneously being a basis for $W$. Let $\phi:V\to V$ be the projection onto the subspace $W$, defined for basis elements  $\phi(v_{i}) = v_{i}$ if $1\leq i \leq l$ and $0$ otherwise. Now define a new linear map $\psi:V \to V$ as follows:
	\[\psi(v) = \frac{1}{|G|}\sum_{g\in G} g^{-1} \, \phi(g \,v).\]
	The map $\psi$ is a homomorphism of $G$-modules. Indeed for any $h\in G$,
	\[h\,\psi(v) = \frac{1}{|G|}\sum_{g\in G}hg^{-1} \, \phi ( g \,v) = \frac{1}{|G|}\sum_{g\in G}g^{-1} \, \phi ( gh\,v) = \psi(hv).\]
	The image of $\psi$ must be contained in the submodule $W$ because of the presence of the projection $\phi$. In fact it is equal to $W$, taking $w\in W$ we have,
	\[ \psi(w) = \frac{1}{|G|}\sum_{g \in G} g^{-1}\phi( g w)  =\frac{1}{|G|}\sum_{g \in G} g^{-1} g w = w  .\]
	Set $W^{0} = \Ker(\psi)$, by Lemma \ref{chpt3:lem:kerim} this forms a stable submodule of $V$. 
	To summarise, $\psi$ is a module homomorphism, with $\im(\psi) = W$ and $\Ker(\psi)= W^{0}$ therefore we must have $W \oplus W^{0} =V$, and $W \cap W^{0} = \{0\}$.
\end{proof}

\begin{thm}[Maschke's Theorem]
	\label{chpt3:thm:decomposition} Suppose $V$ is a 
	module of a finite group $G$. Then we may 
	decompose $V$ into a direct sum of simple modules. This means the module $V$ may be decomposed as follows, 
	\begin{eqnarray}
	\label{chpt3:eqn:decomposition}
	V \cong \bigoplus_{M_{i} \in \rep(G)} k_{i} \; M_{i}
	\end{eqnarray}
	where the direct sum is over all simple modules $
	M_{i} \in \rep(G)$ and $k_{i}$ denotes the multiplicity of $M_{i}$ in the direct sum.
\end{thm}

\begin{proof}

	We proceed by induction on dimension. The statement is clearly true for one-dimensional modules. If $\dim(V) >1$ then either $V$ is simple, in which case we are done, or $V$ has a proper submodule $W$. By \ref{chpt3:lem:spacedecomp} we may write $V = W \oplus W^{0}$, for modules $W,W^{0}$ with $\dim(W), \dim(W^{0}) < \dim(V)$. Hence, $W$ and $W^{0}$ can be written as direct sums of simple submodules, which implies the same for $V$.
\end{proof}

Theorem \ref{chpt3:thm:decomposition} is essential 
to the study of modules (and representations) because it tells us that to understand the  
modules of $G$ it is enough to understand all the simple modules. However, to completely describe any module $V$ we also need to know the value of the constants $k_{i}$, before 
finding these multiplicities we should be sure that the decomposition \eqref{chpt3:eqn:decomposition} is unique 
up to isomorphism. This follows from an application of Schur's Lemma.

\begin{lemma}
	\label{chpt3:lem:unique}
	Let $G$ be a finite group. Then for constants $k_{i},l_{i} \in\mathbb{N}^{0}$ 
	we have:
	\[ \bigoplus_{M_{i}\in \rep(G)} k_{i} \; M_{i} \cong \bigoplus_{M_{i} \in 
		\rep(G)} l_{i} \; 	M_{i} \Leftrightarrow k_{i} = l_{i} \textnormal{ for all $i$}.\]
\end{lemma}

\begin{proof}
	If $k_{i}=l_{i}$ for all $i$ then the conclusion is immediate. Suppose that our direct sums are isomorphic with isomorphism $\psi$. If we restrict to a summand $M_{i}$ on the left hand side then we must get an isomorphic copy of $M_{i}$ on the right hand side, by Schur's lemma. Hence, the summand $k_{i}M_{i}$ must land in the summand $l_{i} M_{i}$ on the right. Counting the dimensions on each side now gives $l_{i}= k_{i}$ for all $i$.

\end{proof}

\begin{corollary}
	\label{chpt3:cor:decomposition}
	The decomposition shown in Theorem \ref{chpt3:thm:decomposition} is 
	unique up to isomorphism. 
\end{corollary}

\begin{defn}
	\label{chpt3:def:isocompotent}
	Let $G$ be a finite group, $M_{i}$ a simple module, and $V$ a module with decomposition $V \cong \oplus k_{i}	\, M_{i}$. The \emph{isotypic component} of $M_{i}$ in $V$ is the unique submodule $k_{i} \, M_{i}$.
\end{defn}

Corollary 
\ref{chpt3:cor:decomposition} tells us we can 
decompose any module into its 
simple parts but it does not tell us the values of the constants $k_{i}$. 
Finding efficient ways to know or compute the constants in a decomposition is an important problem in representation theory. To end this section we introduce  the group algebra as a $\mathbb{C}[G]$-module itself, called the \emph{regular module} of $G$.

\begin{defn}[Definition 6.5 \cite{james2001representations}]
	\label{chpt3:def:regularmod}
	The \emph{regular module} $\mathbb{C}[G]$ has the same vector space and action given by the multiplication of the 
	group algebra and it has dimension $|G|$.
	
\end{defn}

\begin{lemma}
	Let $G$ be a finite group.
	The regular module $\mathbb{C}[G]$ has canonical 
	decomposition
	\begin{eqnarray}
	\mathbb{C}[G] \cong \bigoplus_{M_{i} \in \rep(G)} d_{i} \, M_{i} 
	\end{eqnarray} 
	where $d_{i}$ is the dimension of the simple module $M_{i}$. We prove this decomposition in Lemma \ref{chpt3:lem:charreg}, after the introduction of characters.
\end{lemma}

Knowledge of the simple 
modules of the symmetric group $S_{n}$ will be key to our analysis of the one-sided transposition 
shuffle in Chapter \ref{chpt5:chpt}. 
In particular we study the action of the shuffle $P_{n}$  
on the regular module of $S_{n}$.
This allows us to use the decomposition in Definition \ref{chpt3:def:regularmod} to reduce our goal of finding eigenvalues of the regular module to finding 
eigenvalues of the action on the simple modules of $S_{n}$. Before we perform our analysis we shall present a detailed construction of the simple modules of the symmetric group in Chapter \ref{chpt4:chpt}.

\subsubsection{Character Theory}

In the beginning of this section we assumed the field we are working 
over to be $\mathbb{C}$. This assumption lets us condense the information of any 
representation $(\rho,V)$ into a single function $\chi_{\rho}:G \to \mathbb{C}$, 
called the \emph{character} of our representation. 
The characters gives is a simple way to view the information given by any 
representation.

\begin{defn}
	Let $(\rho,V)$ be a representation of a group $G$. Define the 
	\emph{character} 
	$\chi_{\rho}:G\to \mathbb{C}$ of a representation $\rho$ by $\chi_{\rho}(g) = \Tr \, 
	\rho(g)$, where $\Tr(\cdot)$ denotes the trace of a matrix. 
	We say the \emph{dimension} of the character, denoted $d_{\chi}$, is the same as its 
	corresponding representation.
	If $(\rho,V)$ is an irreducible representation then $\chi_{\rho}$ is called an 
	\emph{irreducible character}.
\end{defn}

\begin{lemma}
	Let $\chi$ be a character for the group $G$. We find:
	\begin{enumerate}
		\item $\chi(e) = d_{\chi}$
		\item For any $g \in G$ we have  $\chi(g^{-1}) = \overline{\chi(g)}$, 
		where $\overline{c}$ denotes the complex conjugate of $c \in \mathbb{C}$
		\item For any $g,h \in G$, we have $\chi(hgh^{-1}) = \chi(g)$.
	\end{enumerate}
\end{lemma}	

\begin{proof}

	\begin{enumerate}
		\item We clearly have $\chi(\gpid) = \Tr(\I) = d_{\chi}$
		\item We know that $\rho(g)^{n} = \rho(g^{n}) = \I$ for some $n$, so the eigenvalues of $\rho(g)$ must be roots of unity, denoted $\xi_{i}$. Then we have
		\[\overline{\chi(g)} = \overline{(\Tr \, \rho(g))} = \sum_{i} \overline{\xi_{i}} = \sum_{i}\xi_{i}^{-1} = \Tr \left(\rho(g)^{-1}\right) = \Tr \, \rho(g^{-1}) = \chi(g^{-1}) \]
		\item Follows from the property $\Tr(AB) = \Tr(BA)$ of the trace .
	\end{enumerate}
\end{proof}

The last of these properties tells us that characters are constant on the 
conjugacy classes of any group $G$. Representations and characters have a close 
connection to the conjugacy classes of the group they are defined over. In fact 
the number of irreducible representations of $G$ is exactly the number of conjugacy classes of $G$. To establish this we introduce an inner product on the space of class functions of $G$.

\begin{defn}
	Let $G$ be a finite group. A \emph{class function} for $G$ is a function $\phi:G\to\mathbb{C}$ that is constant on the conjugacy classes of $G$, i.e., for all $g,h\in G$ we have $\phi(g) = \phi(h^{-1} gh)$.  
	Let $\chi$, $\phi$ be class functions on $G$. Define the inner product $\langle \chi \, | \, \phi \rangle$ as follows:
	\begin{eqnarray}
	\langle \chi \, | \, \phi \rangle = \frac{1}{|G|}\sum_{g\in G} \chi(g) \, \overline{\phi(g)}.
	\end{eqnarray}
\end{defn}

\begin{thm}[See Chapter 15 \cite{james2001representations}]
	Let $\chi_{i}$ denote the character of the irreducible representation $\rho_{i} \in \rep(G)$. The irreducible characters are orthonormal, that is:
	\[
	\langle \chi_{i} \, | \, \chi_{j} \rangle =\begin{cases}
	1 & \textnormal{ if } i =j \\
	0 & \textnormal{ if } i\neq j
	\end{cases}.\]
	Furthermore, the irreducible characters form a basis for the class functions of $G$.
\end{thm}

\begin{corollary}
	\label{chpt3:lem:dimsum}
	Let $G$ be a finite group. The number of irreducible representations of $G$ is exactly the number of conjugacy classes of $G$. 
\end{corollary}
\begin{proof}
	Let $n$ be the number of conjugacy classes of $G$. The space of class functions on $G$ is spanned by exactly $n$ functions which take value $1$ on a single conjugacy class and zero otherwise. Therefore, our basis of irreducible characters must be formed from $n$ characters.
\end{proof}

Using the inner product on class functions and the irreducible characters $\chi_{i}$ we now give a second proof of the uniqueness of decomposition \eqref{chpt3:eqn:decomposition}.

\begin{lemma}
	\label{chpt3:lem:splitchar}
	Let $(\rho,V),(\varphi,W)$ be representations of a finite group $G$. Then we have $\chi_{\rho \oplus \varphi} = \chi_{\rho} + \chi_{\varphi}$.
\end{lemma}

\begin{proof}
	We may choose a basis of the space $V\oplus W$ such that the group homomorphism  has form $(\rho \oplus \varphi) (g) = \left({\begin{array}{cc} \rho(g) & 0 \\	0 & \varphi(g)	\end{array}} \right)$, then take traces.
\end{proof}

\begin{corollary}
	Let $(\rho,V)$ be a representation of a finite group $G$, with decomposition $\oplus_{i} k_{i} \rho_{i}$ into irreducible representations. Then $\langle \chi_{\rho} \, | \, \chi_{i} \rangle = k_{i}$.
\end{corollary}
\begin{proof}
	We split the character $\chi_{\rho}$ into a sum of its irreducible characters via Lemma \ref{chpt3:lem:splitchar}. Taking the inner product, all non $\chi_{i}$ terms disappear by the orthogonality of irreducible characters. Note that the character $\chi_{\rho}$ does not depend on our decomposition, thus the constants $k_{i}$ must be unique.
\end{proof}

Finally we may prove the decomposition of the regular module stated in Definition \ref{chpt3:def:regularmod}. As a consequence we find that $\sum_{\rho \in \rep(G)} d_{\rho}^{2} = |G|$.

\begin{lemma}
	\label{chpt3:lem:charreg}
	Let $(\rho,\mathbb{C}[G])$ be the regular representation for a finite group $G$. Let $\chi_{i}$ be an irreducible character with dimension $d_{i}$. Then $\langle \chi_{\rho} \,|\, \chi_{i} \rangle = d_{i}$, and $ |G| = \chi_{\rho}(\gpid) = \sum_{i} d_{i} \,\chi_{i}(\gpid) = \sum_{\rho \in \rep(G)} d_{\rho}^{2}$.
\end{lemma}
\begin{proof}
	By definition of the regular module for every $g,h \in G$ we have $\rho(g)(h) = gh$,  therefore, $\rho(g)(h) \neq h$ unless $g =\gpid$. This implies that if $g\neq \gpid$ then all diagonals of $\rho(g)$ must be $0$, and thus $\chi_{\rho}(g)=0$.
	Computing the inner product gives us,
	\[\langle \chi_{\rho} \,|\, \chi_{i} \rangle = \frac{1}{|G|}\sum_{g\in G} \chi_{i}(g) \, \overline{\chi_{\rho}(g)}   = \frac{1}{|G|} \chi_{i}(\gpid) \chi_{\rho}(\gpid) = d_{i}  \qedhere.\]
\end{proof}

\subsection{Discrete Fourier Transforms}
\label{chpt3:sec:fourier}
The discrete Fourier 
transform links a $\mathbb{C}$-valued function $P$ on a group $G$ and a representation $\rho$ of $G$ into a 
single algebraic object. Understanding Fourier transforms and how they 
interact with convolutions is a key step in linking the worlds of probability and algebra. Our goal in this 
section is to demonstrate the key properties of the transform, which are used in the next section to bound total 
variation distance of a random walk. The results below and their proofs come from \cite[Chapter 2]{Diaconis1988}

\begin{defn}
	Let $(\rho,V)$ be a representation of $G$, and $P$ a function on $G$. The 
	\emph{Fourier transform} of $P$ at $\rho$ is:
	\begin{eqnarray}
	\label{chpt3:eqn:convolution}
	\hat{P}(\rho) = \sum_{g \in G} P(g)\rho(g).
	\end{eqnarray}
	Note that $\hat{P}(\rho)$ may be viewed as a mapping from $V\to V$.
\end{defn}

It may seem opaque at the moment why we want to consider taking Fourier 
transforms of functions $P$. It turns out that Fourier transforms have 
properties which help simplify the analysis of the convolution $P^{t}$. One of the most immediate results is that the discrete Fourier transform allows us to split up convolutions.

\begin{lemma}
	Let $P,Q$ be functions on  a finite group $G$, and $\rho$ a representation. 
	Then
	\begin{eqnarray}
	\widehat{P\star Q}(\rho) = \hat{P}(\rho)\hat{Q}(\rho).
	\end{eqnarray}
\end{lemma}

\begin{proof}
	Proceeding from left to right we have:
	\begin{eqnarray*}
		\widehat{P\star Q}(\rho) & = & \sum_{g\in G} (P\star Q)(g) \rho(g) \\
		& = & \sum_{g \in G}\sum_{h\in G} P(gh^{-1})Q(h) \rho(g) \\
		& = & \sum_{h \in G} P(gh^{-1}) \rho(gh^{-1}) \sum_{g \in G} Q(h) 
		\rho(h) = \hat{P}(\rho)\hat{Q}(\rho). \qedhere
	\end{eqnarray*}
\end{proof}

This simple lemma is key to the usefulness of the discrete Fourier transform, applying it repeatedly we can see that $\widehat{P^{t}} = 
\hat{P}^{t}$. 
The second property of the Fourier transform which we exploit is the Fourier Inversion Theorem, which allows us to 
recover $P$ from the information provided by $\hat{P}(\rho)$ on all the irreducible representations of $G$.

\begin{lemma}[Fourier Inversion Theorem - See Chapter 2C \cite{Diaconis1988}]
	Let $P$ be a function on  a finite group $G$. Then
	\begin{eqnarray}
	\label{chpt3:eqn:fourinv}
	P(g) = \frac{1}{|G|}\sum_{\rho \in \rep(G)} d_{\rho} 
	\Tr\left(\rho(g^{-1})\hat{P}(\rho)\right).
	\end{eqnarray}
\end{lemma}

The Fourier inversion Theorem is important because it allows us to go back and 
forth between the $\mathbb{C}$-valued function $P$ and the Fourier transform $\hat{P}$.  However, to be able to use 
this technique effectively we also need the details of the irreducible 
representations of our group $G$. Using discrete Fourier transforms we may reduce the study of $P^{t}$ to $P$. This connection was one of the crucial insights which permitted the study of random walks on groups via their representation theory.
Let us now state some simple results following from the Fourier inversion 
Theorem.

\begin{corollary}[Plancherel Theorem]
	Let $P, Q$ be functions on a finite group $G$. Then
	\begin{eqnarray}
	\label{chpt3:eqn:plancherel}
	\sum_{g \in G} P(g)Q(g^{-1}) = \frac{1}{|G|}\sum_{\rho \in \rep(G)} 
	d_{\rho} 
	\textnormal{Tr}\left(\hat{P}(\rho) \hat{Q}(\rho)\right)   .
	\end{eqnarray}
\end{corollary}

\begin{proof}
	Both sides of \eqref{chpt3:eqn:plancherel} are linear in $Q$, so taking 
	$Q(h) = \delta_{g,h}$ we just have to show
	\[P(g^{-1}) = \frac{1}{|G|}\sum_{\rho \in \rep(G)} d_{\rho} 
	\textnormal{Tr}\left(\rho(g)\hat{P}(\rho)\right).\]
	This is just an application of the Fourier Inversion Theorem. Further 
	Details may be found in \cite[Chapter 2C]{Diaconis1988}
\end{proof}

\begin{corollary}
	\label{chpt3:cor:plan2}
	Let $P$ be a function on a finite group $G$. Then
	\begin{eqnarray}
	\label{chpt3:eqn:plancherel2}
	\sum_{g \in G} P(g)P(g) = \frac{1}{|G|}\sum_{\rho \in \rep(G)} d_{\rho} 
	\textnormal{Tr}\left(\hat{P}(\rho)\hat{P}(\rho)^{\star}\right)   .
	\end{eqnarray}
\end{corollary}

\begin{proof}
	Apply Plancherel's Theorem with $Q$ defined as $Q(g) = P(g^{-1})$. 
\end{proof}

\begin{corollary}
	Let $P$ be a probability distribution on  a finite group $G$, and $\pi$ the uniform 
	distribution on $G$. Then
	\begin{eqnarray}
	\label{chpt3:eqn:plancherel3}
	\sum_{g \in G} P(g)\pi(g) = \frac{1}{|G|}.
	\end{eqnarray}
\end{corollary}

\begin{proof}
	Following from Plancherel's Theorem we have
	\begin{eqnarray*}
		\sum_{g \in G} P(g)\pi(g) & = &\frac{1}{|G|}\sum_{\rho \in \rep(G)} 
		d_{\rho_{i}} 
		\textnormal{Tr}\left(\hat{P}(\rho)\hat{\pi}(\rho)\right)\\
		& = &\frac{1}{|G|}\sum_{\rho \in \rep(G)} d_{\rho} 
		\textnormal{Tr}\left(\hat{\pi}(\rho)\right) = \frac{1}{|G|}\sum_{\rho \in 
			\rep(G)} 
		d_{\rho} \frac{d_{\rho}}{|G|} = \frac{1}{|G|}. 
	\end{eqnarray*}	
	The second equality above comes from merging the Fourier transforms into 
	one and applying Lemma \ref{chpt2:lem:uniform}.
\end{proof}

The trace of Fourier transforms present in the Fourier inversion and Plancherel's Theorem turn the representation present in $\hat{P}(\rho),\hat{Q}(\rho)$ into the associated character 
$\chi_{\rho}$. We have seen that characters are constant on 
conjugacy classes of $G$, so if also we impose this condition on $P$ we may simplify the trace of $\hat{P}(\rho)$. 

\begin{lemma}
	\label{chpt3:lem:trace}
	Let $C_{i}$ denote the conjugacy classes of a finite group $G$, with 
	representatives 
	$g_{i}$. Let $P$ be a probability on $G$ which is constant on conjugacy classes of $G$. 
	Then we have $\hat{P}(\rho) = \beta \, \I$ with
	\begin{eqnarray}
	\beta = \sum_{i} |C_{i}| P(g_{i}) 	\frac{\chi_{\rho}(g_{i})}{d_{\rho}}
	\end{eqnarray} 
\end{lemma}

\begin{proof}
	It is enough to show that 	$\hat{P}:V \to V$ is a module homomorphism. This follows from the following calculation: for any $g \in G$ we ahve
	$\rho(g)\hat{P}(\rho)\rho(g^{-1}) = \sum_{h} P(h) \rho(ghg^{-1})  = \sum_{h} P(h) \rho(h) = \hat{P}(\rho)$. Therefore, applying Schur's Lemma gives $\hat{P}(\rho) = 
	\beta \,\I$, and taking traces gives the required value of $\beta$.

\end{proof}

\subsection{The Upper Bound Lemma}
\label{chpt3:section:upperboundlemma}

We are now in a position to state one of the most important results on the 
mixing times of random walks on groups. The upper bound lemma links the worlds 
of probability and algebra in a surprising way, stating that the total variation distance of a random walk on a group $G$ may be upper bounded via Fourier analysis on
the driving probability of the walk. If the 
Fourier transforms on irreducible representation are well understood this reduces the issue of analysing 
$\lVert P^{t}-\pi \rVert_{\textnormal{\tiny TV}}$ to bounding a summation of numbers in $\mathbb{C}$ (often just $\mathbb{R}$). The results and proof presented below are taken from \cite[Chapter 3B]{Diaconis1988}

\begin{lemma}[Upper Bound Lemma \cite{diaconis1981generating}]
	\label{chpt3:lem:upperboundlemma}
	Let $P$ be a probability on a finite group $G$, $\pi$ the uniform 
	distribution on $G$. Then
	\begin{eqnarray}
	\label{chpt3:eqn:upperbound}
	\lVert P - \pi \rVert_{\textnormal{\tiny TV}}^{2} \leq 
	\frac{1}{4}\sum_{\substack{\rho \in \rep(G) \\ \rho \neq \Triv}} 
	d_{\rho} \textnormal{Tr}\left( 
	\hat{P}(\rho)\overline{\hat{P}(\rho)}\right)
	\end{eqnarray}
	
\end{lemma}

\begin{proof}
	Following from the definition of total variation distance we find:
	\begin{eqnarray}
	4\lVert P - \pi \rVert_{\textnormal{\tiny TV}}^{2} & = &\left(\sum_{g\in G} 
	|P(g) -\pi(g)|\right)^{2} \nonumber\\
	& \leq & |G|\sum_{g \in G} |P(g) -\pi(g)|^{2} \label{chpt3:eqn:proofupcs}\\
	& = & |G|\sum_{g \in G} P(g)P(g) - 2 |G|\sum_{g \in G} P(g)\pi(g) 
	+|G|\sum_{g \in G} \pi(g)\pi(g) \nonumber \\
	& = & \sum_{\substack{\rho \in \rep(G) \\ \rho \neq \Triv}} 
	d_{\rho} \textnormal{Tr}\left( 	\hat{P}(\rho)\overline{\hat{P}(\rho)}\right)
	\end{eqnarray}
	
	Step \eqref{chpt3:eqn:proofupcs} is due to the Cauchty-Schwarz inequality, 
	and the final equality is Plancherel's Theorem and its Corollaries.
\end{proof}

\begin{lemma}[Lower Bound Lemma \cite{Diaconis1988}]
	\label{chpt2:lem:upperboundlemma}
	Let $P$ be a probability on a finite group $G$, $\pi$ the uniform distribution on $G$.
	Then
	\begin{eqnarray}
	\label{chpt3:eqn:lowerbound} 
	\frac{1}{4|G|}\sum_{\substack{\rho \in \rep(G) \\ \rho \neq \Triv}} 
	d_{\rho} \textnormal{Tr}\left( 
	\hat{P}(\rho)\overline{\hat{P}(\rho)}\right) \leq \lVert P - \pi \rVert_{\textnormal{\tiny TV}}^{2}
	\end{eqnarray}
\end{lemma}

\begin{proof}
	Following from the definition of total variation distance we find:
	\begin{eqnarray*}
		4\lVert P - \pi \rVert_{\textnormal{\tiny TV}}^{2} & = &\left(\sum_{g\in G} 
		|P(g) -\pi(g)|\right)^{2} 
		\geq  \sum_{g \in G} |P(g) -\pi(g)|^{2} 
		=  \frac{1}{|G|}\sum_{\substack{\rho \in \rep(G) \\ \rho \neq \Triv}} 
		d_{\rho} \textnormal{Tr}\left( 
		\hat{P}(\rho)\overline{\hat{P}(\rho)}\right)
	\end{eqnarray*}
	The first inequality is a general one, 
	and the final equality is Plancherel's Theorem and its Corollaries.
\end{proof}

\begin{corollary}
	\label{chpt3:cor:uplowbound}
	Let $P$ be a probability on a finite group $G$, and $\pi$ the uniform 
	distribution on $G$. Then
	\begin{eqnarray}
	\label{chpt3:eqn:lowerbound}
	\frac{1}{4|G|}\sum_{\substack{\rho \in \rep(G) \\ \rho \neq \Triv}} 
	d_{\rho} \textnormal{Tr}\left( 
	\hat{P}(\rho)^{t}\overline{\hat{P}(\rho)^{t}}\right) \leq \lVert 
	P^{t} - \pi \rVert_{\textnormal{\tiny TV}}^{2} \leq 
	\frac{1}{4}\sum_{\substack{\rho \in \rep(G) \\ \rho \neq \Triv}} 
	d_{\rho} \textnormal{Tr}\left( 
	\hat{P}(\rho)^{t}\overline{\hat{P}(\rho)^{t}}\right)
	\end{eqnarray}
\end{corollary}

\paragraph{}
The first use of the Upper Bound Lemma appeared in the seminal paper of 
Diaconis and Shahshahani, where it was used to analyse the mixing time of the 
random transposition shuffle \cite[Lemma 14]{diaconis1981generating}. We review this 
classic argument in detail in Chapter \ref{chpt4:chpt} to demonstrate the power of the upper 
bound lemma. 
Since its first appearance the 
upper bound 
lemma has seen frequent use in proving cutoff for many random walks on groups. 
It 
gives an explicit link between the areas of mixing time 
of Markov chains and representation theory, reducing the study of $\lVert 
P^{t} - \pi \rVert_{\textnormal{\tiny TV}}$ to understanding  
$\hat{P}(\rho)$. The use of the Cauchy-Schwarz inequality in the Upper Bound Lemma 
may seem like a na\"{i}ve way to bound this complicated summation, however we shall see that it frequently gives sharp bounds of the correct mixing time. The upper and lower bounds presented  in Corollary 
\ref{chpt3:cor:uplowbound} are of different orders, this means that we can not get the exact mixing time for a random walk just using representation theory.

In this thesis are main results will be focused on random walks on the 
symmetric group $S_{n}$. To get a better understanding of the techniques 
introduced in this chapter we dedicated the next chapter to exploring the symmetric group fully.  We concentrate the tools previously seen to specific examples which have guided the analysis of the one-sided transposition shuffle, these are the random transposition shuffle, top-to-random shuffle, 
and random-to-random shuffle.

	\chapter{The Symmetric Group and its Modules}
\label{chpt4:chpt}

\paragraph{}
In this chapter we review all the information needed about the symmetric group in order to analyse the mixing times of random walks on $S_{n}$. The first section is dedicated to the  the structure of the group. 	The second section presents a detailed construction of the simple modules of the symmetric group, by way of permutation modules. We then see how to view these modules as vector spaces spanned by words of length $n$. We end this chapter with detailed examples of shuffles which have guided the analysis of the one-sided transposition shuffle. These include cutoff results for the random transposition shuffle and the top-to-random shuffle, as well as an exploration into lifting eigenvectors for the random-to-random shuffle. A detailed account of the symmetric group and its representations may be found in \emph{The Symmetric Group} by Sagan \cite{sagan2013symmetric} or \emph{The Representation Theory of the Symmetric Groups} by James \cite{james1978representation}.

\section{The Symmetric Group}
\label{chpt4:sec:symgp}

\paragraph{}
The symmetric group on $n$ objects, denoted $S_{n}$, is defined as the 
group of all bijections $\sigma:[n]\to[n]$, where $[n] = \{1,\ldots, n\}$. There are exactly $n!$ of 
these  bijections. Given two bijections 
$\sigma,\eta \in S_{n}$ define the product $\sigma \eta$ as the composition of functions from right to left.
Another common way to think about the symmetric group 
$S_{n}$ is the arrangements of a deck of $n$ cards labelled, $1,\ldots, n$, from bottom to top. If we have $n$ 
positions in a deck of cards labelled $1$ to $n$, and $n$ cards labelled 
$1$ to $n$ then we may 
view the bijection $\sigma$ as telling us what position each card is in, i.e., card $i$ is in position $\sigma(i)$. The identity element $e$ has every card in 
its labelled position, i.e., for all $i \in [n]$ we have $e(i) = i$. 
Viewing the symmetric group as a deck of cards is often 
useful for random walks on $S_{n}$ because it allows us to formulate our random walks in the expressive terms of shuffling a deck of cards.

We express elements of the symmetric group using \emph{cycle notation}, for example,
\[\sigma = (1\,2\,5)(4\,6)(3) = (1\,2\,5)(4\,6) .\]
We read the cycle $(125)$ as $1$ maps to $2$ (i.e. $\sigma(1) =2$), $2$ maps to $5$, and $5$ maps to $1$. The \emph{cycle structure} of a permutation is the tuple of the lengths of its cycles arranged in non-increasing order, example the cycle structure of $\sigma$ defined above is $(3,2,1)$.  Note that the sum of cycle lengths for $\sigma \in S_{n}$ is always $n$, therefore the cycle structure for any permutation is a \emph{partition} of $n$.

\begin{defn}
	A \emph{partition} of $n$ is a tuple of positive integers $\lambda = (\lambda_{1}, \dots, 
	\lambda_{r})$ such that, $\sum_{i=1}^{r} \lambda_{i} =n$, and $\lambda_{1} \geq \dots 
	\geq \lambda_{r}$, we denote this by $\lambda \vdash n$. We call $n$ the \emph{size} of the partition and $r$ the \emph{length} of the partition, denoted $|\lambda|$ and $l(\lambda)$ respectively. If a partition contains repeated digits we may write them as a power for 
	brevity, for example we denote the 	partition $(1,\dots,1) = (1^{n})$, similarly 	$(3,3,2,2,2,1) = (3^{2},2^{3}, 1)$. 
\end{defn}

Cycle notation can help us identify the conjugacy classes of the 
symmetric group. To know if two elements are conjugate in $S_{n}$ we need only compare their cycle structures. 
\begin{lemma}[See Chapter 1 \cite{sagan2013symmetric}]
	\label{chpt4:lem:conjclasses}
	Two elements in $S_{n}$ are conjugate if and only if they have the same 
	cycle structure. Hence, the conjugacy classes of $S_{n}$ are labelled by partitions of $n$.
\end{lemma}
The identity has cycle structure $(1^{n})$, as every element belongs to its own trivial cycle. One particularly important conjugacy class of $S_{n}$ is the class of 
transpositions, which is formed from all two cycles $(i \,j)$ for $i < j$, these elements have cycle type $(2,1^{n-2})$. 
We call the transposition $(i\, j)$ an \emph{adjacent transposition} if $j = i+1$. The conjugacy class of transpositions generates the whole symmetric group, in fact we only need the $n-1$ adjacent transpositions for this, $\langle  (i, i+1) \, | \, 1\leq i \leq n-1 \rangle = S_{n}$. 	We may 
decompose any element in $S_{n}$ into a product of transpositions, and in any decomposition the number of transpositions required remains constant modulo $2$. This leads us to make the following definition.
\begin{defn}
	Let $\sigma \in S_{n}$, and decompose our permutation as $\sigma = \tau_{k} \, \tau_{k-1}\,  \dots \, 
	\tau_{1}$, where $\tau_{i}$ is a transposition. Define the \emph{sign} 
	function for $S_{n}$, denoted $\sign:S_{n} \to \{1,-1\}$, 
	as follows
	\[\sign(\sigma) = (-1)^{k}.\]
	This is a well defined function, i.e. independent of the decomposition 
	of $\sigma$ into transpositions. An element of $\sigma \in S_{n}$ is called \emph{odd} if $\sign(\sigma) =-1$, and called \emph{even} if $\sign(\sigma) =1$. The sign function is also multiplicative, i.e., for $\sigma, \eta \in S_{n}$ we have	$\sign(\sigma \eta) = \sign(\sigma) \sign(\eta)$. Therefore, the set of even permutations defines a subgroup of $S_{n}$ called the \emph{alternating group} and denoted $A_{n}$.
\end{defn}

Random walks supported on the conjugacy class of transpositions have been well studied since the random transposition shuffle, variants include the semi-random transposition shuffles \cite{Mossel2004}, adjacent transposition shuffle \cite{lacoin2016mixing} and biased transposition shuffle\cite{BernsteinBaised}.
In Chapter \ref{chpt5:chpt} we present a novel modification to general transposition shuffles called the one-sided transposition shuffle where the probability of applying transposition $(i\,j)$ depends only on the position $j$.	
We now move on to describe the module structure of the symmetric group.

\section{The Structure of Modules for The Symmetric Group}
\label{chpt4:sec:modules}
Denote the group algebra of the symmetric group as $\mathfrak{S}_{n} := \mathbb{C}[S_{n}]$. All the 
simple modules of $\mathfrak{S}_{n}$ are indexed  by partitions of $ \lambda \vdash n$, we denote the simple module indexed by $\lambda$ as $S^{\lambda}$, these are also called \emph{Specht} modules.  Understanding the simple modules of $\mathfrak{S}_{n}$ will be important to our analysis of the one-sided transposition shuffle presented in 
Chapter \ref{chpt5:chpt}. 
Before we construct modules of the symmetric group, we need to recall some facts about partitions, Young diagrams and Young tableaux.

\subsection{Young Diagrams}
\label{chpt4:subsec:youngdigrams}

Every partition $\lambda$ has an 
associated 
\emph{Young diagram}, 
made by forming a left adjusted stack of boxes with rows labelled downwards 
and 
row $i$ having 
$\lambda_{i}$ boxes.  We often blur the distinction between a partition and 
its Young diagram, e.g. $(3,2) = 	
\ytableausetup{mathmode,baseline,aligntableaux=center,boxsize=0.6em}
\ydiagram{3,2}$. 
We may refer to the boxes of a diagram $\lambda$ by using 
coordinates $(i,j)$ to mean the box in the $i^{th}$ row and $j^{th}$ 
column. Define the \emph{diagonal index} of the box $(i,j)$ in partition 
$\lambda$ to be the value $j-i$. Define the \emph{diagonal sum} of $\lambda$ to be $\D(\lambda) =\sum_{(i,j) \in \lambda} (j-i)$, this is the sum over all the diagonal indexes of $\lambda$.

Given a partition $\lambda\vdash n$, we may form the 
\emph{transpose} of $\lambda$, denoted $\lambda^{\prime}$, by swapping rows 
and	columns in the Young diagram, e.g. $(3,2)^{\prime} = (2,2,1)$. We have 
$\lambda \vdash n$ if and only if 
$\lambda^{\prime} \vdash n$. Define a partial order on partitions of $n$ called the \emph{dominance order}: 
in terms of Young diagrams, for two partitions $\mu,\lambda \vdash n$, we say $\lambda$ \emph{dominates} $\mu$ 
if we can form $\mu$ by moving boxes of 
$\lambda$ down and to the left, we denote this by $\lambda \trianglerighteq \mu$. Equivalently $\lambda \trianglerighteq \mu$  if and only if $\sum_{i=1}^{j} \lambda_{i} \geq \sum_{i=1}^{j} \mu_{i}$ for all choices of $j$.
Furthermore, we have $\lambda \trianglerighteq \mu$ if and only if $\mu'  
\trianglerighteq \lambda'$, see \cite[Lemma 1.4.11]{james_1984}. 

\begin{example}
	Let $n=8$ and consider the partitions $(3,2,2,1)$ and $(2^{4})$. We may see that $(3,2,2,1) \trianglerighteq (2^{4})$, and  $(2^{4})^{\prime} = (4,4) \trianglerighteq(4,3,1)= (3,2,1,1)^{\prime}$. Now consider partitions $(5,1^{3})$ and $(4^{2})$, neither partition here dominates the other. Hence, the dominance order on partitions is not necessarily a total ordering.

\end{example}

Given two partitions $\mu,\lambda$ of different sizes, we 
write $\mu \subseteq \lambda$ if $\mu$ is fully contained in $\lambda$ when 
we 
align the Young diagrams of $\mu$ and $\lambda$ at the top left corners;
equivalently, if we write $\lambda = (\lambda_1,\dots,\lambda_r)$ and 
$\mu=(\mu_1,\dots,\mu_s)$, 
this means that $s\leq r$ and $\mu_i\leq \lambda_i$ for each $1\leq i \leq 
s$. For example $(3,2) \subseteq (4,3)$ -- this is simpler to see from the 
corresponding Young diagrams, 
$\ytableausetup{mathmode,baseline,aligntableaux=center,boxsize=0.6em}
\ydiagram{3,2} \subseteq \ydiagram{4,3}$.
If we have two partitions $\mu,\lambda$ such that $\mu \subseteq \lambda$ we 
may define the \emph{skew diagram} $\lambda/\mu$ as the diagram containing all  
boxes which are in $\lambda$ but not in $\mu$. For example the skew diagram 
of $ (4,3)/(3,2)$ is 
$\ytableausetup{mathmode,baseline,aligntableaux=center,boxsize=0.6em} 
\ydiagram{1+1,0+1}$. A skew diagram $\lambda / \mu$ is called a horizontal 
strip if it has at most one box per column, e.g. $(4,3)/(3,2)$ is a 
horizontal strip but $(4,3)/(2,2) = \ytableausetup{mathmode,baseline,aligntableaux=center,boxsize=0.6em} 
\ydiagram{2,1}$ is not.
\paragraph{}
Given a partition of $n$ we may turn it into a partition of $n+1$ by adding a box to the corresponding Young diagram. In order to formalise this we view partitions as $n$-tuples, filling the end with zeros if necessary.	Given an $n$-tuple $\lambda=(\lambda_1,\ldots,\lambda_n)$ of non-negative 		integers summing to $n$ and an element $i \in [n+1]$, we form an $(n+1)$-tuple denoted $\lambda+e_i$ by first adding a zero to the end 	of 	$\lambda$ and then adding $1$ to this $(n+1)$-tuple in position $i$. Then $\lambda+e_i$ is an $(n+1)$-tuple of non-negative integers summing to $n+1$, e.g. $(1,1,1) + e_{3} = (1,1,1,0) + (0,0,1,0) = (1,1,2,0)$.
In terms of Young diagrams $\lambda +e_{i}$ represents adding a box to $\lambda$ on row $i$.

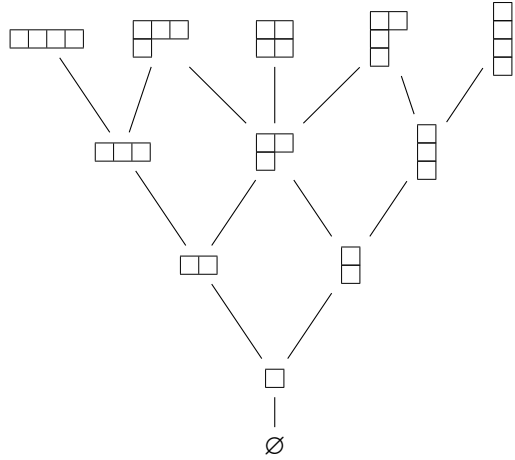
\begin{figure}[H]
	\[
	\begin{tikzpicture}
	\node (realstart) at (0,-0.9) {$\emptyset$};
	
	\node (start) at (0,0) {\tiny\yng(1)};
	
	\node (2) at (-1,1.5) {\tiny\yng(2)};
	\node (11) at (1,1.5) {\tiny\yng(1,1)};
	
	\node (3) at (-2,3) {\tiny\yng(3)};
	\node (21) at (0,3) {\tiny\yng(2,1)};
	\node (111)  at (2,3)  {\tiny\yng(1,1,1)};
	
	\node (4) at (-3,4.5){\tiny\yng(4)};
	\node (31) at(-1.5,4.5) {\tiny\yng(3,1)};
	\node (22) at(0,4.5) {\tiny\yng(2,2)};
	\node (211) at (1.5,4.5) {\tiny\yng(2,1,1)};
	\node(1111) at (3,4.5) {\tiny\yng(1,1,1,1)};

	\draw[-] (realstart) -- (start);
	\draw[-] (start) -- (2);
	\draw[-] (start) -- (11);
	
	\draw[-] (2) -- (21);
	\draw[-] (11) -- (111);
	\draw[-] (11) -- (21);
	\draw[-] (2) -- (3);
	
	\draw[-] (111) -- (211);
	\draw[-] (21) -- (211);
	\draw[-] (111) -- (1111);
	\draw[-] (21) -- (22);
	\draw[-] (21) -- (31);
	\draw[-] (3) -- (31);
	\draw[-] (3) -- (4);
	
	\end{tikzpicture}
	\]
	\caption[Young's lattice for partitions of size $n \in \{0,1,2,3,4\}$]{Young's lattice for partitions of size $n \in \{0,1,2,3,4\}$.}
	\label{chpt4:fig:lattice}
\end{figure}

Note that if $\mu \subseteq \lambda$ we may add boxes to $\mu$ to form $\lambda$, the boxes we have to add are exactly those contained in the skew diagram $\lambda /\mu$, e.g. for $(3,2) \subseteq (4,3)$ we may see that $(3,2) +e_{1} +e_{2} =(4,3)$. 
If we restrict our attention to choices of box $e_{i}$ that result in another partition we uncover a structure on Young diagrams called \emph{Young's lattice}, seen in Figure \ref{chpt4:fig:lattice}.
We allow $n=0$ as a special case with empty partition $(0)$ and corresponding Young diagram $\emptyset$.

Each path in Figure \ref{chpt4:fig:lattice} represents 
the placement or removal of a box to form a new partition. Young's lattice is an important structure as it allows us to link the partitions of $n$ to those of $n+1$. This structure play an important role in the   linking the modules of $S_{n}$ and $S_{n+1}$ in Section \ref{chpt4:subsec:branching}.

\subsection{Young Tableaux}
\label{chpt4:subsec:youngtableau}
Given a partition $\lambda\vdash n$, we may form a \emph{Young tableau} (alternatively a \emph{$\lambda$-tableau}) $T$ 
by putting numbers $1,\ldots, n$ into the boxes of (the Young diagram of) $\lambda$, such that each number only appears once. The set of Young tableaux of shape $\lambda$ is denoted $\YT(\lambda)$.	We says a value $m$ \emph{occurs} in the Young tableau $T$ if $m$ is present in a box of $T$. A \emph{standard Young tableau} $T$ is a Young tableau where the values in the boxes of $T$ 
are increasing across rows and down columns. 
The set of standard Young tableaux of shape $\lambda$  is 
denoted by $\SYT(\lambda)$. The size of the set $\SYT(\lambda)$ is
called the \emph{dimension} of $\lambda$, denoted $d_{\lambda}$. 
For a tableau $T$, form the transpose of $T$, denoted $T^{\prime}$, by 
swapping rows and columns while preserving the value in each box. A tableau $T$ has $T \in 
\SYT(\lambda)$ if and only if $T^{\prime} \in \SYT(\lambda^{\prime})$.
If a tableau $T$ has a box in position $(i,j)$, we let $T(i,j)$ denote the value in that box $(i,j)$; otherwise $T(i,j)$ is undefined. The partition $(3,2)$ has $5$ standard 
Young tableaux therefore $d_{(3,2)} = 5$, these are given below:
\[\young(123,45) \hspace{2cm} \young(124,35) \hspace{2cm} \young(125,34)\]
\[\young(134,25) \hspace{2cm} \young(135,24)\]
Each standard Young tableau corresponds to one path up Young's lattice, starting at $\emptyset$ and ending at $\lambda$. To form this correspondence take a standard Young tableau $T\in \SYT(\lambda)$ and form a path up Young's lattice  from $\emptyset$ to $\lambda$ by adding the boxes to $\emptyset$ in the order given by the entries in $T$. 
A standard Young tableau $T \in \SYT(\lambda)$ is called a 
\emph{desarrangement} tableau if there is no box $(1,2) \in \lambda$ and 
$|\lambda|$ is even, or the value $T(1,2)$ is odd. We denote the number of 
desarrangement tableaux of a shape $\lambda$ as $d^{\lambda}$. For example, in the collection of standard Young tableaux of shape $(3,2)$ displayed above, those in the bottom row are  desarrangement tableaux.

Finally we end with a result that links the tableaux of different shapes to the dominance order of partitions.

\begin{lemma}[Lemma 2.2.4 \cite{sagan2013symmetric}]
	\label{chpt4:lem:tabdom}
	Let $\lambda,\mu \vdash n$. Take $T_{\lambda}$ to be a $\lambda$-tableau and $T_{\mu}$ to be a $\mu$-tableau.
	Suppose that for each index $i$, the elements in row $i$ of $T_{\mu}$ are all in different columns of $T_{\lambda}$. Then $\lambda \trianglerighteq \mu$.
\end{lemma}
\begin{proof}
	
	By our hypothesis the elements in row $1$ of $T_{\mu}$ are in different columns of $T_{\lambda}$. We may rearrange the columns of $T_{\lambda}$ so that these elements all appear in row $1$ of $T_{\lambda}$. Continuing this we may sort the columns of $T_{\lambda}$ so that the elements of rows $1,2,\ldots, j$ of $T_{\mu}$ all occur in the first $j$ rows of $T_{\lambda}$. Then,
	\begin{eqnarray*}
		\sum_{i=1}^{j} \lambda_{i} \; = &\textnormal{ number of elements in the first $j$ rows of $T_{\lambda}$}& \\
		\geq  & \textnormal{ number of elements in the first $j$ rows of $T_{\mu}$}& = \; \sum_{i=1}^{j} \mu_{i}. \qedhere
	\end{eqnarray*}
\end{proof}

\subsection{Permutation Modules}
On our way to defining the simple modules $S^{\lambda}$ for 
$\mathfrak{S}_{n}$ we first need to define the \emph{permutation modules}, denoted
$M^{\lambda}$, which contain the simple modules as 
submodules. To define the permutation module $M^{\lambda}$ we use the natural action of the symmetric group on the Young tableaux of shape $\lambda$.

\begin{defn}
	\label{chpt4:def:tabaction}
	Let $T$ be a Young tableau of $n$, and $\sigma \in S_{n}$. Define an 
	action of $\sigma$ on $T$, by applying $\sigma$ to the values in $T$ 
	box wise. Thus $\sigma \,T$ defines a new Young tableau of the same shape 
	as $T$. For example, given $T = 
	\ytableausetup{mathmode,baseline,aligntableaux=center,boxsize=1em}
	\begin{ytableau} 5 & 4 & 1 \\
	2 & 3 \end{ytableau}$ and $\sigma\in S_{n}$, 
	the product 
	\[\sigma \, T = 
	\ytableausetup{mathmode,baseline,aligntableaux=center,boxsize=2.2em}
	\begin{ytableau}\sigma(5) &\sigma(4) & \sigma(1) \\
	\sigma(2) & \sigma(3) \end{ytableau}.\]

\end{defn}

To form the permutation module $M^{\lambda}$ we first establish an important equivalence relation on the tableaux of shape $\lambda$. Using the action in Definition 
\ref{chpt4:def:tabaction} we may define the groups of row and column 	
permutations of a tableau $T$.
\begin{defn}
	\label{chpt4:def:rowequiv}
	Let $T$ be a tableau of size $n$. Define the subgroup of \emph{row permutations of T}, denoted $R_{T}$, as the set of elements of $S_{n}$ which preserves each row of $T$.  Define the subgroup of \emph{column permutations of T}, denoted $C_{T}$, as the set of elements of $S_{n}$ which preserves each columns of $T$. For example, if 
	$T$ is the tableau given in Definition \ref{chpt4:def:tabaction} then 
	$R_{T} \cong S_{3} \times S_{2}$ as we are 
	allowed to swap the numbers $5,4,1$ freely and $2,3$ freely but not 
	swap between them, similarly $C_{T} \cong S_{2} \times S_{2} \times S_{1}$.
\end{defn}

\begin{defn}
	Let $T_{1}, T_{2}$ be two tableaux of the same Young diagram $\lambda$. 
	We say $T_{1}$ is \emph{row equivalent} to 
	$T_{2}$, denoted $T_{1} \sim_{\textnormal{R}} T_{2}$, if there exists 
	$\sigma \in R_{T_{1}}$ such that $\sigma T_{1} = T_{2}$. Similarly we 
	say $T_{1}$ is \emph{column equivalent} to $T_{2}$, denoted $T_{1} 
	\sim_{\textnormal{C}} T_{2}$, if there exists $\sigma \in C_{T_{1}}$ 
	such that $\sigma T_{1} = T_{2}$.
\end{defn}

\begin{example}
	Let $T = \ytableausetup{mathmode,baseline,aligntableaux=center,boxsize=1em}
	\begin{ytableau} 5 & 4 & 1 \\	2 & 3 \end{ytableau}$, we know that $|R_{T}| = 12$, therefore there are $12$ row equivalent tableaux. $6$ of these are found by permuting the first row of $T$, we list them below:
	\[
	\young(145,23) \hspace{1cm} \young(154,23)
	\hspace{1cm}  
	\young(415,23) \hspace{1cm} \young(451,23) \hspace{1cm}  
	\young(514,23) \hspace{1cm}  
	\young(541,23) .\]
	The remaining $6$ row equivalent tableaux can be found by swapping $2$ and $3$ in the second row.  		
\end{example}

\begin{lemma}
	The relation $\sim_{\textnormal{R}}$ defines an 
	equivalence relation on the tableaux of shape $\lambda$. 
\end{lemma}
\begin{proof}
	This follows from $R_{T}$ forming a subgroup of $S_{n}$.
\end{proof}

\begin{defn}
	The row equivalence class of a tableau $T$ is denoted by $\{T\}$ and called a \emph{tabloid}.
\end{defn}

Our action of $S_{n}$ on tableaux extends naturally to an action 
on 
tabloids given by $\sigma \{T\} = \{\sigma T\}$. We are now in a position 
to define our permutation modules using tabloids 
of a given partition. 

\begin{defn}
	\label{chpt4:def:permmod}
	Let $\lambda \vdash n$. The \emph{permutation module} for $\mathfrak{S}_{n}$, denoted $M^{\lambda}$, is the following vector space
	\[M^{\lambda} = \langle \{T\} \, | \, T \textnormal{ is a tableau of shape } 
	\lambda \rangle\]
	with action of $\mathfrak{S}_{n}$ given by extending the action of the symmetric group on tabloids (Definition \ref{chpt4:def:tabaction}) linearly.
\end{defn}

There are $5!$ Young tableaux of shape $(3,2)$ (in fact of any shape $\lambda \vdash 5$), but the module $M^{(3,2)}$ has a basis of $10$ tabloids, because there are only two rows to split the five numbers $\{1,2,3,4,5\}$ across, giving 
${5\choose 3}$ choices for the first row, with the second row containing the leftover numbers. Note that for the permutation module we may have basis 
elements that are not in the equivalence class of a standard Young tableau. 
E.g., the module $M^{(3,2)}$ is spanned by the 
following 10 tabloids:
\[\ytableausetup{mathmode,baseline,aligntableaux=center,boxsize=1em} 
\left\lbrace \begin{ytableau}1 & 2 & 3 \\4 & 5 
\end{ytableau}  
\right\rbrace 
\hspace{1.4cm} 
\left\lbrace 
\begin{ytableau}1 & 2 & 4 \\3 & 5 
\end{ytableau}  
\right\rbrace \hspace{1.4cm}
\left\lbrace 
\begin{ytableau}1 & 2 & 5 \\3 & 4 
\end{ytableau}  
\right\rbrace \hspace{1.4cm}
\left\lbrace 
\begin{ytableau}1 & 3 & 4 \\2 & 5 
\end{ytableau}  
\right\rbrace \hspace{1.4cm}
\left\lbrace 
\begin{ytableau}1 & 3 & 5 \\2 & 4 
\end{ytableau}  
\right\rbrace\]
\[\ytableausetup{mathmode,baseline,aligntableaux=center,boxsize=1em} 
\left\lbrace \begin{ytableau}1 & 4 & 5 \\ 2& 3 
\end{ytableau}  
\right\rbrace 
\hspace{1.4cm} 
\left\lbrace 
\begin{ytableau}2 & 3 & 4 \\1 & 5 
\end{ytableau}  
\right\rbrace \hspace{1.4cm}
\left\lbrace 
\begin{ytableau}2 & 3 & 5 \\1 & 4 
\end{ytableau}  
\right\rbrace \hspace{1.4cm}
\left\lbrace 
\begin{ytableau}2 & 4 & 5 \\1 & 3 
\end{ytableau}  
\right\rbrace \hspace{1.4cm}
\left\lbrace 
\begin{ytableau}3 & 4 & 5 \\1 & 2 
\end{ytableau}  
\right\rbrace\]
Only the first row above are the tabloid classes of standard Young 
tableaux.

\subsection{Specht Modules}
\label{chpt4:subsec:spechtmod}

Every $M^{\lambda}$ contains a unique copy of the simple module 
$S^{\lambda}$. We construct the Specht module $S^{\lambda}$ as a stable subspace of $M^{\lambda}$ before going on to show it is simple. We form a basis for the Specht module $S^{\lambda}$ by taking taking linear combinations of tabloids, we call these new elements \emph{polytabloids}.

\begin{defn}
	Given a Young tableau $T$ form the element $k_{T} \in \mathfrak{S}_{n}$ as the following sum:
	\begin{eqnarray}
	k_{T} = \sum_{\sigma \in C_{T}} \textnormal{sgn}(\sigma) \, \sigma \label{chpt4:eqn:polytabsumelement}
	\end{eqnarray}
	The \emph{polytabloid} generated by tableau $T$, denoted $e_{T}$, is given by the element $e_{T} = k_{T} \, \{T\}$.
	Every polytabloid is an element of the permutation module $M^{\lambda}$. We say the polytabloid $e_{T_{1}}$ contains the tabloid $\{T_{2}\}$ if it appears with a non-zero coefficient in $e_{T_{1}}$
\end{defn}

\begin{lemma}[Lemma 2.3.2 \cite{sagan2013symmetric}]
	\label{chpt4:lem:speccyclic}
	Let $T$ be a tableau of size $n$, and $\eta$ is an element in $S_{n}$. Then 
	$\eta \, e_{T} = e_{\eta \, T}$. 
\end{lemma}	
\begin{proof}
	\begin{eqnarray*}
		\eta e_{T}  =  \eta \left(\sum_{\sigma \in C_{T}} 
		\textnormal{sgn}(\sigma)  \sigma \, \right)	\eta^{-1} \eta \{T\}
		= 
		\sum_{\sigma \in C_{T}} 
		\textnormal{sgn}(\sigma) \eta \sigma \eta^{-1} \, \{ \eta T\} = \sum_{ 
			\sigma \in C_{\eta T}} 
		\textnormal{sgn}(\sigma) \sigma\, \{ \eta T\}. \qedhere
	\end{eqnarray*}
\end{proof}

\begin{defn}
	Let $\lambda \vdash n$, the \emph{Specht module} $S^{\lambda}$ is defined to be 
	the following vector space
	\[S^{\lambda} = \langle e_{T} \, | \, T \in \YT(\lambda)\rangle \]
	with the natural action of $\mathfrak{S}_{n}$ on tabloids. In fact we may restrict this spanning set to a basis by considering only  polytabloids formed from standard Young tableaux, 
	\[S^{\lambda} = \langle e_{T} \, | \, T \in \SYT(\lambda)\rangle \]
	see \cite[Theorem 2.6.2]{sagan2013symmetric} for a proof of this. Thus, the Specht module $S^{\lambda}$ has dimension $d_{\lambda}$.  The Specht modules are also cyclic,  that is $S^{\lambda}$ is generated as an $\mathfrak{S}_{n}$-module by any one polytabloid.
\end{defn}

\begin{example}
	\label{chpt4:ex:specht}
	Let $\lambda = (3,1)$, up to row equivalence we have $4$ Young tableaux of this shape, we label them by which number appears in the second row
	\[\ytableausetup{mathmode,baseline,aligntableaux=center,boxsize=1em}
	T_{4} = \begin{ytableau}1 & 2 & 3 \\4 \end{ytableau} \hspace{1.5cm} 
	T_{3} 
	=\begin{ytableau}1 & 2 & 4 \\3 \end{ytableau}  \hspace{1.5cm} 
	T_{2} = \begin{ytableau}1 & 3 & 4 \\2 \end{ytableau} \hspace{1.5cm} 
	T_{1} = \begin{ytableau}2 & 3 & 4 \\1 \end{ytableau} . \]
	The permutation module $M^{\lambda}$ is four dimensional and spanned by the 
	four tabloids above.  We can see 
	that $T_{4},T_{3},T_{2}$ are all standard Young tableau of shape 
	$\lambda$, and so the dimension of $S^{(3,1)}$ is $3$.
	To form the polytabloids that make up $S^{\lambda}$ we need to 
	use the group of column permutations of each tableau. In this 
	example $C_{T} \cong S_{2}$ with each tableau having a single 
	transposition which may be applied to the first column. Therefore,
	\[\ytableausetup{mathmode,baseline,aligntableaux=center,boxsize=1em} 
	e_{T_{4}} = \left\lbrace \begin{ytableau}1 & 2 & 3 \\4 \end{ytableau} 
	\right\rbrace - (1 \, 4) \left\lbrace \begin{ytableau}1 & 2 & 3 \\4 
	\end{ytableau} 
	\right\rbrace =\left\lbrace \begin{ytableau}1 & 2 & 3 \\4 
	\end{ytableau} 
	\right\rbrace -\left\lbrace \begin{ytableau}4 & 2 & 3 \\1 
	\end{ytableau}  
	\right\rbrace = \{T_{4}\} - \{T_{1}\}.\]
	Performing the same calculation for the tableaux $T_{3}$ and $T_{2}$, 
	we may write a basis for $S^{\lambda}$ as follows:
	\[S^{(3,1)} = \langle e_{T_{4}}, e_{T_{3}}, e_{T_{2}} \rangle = \langle \{T_{4}\} -\{T_{1}\}, \{T_{3}\} -\{T_{1}\}, 
	\{T_{2}\} -\{T_{1}\} \rangle .\]
	
	Furthermore, the Specht module $S^{(3,1)}$ may be generated by a single polytabloid $e_{T_{4}}$, 
	\begin{eqnarray*}
		(3 \, 4) e_{T_{4}} = (3 \, 4) \{T_{4}\} - (3 \, 4) \{T_{1}\} = \{ 
		T_{3}\} - \{T_{1}\} = e_{T_{3}} \\
		(2 \, 4) e_{T_{4}} = (2 \, 4) \{T_{4}\} - (2 \, 4) \{T_{1}\} = \{ 
		T_{2}\} - \{T_{1}\} = e_{T_{2}}  
	\end{eqnarray*} 
	as Lemma  \ref{chpt4:lem:speccyclic} asserts.
\end{example}

We know the conjugacy classes of $S_{n}$ are labelled by partitions of $n$, and we have found exactly one Specht module for each partition. To be sure $\{S^{\lambda} \, | \, \lambda \vdash n\}$ forms a complete set of simple modules we have to check that they are in fact simple, and all distinct from one another. To do this we study how $k_{T}$ acts on tableaux of different shapes. The following arguments were presented by James in  \cite{JamesSym}, and further details of them can be found in Sections 2.4 and 2.5 in \cite{sagan2013symmetric}

\begin{lemma}
	\label{chpt4:lem:samecol}
	Let $\lambda,\mu \vdash n$. Take $T_{\lambda}$ to be a $\lambda$-tableau and $T_{\mu}$ to be a $\mu$-tableau. Suppose there exists $a,b$ which occur in the same row of $T_{\mu}$ and same column of $T_{\lambda}$. Then $k_{T_{\lambda}} \{T_{\mu}\} =0$.
\end{lemma}
\begin{proof}
	From our hypothesis we have $(\gpid - (a \, b)) \{T_{\mu}\} =0$. Take signed coset representatives $\{\sigma_{1},\ldots, \sigma_{k}\}$ for the subgroup $\{\gpid, (a \, b)\} \subseteq C_{T_{\lambda}}$. Then $K_{T_{\lambda}}\, \{T_{\mu}\} = \left( \sum_{\sigma \in C_{T_{\lambda}}} \sign(\sigma) \, \sigma \right) \{T_{\mu}\}=\left(\sum_{i=1}^{k}\sigma_{i}(\gpid - (a \, b) ) \right) \{T_{\mu}\} =0$.
\end{proof}

Applying Lemma \ref{chpt4:lem:samecol} in combination with Lemma \ref{chpt4:lem:kdom}, we are able to use the value of $k_{T_{\lambda}}\{T_{\mu}\}$ to tell us information about $\lambda$ and $\mu$.

\begin{lemma}
	\label{chpt4:lem:kdom}
	Let $\lambda,\mu \vdash n$. Take $T_{\lambda}$ to be a $\lambda$-tableau and $T_{\mu}$ to be a $\mu$-tableau. If $k_{T_{\lambda}} \{T_{\mu}\} \neq 0$, then $\lambda \trianglerighteq \mu$. Furthermore if $\lambda = \mu$ then $k_{T_{\lambda}} \{T_{\mu}\} = \pm e_{T_{\lambda}}$.
\end{lemma}

\begin{proof}
	If $k_{T_{\lambda}} \{T_{\mu}\} \neq 0$ then for all $a,b$ in the same row of $\mu$ they are in different columns of $\lambda$. Therefore, by Lemma \ref{chpt4:lem:tabdom} we know $\lambda \trianglerighteq \mu$. If $\lambda = \mu$, then we must have our tabloids being column equivalent otherwise we violate $k_{T_{\lambda}}\{T_{\mu}\}\neq 0$ (see \cite[Corollary 2.4.2]{sagan2013symmetric}). Hence, there must exist $\eta \in  C_{T_{\lambda}}$ such that $\{T_{\mu}\} = \eta \{T_{\lambda}\}$, therefore
	\[k_{T_{\lambda}} \{T_{\mu}\} = \left( \sum_{\sigma \in C_{T}} \textnormal{sgn}(\sigma) \, \sigma \, \eta \right)  \{T_{\lambda}\} = \left( \sum_{\sigma \in C_{T}} \sign(\eta^{-1})\sign(\sigma) \, \sigma  \right) \{T_{\lambda}\} = \sign(\eta) e_{T_{\lambda}}. \qedhere\]
\end{proof}

\begin{corollary}
	\label{chpt4:cor:kmapmulti}
	Let $v \in M^{\mu}$, and take $T_{\mu}$ a $\mu$-tableau. Then $k_{T_{\mu}} \,v$ is a multiple of $e_{T_{\mu}}$.
\end{corollary}
\begin{proof}
	Write $v = \sum_{i} c_{i} \{T_{i}\}$ where $T_{i}$ are $\mu$-tableaux. By Lemma \ref{chpt4:lem:kdom} each summand $k_{T_{\mu}} \{T_{i}\}$ is either a multiple of $e_{T_{\mu}}$ or $0$.
\end{proof}	

We are now in a position to prove that the Specht module $S^{\lambda}$ is simple, and that $S^{\lambda} \cong S^{\mu}$ if and only if $\lambda = \mu$. These facts together tell us that the Specht modules form a complete set of simple modules for $\mathfrak{S}_{n}$.

\begin{thm}[Submodule Theorem]
	\label{chpt4:thm:submodthm}
	Let $V \subseteq M^{\lambda}$ be a submodule. Then $V \supseteq S^{\lambda}$ or $V \subseteq (S^{\lambda})^{0}$. Therefore, the Specht module $S^{\lambda}$ is simple.
\end{thm}

\begin{proof}
	Take $v \in V$, and a $\lambda$-tableau $T$. By Corollary \ref{chpt4:cor:kmapmulti} we know that $k_{T} v = c \, e_{T}$ for some $c\in \mathbb{C}$. We consider two cases:  $c\neq 0$ for some $T$, and $c=0$ always. Suppose that there exists $v$ and $T$ with $ k_{T} v = c \, e_{T} \neq 0$. Then we have $ c^{-1} k_{T} v = e_{T} \in V$. Hence, by Lemma \ref{chpt4:lem:speccyclic} we may generate $S^{\lambda}$ inside of $V$.
	
	Now suppose we always have $k_{T} v = 0$. Consider the inner product on $M^{\lambda}$ defined on tabloids by $\langle \{T_{1}\}, \{T_{2}\} \rangle = \delta_{\{T_{1}\},\{T_{2}\}}$, this inner product is $S_{n}$-invariant.
	We now find,
	\[\langle v, e_{T} \rangle = \langle v, \sum_{\sigma \in C_{T}} \sign(\sigma) \, \sigma \, \{T\} \rangle = \langle \sum_{\sigma \in C_{T}}\sign(\sigma)  \,\sigma^{-1} \,v, \{T\} \rangle  =\langle k_{T} \,v,  \{T\} \rangle = \langle 0 , \{T\} \rangle = 0 .\]
	A single polytabloid $e_{T}$ spans $S^{\lambda}$, therefore $v \notin S^{\lambda} \Rightarrow v \in (S^{\lambda})^{0}$.
\end{proof}

\begin{thm}
	\label{chpt4:thm:spechtmaps}
	Let $\lambda,\mu\vdash n$. Suppose we have a non-zero homomorphism $\psi:S^{\lambda} \to M^{\mu}$. Then $\lambda \trianglerighteq \mu$ and if $\lambda = \mu$ then $\psi$ is multiplication by a scalar.
\end{thm}
\begin{proof}
	Take a basis vector $e_{T} \in S^{\lambda}$ such that $\psi(e_{T}) \neq 0$. Extend the homomorphism $\psi$ to a homomorphism $\psi:M^{\lambda} \to M^{\lambda}$ by setting it to be zero on the complement of $S^{\lambda}$. Then,
	\[0 \neq \psi(e_{T}) = k_{T}\psi(\{T\}) = k_{T} \left( \sum_{i} c_{i} \{T_{i}\}\right)\]
	where the $T_{i}$ are $\mu$-tableaux, and we must have at least one $c_{i}$ being non-zero. Hence, by Lemma \ref{chpt4:lem:kdom} we have $\lambda \trianglerighteq \mu$.	If $\lambda = \mu$, then we know $\psi(e_{T}) = c \, e_{T}$ for some constant $c\in\mathbb{C}$, and for any $\sigma \in S_{n}$ we have
	\[\psi(e_{\sigma\,T}) = \psi( \sigma e_{T}) = \sigma \psi (e_{T}) = c\cdot \sigma e_{T} = c \cdot e_{\sigma \, T} \qedhere .\]
\end{proof}

\begin{corollary}
	\label{chpt4:cor:speciso}
	Let $\lambda,\mu \vdash n$. Then $S^{\lambda} \cong S^{\mu}$ if and only if $\mu = \lambda$.
\end{corollary}

\begin{proof}
	If $\mu =\lambda$ then the conclusion is immediate. Suppose $S^{\lambda} \cong S^{\mu}$, then there exist non-zero homomorphisms $\psi: S^{\lambda} \to M^{\mu}$, and $\varphi:S^{\mu} \to M^{\lambda}$. Therefore, by Theorem \ref{chpt4:thm:spechtmaps}, $\lambda \trianglerighteq \mu$ and $\mu \trianglerighteq \lambda$, which implies $\lambda = \mu$.
\end{proof}

\begin{corollary}[Theorem 2.4.6 \cite{sagan2013symmetric}]
	The Specht modules $S^{\lambda}$ for $\lambda\vdash n$ form a complete 
	set of non-isomorphic simple modules for $S_{n}$. 
\end{corollary}

\begin{proof}
	We have found the correct number of non-isomorphic simple modules.
\end{proof}

Now we have a complete set of simple modules for the symmetric groups we proceed to  decompose each permutation module into its simple parts. The following result is sometimes known as \emph{Young's Rule}.

\begin{lemma}[Young's Rule]
	\label{chpt4:lem:youngrule}
	For $\mu \vdash n$ we have,
	\[M^{\mu} \cong \bigoplus_{\lambda \,\trianglerighteq  \,\mu} 
	K_{\lambda,\mu} S^{\lambda},\]
	where $K_{\lambda,\mu}S^\mu$ denotes a direct 
	sum of $K_{\lambda,\mu}$ copies of $S^\mu$.
	The coefficients $K_{\lambda,\mu} \in \mathbb{N}$ are called \emph{Kostka 
		numbers}, and for 
	all $\lambda \vdash n$ we know $K_{\lambda,\lambda} = 1$.
\end{lemma}

\begin{proof}
	If $S^{\lambda}$ appears with a non-zero coefficient in the decomposition of $M^{\mu}$ then we clearly have a non-zero homomorphism $\psi:S^{\lambda} \to M^{\mu}$, therefore $\lambda \trianglerighteq \mu$. To establish that all $\lambda \trianglerighteq \mu$ appear in our decomposition requires information about semi-standard Young tableaux which will not feature elsewhere in this thesis, so we leave the details which may be found in Section 2.10 \cite{sagan2013symmetric}.	If $\lambda=\mu$, we know any morphism $\psi:S^{\lambda} \to M^{\lambda}$ is multiplication by a scalar hence there is only one copy in of $S^{\lambda}$ in $M^{\lambda}$.
\end{proof}

Lemma \ref{chpt4:lem:youngrule} shows a instance of Maschke's Theorem 
(Theorem \ref{chpt3:thm:decomposition}).
Notably this theorem does not tell us the value of the Kostka numbers but 
importantly it tells us that $S^{\mu}$ appears as a submodule of 
$M^{\lambda}$ if and only if $\mu \trianglerighteq \lambda$.
Consider the permutation module $M^{(3,1)}$ from Example \ref{chpt4:ex:specht}, we know that $S^{(3,1)}$ appears as one composition factor. This leaves a $1$ dimensional submodule left to find, and Young's rule tells us it must be the Specht module $S^{(4)}$, as the only partition of $4$ which dominates $(3,1)$ is $(4)$.

We have previously seen that 
the regular module for any group has a canonical decomposition into the simple modules for that group. For the symmetric group the regular module 
$\mathfrak{S}_{n}$ may be seen to be the permutation module $M^{(1^{n})}$. 

\begin{lemma}
	\label{chpt4:lem:regulardecomp}
	The permutation module $M^{(1^{n})} \cong \mathfrak{S}_{n}$ as modules. 
	Therefore, $M^{(1^{n})}$ has canonical decomposition
	\[M^{(1^{n})} \cong \bigoplus_{\lambda \vdash n} d_{\lambda} 
	S^{\lambda}  \text{ as $\mathfrak{S}_{n}$-modules}.\]
	This decomposition satisfies Lemma \ref{chpt4:lem:youngrule}, as every 
	partition of $n$ dominates $(1^{n})$.
\end{lemma}

\begin{proof}
	The module $M^{(1^{n})}$ is spanned by $n!$ tabloids. To define a linear map $\psi: 
	\mathfrak{S}_{n} \to M^{(1^{n})}$ it is enough to define it on each permutation $\sigma \in S_{n}$. To do this set,
	\[\psi(\sigma)  = 
	\ytableausetup{mathmode,baseline,aligntableaux=center,boxsize=2.2em}
	\left\lbrace \begin{ytableau} \sigma(1) \\ \dots  \\ 
	\sigma(n)\end{ytableau}
	\right\rbrace .\]
	The map $\psi$ respects the action of both modules and is an 
	isomorphism of vector spaces, therefore $\mathfrak{S}_{n} \cong 
	M^{(1^{n})}$ as $\mathfrak{S}_{n}$-modules.
\end{proof}

The decomposition present in Lemma \ref{chpt4:lem:regulardecomp} is 
important for the analysis of the random-to-random shuffle and the one-sided transposition shuffle. It allows us to focus on the Specht modules of $S_{n}$ as opposed to the larger space $\mathfrak{S}_{n}$.

\subsection{Branching Rules for Specht Modules}
\label{chpt4:subsec:branching}

Every symmetric group has a natural embedding into symmetric groups of a greater size.	
To embed $S_{m} \hookrightarrow S_{n}$ for $m<n$, we extend every permutation $\sigma \in S_{m}$ to a permutation in $S_{n}$ by choosing it to fix all elements of $[m] \setminus [n]$. The group algebras of the symmetric groups and subsequently $\mathfrak{S}_{n}$-modules inherit this recursive structure. A module $S^{\lambda}$ for $\lambda\vdash n$ may be viewed as an $\mathfrak{S}_{n-1}$ module by the action of $S_{n-1}$ inside of $S_{n}$. Conversely we may also turn $\mathfrak{S}_{n}$-modules into $\mathfrak{S}_{n+1}$-modules by a process known as induction. We now introduce the definition of restricted and induced modules.

\begin{defn}
	\label{chpt4:def:restrictandinduce}
	Let $H \subseteq G$ be groups. Suppose $W$ is a $H$-module and $V$ is a $G$-module. The \emph{restriction} $\res_{H}^{G} V$ of $V$ consists of the vector space $V$ viewed as a module for $H$ by restricting the action of $G$ to the subgroup $H$.

	Let $\{g_{1},\ldots, g_{n}\}$ be a set of coset representatives for $H$ in $G$. Then the \emph{induced module}  of $W$ to a $G$-module, denoted $\induce_{H}^{G}W$, has vector space $\oplus_{i} \, \left( \{g_{i} \} \oplus W\right)$ and action given by
	\[g \,\left( \sum_{i=1}^{n}(g_{i}, \, w_{i}) \right) = \sum_{i=1}^{n} (g_{r(i)}, \, h_{i} \, w_{i} ) \]
	where $r_{i}, h_{i}$ are the unique values satisfying $g \cdot g_{i} = g_{r(i)} \cdot h_{i}$. This definition may be extended to include restricting/inducing representations and characters as well, see \cite[Chapter 21]{james2001representations}.
\end{defn}

There is a close link between the restriction and induction of modules, this is shown off best in a theorem known as \emph{Frobenius reciprocity} (see  \cite[Theorem 21.16]{james2001representations}).

\begin{thm}[Frobenius Reciprocity]
	Let $H\subseteq G$ be groups, and suppose $\varphi$ and $\chi$ are characters of $H$ and $G$ respectively. Then
	\[\langle\induce_{H}^{G}\varphi, \chi \rangle_{G} = \langle \varphi, \res_{H}^{G}\chi \rangle_{H},\]
	where the inner product for characters is computed over $G$ and $H$ respectively.
\end{thm}

The decomposition of an induced module $\induce_{H}^{G} V$ into the simple modules of $G$, or the restricted module $\res_{H}^{G} W$ into simple modules of $H$, is an important question in representation theory.
Results which allow us to relate the of the modules of $G$ in terms of those for $H$ and vice versa are called \emph{branching rules}. The natural recursive structure of the symmetric group allows us an answer to this question for $S_{n-1} \subseteq S_{n}$. A proof of the following result can be found in \cite[Theorem 2.8.3]{sagan2013symmetric}

\begin{thm}[Branching rules for $S_{n}$]
	\label{chpt4:thm:branching}
	Let $n \geq 1$, and $\lambda \vdash n$. The \emph{branching rules} for the simple module of the symmetric group are as follows:
	\begin{eqnarray}
	\res_{S_{n-1}}^{S_{n}} S^{\lambda} & \cong & \bigoplus_{\substack{\mu \vdash n-1 \\ \mu \subseteq \lambda}} S^{\mu} \textnormal{ as $\mathfrak{S}_{n-1}$-modules}\label{chpt4:eqn:restrict}\\
	\induce_{S_{n}}^{S_{n+1}} S^{\lambda} & \cong &\bigoplus_{\substack{\mu \vdash n+1 \\ \lambda \subseteq \mu }} S^{\mu} \textnormal{ as $\mathfrak{S}_{n+1}$-modules} \label{chpt4:eqn:induce}
	\end{eqnarray}
\end{thm}	
The branching rules for the symmetric group are closely related to Young's lattice. The direct sum of \eqref{chpt4:eqn:restrict} could be rephrased as take a direct sum of all Specht modules found by removing a box of $\lambda$, similarly the direct sum \eqref{chpt4:eqn:induce} may be thought as taking a sum of all Specht modules found by adding a box to $\lambda$. The branching rules allow us to describe the restriction or induction of any $\mathfrak{S}_{n}$-module, all we must do is decompose it into its simple summands and then apply Theorem \ref{chpt4:thm:branching} to each part individually. The recursive structure of $\mathfrak{S}_{n}$-modules is key to the study of the random-to-random shuffle and one-sided transposition shuffle.

\subsection{Switching to Words }
\label{chpt4:subsec:words}
The notation of tabloids is cumbersome, we therefore introduce a one-to-one correspondence between certain words of length $n$ and tabloids of size $n$, allowing us to describe our permutation and Specht modules in more succinct notation. This notation is particularly useful when studying the action of $\mathfrak{S}_{n}$ on these modules.

\paragraph{}
Given $n\in \mathbb{N}^{0}$ we denote by $W^n$ the set of words of length $n$ 
with letters in $[n]$, where by a word of length $n$ we simply mean a 
string $w=w_1\, w_2 \,\ldots \,
w_n$ with $w_{i} \in [n]$ for all $i$. For $n=0$ we allow a special case where $W^{0}$ is comprised solely of the empty word denoted $\omega$. Note that in forming words we regard 
the elements of $[n]$ as distinct symbols. Later on it will be notationally convenient to have our words comprised of positive integers, for example $W^{2} = \{11,\, 
12, \,  21, \, 22\}$. The size of the set $W^{n}$ is $n^{n}$.

There is a natural action of the symmetric group $S_n$ on $W^n$. For a 
word $w=w_1\, w_2\, \ldots \, w_n \in W^n$ and an element $\sigma 
\in S_n$, we let 
$\sigma \, w := w_{\sigma^{-1}(1)}\, w_{\sigma^{-1}(2)}\, \ldots 
\, w_{\sigma^{-1}(n)} \in W^n$. 
We emphasise that this is the action of $S_n$ on words by place 
permutations, it is not the action of $S_n$ acting on the individual 
letters that comprise a word, 
e.g. if $\sigma = (123)$ then $\sigma \, (2\,3\,2) = 223  \neq 313$.
Let $M^{n}$ be the vector space over the field $\mathbb{C}$ with basis of 
words in $W^n$. The action of $S_{n}$ on words in $W^{n}$ extends linearly 
to an action of $\mathfrak{S}_{n}$ on the vector space $M^{n}$. Thus, 
$M^{n}$ is an $n^{n}$-dimensional module for the group algebra 
$\mathfrak{S}_{n}$. For elements of $M^{n}$ we use the notation $\cdot$ to separate the complex coefficients from the words in $W^{n}$, e.g. $2 \cdot 232 + 4i \cdot 213$.

To each word $w \in W^n$ we can associate an $n$-tuple of non-negative 
integers, called its \emph{evaluation}, denoted $\eval(w)$, as follows. For $1\leq i 
\leq n$, let $\eval_i(w)$ count the number of occurrences of the symbol $i$ 
in the word $w$, and then let $\eval(w):= (\eval_1(w),\ldots,\eval_n(w))$.
Note that $\sum_{i=1}^n \eval_i(w) = n$ for any word $w$ in $W^n$.
For example taking $232 \in W^{3}$, its evaluation is 
$\eval(2\,3\,2) = (0,2,1)$. If in addition $\eval(w)$ is a 
non-increasing sequence of integers, then we identify $\eval(w)$ with the 
corresponding partition of $n$, ignoring possible ending zeros at the end 
of the partition. For example $ w = 1231 \in W^{4}$ has 
evaluation $\eval(w) = (2,1,1,0)$ and we associate it to the partition $(2,1,1)$.
Note that the evaluation of a word is unchanged by the action of $S_{n}$, 
i.e. for any $\sigma \in S^{n}$ and $w \in W^{n}$, we have  $\eval(w) 
= \eval(\sigma \, w)$. Thus we may find a stable subspace of $M^{n}$ 
corresponding to the words with a given evaluation.

\begin{defn}
	\label{chpt4:def:permwords}
	Let $\nu$ be a $n$-tuple of non-negative integers. Define the module 
	$M^{\nu}$ as the following vector space
	\[M^{\nu} = \langle w \in W^{n} \, | \, \eval(w) = \nu \rangle \subseteq M^{n}.\]
\end{defn}

For a partition $\lambda \vdash n$ the module $M^{\lambda}$ from Definition 
\ref{chpt4:def:permwords} is exactly the permutation module we defined in 
Definition \ref{chpt4:def:permmod}. To establish this we create a 
correspondence between words of evaluation $\lambda$ and Young tableau of shape $\lambda$. For a partition $\lambda\vdash n$, a Young tableaux $T$ of shape $\lambda$ naturally corresponds to a word in $W^n$. 
\begin{defn}
	Let $\lambda \vdash n$. Define a map $w:\YT(\lambda) \rightarrow W^{n}$ as follows: for each tableau $T$ of shape $\lambda$, let $w(T) = w_{1} \,\ldots \, w_{n}$ be the word with $w_{T(i,j)} = i$ for each box $(i,j)$ in $T$.	Equivalently, the numerical entries in the $i^{\rm th}$ row of $T$ tell us in which positions to put the symbol $i$ in the word $w(T)$. Importantly the words formed by tableau of shape $\lambda$ have evaluation $w(T) = \lambda$ (possibly ignoring some zeroes).
\end{defn}

The map $w$ respects the action of $S_{n}$, i.e. for any $\sigma \in S_{n}$, we have $\sigma \, w(T) = w( \sigma T)$. Also for any $w \in M^{n}$ with $\eval(w)=\lambda$, there exists a tableau $T\in\YT(\lambda)$ such that $w(T) =w$. Thus, the map is surjective on words of the correct evaluation. 

\begin{example}
	Let $\lambda = (3,2)$, and $T =  
	\ytableausetup{mathmode,baseline,aligntableaux=center,boxsize=1em}
	\begin{ytableau} 5 & 4 & 1 \\2 & 3 \end{ytableau}$. The word in $W^{n}$ 
	corresponding to $T$ is $w(T) = 1 2 2  1 1$. Now 
	take the element $(123) \in 
	S_{n}$, then the action on tableaux (equivalently words) gives us: $(123) \, w(T) = 2 
	1 2 1 1 = w( (123) \, T)$. Note that the map $w$ is not injective, for example
	\[\ytableausetup{mathmode,baseline,aligntableaux=center,boxsize=1em} 
	w\left( 	\begin{ytableau} 5 & 4 & 1 \\2 & 3 \end{ytableau}\right) = 
	1 2 2 1  1 = w\left( 	\begin{ytableau} 1 & 4 
	& 5 	\\3 & 2 \end{ytableau}\right).\]

\end{example}
To 
link Definitions \ref{chpt4:def:permmod} and \ref{chpt4:def:permwords} we 
need to find a bijection from tabloids of shape $\lambda$ to words of evaluation $\lambda$, to do this we observe that two tableaux have the same corresponding word if and only if they 
belong to the same  class of tabloids.

\begin{lemma}
	\label{chpt4:lem:maptowords}
	Let $T_{1},T_{2}$ be two tableaux of shape $\lambda$. Then
	\[w(T_{1}) = w(T_{2}) \Leftrightarrow T_{1} \sim_{R} T_{2} \Leftrightarrow 
	\{T_{1}\} = \{T_{2}\}.\]
\end{lemma}
\begin{proof}
	This is clear from the description of $w(T)$.
\end{proof}

Following from this we can see that the map $w(T)$ on all tableaux induces a bijection between tabloids and words. Thus we can 
quickly establish our permutation modules as vector spaces over words.
\begin{lemma}
	Let $\lambda \vdash n$, the permutation module $M^{\lambda}$ for 
	$\mathfrak{S}_{n}$ may be seen as a vector 
	space over the following bases:
	\begin{eqnarray*}
		M^{\lambda} & = & \langle \{T\} \,|\, T \textnormal{ is a tableau of 
			shape 
		} 	\lambda \rangle \\
		& \cong & \langle w(T) \, | \,T \textnormal{ is a tableau of shape } 	
		\lambda \rangle \\ 
		& \cong & \langle w \in M^{n} \, | \,\eval(w) = \lambda \rangle
	\end{eqnarray*}  
\end{lemma}
\begin{proof}
	The map $w(T)$ forms a bijection between tabloids of shape $\lambda$ and 
	words of evaluation $\lambda$, which respects the action of 
	$\mathfrak{S}_{n}$.	
\end{proof}
Thus to find elements of the module $M^{\lambda}$ we now only have to think 
of words with a fixed evaluation. The regular module $M^{(1^{n})}$ may now be viewed as the module of words of length $n$ which contain every letter of $[n]$, e.g.
\begin{eqnarray}	
M^{(1^{3})} \cong \langle 123, \, 132, \, 213, \, 231, \, 312, \, 321 \rangle.
\end{eqnarray}
We may also establish a basis for the Specht module 
$S^{\lambda}$  using words, by replacing each tabloid by its 
corresponding word. E.g. rewriting $S^{(3,1)}$ from Example \ref{chpt4:ex:specht},
\begin{eqnarray}	
S^{(3,1)} \cong \langle 1112 -2111, \, 1121 -2111, \, 1211 -2111 \rangle.
\end{eqnarray}

\section{Random Walks on The Symmetric Group}

We now use the techniques introduced in Chapter \ref{chpt1:chpt} to study random walks on the symmetric group, otherwise know as \emph{shuffles}. 
We explore three different shuffles which have inspired the analysis of the one-sided transposition shuffle. These illustrate both the algebraic and probabilistic techniques needed to prove cutoff results on mixing time. We begin with the classic random transposition shuffle, where we present Diaconis and Shahshahani's argument using the upper bound lemma (Lemma \ref{chpt2:lem:upperboundlemma}). 
Next we study the top-to-random shuffle where we use strong stationary times to get an effective bound on our mixing time. 
Finally we end the section with a discussion of the random-to-random shuffle, for which the existence of a cutoff was recently proven. The random-to-random shuffle is our first look at the technique of lifting eigenvectors. We give explicit examples of this procedure but leave the full algebraic details for our analysis of the one-sided transposition shuffle.

\subsection{The Random Transposition Shuffle}
\label{chpt4:subsec:rt}

Our first significant example is the random transposition shuffle. This 
was first analysed by Diaconis and Shahshahani in 1981 
\cite{diaconis1981generating}. In this seminal paper they proved the upper 
bound lemma and used it to find tight bounds on the mixing time of the 
random transposition shuffle.	Since this important work there has been much research into the behaviour of the random transposition shuffle. Berestycki has studied the hyperbolic geometry of the random transposition shuffle when formulated as a random walk on the Cayley graph of $S_{n}$ generated by transpositions \cite{berestycki2006hyperbolic}. Furthermore, Berestycki and Durrett showed that under a continuous time random transposition shuffle the number of transposition required to return to the identity undergoes a phase transition around time $n/2$ \cite{berestycki2006phase}. A natural extension of the random transposition shuffle is the random $k$-cycles shuffle where at each step we of the shuffle we apply a uniformly chosen $k$-cycle. In 2011, Berestycki, Schramm and Zeitouni proved that the random $k$-cycles shuffle has mixing time $(n/k)\log n$ \cite{Berestycki2011}.

Lately there has been an effort to study the cutoff (or limit) profile of random walks on $S_{n}$, that is the exact behaviour of $\lVert \textnormal{P}_{n}^{t_{n}  +cw_{n}} - \pi_{n} \rVert_{\tiny \textnormal{TV}}$ as $n\to\infty$ where $\textnormal{P}_{n}$ exhibits a cutoff at time $t_{n}$. 
Naturally the first cutoff profile we want to understand is that of the random transposition shuffle and in 2016 Berestycki posed this exact question at an AIM workshop \cite{Bconj}. This conjecture was recently settled in a breakthrough paper by Teyssier where he proved that 
\[\lim_{n\to\infty}\lVert \RT^{(n/2)\log n +cn} - \pi \rVert_{\tiny \textnormal{TV}} =  \lVert \textnormal{Po}(1 + e^{-2c}) - \textnormal{Po}(1)\rVert_{\tiny \textnormal{TV}} \]
where $\textnormal{Po}(r)$ represents the Poisson distribution with rate $r$ \cite{teyssier2020limit}. In order to prove the limit profile of the random transposition shuffle Teyssier derived a improvement of Diaconis' upper bound lemma (Lemma \ref{chpt3:lem:upperboundlemma}). This improvement of the upper bound lemma has opened the path to study limit profiles of many other random walks, including the previously mentioned random $k$-cycles shuffle \cite{nestoridi2020limit, freslon2020cutoff}. 
This collection of work demonstrates how important the random transposition shuffle is to the study of random walks on groups and how it is still influencing the field even today.

\paragraph{}
Throughout this section we use 
information about the irreducible representations of $S_{n}$ which we 
proved in Section \ref{chpt4:sec:modules}. The argument we present here is 
given by Diaconis in his book \emph{Group Representations in Probability 
	and Statistics} \cite[Chapter 3D, Theorem 5]{Diaconis1988}.\\

\paragraph{}
The random transposition shuffle is a random walk on $S_{n}$ given by the following procedure: at time $t$, choose position $i$ uniformly at random with your left hand, and independently choose position $j$ uniformly at random with your right hand, then swap the cards at the positions. Swapping the cards at positions $i$ and $j$ amounts to applying the transposition $(i \, j)$.

\begin{defn}
	\label{chpt4:def:RT}
	The \emph{random transposition shuffle} is the  
	random walk on 	$S_n$ generated by the following probability distribution:
	\[\RT(\sigma) = \begin{cases}
	1/n & \text{ if } \sigma = e \\
	2/n^{2} & \text{ if } \sigma = (i \, j) \text{ with } i < j \\ 
	0 & \text{ otherwise }
	\end{cases}.\]
\end{defn}

We now state the cutoff result for the random transposition shuffle which we 
work towards throughout this section.

\begin{thm}
	\label{chpt4:thm:RTcutoff}
	Let $t_{n} = (n/2)\log n$, 	and $w_{n} = n$. The random transposition shuffle $\RT$ satisfies the following bounds for any $c>0$:
	\begin{eqnarray}
	\limsup_{n\to \infty}\lVert \RT^{t_{n} +cw_{n}} -\pi_{n} \rVert_{\tiny 
		\textnormal{TV}} & \leq &A e^{-2c} \textnormal{ for a universal 
		constant } A\\
	\liminf_{n\to\infty}\lVert \RT^{t_{n} -cw_{n}} -\pi_{n} \rVert_{\tiny 
		\textnormal{TV}} & \geq & \frac{1}{e} - e^{-e^{2c}} 
	\end{eqnarray}
	Thus, the random transposition shuffle exhibits a cutoff in total variation distance at time $(n/2)\log n$ with a window of size $n$.
\end{thm}

\subsubsection{Upper Bound}
We begin by establishing the upper bound present in Theorem \ref{chpt4:thm:RTcutoff}.  
The random transposition shuffle is constant on the conjugacy classes of $S_{n}$, thus 
making it amenable to analysis via discrete Fourier transforms. 
We have seen that the simple modules of $S_{n}$ are labelled by partitions of $n$.  Let $\lambda \vdash n$ be a partition, and let $\chi_{\lambda}$ denote the character of the corresponding irreducible representation.

The random transposition shuffle only takes non-zero values of the conjugacy class of the identity and of transpositions.
The conjugacy class of the identity has size $1$, and the class of transpositions has size ${n 
	\choose 2} = (n(n-1))/2$. Applying 
Lemma \ref{chpt3:lem:trace} we find the Fourier transform 
of the random transposition shuffle at the irreducible representation $\lambda$ to be,
\[\widehat{\RT}(\lambda) = \left(\sum_{i} |C_{i}| \RT(g_{i}) 	
\frac{\chi_{\lambda}(g_{i})}{d_{\lambda}}\right)\cdot \I = 
\left(\frac{1}{n} + \frac{n-1}{n}\frac{\chi_{\lambda}(\tau)}{d_{\lambda}} 
\right)\cdot \I\]
where the first sum is over the conjugacy classes of $S_{n}$ with $g_{i}$ a class representative, and $\tau$ in the second equation is any transposition.	A straightforward application of the Upper Bound Lemma (Lemma \ref{chpt3:lem:upperboundlemma}) now gives us
\begin{eqnarray}
\label{chpt4:eqn:RTupperboundlemma}
4\lVert \RT^{t} -\pi_{n} \rVert_{\tiny \textnormal{TV}}^{2} \leq \sum_{\substack{\lambda \vdash n \\ \lambda \neq (n)}} 
d_{\lambda} \textnormal{Tr}\left( 
\widehat{\RT}(\lambda)^{t}\overline{\widehat{\RT}(\lambda)^{t}}\right) = 
\sum_{\substack{\lambda \vdash n\\\lambda \neq (n)}} 
d_{\lambda}^{2}\left(\frac{1}{n} + 
\frac{n-1}{n}\frac{\chi_{\lambda}(\tau)}{d_{\lambda}} \right)^{2t}.
\end{eqnarray}
To bound this sum we need to understand the irreducible characters of 
$S_{n}$, particularly their values on the conjugacy class of transpositions.

Below we give a combinatorial formula for the value 
of $\chi_{\lambda}(\tau) / d_{\lambda}$ in terms of the Young diagram 
$\lambda$, see \cite[Chapter 3D, Fact 2]{Diaconis1988} for details.
\begin{lemma}
	Let $\lambda \vdash n$, and $\tau$ be any transposition. Then
	\[ \frac{n(n-1)}{2} \frac{\chi_{\lambda}(\tau)}{d_{\lambda}} = \D(\lambda) 
	\]
	where $\D(\lambda) = 
	\sum_{(i,j) \in \lambda} (j-i)$. 
\end{lemma} 

Using this equality we may rewrite the sum 
\eqref{chpt4:eqn:RTupperboundlemma} as 
\begin{eqnarray}
\label{chpt4:eqn:RTdiag}
\sum_{\substack{\lambda \vdash n\\\lambda \neq (n)}} 
d_{\lambda}^{2}\left(\frac{1}{n} + 
\frac{n-1}{n}\frac{\chi_{\lambda}(\tau)}{d_{\lambda}} \right)^{2t} = 	\sum_{\substack{\lambda \vdash n\\\lambda \neq (n)}} 
d_{\lambda}^{2}\left( \frac{n + 2\D(\lambda)}{n^{2}}\right)^{2t} \label{chpt4:eqn:rtfirstbound1}
.
\end{eqnarray}

The combinatorics of Young diagrams allows us to find and prove 
bounds on the terms present in \eqref{chpt4:eqn:RTdiag}.
We now state several facts which have intuitive proofs using our knowledge of Young diagrams.

\begin{lemma}
	\label{chpt4:lem:diagbounds}
	Let $\lambda,\mu$ be partitions of $n$,  such that $\lambda 
	\trianglerighteq \mu$. Then we have:
	\begin{enumerate}
		\item $\D(\lambda) \geq \D(\mu)$.
		\item $\D(\lambda^{\prime}) = - \D(\lambda)$.
		\item The following bound holds,
		\begin{eqnarray}
		2\D(\lambda) \leq \begin{cases}
		(n-1)n - 2(n-\lambda_{1})(\lambda_{1}+1) & 
		\textnormal{ 
			if } \lambda_{1} \geq \frac{n}{2}\\
		(\lambda_{1}-1)n & \textnormal{ for all } \lambda  
		\end{cases}
		\end{eqnarray}
	\end{enumerate}
\end{lemma}

\begin{proof}
	For our first assertion suppose that $\mu$ is one step below $\lambda$, 
	that is we may form $\mu$ by moving one box of $\lambda$ down and to 
	the left. When we move this box its diagonal index (column - row) must decrease, 
	therefore $\D(\lambda) \geq \D(\mu)$. If $\mu$ is more than one step 
	down from $\lambda$ we may apply this result inductively to find 
	$\D(\lambda) \geq \D(\mu)$.
	The second assertion follows from the construction of the transpose  $(i,j) \in \lambda^{\prime} 
	\Leftrightarrow (j,i) \in \lambda$.
	For the last property take any $\lambda \vdash n$, in general we 
	know that the first row of $\lambda$ gives 
	the biggest contribution to its diagonal index sum. Therefore, by taking 
	this value and multiplying by the number of rows in $\lambda$ we may 
	form the following bound
	\[\D(\lambda) = \sum_{(i,j) \in \lambda} (j-i) \leq 
	\frac{n}{\lambda_{1}} \sum_{j=1}^{\lambda_{1}} (j-1) = 
	\frac{(\lambda_{1}-1)n}{2} .\]
	Now suppose $\lambda_{1} \geq n/2$, then we may see that our partition 
	is dominated by $(\lambda_{1}, 
	n-\lambda_{1}) \trianglerighteq \lambda$, by fixing the first row and move all other boxes up and to the right to join the second row.  Therefore, by property 2 we 
	have:
	\[\D(\lambda) \leq \D(\lambda_{1}, 
	n-\lambda_{1}) = \sum_{j=1}^{\lambda_{1}} (j-1) + 
	\sum_{j=1}^{n-\lambda_{1}} (j-2) = \frac{n(n-1)}{2} - 
	(n-\lambda_{1})(\lambda_{1}+1). \qedhere \]
\end{proof}
To simplify our upper bound we also need a bound on the dimension sum of partitions with fixed first row $\lambda_{1}$. Details of the following result may be found in \cite[Chapter 3D]{Diaconis1988}.
\begin{lemma}
	\label{chpt4:lem:youngdimbound}
	Let $\lambda \vdash n$, then
	\begin{eqnarray}
	\sum_{\substack{\lambda \vdash n \\ \lambda_{1} = n-k}} d_{\lambda}^{2} 
	\leq {n\choose k}^{2} \sum_{\substack{\lambda \vdash n \\ \lambda_{1} = 
			n-k}} d_{\lambda / \lambda_{1}}^{2} \leq {n\choose k}^{2} k!
	\label{chpt4:eqn:youngdimbound}
	\end{eqnarray}
\end{lemma}

Using the results we have just proved we may simplify the summation 
\eqref{chpt4:eqn:RTdiag}. 
We first double up our sum by only concentrating on those partitions 
which have positive diagonal sum. Note that we have to deal with the 
partition $(1^{n})$ as a special case because its transpose $(1^{n})^{\prime} = (n)$ is not 
presented in the sum \eqref{chpt4:eqn:RTdiag}. We then parametrise our sum based on the value of $\lambda_{1}$ splitting our sum around the value of $\lambda_{1} \lessgtr 3n/4$. Lemma \ref{chpt4:lem:diagbounds} helps to bound the terms with given $\lambda_{1}$ as follows:
\begin{eqnarray}
\frac{n + 2\D(\lambda)}{n^{2}} \leq  \begin{cases} 1 - \frac{2(\lambda_{1}+1)(n-\lambda_{1})}{n^{2}} & \text{ if } \lambda_{1} \geq 3n/4 \\
\frac{\lambda_{1}}{n} & \text{ for all } \lambda
\end{cases} \label{chpt4:eqn:firstrowbound}.
\end{eqnarray}
Subsequently Lemma \ref{chpt4:lem:youngdimbound} may be used to bound the multiplicities of the new terms. Performing these steps gives us,

\begin{eqnarray}
\eqref{chpt4:eqn:RTdiag} & =   & \sum_{\substack{\lambda \vdash n \\\lambda \neq (n) \\ 
		\D(\lambda) \geq 
		0}} d_{\lambda}^{2}\left(\frac{n + 2\D(\lambda)}{n^{2}} \right)^{2t} + 
\sum_{\substack{\lambda \vdash n \\\lambda \neq (1^{n}) \\ \D(\lambda) < 
		0}} 
d_{\lambda}^{2}\left(\frac{n + 2\D(\lambda)}{n^{2}} \right)^{2t} + \left( 
\frac{n +2\D((1^{n}))}{n^{2}}\right)^{2t} \nonumber \\
& \leq  & \sum_{\substack{\lambda \vdash n \\\lambda \neq (n) \\ 
		\D(\lambda) \geq 0}} 
d_{\lambda}^{2}\left(\frac{n + 2\D(\lambda)}{n^{2}} \right)^{2t} + 
\sum_{\substack{\lambda \vdash n \\\lambda \neq (1^{n}) \\ \D(\lambda) < 0}}
d_{\lambda^{\prime}}^{2}\left(\frac{n + 2\D(\lambda^{\prime})}{n^{2}} 
\right)^{2t}  + \left( 
\frac{n -n(n-1)}{n^{2}} \right)^{2t} \nonumber\\
& = & 2\sum_{\substack{\lambda \vdash n \\\lambda \neq (n) \\ 
		\D(\lambda) \geq 0}} 
d_{\lambda}^{2}\left(\frac{n + 2\D(\lambda)}{n^{2}} \right)^{2t} + \left( 1 - \frac{2}{n}\right)^{2t} \label{chpt4:eqn:RTfirstsum} \\
& \leq & 2\sum_{k=1}^{n/4}\left( 1 - \frac{2(n-k+1)k)}{n^{2}} 
\right)^{2t} \sum_{\substack{\lambda \vdash n \\ \lambda_{1} =n-k}} 
d_{\lambda}^{2}+ 2\sum_{k>n/4}^{n-2}\left(1 - \frac{k}{n} 
\right)^{2t} 
\sum_{\substack{\lambda \vdash n \\ \lambda_{1} =n-k}} 
d_{\lambda}^{2}  + \left( 1 - \frac{2}{n}\right)^{2t} \nonumber\\
& \leq & 2\sum_{k=1}^{n/4}\left( 1 - \frac{2(n-k+1)k)}{n^{2}} 
\right)^{2t}{n \choose k}^{2} k!+ 2\sum_{k> n/4}^{n-2}\left(1 - \frac{k}{n} \right)^{2t} {n \choose k}^{2} 
k!   + \left( 1 - \frac{2}{n} \right)^{2t} \label{chpt4:eqn:RTbothsums}
\end{eqnarray}
After our simplifications we are left with summations in  
$n,k,$ and $t$. Finally we are in a position to analyse 
the behaviour of these sums around time $t = (n/2)\log n +cn$. The singular term given by the partition $(1^{n})$ may be easily dealt with at this time, 
\[\limsup_{n\to \infty} \left( \frac{n -n(n-1)}{n^{2}}\right)^{2t} 
=\limsup_{n\to 
	\infty} \left(1 - 
\frac{2}{n}\right)^{n\log n +2 cn} = 0.\]
The first summation of \eqref{chpt4:eqn:RTbothsums} contains the largest values of the Fourier transforms $\widehat{\RT}(\lambda)$. These values control the mixing time of the random transposition shuffle. The first term, $k=1$, in the summation is given by
\begin{eqnarray}
n^{2}\left(1 - \frac{2}{n}\right)^{2t} \label{chpt4:eqn:rtsecondbiggesteig}
\end{eqnarray}
this term is tightly bounded by $1$ at time $t =(n/2)\log n$ as $n\to\infty$. If we add an additional window of $cn$ to give time $t =(n/2)\log + cn$ we find the first term bounded by an exponential decay,
\[n^{2}\left(1 - \frac{2}{n}\right)^{2\left((n/2)\log n +cn\right)} \leq n^{2} e^{-2\log n - 4c} \leq e^{-4c}.\]
We may show that the ratio of subsequent terms in the first sum is less than $1$ for $n\geq 17$. This allows us to bound the first summation via a geometric series and conclude that it may be bounded by $Ae^{-4c}$ for a universal constant $A$ at the time $(n/2)\log n + c n$ for $n$ sufficiently large. 
Now we look to bound the second summation in \eqref{chpt4:eqn:RTbothsums} at time $(n/2)\log n +cn$.
The terms in the second summation can be shown to be decreasing in $k$ for $n$ sufficiently large and thus we may bound the sum by $3n/4$ times the first term. Using Stirling's formula the first term can be seen to be tending to $0$ at time $(n/2)\log n + cn$ as $n\to\infty$. Further details of this analysis can be found in \cite{Diaconis1988}. Putting these bounds together we are now in a position to prove our upper bound on mixing time for the random transposition shuffle.

\begin{proof}[Proof of the upper bound in Theorem \ref{chpt4:thm:RTcutoff}]
	We reduced the analysis of total variation distance to the bound of summation \eqref{chpt4:eqn:RTbothsums}. Taking the limit as $n\to\infty$ we may establish the following bound on total variation distance at time $(n/2)\log n +cn$,
	\begin{eqnarray*}
		\limsup_{n\to \infty} \; 4\lVert \RT^{t_{n} + cw_{n}} -\pi_{n} 
		\rVert_{\tiny \textnormal{TV}}^{2} & \leq  &\limsup_{n\to \infty}
		\sum_{\substack{\lambda \vdash n\\\lambda \neq (n)}} 
		d_{\lambda}^{2}\left(\frac{1}{n} + 
		\frac{n-1}{n}\frac{\chi_{\lambda}(\tau)}{d_{\lambda}} \right)^{2(t_{n} + cw_{n})} \\ 
		&\leq & A e^{-4c} \, \textnormal{ for a universal constant } A. \qedhere 
	\end{eqnarray*}
\end{proof}

\subsubsection{Lower Bound}
We now prove the lower bound of Theorem \ref{chpt4:thm:RTcutoff}. Let $
F_{n} = \{\sigma \in S_{n} \,| \, \sigma $ has at least one 
fixed point$\}$. Using the definition of total variation distance we may 
find
\[\lVert \RT - \pi_{n} \rVert_{\tiny \textnormal{TV}} \geq | \RT(F_{n}) - 
\pi_{n}(F_{n}) | .\]
It is a classical result \cite{feller2008introduction} that as $n\rightarrow \infty$ the chance that a permutation picked uniformly at random has no fixed points tends to $1/e$. Therefore,
\[\lim_{n\to\infty} \pi_{n}(F_{n}) = 1 - \frac{1}{e}.\]
Next we show that under the random transposition measure we have a 
high probability of being in set $F_{n}$ at time $(n/2)\log n -cn$. 
Consider the set of random variables $\{R^{i},L^{i}\}_{i\in\mathbb{N}}$ which track 
which positions our hands have picked at time $i$ in the random transposition shuffle. Clearly if the random 
variables $R^{i},L^{i}$, have not chosen every position by time $t$ then 
the permutation 
we find ourselves at after $t$ steps of random transposition shuffle must have at least one 
fixed point. Let $J_{n}^{t}$ be the event $\{\cup_{i=1}^{t}\{R^{i},L^{i}\} 
\subsetneq [n]\}$, then by the reasoning above we have $J_{n}^{t}\subseteq 
F_{n}$, and hence $\mathbb{P}(J_{n}^{t}) \leq \RT^{t}(F_{n})$. The 
probability $\mathbb{P}(J_{n}^{t})$ may be computed using a coupon collector's 
argument. In the uniform coupon collector's problem, the expected time to 
collect all coupons is $n\log n$. However, we are performing two independent 
trials at each step, leading our expected time to choose every card to quicken to $(n/2)\log n$. More careful arguments in Feller 
\cite[Section IV.2]{feller2008introduction} show that,
\[\mathbb{P}(J_{n}^{t}) = 1 - e^{-ne^{-2t/n}} + o(1) \textnormal{, as } n\to \infty .\]
\begin{proof}[Proof of the lower bound in Theorem \ref{chpt4:thm:RTcutoff}]
	Let $t_{n} = (n/2)\log n$ and $w_{n} =n$. Putting together our bounds on the probability of $F_{n}$ under our different distributions we may complete our lower bound on total variation distance as follows,
	\[\liminf_{n\to\infty}\lVert \RT^{t_{n} -cw_{n}} - \pi_{n} \rVert_{\tiny \textnormal{TV}} \geq \liminf_{n\to\infty} \, \RT^{t_{n} -cw_{n}}(F_{n}) - \pi(F_{n})\geq 	\frac{1}{e} - e^{-e^{2c}} .\]
	This completes the proof of a total variation cutoff for the random transposition shuffle.
\end{proof}

\paragraph{}
The analysis of the random transposition shuffle is exemplary in a few important ways.  We always look to reduce the upper bound on total variation distance to a collection of summations like those in \eqref{chpt4:eqn:RTbothsums}. At this point we may perform analysis of the summations to find a value of $t$ for which they are bounded. The mixing time of random walks on groups are thought to be controlled by their second biggest eigenvalue. The first sum \eqref{chpt4:eqn:RTbothsums} contains the largest eigenvalues of the random transposition shuffle, including the second biggest eigenvalue \eqref{chpt4:eqn:rtsecondbiggesteig}, and is closely bounded at time $(n/2)\log n$, whereas the second summation disappears at time $(n/2)\log n$ as $n\to\infty$. Therefore, the first summation in  \eqref{chpt4:eqn:RTbothsums} can be thought of as determining the mixing time for the random transposition shuffle.

The lower bound argument presented for the random transposition shuffle  shows how versatile the simple bound in Lemma \ref{chpt2:lem:lowerbound} can be when applied to a careful choice of set. Often we like to reduce a lower bound in total variation distance to other probabilistic problems, in this instance a coupon collector's problem. This gives us another insight into why a cutoff occurs, hard to reach sets stop the total variation decaying before they are hit, keeping it close to 1 until we pass the critical time of $(n/2)\log n$.

\subsection{The Top-to-Random Shuffle}
\label{chpt4:subsec:t2r}

\paragraph{}
The top-to-random shuffle is defined by the following procedure: choose 
a position of the deck uniformly at random, then insert the top card into 
this position. The top-to-random shuffle is an example of a random walk which 
is not constant on the conjugacy classes of $S_{n}$. This means the 
technique of discrete Fourier transforms does not work as effectively for the top-to-random shuffle, which leads us to use the method of strong stationary times instead.

\begin{defn}
	The \emph{top-to-random shuffle} has driving probability $\ttr$ defined by:
	\[\ttr(\sigma) = \begin{cases}
	1/n & \textnormal{ if } \sigma = (k\hspace{0.1cm} k-1 \, 
	\ldots \, 1)\textnormal{ for some } k\in [n] \\ 
	0 & \textnormal{ otherwise }
	\end{cases}
	.\]
	The elements $(k\, k-1 \, \dots \, 1)$ represent cycling the top $k$ 
	cards of the deck one space up the deck with the card in position $k$ moving to position $k-1$, and finally the top card of the deck moving to position $k$.
\end{defn}
We may clearly see that the top-to-random shuffle is not constant on the 
conjugacy classes of $S_{n}$, for a quick counter example notice that $\ttr( 
(2 \, 3)) \neq \ttr( (1 \, 2))$ for $n\geq 3$. In this section we prove the following total variation cutoff 
for the top-to-random shuffle, the arguments we present are taken from 
\cite[Chapter 4A, Theorem 1]{Diaconis1988} and \cite[Section 6.5.3]{Levin2017}.

\begin{thm}
	\label{chpt4:thm:ttrcutoff}
	Let $t_{n} = n\log n$ and $w_{n} = n$. The top-to-random shuffle $\ttr$ satisfies the following bounds for any $c>0$:
	\begin{eqnarray}
	\limsup_{n\to \infty} \lVert \ttr^{t_{n} + cw_{n}} - \pi_{n} 
	\rVert_{\tiny \textnormal{TV}} & \leq & e^{-c}
	\label{chpt4:eqn:ttrupper}\\ 
	\liminf_{n\to \infty} \lVert \ttr^{t_{n} - cw_{n}} 
	- \pi_{n} 
	\rVert_{\tiny \textnormal{TV}} & \geq &  1 - \frac{2}{e^{c-1}}
	\label{chpt4:eqn:ttrlower}
	\end{eqnarray}
	
	Thus, the top-to-random shuffle exhibits a cutoff in total variation distance at time $n\log n$ with a 
	window of size $n$.
\end{thm}

\subsubsection{Upper Bound}
We start by proving the upper bound \eqref{chpt4:eqn:ttrupper}, but for separation distance -- we can then appeal to
Lemma \ref{chpt2:lem:sepupperbound} to give a bound on total variation 
distance. Recall that our deck of cards is labelled $\{1,\ldots, n\}$ from bottom to top. As mentioned above to find a bound on separation distance we 
use a strong stationary time for the top-to-random shuffle.
To this end define $T$ to be the first time card $1$ is randomly 
inserted into the deck when following the top-ro-random shuffle. We prove that $T$ is a strong stationary time  and then reduce the analysis of $\mathbb{P}(T>n\log n +cn)$ to the uniform coupon collector's problem. 

\begin{lemma}
	\label{chpt4:lem:sst}
	Let $T$ be the first time card $1$ is inserted into the deck during the top-to-random shuffle. Then $T$ is a strong stationary time for the top-to-random shuffle.
\end{lemma}
\begin{proof}
	We start our shuffle from the identity permutation with card $1$ at the 
	bottom of the deck and card $n$ at the top. Define a sequence of stopping 
	times $T_{i}$, as the first time $i$ cards have been placed under the bottom card of the deck, i.e. card $1$. Notice that $T_{n-1}$ is the first time card $1$ has reached the top of our deck, define the stopping time $T_{n}$ as the time when card $1$ is inserted into the deck uniformly at random. Our proposed strong stationary time is $T= T_{n} = T_{n-1}+1$. By induction we prove that at time $T_{i}$ the bottom $i$ cards of the deck are distributed according to a uniformly chosen permutation of $S_{i}$. Clearly at time $T_{1}$ the first card to be inserted below card $1$ is in a uniformly chosen permutation of $S_{1}$. Assume our hypothesis to be true for time 
	$T_{i}$, so the bottom $i$ cards of our deck are arranged in a uniformly chosen permutation of $S_{i}$. At time $T_{i+1}$ we have moved a new card into one of 
	$i+1$ positions below card $1$, thus forming a permutation of $S_{i+1}$ 
	amongst the bottom $i+1$ cards. To see that this permutation is 
	uniformly chosen from $S_{i+1}$ we must recognise that the new card has equal probability of being inserted into of any of the $i+1$ positions below card $1$. Pairing this fact with our inductive hypothesis means that all permutations of the bottom $i+1$ cards are equiprobable. Furthermore, the arrangement of the bottom $i+1$ cards is independent of the value of $T_{i+1}$.
	At time $T = T_{n}$ the entire deck is distributed according to a uniformly chosen permutation of $S_{n}$, and so $T$ is a strong stationary time.
\end{proof}

\begin{lemma}
	\label{chpt4:lem:ttrtimebound}
	The time $T$ satisfies the following bound: 
	\begin{eqnarray}
	\mathbb{P}(T> n\log n +cn) \leq e^{-c}
	\end{eqnarray} 
\end{lemma}
\begin{proof}
	Define a sequence of stopping times $T_{i}$ as in Lemma \ref{chpt4:lem:sst}. Our first observation is that
	\[T  = T_{n}= (T_{n} - T_{n-1}) + (T_{n-1}- T_{n-2}) + \dots + (T_{2} -T_{1}) + 
	T_{1}.\]
	The increments $T_{i} - T_{i-1}$ are independent and geometrically 
	distributed with parameter $i/n$. The random variable $T$ shares the same distribution as the uniform coupon collector's problem which may be described as follows:
	Suppose at each time step we choose a card from our deck of $n$ cards uniformly at random, if this 
	is the first time we have seen this card we say we \emph{collect} this 
	card. The coupon collector's problem asks how long does it takes to 
	collect every card of the deck?

	Let 
	$\mathcal{T}_{i}$ be the first time we have collected $i$ cards from the deck, 
	clearly we have
	\[\mathcal{T}_{n} = (\mathcal{T}_{n}-\mathcal{T}_{n-1}) + \dots + 
	(\mathcal{T}_{2}-\mathcal{T}_{1}) + \mathcal{T}_{1}.\]
	The increments $\mathcal{T}_{i} - \mathcal{T}_{i-1}$ are independent 
	and are geometrically distributed with parameter $(n-i)/n$.
	Therefore, $T_{n-i} - T_{n-i-1} \sim_{D} \mathcal{T}_{i} - 
	\mathcal{T}_{i-1}$ and so the random variables $T_{n} \sim_{D} \mathcal{T}_{n}$ 
	follow the same law.
	In the uniform coupon collector's problem let $C_{i}^{t}$ be the event 
	we have not collected the card $i$ by time $t$. Using the events $C_{i}^{t}$ we 
	may find a simple bound on $\mathcal{T}_{n}$ (and respectively $T$) at 
	the time $ t= n\log n +cn$,
	\begin{eqnarray}
	\mathbb{P}(T_{n} >t) = 
	\mathbb{P}(\mathcal{T}_{n} >t)  \leq  \mathbb{P} 
	\left(\bigcup_{i=1}^{n} C_{i}^{t} \right) 
	\leq  \sum_{i=1}^{n} \mathbb{P}(C_{i}^{t}) \leq \sum_{i=1}^{n} 
	\left( 1- \frac{1}{n}\right)^{n\log n +cn} \leq e^{-c}. \label{chpt4:eqn:classiccoupon}
	\end{eqnarray}
\end{proof}

\begin{proof}[Proof of Upper Bound in Theorem \ref{chpt4:thm:ttrcutoff}]
	Using the strong stationary time $T$, together with Lemma \ref{chpt4:lem:ttrtimebound}, we proven the upper bound \eqref{chpt4:eqn:ttrupper} as follows,
	\begin{eqnarray*}
		\limsup_{n\to\infty} \,\lVert \ttr^{n \log n +cn} -\pi \rVert_{\tiny 
			\textnormal{TV}} & \leq &\limsup_{n\to\infty} \, \lVert 
		\ttr^{n \log n +cn} -\pi_{n} \rVert_{\textnormal{Sep}}\\ 
		&\leq & \limsup_{n\to\infty}\, \mathbb{P}(T > 
		n\log n +cn ) \; \leq \; e^{-c} \qedhere
	\end{eqnarray*}
\end{proof}

\subsubsection{Lower Bound}
To prove the lower bound present in
Theorem \ref{chpt4:thm:ttrcutoff} we again appeal to Lemma 
\ref{chpt2:lem:lowerbound}. For $k\geq 2$, let  $I_{k}$ be the set of 
permutations our deck
such that the bottom $k$ cards remain in their original relative order, in other words $\sigma \in T_{k}$ if and only if $\sigma(1) < \sigma(2) < \dots < \sigma(k)$. 
Under the uniform distribution we have ${n\choose k}$ ways to arrange the 
bottom $k$ cards while retaining their relative order, after which we may 
place the remaining $n-k$ cards in any order giving us a factor of $(n-k)!$. 
Therefore, 
\begin{eqnarray}
\pi_{n}(I_{k}) = \frac{(n-k)!}{n!}{n\choose k } = \frac{1}{k!} \leq \frac{1}{k-1}
\label{chpt4:eqn:ttrlow1}.
\end{eqnarray}

To bound the probability of being in set $I_{k}$ under the top-to-random 
shuffle we make use of the random variables $T_{i}$ introduced in the proof 
of Lemma \ref{chpt4:lem:ttrtimebound}. The time $T - T_{k-1}$ is the first 
time the card $k$ has come to the top of our deck and been reinserted into the deck. 
Before this time the bottom $k$ cards always remain in the same relative order, therefore 
$\ttr^{t}(I_{k}) \geq \mathbb{P}(T- T_{k-1} > t)$. Using the decomposition 
of $T$ in Lemma \ref{chpt4:lem:ttrtimebound} we may rewrite the random 
variable into independent geometrically distributed increments $T - T_{k-1}= \sum_{i=k}^{n} (T_{i} -T_{i-1})$. The expectation and variance of the random 
variable $T - T_{k-1}$ is as follows:	
\begin{eqnarray}
\E[T - T_{k-1}] & = & \sum_{i=k}^{n} \frac{n}{i} \geq n (\log n -\log k)\\
\Var[T - T_{i-1}] & \leq & \sum_{i=k}^{\infty} \frac{n^{2}}{i(i-1)} \leq \frac{n^{2}}{k-1}.
\end{eqnarray}

Applying Chebyshev's inequality to the random variable $T - T_{k-1}$ we find,
\begin{eqnarray}
\mathbb{P}(T- T_{k-1} \leq n\log n - cn)  \leq \frac{1}{k-1} \label{chpt4:eqn:cheby}
\end{eqnarray}
provided that $c\geq \log k +1$. We are now in a position to prove the lower bound on the mixing time of the top-to-random shuffle.

\begin{proof}[Proof of Lower Bound in Theorem \ref{chpt4:thm:ttrcutoff}]
	Let $I_{k}$ be the set defined above and set $k-1= \lceil e^{c-1} \rceil$. Provided $n \geq e^{c-1} + 1$ we may use our bounds \eqref{chpt4:eqn:ttrlow1}, \eqref{chpt4:eqn:cheby}, to form the following lower bound on total variation distance,
	\[\lVert \ttr^{t_{n} -cw_{n}} - \pi_{n} 
	\rVert_{\tiny \textnormal{TV}} \geq  \ttr^{t_{n} - cw_{n}}(I_{k}) - \pi_{n}(I_{k})  \geq 1 -\frac{2}{e^{c-1}}.\]
	Taking the $\liminf$ as $n \to \infty$ gives the desired result.

\end{proof}

\paragraph{}
The top-to-random shuffle showcases the advantages of strong stationary times. 

In this case
analysis via discrete Fourier transforms would 
be difficult owing to its support being a mix of elements from different 
conjugacy classes. Instead a strong stationary time leads 
to a simple and effective upper bound on its mixing time, reducing the 
analysis to the well studied uniform coupon collector's problem. 
The strong stationary time $T$ was the limit of a sequence of stopping times, which build up a uniform set from $S_{1}$ to $S_{n}$. 

In Chapter \ref{chpt5:chpt} we present a novel strong stationary time for the one-sided transposition shuffle which does not rely on building up a sequence of uniform subgroups. The lower bound for the top-to-random shuffle further showcases the usefulness of Lemma \ref{chpt2:lem:lowerbound} in reducing the analysis of total variation distance to discrete random variables.

\subsection{The Random-to-Random Shuffle}
\label{chpt4:subsec:rtr}
The last shuffle we study in this chapter is the random-to-random shuffle. The random-to-random shuffle is described by the following procedure: at each step pick a card uniformly at random and take it from the deck, then pick a position uniformly at random and insert the card back into the deck at this position. 	This 
shuffle has a long history in the literature of mixing times of random 
walks. 	First introduced by Diaconis and Saloff-Coste 
\cite{diaconis1993comparison} in 1993, together they proved its mixing time was of  order  $n\log n$, and conjectured that a total variation cutoff occurs at time $(3/4)n\log n$.	This bound was improved by Uyemura-Reyes 
\cite{uyemura2002random} who showed the mixing time to be in the range $[(n/2)\log n, \,4n\log n]$. In 2012, Subag \cite{subag2013lower} was able to prove a lower bound of $(3/4)n\log n -(1/4) n\log \log n$ on the mixing time and in 2017, Bernstein and Nestoridi \cite{Bernstein2017} proved a matching upper bound of $(3/4)n\log n -(1/4)n\log \log n$ on the mixing time. Together these results prove the existence of a total variation cutoff for the random-to-random shuffle.

The recent work of Bernstein and Nestoridi \cite{Bernstein2017} relied on knowledge of the eigenvalues for the 
random-to-random shuffle. These were computed by Dieker and Saliola 
\cite{dieker2018spectral}, in a breakthrough publication where they 
invented the technique of \emph{lifting eigenvectors and eigenvalues}.  We 
dedicate this subsection to an exploration of the technique of lifting eigenvectors for the the random-to-random shuffle.
This technique uses the branching structure of the symmetric group to turn eigenvalues of the random-to-random shuffle on $n$ cards into eigenvalues of the random-to-random shuffle on $n+1$ cards.

\paragraph{}
To describe the random-to-random shuffle we introduce the concept of 
\emph{symmetrizing} a random walk. Consider a Markov chain $\{X^{t}\}_{t\geq 0}$ on 
$S_{n}$ driven by a probability $P$. By viewing our Markov chain backwards in time we may create a new a Markov chain $\{Y^{t}\}_{t\geq 0}$ driven by probability $P^{-1}(\sigma) := P(\sigma^{-1})$, this is called the \emph{time reversal process}.
Note that for a reversible random walk on a group driven by $P$ we necessarily have that $P^{-1}=P$, this follows straight from the detailed balance equations \eqref{chpt2:eqn:detail}. From any shuffle $P$ we may create a reversible random walk called the \emph{symmetrization} of $P$,  by taking the convolution of $P$ with its time reversal process, i.e., $P^{-1} \star P$.

Consider the top-to-random shuffle, its time reversal process does the following: it picks a card uniformly at random and places it back on top of the deck, the shuffle described by this process is called the \emph{random-to-top shuffle}, denoted $\rtt$. Imagine that we perform one random-to-top shuffle and then one top-to-random shuffle, the two steps together define a process whereby we pick a card uniformly at random and insert it back into the deck uniformly at random, we know this as the random-to-random shuffle.

\begin{defn}
	The \emph{random-to-random shuffle} on $S_{n}$ picks a card 
	uniformly at random and places it back into the deck uniformly at 
	random. It is defined by the 	following probability distribution:
	\begin{eqnarray}
	\rtr(\sigma) = \begin{cases}
	1/n & \textnormal{ if } \sigma = e \\
	1/n^{2} & \textnormal{ if } \sigma = (i \hspace{0.2cm} i+1 \hspace{0.2cm} \dots \hspace{0.2cm} j-1 \hspace{0.2cm} j) \text{ for } 1\leq i < j \leq n  \\
	1/n^{2} & \textnormal{ if } \sigma = (i \hspace{0.2cm} i-1 \hspace{0.2cm} \dots \hspace{0.2cm} j+1 \hspace{0.2cm} j) \text{ for } 1\leq j < i \leq n \\
	0 & \textnormal{ otherwise }
	\end{cases}.
	\end{eqnarray}
	The first non-identity permutations above correspond to taking a card from position $j$ and moving it down the deck to position $i$. The second non-identity permutation correspond to taking a card from position $j$ and moving it up the deck to position $i$.	We may also define the random-to-random shuffle as the \emph{symmetrization} of the top-to-random (or equivalently random-to-top) shuffle, $\rtr = \ttr \star \rtt$. The random-to-random shuffle defines an aperiodic, transitive, reversible Markov chain.
\end{defn}

The random-to-random shuffle is not constant on the conjugacy classes of 
$S_{n}$, again this may see by considering two transpositions, $\rtr ((1 
\, 2) ) \neq \rtr((1 \, n))$ for all $n \geq 3$. 
This makes the discrete Fourier transforms of probability $\rtr$ difficult to compute, leading us to require a different method to study this shuffle.
The random-to-random shuffle is reversible, which means if we can compute its eigenvalues we may apply the classical $\ell^{2}$ 
bound given in Theorem \ref{chpt2:thm:classicL2} to find an upper bound on the total variation distance of the random-to-random shuffle and uniform distribution. 
Dieker and Saliola provide an interesting combinatorial description (involving Young diagrams) of all the eigenvalues for the random-to-random shuffle which follows directly from the procedure of lifting eigenvectors.

\begin{thm}[Theorem 5 \cite{dieker2018spectral}]
	\label{chpt4:thm:rtreigenvalues}
	The eigenvalues of the random-to-random shuffle on $n$ cards are indexed by pairs of partitions
	$(\lambda, \mu)$ with $\lambda \vdash n$ and  $\lambda/\mu$ a horizontal strip. The 
	eigenvalue $\eig(\lambda,\mu)$ corresponding to $(\lambda,\mu)$ occurs with multiplicity $d_{\lambda}d^{\mu}$ and has value,
	\[\eig(\lambda,\mu) = \frac{1}{n^{2}}\left( {|\lambda|+1 \choose 2} - 
	{|\mu|+1\choose 2} + \D(\lambda) - \D(\mu) \right).\]
\end{thm}
The combinatorial description of the eigenvalues of the random-to-random shuffle hides away the algebraic techniques used to find them. 
The key idea behind Deiker and Saliola's work is to exploit the module structure of $S_{n}$ in order to create new eigenvectors from known ones. 
The first step is to change from studying the random-to-random shuffle as a probability to instead studying it as an element of the group algebra $\mathfrak{S}_{n}$.

\begin{defn}
	The random-to-random shuffle on $n$ cards can be viewed as the following element of $\mathfrak{S}_{n}$:
	\begin{eqnarray}
	\artr & := &n^{2} \sum_{\sigma \in S_{n}} \rtr(\sigma) \, \sigma 
	\nonumber\\
	& = &n \, 
	e + \sum_{1 \leq i <j \leq n} (i \hspace{0.2cm} i+1 \hspace{0.2cm} 
	\dots 
	\hspace{0.2cm} j-1 \hspace{0.2cm} j) + \sum_{1 
		\leq j <i \leq n} (i \hspace{0.2cm} i-1 \hspace{0.2cm} \dots 
	\hspace{0.2cm} j+1 \hspace{0.2cm} j)
	\end{eqnarray}
	We call this element the \emph{algebraic random-to-random shuffle}.
	Note that to create $\artr$ we have scaled our probability by a factor of $n^{2}$.
\end{defn}

This new viewpoint allows the random-to-random shuffle to act on modules of $S_{n}$ via the element $\artr$.

In particular we are interested in the action of $\artr$ on the regular module of $\mathfrak{S}_{n}$, because the eigenvectors and eigenvalues of this action are exactly the eigenvectors and eigenvalues for the probability $\rtr$.

\begin{lemma}
	\label{chpt4:lem:eigensame}
	A probability distribution $\nu:S_{n} \to [0,1]$ is an eigenvector for $\rtr$ with eigenvalue 
	$\epsilon$, if and only if $v = \sum_{\eta \in S_{n}} \nu(\eta) \, \eta  \, \in \mathfrak{S}_{n}$ is an eigenvector for $\artr$ with 
	eigenvalue $n^{2} \epsilon$.
\end{lemma} 

\begin{proof}
	Suppose $\nu$ is an eigenvector for $\rtr$ with eigenvalue $\epsilon$, this means that $P \star \nu  = \epsilon \; \nu$. Therefore, $\epsilon 
	\; \nu(\sigma) = \sum_{\eta} \rtr(\sigma \eta^{-1}) \nu(\eta)$.
	Now consider the action 
	\begin{eqnarray*}
		\artr \left( \sum_{\eta \in S_{n}} \nu(\eta)  \; 
		\eta\right) & = & n^{2}\sum_{\sigma \in S_{n}} \sum_{\eta\in S_{n}}  
		\rtr(\sigma) \nu(\eta) \, \sigma \eta\\ 
		& = &n^{2}\sum_{\tau\eta^{-1} \in 
			S_{n}} \sum_{\eta \in S_{n}}  \rtr(\tau \eta^{-1}) \nu(\eta) \, 
		\tau\\
		& = & 
		n^{2} \sum_{\tau \eta^{-1} \in S_{n}} \epsilon \, \nu(\tau) \, \tau = 
		n^{2}\epsilon \left(\sum_{\eta \in S_{n}} 
		\nu(\eta) \, \eta \right) .
	\end{eqnarray*}
	In the second equally we have relabelled $\sigma \eta = \tau$, and 
	the last equality we relabel the sum to match our original vector.
\end{proof}

Following Lemma \ref{chpt4:lem:eigensame} we may focus on the action of 
$\artr$ on the regular module with the understanding that all of the eigenvectors and eigenvalues of the random-to-random shuffle may be recovered. Recall that the regular module has the following decomposition
\begin{eqnarray}
\mathfrak{S}_{n} \cong \bigoplus_{\lambda \vdash n} d_{\lambda} 
S^{\lambda} \textnormal{ as $\mathfrak{S}_{n}$-modules}. \label{chpt4:eqn:spechtdecomp}
\end{eqnarray}

Note that because $\artr$ is an element of our group algebra it acts on the module $\mathfrak{S}_{n}$ and stabilises the  submodules $S^{\lambda}$. 
Therefore, the decomposition \eqref{chpt4:eqn:spechtdecomp} means we can reduce the problem of 
finding 
eigenvectors for the random-to-random shuffle on $\mathfrak{S}_{n}$ to the 
problem of finding eigenvectors belonging to the individual Specht modules 
$S^\lambda$.
Moreover, since the shuffle is acting as an element of the group algebra, 
we are then free to study its action on the single copy of $S^\lambda$ 
inside $M^\lambda$ to solve this problem rather than having to 
stick to the copies of $S^\lambda$ which appear inside $\mathfrak{S}_{n}$. 
In work that follows we present the eigenvectors of $\artr$ as 
elements of the Specht modules comprised of words as described in Section \ref{chpt4:subsec:words}.

\paragraph{}
The Specht modules of $\mathfrak{S}_{n}$ and $\mathfrak{S}_{n+1}$ are closely related by the branching structure in Theorem \ref{chpt4:thm:branching}, and by Young's lattice shown in Figure \ref{chpt4:fig:lattice}. The random-to-random shuffle also exhibits a recursive structure which can be seen by taking the difference of the shuffle on $n+1$ cards and the shuffle on $n$-cards,
\begin{eqnarray*}
	\artrp - \artr = \gpid + \sum_{\substack{1\leq i \leq n \\ j=n+1}} (i \hspace{0.2cm} i+1 \hspace{0.2cm} 
	\dots 
	\hspace{0.2cm} n \hspace{0.2cm} n+1) + \sum_{\substack{1\leq j \leq n \\ i=n+1}}	 (n+1 \hspace{0.2cm} n \hspace{0.2cm} \dots 
	\hspace{0.2cm} j+1 \hspace{0.2cm} j) .
\end{eqnarray*}
The permutations that do not involve the new card in position $n+1$ have disappeared in our comparison of the two elements. The key idea is that moving the card in position $j$ to position $i<j$ does not affect any of the cards above $j$, thus it does not matter how many cards above $j$ our deck contains. This recursive structure is at the heart of why lifting eigenvectors works and we recover a similar recursive structure for the one-sided transposition shuffle (see equation \eqref{chpt5:eqn:lrdifference}). The technique of lifting eigenvectors uses the branching structure for Specht modules and the recursive structure of the random-to-random shuffle, to turn eigenvectors of $\artr$ into eigenvectors of $\artrp$. In particular we lift eigenvectors belonging to the Specht modules $S^{\mu}$ into eigenvectors belonging to Specht modules $S^{\lambda}$ where $\mu \subset \lambda$, i.e., for $\lambda$ we can create from $\mu$ by adding boxes.

The following theorem summarises the main result of lifting eigenvectors for the random-to-random shuffle.

\begin{thm}[Theorem 21 \cite{dieker2018spectral}]
	\label{chpt4:thm:rtrlifting}
	Let $\lambda \vdash n$, and $1\leq i \leq l(\lambda)+1$ be such that $\lambda+e_{i} \vdash n+1$. Then there exists a linear map $\mathcal{L}^{\lambda}_{i}:S^{\lambda} \to S^{\lambda+e_{i}}$, such that,
	$\mathcal{L}^{\lambda}_{i}$ maps eigenvectors of the random-to-random shuffle on $n$ cards belonging to $S^{\lambda}$, to eigenvectors of the random-to-random shuffle on $n+1$ cards belonging to $S^{\lambda+e_{i}}$. In particular, if $v \in S^{\lambda}$ is an eigenvector of $\artr$ with eigenvalue $\epsilon$, then $\mathcal{L}^{\lambda}_{i}(v)$, if \textbf{non-zero}, is an eigenvector of $\artrp$ with eigenvalue $\epsilon + (n+1) + (\lambda_{i} +1) -i$.
\end{thm}

Using Theorem \ref{chpt4:thm:rtrlifting} we may find eigenvectors of a Specht module $S^{\lambda}$ from by lifting eigenvectors from below, however we must be careful that the map $\mathcal{L}^{\lambda}_{i}$ we apply does not kill the eigenvectors by mapping them to $0$. The non-injectivity of the lifting operators stops all the eigenvectors belonging to a Specht module $S^{\lambda}$ being described by lifting. However,
Deiker and Saliola managed to recover all the eigenvalues belonging to $S^{\lambda}$ by proving that the eigenvectors we fail to find by lifting belong to the kernel of $\artrp$ and thus have eigenvalue $0$. This was taken into account when formulating Theorem \ref{chpt4:thm:rtreigenvalues}. The lifting operators $\mathcal{L}^{\lambda}_{i}$ may be fully described in terms of linear operators on the space of words $M^{n}$ (see  Section \ref{chpt4:subsec:words}), which we now introduce.

\begin{defn}
	\label{chpt4:def:shuffleandswitch}
	Define two linear operators on the space $M^{n}$ spanned by words, to do 
	so, it is enough to define the effect on any given word. Let $w = w_{1} 
	\dots w_{n}$ in	$W^{n}$.
	\begin{enumerate}
		\item Let  $a \in [n+1]$.
		Define the \emph{shuffling operator}, denoted $\rtrsh_{i}: M^{n} \to M^{n+1}$ as 
		the following linear map
		\[\rtrsh_{a}(w) = \sum_{j=0}^{n+1} w_{1} \, \dots \, w_{j} \, a \, 
		w_{j+1} \, \dots \, w_{n}.\]
		The operator $\rtrsh_{a}$ is a linear combination of all words that are 
		formed by inserting $a$ into one position. The shuffling operator is	a linear map from $M^{\lambda}$ (or $S^{\lambda}$) to $M^{\lambda+e_{a}}$.
		
		\item Let $a,b\in[n]$. Define the \emph{switching operator}, denoted
		$\Theta_{b,a}:M^{n}\rightarrow 	M^{n}$ as follows:
		\[\Theta_{b,a}(w) := \sum_{\substack{1 \leq k \leq n \\ w_{k} = b}} 
		w_{1}\cdot \ldots \cdot w_{k-1} \cdot a \cdot w_{k+1}\cdot \ldots \cdot 
		w_{n}.\]
		The operator $\Theta_{b,a}$ forms a linear sum of all words created 
		from $w$ by replacing one occurrence of $b$ with $a$. If we restrict to domain $M^{\lambda}$ (or $S^{\lambda}$) then we end in the module $M^{\mu}$ with $\lambda +e_{a} = \mu + e_{b}$. 
	\end{enumerate}
\end{defn}
\begin{rem}
	We should define the shuffling and switching operators for all possible domains $M^{\lambda}$, however to do so would burden us with more notation, when it will always be clear from context which domain and codomain we are considering.
\end{rem}

\begin{example}
	Let $w = 1231 \in M^{4}$. Then  $\rtrsh_{1}(w) = 11231 + 11231 + 12131 + 
	12311 +12311 = 2 \cdot 11231 + 12131 + 2 \cdot 12311$, this is a element 
	belongs to the vector space $M^{5}$. We may see that $1231 \in M^{(2,1,1)} \subseteq M^{4}$ and $\rtrsh_{1}(1231) \in M^{(3,1,1)} \subseteq M^{5}$.
	An example of our switching operator is 
	$\Theta_{3,1}(w) = 1211$, and $\Theta_{1,3} (w) = 3231 + 1233$. Our images for the switching operator are $\Theta_{3,1}(w) \in M^{(3,1)}$ and $\Theta_{1,3}(w) \in M^{(1,1,2)}$.
	Another characteristic of the switching operator is that for any $w \in M^{n}$ and $i\in[n]$ we have $\Theta_{i,i}(w) =  \eval_i(w) \, w$ .
\end{example}

We may construct the lifting operators $\mathcal{L}_{i}$ as linear combinations of the shuffling and switching operators defined in Definition \ref{chpt4:def:shuffleandswitch}. Moreover, we may tell exactly when there is no possible lifting from module $S^{\mu}$ to $S^{\lambda}$ by analysis of the partitions $\mu$ and $\lambda$.

\begin{lemma}[Theorem 21 \cite{dieker2018spectral}]
	\label{chpt4:lem:liftingoperators}
	Let $\lambda \vdash n$. The lifting operators as defined in Theorem 
	\ref{chpt4:thm:rtrlifting} are a linear combination of the shuffling and switching operators, explicitly they are:	
	\begin{eqnarray}
	\mathcal{L}_{i}^{\lambda} = \sum_{1\leq b_{1} <\ldots <b_{m} <b_{m+1}=i} \left( \prod_{j=1}^{m} \frac{1}{(\lambda_{i} -i) - (\lambda_{b_{j}} -b_{j})}\Theta_{b_{j},b_{j+1}} \right) \circ \rtrsh_{b_{1}}.	
	\end{eqnarray}
	For $i=1,2$ the above equation simplifies to:
	\begin{eqnarray}
	\mathcal{L}^{\lambda}_{1} & = & \rtrsh_{1} \\
	\mathcal{L}^{\lambda}_{2} & = & \rtrsh_{2} + \frac{1}{(\lambda_{2}-2) - 
		(\lambda_{1}-1)}\Theta_{1,2} \circ \rtrsh_{1} 
	\end{eqnarray}
	Furthermore, let $\mu \subset \lambda$, such that $\mu + \sum_{j=1}^{k} e_{i_{j}} = \lambda$ where $i_{1} \leq \ldots \leq i_{k}$. Define the lifting operator from $\mu$ to $\lambda$, denoted $\mathcal{L}^{\lambda/\mu}$, as the map $\mathcal{L}^{\lambda/\mu} = \prod_{i=1}^{k} \mathcal{L}_{i_{j}}$, that is lift $\mu$ to $\lambda$ by forming each row in order.
\end{lemma}

\begin{thm}[Theorem 26 \cite{dieker2018spectral}]
	\label{chpt4:thm:detailedlifting}
	Suppose $\lambda \vdash n$. Every non-kernel eigenvector of $\artr$ belonging to $S^{\lambda}$ may be found by lifting a eigenvector in the kernel of $\artrn_{|\mu|}$ belonging $S^{\mu}$ for some $\mu \subset \lambda$ with $\lambda / \mu$ a horizontal strip. In other words if $v \in \mu$ is a eigenvector in the kernel of $\artrn_{|\mu|}$, then $\mathcal{L}^{\lambda /\mu} (v) \neq 0$ if and only if $\lambda/\mu$ is a horizontal strip. Thus, all eigenvalues for the module $S^{\lambda}$ may be recovered as they are the eigenvalue of a lifted eigenvector, or $0$ otherwise.
\end{thm}

We end this chapter with explicit examples of lifting eigenvectors of the random-to-random shuffle to Specht modules of $\mathfrak{S}_{1}$,  $\mathfrak{S}_{2}$, and $\mathfrak{S}_{3}$, using the lifting operators $\mathcal{L}^{\lambda}_{1},\mathcal{L}^{\lambda}_{2}$. We start with the only Specht module of $\mathfrak{S}_{1}$ which is $S^{(1)} = \langle 1 \rangle$.
Define the trivial partition $(0)$ to have a trivial eigenvector of the empty word, denoted $\omega$, with eigenvalue $0$. The only partition strictly contained in the partition $(1)$ is the empty partition $(0)$ and so we the only possible lifting is 
\[\mathcal{L}^{(1)/(0)}(\omega) = \mathcal{L}^{(0)}_{1}(\omega) =1.\] 
We easily see that  $\mathcal{L}^{(1)/(0)}(\omega)$ is the only eigenvector of $\artrn_{1}$ belonging to $S^{(1)}$ with eigenvalue $1$.

We now compute the eigenvectors for $\artrn_{2}$ belonging to the Specht modules of $\mathfrak{S}_{2}$. For the simple modules $S^{(2)} = \langle 11 \rangle$ and $S^{(1,1)} = \langle 12-21 \rangle$ we must start our lifting in the module $S^{(0)}$ because the we saw above that $S^{(1)}$ has no eigenvector in the kernel of $\artrn_{1}$ thus it can not be used for lifting. The skew diagram of $(2)/(0)$ is a horizontal strip, therefore 
\[\mathcal{L}^{(2)/(0)} (\omega) = \mathcal{L}_{1}^{(1)} \mathcal{L}_{1}^{(0)} (\omega) = \mathcal{L}_{1}^{(1)} (1) = 2\cdot 11\] 
is an eigenvector of $\artrn_{2}$. The skew diagram of $(1,1)/(0)$ is not a horizontal strip therefore when we applying the lifting operator we find 
\[\mathcal{L}^{(1,1)/(0)} (\omega) = \mathcal{L}_{2}^{(1)} \mathcal{L}_{1}^{(0)}(\omega) =  (12 +21 ) - (12 +21) = 0\]
as Theorem \ref{chpt4:thm:detailedlifting} asserts. The eigenvector of $\artrn_{2}$ belonging to $S^{(1,1)}$ is actually the basis element $12 -21$ which has corresponding eigenvalue $0$.

Our final example is lifting eigenvectors to the 2 dimensional Specht module $S^{(2,1)} = \langle 112 -211, 121 -211 \rangle$. Consider $\mu \subset (2,1)$, from our above work we know that the only modules $S^{\mu}$ which contain eigenvectors in the kernel of $\artrn_{|\mu|}$ are $\mu = (0)$ or $\mu= (1,1)$. The skew diagram of $(2,1)/(0)$ is not a horizontal strip so lifting $\omega$ does not recover an eigenvector. The skew diagram $(2,1)/(1,1)$ is a horizontal strip and so lifting the eigenvector $12-21$ we find
\begin{eqnarray}
\mathcal{L}_{1}^{(1,1)}(12-21) & = &2\cdot (112 - 211)  \nonumber
\end{eqnarray}
which belongs to $S^{(2,1)}$ and is an eigenvector of $\artrn_{3}$ with eigenvalue $4$. To find the remaining eigenvector for $S^{(2,1)}$ it is enough to look at the orthogonal complement of the subspace $\langle 112-211 \rangle \subset S^{(2,1)}$, this gives an eigenvector $112 - 2\cdot 121 +211$ with eigenvalue $0$.

\paragraph{}	
This ends our discussion of the random-to-random shuffle.  Since Dieker and 
Saliola's work Lafreni{\`e}re has shown that similar techniques can be applied to more general symmetrized shuffling operators \cite{Lafreni2018}. In the next chapter we show that the technique of lifting eigenvectors can be used to recover the spectrum of a variety of transposition shuffles, including the one-sided transposition shuffle . The lifting we present in Chapter \ref{chpt5:chpt} differs in key ways from Deiker and Saliola's \cite{dieker2018spectral} because we are analysing very different random walks. The random-to-random shuffle is a uniform measure on its non-identity support which covers many conjugacy class of $S_{n}$, whereas the one-sided transposition shuffle is not uniform on its single conjugacy class of transpositions.
These differences lead to changes in the underlying algebra of lifting eigenvectors and subsequentially a change in the required lifting operators. The biggest difference for the one-sided transposition shuffle is the corresponding lifting operators are injective and thus never kill our eigenvectors (see Theorem \ref{mapsnew}). This means that \emph{all} the eigenvectors as well as the eigenvalues of the one-sided transposition shuffle belonging to $S^{\lambda}$ may be recovered by lifting eigenvectors from modules $S^{\mu}$ with $\mu \subset \lambda$. In Chapter \ref{chpt5:chpt} we give the full details of this technique for the one-sided transposition shuffle.  Later in Chapter \ref{chpt6:chpt} we extend the technique of lifting eigenvectors to the random transposition shuffle and one-sided transposition shuffle on the hyperoctohedral group.

	\chapter{The One-sided Transposition Shuffle}
\label{chpt5:chpt}

In this chapter we introduce a new class of shuffles called 
\emph{one-sided transposition shuffles}: these have the defining property 
that at step $i$ the right hand's position ($\righthand$) is chosen 
according to a distribution supported on $[n]$, and given the value of 
$\righthand$ the distribution of the left hand's position ($\lefthand$) is supported on the set $\{1,\dots,\righthand\}$. In this chapter we 	restrict ourselves to the case where the left hand is chosen uniformly from $\{1,\ldots,\righthand\}$. We begin by focusing on the situation where our right hand is chosen uniformly from its possible range, we  call this the \emph{(unbiased) one-sided transposition shuffle}. Afterwards we generalise our techniques  to allow our right hand to be driven by a non-uniform distribution, we call these the \emph{biased one-sided transposition shuffles}.

\section{Main Results}
\label{chpt5:sec:mainresults}

In order to state our main results we introduce the (unbiased) one-sided transposition shuffle formally as follows.

\begin{defn} 
	\label{LRprob}
	The \emph{(unbiased) one-sided transposition shuffle} is the  
	random walk on 
	$S_n$ generated by the following distribution on the conjugacy 
	class of transpositions:
	\[\RL_{n}(\tau) = 
	\begin{cases}
	\frac{1}{n}\cdot\frac{1}{j} &  \text{if } \tau = (i\,j) \text{ for some 
	} 
	1\leq i 
	\leq j 
	\leq n\\
	0 & \text{otherwise.}
	\end{cases}
	\]
	
	We use the convention that all permutations $(i \,i)$ are equal to 
	the identity element $\gpid$, and therefore 
	$\RL_{n}(\gpid) = \frac{1}{n} (1 + \frac{1}{2} + \dots + 
	\frac{1}{n}) = H_{n}/n$, where $H_{k}$ denotes the 
	$k^{\textnormal{th}}$ 
	harmonic number.	
\end{defn}
This shuffle is clearly  reversible, transitive, and has stationary 
distribution equal to the uniform distribution on $S_{n}$, denoted 
$\pi_{n}$. However, the shuffle is not constant on the conjugacy class of transpositions, unlike the previously seen random transposition shuffle. We look to study the total variation mixing time of the family 
of one-sided transposition shuffles $\{\LR_{n}\}_{n\in\mathbb{N}}$. Recall 
that mixing time is defined as follow:
\[\mt(\varepsilon) = \min\{t : \lVert \LR_{n}^{t} - 
\pi_{n} \rVert_{\textnormal{TV}} < \varepsilon\} .\]

Existence of a cutoff at time $t_{n}$ implies that 
$t_n^{\textnormal{mix}}(\varepsilon) \sim t_n$ for all 
$\varepsilon\in(0,1)$ (see Definition \ref{chpt2:def:cutoff}). The main conclusion of our work 
is that the one-sided transposition shuffle exhibits a cutoff at time $t_n 
= n\log n$.
\begin{thm}
	\label{chpt5:thm:cutoff}
	The one-sided transposition shuffle $\RL_{n}$ satisfies the following bounds for any $c_{1}>0$ and $c_{2}>2$:
	\begin{eqnarray}
	\limsup_{n\rightarrow \infty} \,\lVert \RL_{n}^{n\log n+c_{1}n} - 
	\pi_{n} \rVert_{\textnormal{TV}} & \leq & 
	\sqrt 2e^{-c_{1}} \,, \label{eqn:main_UB}\\
	\text{and} \quad \liminf_{n\to\infty}\,\lVert \RL_{n}^{n\log n - 
		n\log\log n - c_{2}n} - 
	\pi_{n}
	\rVert_{\textnormal{TV}} 
	& \geq & 1 - \frac{\pi^{2}}{6(c_{2}-2)^{2}}\,.\label{eqn:main_LB}
	\end{eqnarray}
	Thus, the one-sided transposition shuffle exhibits a cutoff at time $n\log n$ with a window of order $n\log\log n$.
\end{thm}

The lower bound on the total variation distance in \eqref{eqn:main_LB} will be obtained 
via a coupling argument which allows us to compare the one-sided 
transposition shuffle to a variation of a coupon collector's problem. 
To establish the upper bound on total variation distance we make use of the classical $\ell^{2}$ bound given in Theorem \ref{chpt2:thm:classicL2}. In order to use this result we compute the eigenvalues of the one-sided 
transposition shuffle.

To analyse the spectrum of the one-sided 
transposition shuffle we make use of the technique of lifting eigenvectors.	We make several non-trivial changes to the technique presented in Section \ref{chpt4:subsec:rtr} in order to employ it in the analysis of transposition shuffles: we believe that this is the first time such a technique has been shown to be applicable to non-symmetrized shuffles or to a transposition shuffle. In Section \ref{chpt5:sec:algebra} we describe an explicit method for obtaining the eigenvectors of $\LR_{n+1}$ 
from those of $\LR_n$.	The key to our method is to show that each eigenvalue of $\LR_n$ corresponds to a 
standard Young tableau, and may be computed explicitly from the 
entries in the tableau. 
We state here the main 
result which we aim towards with our analysis. 

\begin{thm}
	\label{chpt5:thm:liftingeig}
	The eigenvalues of $\RL_{n}$ are labelled by standard Young tableaux of 
	size $n$, and the eigenvalue represented by a tableau of shape 
	$\lambda$ appears 
	$d_{\lambda}$ times, where $d_\lambda$ is the dimension of $\lambda$.	For a standard Young tableau $T$ of shape $\lambda$ the eigenvalue 
	corresponding to 
	$T$ is given by
	\begin{eqnarray}
	\label{eigenvaluesum}
	\eig(T) = \frac{1}{n}\sum_{\substack{\textnormal{boxes}\\ (i,j)}} 
	\frac{j-i+1}{T(i,j)}\,,
	\end{eqnarray}
	where the sum is performed over all boxes $(i,j)$ in $T$.
\end{thm}

The organisation of this chapter is as follows.
Section \ref{chpt5:sec:algebra} will be dedicated to the proof of Theorem 
\ref{chpt5:thm:liftingeig}. We give full details of the lifting 
procedure, highlighting original contributions to the method which 
allow the recovery of the eigenvalues for the one-sided transposition shuffle. In 
Section~\ref{chpt5:sec:eigenvalues} we first explore some important 
properties of 
the eigenvalues for $\LR_{n}$,  and then  use these to prove the upper 
bound on the mixing time given by Theorem~\ref{chpt5:thm:cutoff}. The 
corresponding lower bound will be proved in Section~\ref{chpt5:sec:lowerbound}, using 
entirely probabilistic arguments. In Section~\ref{chpt5:sec:biased} 
we  consider a generalisation of the one-sided transposition shuffle, in 
which 
$\righthand$ is chosen according to a non-uniform distribution: we show 
that the 
algebraic technique developed for $\RL_{n}$ holds in this more general 
setting, and that the well-known mixing time for the random 
transposition shuffle may be recovered in this way. Finally in Section \ref{chpt5:sec:separation} we show that the unbiased one-sided transposition shuffle exhibits a cutoff in separation distance at $n\log n$, the same time as our cutoff in total variation distance.

\section{Lifting Eigenvectors for Transposition Shuffles}
\label{chpt5:sec:algebra}

In this section we explore the technique of lifting eigenvalues for transposition shuffles. 
Our analysis follows a similar path to Dieker and Saliola's but with several novel changes which allow us to describe every eigenvector of the one-sided transposition shuffle by lifting.
We end this section showing how lifting may be used to recover all the eigenvalues for the random transposition shuffle which were previously stated in Section \ref{chpt4:subsec:rt}.

\subsection{Lifting Eigenvectors for the One-sided Transposition Shuffle}
\label{chpt5:subsec:algebra}

Recall that the permutation modules $M^{\lambda}$ of $\mathfrak{S}_{n}$ are 
formed of a basis of words in the alphabet $[n]$ with length $n$ and evaluation 
$\lambda \vdash n$. The simple modules $S^{\lambda}$ of 
$\mathfrak{S}_{n}$ are 	called Specht modules and labelled by $\lambda$ a 
partition of $n$.  Within each permutation module $M^{\lambda}$ we find one 
copy of the Specht module $S^{\lambda}$. The 
permutation module $M^{(1^{n})}$ is isomorphic to the regular module for 
$\mathfrak{S}_{n}$, and has decomposition,
\begin{eqnarray}
M^{(1^{n})} \cong \bigoplus_{\lambda \vdash n} d_{\lambda} 
S^{\lambda}  \text{ as $\mathfrak{S}_{n}$-modules}. \label{chpt5:eqn:classiciso}
\end{eqnarray}

To model our shuffle $\LR_{n}$ acting on the space $\mathfrak{S}_{n}$ 
we need to turn it into a linear operator. In fact we turn it into an 
element of our group algebra $\mathfrak{S}_{n}$.

\begin{defn}
	\label{RLop}
	Let $n \in \mathbb{N}$.
	The one-sided transposition shuffle on $n$ cards may be viewed as the 
	following element of the group algebra $\mathfrak{S}_n$.
	\[\sum_{1 \leq i \leq j \leq n}\RL_{n}((i\,j))  (i\,j)= \sum_{1 \leq i 
		\leq j 
		\leq n}\frac{1}{nj} (i\,j)  .\]
	To simplify our calculations it is convenient to scale this operator by 
	$n$, so we introduce a new element called the \emph{algebraic one-sided transposition shuffle}: 
	\begin{eqnarray}
	\ALR_n:= \sum_{1 \leq i \leq j \leq n}\frac{1}{j} (i\,j) . 
	\label{chpt5:eqn:algebrashuffle}
	\end{eqnarray}
	
\end{defn}

\begin{lemma}
	\label{chpt5:lem:eigensame}
	Let $\nu$ be a distribution of $S_{n}$, and define $v = \sum_{\eta \in S_{n}} \nu(\eta) \, \eta $ an element of $\mathfrak{S}_{n}$. Then 
	$\nu$ is an eigenvector for $\LR_{n}$ with eigenvalue 
	$\varepsilon$ if and only if $v$ is an eigenvalue for $\ALR_{n}$ with 
	eigenvalue $n \, \varepsilon$.
\end{lemma} 
\begin{proof}
	The proof follows from the same argument as Lemma
	\ref{chpt4:lem:eigensame}.
\end{proof}

By realising the one-sided transposition shuffle as an element of the group algebra we 
can 
concentrate on finding the eigenvalues of $\ALR_{n}$ acting on 
$\mathfrak{S}_{n}$ or equivalently $M^{(1^{n})}$. 
Furthermore, applying equation \eqref{chpt5:eqn:classiciso}, we can reduce the problem 
of finding eigenvalues for the shuffle on $M^{(1^n)}$ to the problem of 
finding eigenvalues belonging to the individual Specht modules $S^\lambda$ comprised 
of the words formed by polytabloids. To lift our eigenvectors from 
$\mathfrak{S}_{n}$-modules $S^{\lambda}$ to $\mathfrak{S}_{n+1}$ modules 
$S^{\lambda+ e_{i}}$ we introduce linear operators on the vector 
space of words $M^{n}$. The one-sided transposition shuffle admits a recursive structure which is seen when we focus on the difference of $\ALR_{n+1}$ and $\ALR_{n}$,
\begin{eqnarray}
\ALR_{n+1} - \ALR_{n} = \frac{1}{(n+1)}\sum_{1\leq i \leq (n+1)} (i 
\hspace{0.2cm} n+1). \label{chpt5:eqn:lrdifference}
\end{eqnarray}

This signifies that the only difference between shuffle $n+1$ and $n$ cards is the movement of the new card in position $n+1$. This relation is key to lifting eigenvectors as it allows us to model the action of $\ALR_{n+1}$ using $\ALR_{n}$ and special linear operators which we now define.

\begin{defn}
	\label{chpt5:def:linoperators}
	We define two linear operators on the spaces spanned by words. To do 
	so, it 	is enough to define the effect on any given word. Let $w = 
	w_{1} \, \ldots \, w_{n} \in W^{n}$.
	\begin{enumerate}
		\item Let $a \in [n+1]$. Define the \emph{adding 
			operator} $\sh_{a}:M^{n} \rightarrow M^{n+1}$ as follows:
		\begin{eqnarray}
		\sh_{a}(w) := w \, a \label{chpt5:eqn:addingoperator}
		\end{eqnarray}
		The adding operator appends the symbol $a$ to the end of the word. 	
		If we think about the adding operator acting on tabloids $\{T\}$ 
		instead of words, then $\sh_{a}(\{T\})$ is the process of adding a box labelled $n+1$ onto row $a$ of the tabloid $\{T\}$. 
		
		\item Let $a,b\in[n]$ Define the \emph{switching operator} 
		$\Theta_{b,a}:M^{n}\rightarrow 	M^{n}$ as follows:
		\begin{eqnarray}
		\Theta_{b,a}(w) := \sum_{\substack{1 \leq k \leq n \\ w_{k} = b}} 
		w_{1}\, \ldots \, w_{k-1} \, a \, w_{k+1}\, \ldots 
		\, w_{n}. \label{chpt5:eqn:switchingoperator}
		\end{eqnarray}
		The operator $\Theta_{b,a}$ forms a linear sum of all words created 
		from $w$ by replacing one occurrence of $b$ with $a$.
	\end{enumerate}
\end{defn}

We previously defined the switching operator in Definition \ref{chpt4:def:shuffleandswitch}. Again we should define the operators for all possible domains separately but to do so would burden us with more notation: it will always be clear from context which domain and codomain we are considering.
The adding operator is our analogue of the shuffling operator (Section \ref{chpt4:subsec:rtr}) and it allow us to lift eigenvectors of the one-sided transposition shuffle. 

Recall that given a partition $\lambda$ we may add boxes onto certain rows to form a new partition.  Since we are allowing ourselves to blur the distinction between $n$-tuples and partitions of $n$, 
if we add a box on row $i$ the new tuple/partition formed is 
$\lambda+e_{i}$. By the end of our analysis we only need the cases 
where $\lambda+e_{i}$ is a partition. 

\begin{lemma}
	\label{chpt5:lem:linearrestriction}
	Given $a \in [n+1]$ and an $n$-tuple $\lambda$ of non-negative integers 
	summing to $n$, we have
	$$
	\sh_{a}: M^{\lambda} \to M^{\lambda+e_a}.
	$$
	In other words the restriction of $\sh_{a}$ to $M^\lambda$ (or $S^{\lambda}$) has image in 
	$M^{\lambda+e_a}$. 	
	
	Given $a,b \in [n]$ and $n$-tuples $\lambda,\mu$ of non-negative 
	integers summing to $n$ with $\lambda+e_{a} = \mu +e_{b}$, we have
	$$ \Theta_{b,a}:M^{\lambda} \rightarrow M^{\mu},$$
	i.e., the restriction of $\Theta_{b,a}$ to $M^\lambda$ (or $S^{\lambda}$) has image in 
	$M^{\mu}$.
\end{lemma}

Our next result 
establishes the crucial equation upon which all the subsequent results in this section rely.  It relates the shuffle on $n$ cards to that on $n+1$ cards, and gives us the basis of lifting eigenvectors. The following theorem is an analogue of \cite[Theorem 38]{dieker2018spectral}.

\begin{thm}
	\label{chpt5:thm:master}
	Given $n \in \mathbb{N}$, acting on words in $M^{n}$ we have
	\begin{equation}
	\ALR_{n+1} \circ \sh_{a} - \sh_{a} \circ \ALR_{n} = 
	\frac{1}{n+1}\sh_{a} + 
	\frac{1}{n+1}\sum_{1 \leq b \leq n} \sh_{b}\circ\Theta_{b,a} \,. 
	\label{mastereq}
	\end{equation}
\end{thm}

\begin{proof}
	It suffices to prove the result on a generic word in $M^{n}$.
	Let $w=w_1 \, \ldots \, w_n$ be a word of length $n$ and let $a 
	\in [n+1]$.
	Consider the two terms on the left hand side applied to $w$:	
	\begin{eqnarray}
	\label{LH1}
	(\ALR_{n+1} \circ \sh_{a})(w) = 
	\frac{1}{n+1}\sum_{\substack{j=n+1 \\ 
			1 \leq i \leq n+1}} (i\,j)\,  (w\, a) + 
	\sum_{1 \leq i \leq j \leq n}\frac{1}{j} (i\,j) \, (w \, a)
	\end{eqnarray}
	
	\begin{eqnarray}
	\label{LH2}
	(\sh_{a} \circ \ALR_{n})(w) = \left(\sum_{1 \leq i \leq j \leq n
	}\frac{1}{j} (i\,j)(w) \right) \, a.
	\end{eqnarray}
	
	The second summation in (\ref{LH1}) cancels with (\ref{LH2}) 
	because the adjoined $a$ is in the $(n+1)$-th place, therefore it never 
	moves and may be brought outside the sum.
	This leaves us with the following:
	\begin{eqnarray}
	\label{leftover}
	(\ALR_{n+1} \circ \sh_{a} - \sh_{a} \circ \ALR_{n})(w) 
	= \frac{1}{n+1}	\sum_{1 \leq i \leq n+1} (i \hspace{0.2cm} n+1)  
	(w\, a).
	\end{eqnarray}

	If $i=n+1$ we move nothing, giving the term $w\, a = \sh_{a}(w)$. 
	Otherwise we apply the transposition 
	$(i \hspace{0.2cm} n+1)$ to $w \, a$.
	This has the same effect as replacing the $i^{\rm th}$ symbol $w_i$ in 
	$w$ with $a$ and then appending $w_i$ on the end of the new word. 		 
	Since we do this for all symbols in $w$, the net effect is the same as 
	$\sum_{1\leq b \leq n} \sh_b \circ \Theta_{b,a}$ applied to $w$. 
	The operator $\Theta_{b,a}$ systematically 
	finds all occurrences of the letter $b$ in $w$ and replaces with an 
	$a$, 	and then $\sh_b$ puts the $b$ back on the end. Since $w \in W^n$,
	all possibilities are exhausted by letting $b$ range over every 
	possible letter $1\leq b \leq n$. 
	This completes the proof. 
\end{proof}

In terms of shuffling cards, we can interpret \eqref{mastereq} as taking 
into 
account the difference between shuffling a deck and then adding a card 
versus 
adding a card and then shuffling.
If we can understand how the operators $\sh_{a}$ and $\Theta_{b,a}$ behave, 
then this inductively tells us how the shuffle on 
$n+1$ cards behaves using information about the shuffle on $n$ cards, 
vastly simplifying our original problem. We now record a key property of 
the linear maps $\Theta_{a,b}$.

\begin{lemma}[See Section 
	2.9 of \cite{sagan2013symmetric}]\label{thetamorphism}
	The switching operators $\Theta_{b,a}$ are $\mathfrak{S}_{n}$-module morphisms. 
\end{lemma}	

\begin{proof}
	This is clear from the definitions, since $S_n$ is acting by place 
	permutations, the operator $\Theta_{b,a}$ commutes with the action of
	$\mathfrak{S}_{n}$. It amounts to the same thing to replace an 
	occurrence of the symbol $b$ with a symbol $a$ and then permute the 
	word as to first permute the word and then replace the same symbol $b$ 
	in its new position with an $a$. 
\end{proof}

The above result is helpful in understand how our adding operators behave when restricted to Specht modules. 

\begin{lemma}[Lemma 44 of \cite{dieker2018spectral}]
	\label{chpt5:lem:restrict}
	Let $\lambda \vdash n$ be such that $\lambda + e_a=\mu+e_b$  for 
	some $a,b\in[n]$. 	Then $\Theta_{b,a}$ is non-zero on $S^{\lambda}$ if and only if 
	$\lambda$ 
	dominates the non-increasing rearrangement of $\mu$. 
	In particular, if $b>a$ then $\Theta_{b,a}(S^{\lambda}) =0$.
\end{lemma}	

\begin{proof}
	Since $S^\lambda$ is simple and $\Theta_{b,a}$ is a module 
	homomorphism, 
	the image  $\Theta_{b,a}(S^\lambda)$ is $0$ or isomorphic to 
	$S^\lambda$, by Schur's lemma.
	But $\Theta_{b,a}(S^\lambda)$ lies in $M^\mu$ because of the 
	relationship $\lambda+e_a = \mu+e_b$. Let $\nu$ be the non-increasing 
	rearrangement of $\mu$, then $\nu\vdash n$, and  $M^\mu \cong 
	M^{\nu}$, so $M^{\mu}$ has a submodule 
	isomorphic to $S^\lambda$ if and only if $\lambda$ dominates $\nu$, 
	see Lemma~\ref{chpt4:lem:youngrule}. Therefore, if $\lambda$ does not dominate $\nu$ we have $\Theta_{b,a}(S^{\lambda}) =0$. Conversely if $\lambda$ dominates $\nu$ then $\Theta_{b,a}(S^{\lambda}) \cong S^{\lambda}$ (see \cite[Section 2.10]{sagan2013symmetric} for further details).
	
	To finish, note that in terms of diagrams the fact that $\lambda+e_a = 
	\mu+e_b$ corresponds to the fact that we can get from the diagram 
	for $\lambda$ to that for $\mu$ by moving a box from row $b$ to row $a$.
	Hence, under the given hypothesis, we have that $\lambda$ dominates the non-increasing rearrangement of 
	$\mu$ if and only if $b\leq a$.

\end{proof}

\begin{example}
	Consider the Specht module $S^{(3,1)}$. We have seen previously that 
	this is spanned by the following polytabloids represented as words
	\[S^{(3,1)} = \langle 1112 -2111, \, 1121 -2111, \, 1211 - 2111 
	\rangle. \]
	If we apply $\Theta_{2,1}$ we go from partition $(3,1)$ to the 
	partition $(2,2)$. Lemma \ref{chpt5:lem:restrict} tells us that all the 
	elements of $S^{(3,1)}$ belong to the kernel of $\Theta_{2,1}$. 
	Verifying this result we find,
	\begin{eqnarray*}
		\Theta_{2,1} (1112 - 2111) & = & 1111 -1111 = 0 \\
		\Theta_{2,1} (1211 - 2111) & = & 0 \\
		\Theta_{2,1} (1121 - 2111) & = & 0 	
	\end{eqnarray*}
	Instead if we apply $\Theta_{1,2}$ we find a non-zero elements belonging 
	to the module $M^{(2,2)}$. 
	\begin{eqnarray*}
		\Theta_{1,2} (1112 - 2111) & = & 1212 + 1122 -2211 - 2121 \\
		\Theta_{1,2} (1211 - 2111) & = & 1212 +1221 -2121 - 2112 \\
		\Theta_{1,2} (1121 - 2111) & = & 1221 + 1122 - 2211 - 2121. 	
	\end{eqnarray*}

\end{example}

The preceding result shows that when we restrict equation 
\eqref{mastereq} to a Specht module $S^\lambda$ we can change the index of 
the 
summation in the final term on the right hand side, as follows.

\begin{corollary}[Corollary 45 of \cite{dieker2018spectral}]
	\label{cor:spechtrestrict}
	\begin{align}
	(\ALR_{n+1} \circ \sh_{a} - \sh_{a} \circ \ALR_{n})|_{S^{\lambda}} = 
	\frac{1}{n+1}\sh_{a}|_{S^{\lambda}} + 
	\frac{1}{n+1}\sum_{1 \leq b \leq a} 
	\sh_{b}\circ\Theta_{b,a}|_{S^{\lambda}} \label{eqn:restricted}
	\end{align}
\end{corollary}

Having restricted equation (\ref{mastereq}) to the Specht module 
$S^{\lambda}$, we analyse the image in the module $M^{\lambda +e_{a}}$: 
note that it is clear from the left hand side of \eqref{mastereq} that we 
land in $M^{\lambda+e_a}$. The following lemma and its proof follow Lemma 41 of Deiker and Saliola \cite{dieker2018spectral}.

\begin{lemma}[Lemma 41 of \cite{dieker2018spectral}]
	\label{chpt5:lem:lives}
	Suppose $\lambda\vdash n$.
	Then the subspace $\sh_{a}(S^{\lambda})$ is contained in an 
	$\mathfrak{S}_{n+1}$ submodule of $M^{\lambda+e_{a}}$ that is 
	isomorphic to 
	$\oplus_{\mu} S^{\mu}$, where the sum ranges over the partitions $\mu$ 
	obtained 
	from $\lambda$ by adding a box in row $i$ for $i\leq a$.
\end{lemma}

\begin{proof}
	Let $w$ be a word of length $n$, so that $\sh_{a}(w) = w \, a$.
	If the symbol $b$ does not occur in $w$ then
	\[\sh_{a}(w) = \Theta_{b,a}(\sh_{b}(w)).\]
	Let $b =l(\lambda) +1$, so $b$ does not appear in any $w \in 
	M^{\lambda}$, and 
	consider the $\mathfrak{S}_{n+1}$-submodule $N$ of 
	$M^{\lambda + e_{b}}$ generated by the elements $x\, b$ with $x  \in 
	S^{\lambda}$,
	\[N = \langle x \, b : x \in S^{\lambda} \rangle.\]

	The submodule $N$ is isomorphic to 
	$\text{Ind}^{\mathfrak{S}_{n+1}}_{\mathfrak{S}_{n} \times 
		\mathfrak{S}_{1}} 
	(S^{\lambda} \otimes S^{1})$ (this is essentially the definition of how 
	to induce), and using the branching 
	rules on $S_{n}$ this decomposes as a 
	multiplicity free 
	direct sum of Specht modules $S^{\mu}$, where $\mu \vdash n+1$ and 
	$\lambda 
	\subset \mu$ (see \cite[Theorem I.7]{geissinger1978representations}).
	Using the observation at the start of the proof, we obtain
	\[\sh_{a}(S^{\lambda}) = \Theta_{b,a}(\sh_{b}(S^{\lambda})) \subseteq 
	\Theta_{b,a}(\langle \sh_{b}(S^{\lambda}) \rangle) = \Theta_{b,a}(N) 
	\cong \bigoplus_{\substack{\mu \vdash n+1 \\ \lambda \subset \mu}} \Theta_{b,a}(S^{\mu}).\]
	Now note that $\Theta_{b,a}$ sends any word with evaluation 
	$\lambda+e_b$ to a word with evaluation $\lambda+e_a$, and hence
	$\Theta_{b,a}(M^{\lambda +e_{b}}) \subseteq M^{\lambda +e_{a}}$. 
	It follows that all nonzero summands $S^\mu$ appearing on the right 
	hand side occur for $\mu\vdash n+1$ dominating the non-decreasing rearrangement of $\lambda+e_a$, and then
	by Lemma \ref{chpt5:lem:restrict} we can conclude that $\mu$ is obtained from 
	$\lambda$ by adding a cell in row $i$ with $i \leq a$, as required.
\end{proof}

Lemma \ref{chpt5:lem:restrict} tells us a great deal about the image of $\sh_{a}(S^{\lambda})$. Specifically we know that 
$\sh_{a}(S^{\lambda})$ contains exactly one submodule isomorphic to 
$S^{\lambda+e_{a}}$. This means that after applying our adding operator we may project down onto the required simple module to get a linear operator from $S^{\lambda} \to S^{\lambda+e_{a}}$: these will become our lifting operators. To project onto the required simple module we introduce the isotypic projections.

\begin{defn}
	Let $G$ be a finite group, $V$ be a simple $G-$module with 
	corresponding character 
	$\chi_{V}$. Let $M$ be a module for $G$, we may define the 
	\emph{isotypic projection} $\pi^{V}$ from $M$ onto its unique defined $V$-isotypic 
	component (see Definition \ref{chpt3:def:isocompotent}) in the 
	following way:
	\[\pi^{V}:M \rightarrow M, \hspace{0.2cm} \pi^{V}(m)  = \left( 
	\frac{\textnormal{dim}(V)}{|G|}\sum_{g\in G} \overline{\chi_{V}(g)} g 
	\right) m.\] 
	Note that the isotypic projection is given by the action of the element 
	$\frac{\textnormal{dim}(V)}{|G|}\sum_{g\in G} \overline{\chi_{V}(g)} g$ 
	which belongs to the centre of the group algebra $\mathbb{C}[G]$. Thus, 
	$\pi_{V}$ is a morphism of $G-$modules and commutes with any $G-$module 
	endomorphism.
\end{defn}

Using these projections, we can now define our lifting operators, which will beproven to map eigenvectors of $\ALR_{n}$ 
to those of $\ALR_{n+1}$.

\begin{defn}
	\label{chpt5:def:liftingmaps}
	Suppose $\lambda \vdash n$ and $\lambda + e_{a} = \mu \vdash n+1$ are 
	two 
	partitions. 	
	Define the \emph{lifting operator} 
	$$
	\pro_{a}^{\lambda,\mu} := \iso^{\mu} \circ \sh_{a}:S^{\lambda} 
	\rightarrow S^{\mu}\subseteq M^\mu.
	$$ 
	Note that since $\sh_a(S^\lambda)\subseteq M^\mu$ and $M^\mu$ contains 
	a unique copy of $S^\mu$, therefore the $S^{\mu}$ isotypic component in $M^{\mu}$ is just $S^{\mu}$ and $\pro_a^{\lambda,\mu}(S^{\lambda})$ is actually contained in $S^\mu$.
\end{defn}

We next prove that our lifting operators $\pro_{a}^{\lambda,\mu}$ are injective $\mathfrak{S}_{n}$-module morphisms. This differs from the lifting operators $\mathcal{L}_{i}^{\lambda}$ for the random to random shuffle which we have seen are not injective. The difference in lifting operators depends in an essential way on our adding operator $\sh_{a}$. 	Using the lifting operators $\pro_{a}^{\lambda,\mu}$ we are able 
to find \emph{all} the eigenvectors for a module $S^{\mu}$ by  lifted eigenvectors from partitions $\lambda \subset \mu$. 

\begin{corollary}
	\label{nonzero}
	For any $\lambda \vdash n$ and $\lambda+e_{a} \vdash n+1$, there exists 
	some $v \in S^{\lambda}$ with
	\[\pro^{\lambda,\lambda+e_{a}}_{a}(v) \neq 0. \]
\end{corollary}

\begin{proof}
	If $\pro^{\lambda,\lambda+e_{a}}_{a}(S^\lambda) = \{0\}$, then the image 
	$\sh_a(S^\lambda)$ 
	lies in the kernel of the projection 
	$\iso^{\lambda+e_a}:M^{\lambda+e_a}\to S^{\lambda+e_a}$,
	which is an $\mathfrak{S}_{n+1}$-submodule with no component equal to 
	$S^{\lambda+e_a}$.
	Hence the submodule generated by $\sh_a(S^\lambda)$ has no component 
	equal to $S^{\lambda+e_a}$.
	But we previously observed that (with notation as in the proof of Lemma 
	\ref{chpt5:lem:lives})
	\[\langle \sh_{a}(S^{\lambda}) \rangle =  \langle\Theta_{b,a}( 
	\sh_{b}(S^{\lambda}))\rangle=\Theta_{b,a}( 
	\langle\sh_{b}(S^{\lambda})\rangle ) \cong \Theta_{b,a}(N) 
	\cong \bigoplus_{1\leq i \leq a} S^{\lambda + e_{i}}.\]
	Since the right hand side contains $S^{\lambda+e_a}$ as a summand, we 
	have a contradiction.
\end{proof}

We already know the map $\iso^{\mu}$ is an $\mathfrak{S}_{n+1}$ module 
morphism. Let us realise $\mathfrak{S}_{n}$ inside $\mathfrak{S}_{n+1}$ as 
the stabilizer of the $(n+1)^{\rm th}$ position.
Then any $\mathfrak{S}_{n+1}$-module gives rise to an 
$\mathfrak{S}_n$-module by restriction (see Definition \ref{chpt4:def:restrictandinduce}). Hence, $\iso^{\mu}$ is a $\mathfrak{S}_{n}$-module morphism.

\begin{lemma}
	\label{modulemorph}
	The linear operator $\pro_{a}^{\lambda,\lambda+e_{a}}$ is a 
	$\mathfrak{S}_{n}$-module morphism with trivial kernel. Therefore, the lifting operators are injective.

\end{lemma}

\begin{proof}
	Our key observation is that our adding operator commutes with elements of $\mathfrak{S}_{n}$ inside of $\mathfrak{S}_{n+1}$, i.e., $\sh_a(\sigma(v)) = \sigma(\sh_a(v))$ for 
	all $v \in S^{\lambda}$ and $\sigma \in 
	\mathfrak{S}_{n} \subset \mathfrak{S}_{n+1}$. 
	This is obvious, since $\sh_a$ adds an element in the final position 
	which is not affected by $\sigma$. 
	Hence $\pro_a^{\lambda,\lambda+e_a}$ is the composition of two 
	$\mathfrak{S}_n$-module morphisms.
	The final observation follows from Corollary \ref{nonzero} -- since 
	$\pro_a^{\lambda,\lambda+e_a}$
	is a nonzero module morphism with a simple module as its domain, it 
	must be injective by Schur's lemma (Lemma \ref{chpt3:lem:schur}).
\end{proof}

The lifting operators $\pro^{\lambda, \lambda + e_{a}}_{a}$ 
being injective is a key point which simplifies our analysis  
compared to that of the random to random shuffle - in \cite{dieker2018spectral} the lifting operators can kill 
eigenvectors. The next results 
show that $\pro_{a}^{\lambda,\lambda+e_{a}}$ lifts eigenvectors 
of 
$\ALR_{n}$ into those of $\ALR_{n+1}$. 
To establish this we apply our projection $\iso^{\lambda+e_{a}}$ to 	equation~(\ref{eqn:restricted}). Before we do this we state an identity between our adding and switching operators which is critical to the proofs which follow.

\begin{lemma}
	Take $a,b \in [n]$. Then our adding and switching operators satisfy the following equality
	\begin{eqnarray}
	\sh_{b} \circ \Theta_{b,a} = \Theta_{b,a} \circ \sh_{b}  -\sh_{a} .\label{chpt5:eqn:addswtichequally}
	\end{eqnarray}
\end{lemma}
\begin{proof}
	Let $w$ be a word of length $n$. Let $v := \sh_{b}(w)$ so that $v_{n+1} =b$. Consider the affect of switching operator $\Theta_{b,a}$ on the word $v$ of length $n+1$. We find,
	\begin{eqnarray*}
		\Theta_{b,a} \circ \sh_{b}(w)  & =  &\sum_{\substack{1\leq k \leq n+1 \\ v_{k}=b}} v_{1} \ldots v_{k-1} \, a\, v_{k+1} \ldots v_{n+1} \\
		& = & \sh_{a}(w) + \sum_{\substack{1\leq k \leq n \\ v_{k}=b}} v_{1} \ldots v_{k-1} \, a\, v_{k+1} \ldots v_{n} \,b  \; = \;\sh_{a}(w) + \sh_{b}\circ\Theta_{b,a}(w)
	\end{eqnarray*}
	The second equality follows from taking the $k=n+1$ out of the sum. The last equality follows from the fact that $v_{i} = w_{i}$ for all $i \in[n]$. Rearranging this final expression gives us the desired equality.
\end{proof}

We can now state our versions of \cite[Lemma 
48, Theorem 49]{dieker2018spectral};
the proofs follow \emph{mutatis mutandis} from the ones given there (the 
changes 
needed are to the constants in equation (\ref{mastereq})).

\begin{lemma}[Lemma 48 of \cite{dieker2018spectral}]
	\label{chpt5:lem:lifting}
	Let $\lambda \vdash n$, and  $a \in \{1,2,\ldots,l(\lambda)+1\}$. Take $i\in [n]$ such that  $1 \leq i \leq a$ and set $\mu = \lambda +e_{i}$. Then,

	\[	\ALR_{n+1} \circ \pro_{a}^{\lambda,\mu} - 
	\pro_{a}^{\lambda,\mu} \circ \ALR_{n} = \frac{2 
		+\lambda_{a} -a}{n+1}\pro_{a}^{\lambda,\mu} + 
	\frac{1}{n+1}\sum_{i \leq b < a} \Theta_{b,a} \circ 
	\pro_{b}^{\lambda,\mu}.\]
	
\end{lemma}

\begin{proof}
	This follows from the work in \cite{dieker2018spectral} because we 
	have not changed the switching operators $\Theta_{b,a}$.  The values on 
	the right hand side change to reflect our adding operators and our new equation \eqref{mastereq}. We present the modified proof in full for completeness.
	
	Continuing from Corollary \ref{cor:spechtrestrict} we know
	
	\[ \ALR_{n+1} \circ \sh_{a} - \sh_{a} \circ \ALR_{n}\big|_{S^{\lambda}} 
	= 
	\frac{1}{n+1}\sh_{a}\big|_{S^{\lambda}} + \frac{1}{n+1}\sum_{1 \leq b 
		\leq a} \sh_{b}\circ\Theta_{b,a}\big|_{S^{\lambda}}.\]
	
	Apply the isotypic projection $\iso^{\mu}$ to both sides of 
	the equation.	Since $\ALR_{n+1}$ is given by the action of an element 
	of 
	the group algebra of $\mathfrak{S}_{n+1}$ and $\iso^{\mu}$ is an 
	$\mathfrak{S}_{n+1}$-module morphism, these operators commute and so we 
	have
	\begin{eqnarray}
	\ALR_{n+1} \circ \pro_{a}^{\lambda,\mu} - 
	\pro_{a}^{\lambda,\mu} \circ \ALR_{n} = 
	\frac{1}{n+1}\pro_{a}^{\lambda,\mu} +  
	\frac{1}{n+1}\sum_{1 \leq b \leq a}(\iso^{\mu} 
	\circ \sh_{b}\circ\Theta_{b,a})\big|_{S^{\lambda}}. 
	\label{eqn:proj1}
	\end{eqnarray}
	Our adding operator $\sh_{a}$ satisfies the equation
	\eqref{chpt5:eqn:addswtichequally} 
	and so
	\[\iso^{\mu} \circ \left(\sh_{b}\circ\Theta_{b,a} \right) = \iso^{\mu} 
	\circ \left(\Theta_{b,a} \circ \sh_{b} -\sh_{a} \right)
	= \Theta_{b,a} \circ \iso^{\mu} \circ \sh_{b} - 
	\iso^{\mu} \circ \sh_{a} .\]
	The right side side of the equation \eqref{eqn:proj1} now becomes
	\begin{eqnarray}
	\frac{1-a}{n+1}\pro_{a}^{\lambda,\mu} +  
	\frac{1}{n+1}\sum_{1 \leq b \leq a} \Theta_{b,a} \circ 
	\pro_{b}^{\lambda,\mu} . \label{eqn:proj2}
	\end{eqnarray}
	Notice that if $b=a$ then 
	$\Theta_{a,a}(w)$ acts as a scalar by the number of occurrences of the 
	symbol $a$ in $w$. In our case all words in 
	$\pro_{a}^{\lambda,\mu}(S^{\lambda})$ contain $\lambda_{a}+1$ 
	occurrences of $a$. Finally if $b<i$ we know that $\sh_{b}(S^{\lambda})$ does not contain the module $S^{\lambda+e_{i}}$ (by Lemma \ref{chpt5:lem:lives}), so $\pro_{b}^{\lambda,\mu}=0$. Thus, equation \eqref{eqn:proj2} is equal to
	\[	\frac{2+\lambda_{a}-a}{n+1}\pro_{a}^{\lambda,\mu} +  
	\frac{1}{n+1}\sum_{i \leq b < a} \Theta_{b,a} \circ 
	\pro_{b}^{\lambda,\mu} . \label{eqn:proj3} \qedhere \]
\end{proof}

\begin{thm}[Theorem 49 of \cite{dieker2018spectral}]
	\label{chpt5:thm:lift}

	Let $\lambda \vdash n$, and  $a \in \{1,2,\ldots,l(\lambda)+1\}$. Take $i\in [n]$ such that  $1 \leq i \leq a$ and set $\mu = \lambda +e_{i}$. Then,
	\begin{eqnarray}
	\ALR_{n+1} \circ \pro_{a}^{\lambda,\mu} - 
	\pro_{a}^{\lambda,\mu} \circ \ALR_{n} = \frac{(2 + 
		\lambda_{i} -i)}{n+1}\pro_{a}^{\lambda,\mu} .
	\end{eqnarray}
	
	In particular if $v \in S^{\lambda}$ is an eigenvector of $\ALR_{n}$ 
	with 
	eigenvalue $\varepsilon$, then  
	$\pro_{a}^{\lambda,\mu}(v)$ is an eigenvector of 
	$\ALR_{n+1}$ belonging to $S^{\mu}$ with eigenvalue 
	\begin{eqnarray}
	\varepsilon + \frac{2 +\lambda_{i} -i }{n+1}.
	\end{eqnarray}
	
\end{thm}

\begin{proof}
	This proof follows from the work in \cite{dieker2018spectral} with 
	minor changes to reflect the one-sided transposition shuffle. For $i=a$ this follows from Lemma \ref{chpt5:lem:lifting}, this will be the key	to showing 	it 	holds in all cases.
	Let $\mu = \lambda + e_{i}$, we know from Lemma \ref{chpt5:lem:lifting} that lifting via $\pro_{i}^{\lambda,\mu}$ gives us
	\begin{eqnarray}
	\ALR_{n+1} \circ \pro_{i}^{\lambda,\mu} - 
	\pro_{i}^{\lambda,\mu} \circ \ALR_{n} = \frac{2+\lambda_{i} 
		-i}{n+1}\pro_{i}^{\lambda,\mu}  \label{chpt5:eqn:thmlift1}.
	\end{eqnarray}
	Applying the linear operator $\Theta_{i,a}$ to the above equation,
	\begin{eqnarray}
	\ALR_{n+1} \circ \Theta_{i,a} \circ 
	\pro_{i}^{\lambda,\mu} 
	- \Theta_{i,a} \circ \pro_{i}^{\lambda,\mu} \circ \ALR_{n} 
	= 
	\frac{2 + \lambda_{i}-i}{n+1}\Theta_{i,a} \circ 
	\pro_{i}^{\lambda,\mu}  \label{chpt5:eqn:proj4}.
	\end{eqnarray}
	Consider the left hand side of \eqref{chpt5:eqn:proj4} break up the lifting operator into 
	$\pro_{i}^{\lambda,\mu} = \pi^{\mu} \circ \sh_{i}$. The projection $\pi^{\mu}$ commutes with the $\mathfrak{S}_{n+1}$ module morphism $\Theta_{i,a}$. Performing this we obtain the equation below, restricted to $S^{\lambda}$;
	\begin{eqnarray}
	& & \ALR_{n+1} \circ \Theta_{i,a} \circ \pro_{i}^{\lambda,\mu}		
	- \Theta_{i,a} \circ \pro_{i}^{\lambda,\mu} \circ 		
	\ALR_{n} \nonumber\\
	& = & 	
	\ALR_{n+1}  \circ \iso^{\mu} \circ \Theta_{i,a} \circ 
	\sh_{i} -  \iso^{\mu} \circ \Theta_{i,a}  \circ \sh_{i} 
	\circ \ALR_{n} \label{chpt5:eqn:1}
	\end{eqnarray}
	From the identity $\Theta_{i,a} \circ \sh_{i} = \sh_{a} + \sh_{i} 
	\circ 
	\Theta_{i,a}$ (which is a rearrangement of equation \eqref{chpt5:eqn:addswtichequally}), the equation \eqref{chpt5:eqn:1} becomes 		
	\begin{eqnarray*}
		(\ALR_{n+1}  \circ \iso^{\mu} \circ  \sh_{a}  + 
		\ALR_{n+1}  \circ \iso^{\mu} \circ  \sh_{i} \circ 
		\Theta_{i,a})
		-(\iso^{\mu} \circ \sh_{a} \circ \ALR_{n} + 
		\iso^{\mu} \circ \sh_{i} \circ \Theta_{i,a}  \circ 
		\ALR_{n}) 
	\end{eqnarray*}
	Our one-sided transposition shuffle $\ALR_{n}$ commutes with $\Theta_{i,a}$, swapping these around and using our observation \eqref{chpt5:eqn:thmlift1}, we reduce the above equation to:
	\begin{eqnarray}
	& &	(\ALR_{n+1}  \circ \iso^{\mu} \circ \sh_{a}  - 
	\iso^{\mu} \circ \sh_{a} \circ \ALR_{n} )
	+(\ALR_{n+1}  \circ \iso^{\mu} \circ \sh_{i} - 
	\iso^{\mu} \circ \sh_{i} \circ \ALR_{n}) \circ 
	\Theta_{i,a}. \nonumber \\
	& = &  (\ALR_{n+1} \circ \pro_{a}^{\lambda,\mu}  - 
	\pro_{a}^{\lambda,\mu} \circ \ALR_{n} )
	+(\ALR_{n+1}  \circ \pro_{i}^{\lambda,\mu} - 
	\pro_{i}^{\lambda,\mu} \circ \ALR_{n}) \circ 
	\Theta_{i,a}. \nonumber \\
	& = &  (\ALR_{n+1} \circ \pro_{a}^{\lambda,\mu}  - 
	\pro_{a}^{\lambda,\mu} \circ \ALR_{n} )  +  \frac{2+\lambda_{i} 
		-i}{n+1}\pro_{i}^{\lambda,\mu} \circ \Theta_{i,a} .\label{chpt5:eqn:2} 
	\end{eqnarray}
	We are finished manipulating the left hand side of equation \eqref{chpt5:eqn:proj4}, and now focus our attention on the right hand side. Taking the right hand side we split the lifting operator into  $\pro_{i}^{\lambda,\mu} = \pi^{\mu} \circ \sh_{i}$,  and use the properties of $\Theta_{i,a}$ and equation \eqref{chpt5:eqn:addswtichequally} in order to find the following: 
	\begin{eqnarray}
	& & \frac{2 + \lambda_{i}-i}{n+1}\Theta_{i,a} \circ 
	\pro_{i}^{\lambda,\mu} \nonumber\\
	& = &
	\frac{2 +\lambda_{i}-i}{n+1}\iso^{\mu} \circ \Theta_{i,a} \circ \sh_{i} |_{S^{\lambda}}  
	=  \frac{2 +\lambda_{i}-i}{n+1} \iso^{\mu} \circ \sh_{a} |_{S^{\lambda}} + \frac{2 +\lambda_{i}-i}{n+1} \iso^{\mu} \circ \sh_{i} \circ \Theta_{i,a} |_{S^{\lambda}} \nonumber \\
	& = & \frac{2 +\lambda_{i}-i}{n+1} \pro_{a}^{\lambda,\mu} + \frac{2 +\lambda_{i}-i}{n+1} \pro_{i}^{\lambda,\mu} \circ \Theta_{i,a} \label{chpt5:eqn:thmlift5} .
	\end{eqnarray}
	Combining equations \eqref{chpt5:eqn:2} and \eqref{chpt5:eqn:thmlift5} gives the desired result.
\end{proof}

The last theorem tells us exactly how to turn eigenvectors of $\ALR_{n}$ 
into 
those of $\ALR_{n+1}$ and critically it shows how the eigenvalues 
change in value.
An important observation here is the numerator of the change $\frac{2+\lambda_{i} -i}{n+1}$ depends on what box (or letter) we are adding and the denominator depends on the step of our lifting, thus they are independent of one another. The final part of our analysis rests on showing that \emph{all} of the 
eigenvectors in a Specht module $S^{\mu}$ can be retrieved
by lifting from Specht modules $S^\lambda$ with $\mu = 
\lambda+e_a$. 
In fact, we show that these lifted eigenvectors form a basis of $S^\mu$.

\begin{thm}
	\label{mapsnew}
	For any $\mu \vdash n+1$ we may find a basis of eigenvectors of 
	$\ALR_{n+1}$ 
	for the module $S^{\mu}$, 
	by lifting the eigenvectors of $\ALR_{n}$ belonging to the modules 
	$S^{\lambda}$ with $\lambda\vdash n$ and $\lambda \subset \mu$.
\end{thm}

\begin{proof}
	
	We proceed by induction. For $n=1$ we know that the simple modules 
	$S^{(2)}, S^{(1,1)}$ of $\mathfrak{S}_{2}$ 
	are both one dimensional. Therefore, the eigenvector $1 \in S^{(1)}$ 
	when lifted indeed forms a basis for each simple module, this was demonstrated in Example \ref{chpt5:ex:difflift}.
	
	Consider the simple module of $S^{\mu}$ with $\mu \vdash n+1$. We know  from the classic branching rules of $S_{n}$ (Theorem \ref{chpt4:thm:branching}) that the restriction 
	of this module to $\mathfrak{S}_{n}$ is 
	given by 
	\[\textnormal{Res}^{\mathfrak{S}_{n+1}}_{\mathfrak{S}_{n}}(S^{\mu})\cong
	\bigoplus_{\substack{\lambda
			\vdash n \\\lambda \subset \mu}} S^{\lambda}.\]
	Importantly $S^\mu$ is isomorphic as a vector space to $\oplus_{\substack{\lambda
			\vdash n \\\lambda \subset \mu}} S^{\lambda}$.
	Now suppose we have a basis of eigenvectors for every $S^{\lambda}$. By 
	Lemma~\ref{modulemorph} the map $\pro^{\lambda,\mu}(S^{\lambda})$ 
	gives a basis for the submodule $S^{\lambda}$ inside of the vector space of 
	$\textnormal{Res}^{\mathfrak{S}_{n+1}}_{\mathfrak{S}_{n}}(S^{\mu})$. Hence,
	considering all of the lifted eigenvectors from every $S^{\lambda}$ 
	together we find a basis for $S^{\mu}$. By Theorem~\ref{chpt5:thm:lift} the lifted eigenvectors form a basis of 
	eigenvectors for $S^{\mu}$.
\end{proof}

Inductively, for any $\lambda\vdash n$, Theorem \ref{mapsnew} gives us the 
way to find all 
the eigenvectors for $\ALR_{n}$ belonging to the Specht 
module  $S^{\lambda}$: starting at $S^{\emptyset}$ and recursively applying 
lifting 
operators until we 
reach $S^\lambda$ gives us an eigenvector, and all eigenvectors arise 
in 
this way. 
Note that $S^{\emptyset}$ has no eigenvectors attached to it, but we allow the empty word
$\omega$ to be an eigenvector with eigenvalue $0$, and 
$\sh_{a}(\omega) 
=a$. This agrees with the formula in Theorem~\ref{chpt5:thm:lift} because $a$ 
is the only 
eigenvector of $\ALR_{1}$ with eigenvalue $1 = 0 + (2+0-1)/(1)$.
The inductive process of lifting naturally forms one path up Young's 
lattice which starts at $\emptyset$ and ends at $\lambda$. Furthermore, by 
Theorem 
\ref{mapsnew} each unique path we take  $\emptyset \rightarrow \lambda$ 
results 
in a unique eigenvector for $S^{\lambda}$, and all these eigenvectors 
together form  a basis. We now are in a position to prove Theorem \ref{chpt5:thm:liftingeig}.

\begin{proof}[Proof of Theorem \ref{chpt5:thm:liftingeig}]
	Every eigenvector in our constructed basis gives a distinct eigenvalue 
	of $S^{\lambda}$, hence there are $d_{\lambda}$ distinct eigenvalues. 
	These are eigenvalues for the shuffle $\ALR_{n}$, and each one appears 
	$d_{\lambda}$ times due to the isomorphism in 
	equation~\eqref{chpt5:eqn:classiciso}. 
	Overall we have found $\sum_{\lambda \vdash n} 
	d_{\lambda}^{2} = n!$	eigenvalues and thus have a complete set.	Given a standard tableau $T \in \SYT(\lambda)$, we build up the tableau following its 
	labelling and keeping 
	track of the changes in eigenvalue given by Theorem~\ref{chpt5:thm:lift}. 
	When box 
	$(i,j)$ is added to $T$ we get a change in eigenvalue of 
	$\frac{2+\lambda_{i}-i}{n+1} = 	\frac{2 + (j-1) -i}{T(i,j)} = \frac{j-i+1}{T(i,j)}$. After summing these	changes for all boxes $(i,j)$ in $T$ we divide by $n$ to normalise the eigenvalue, recovering equation \ref{eigenvaluesum}.
\end{proof}

We have given an explicit description of how to compute the eigenvalues of the shuffle $\ALR_{n}$. We now state a description of maps $\pro_{a}^{\lambda,\lambda+e_{a}}$ in terms of adding and switching operators, these allow explicit computation of any eigenvector if required. The proof of Theorem \ref{chpt5:lem:liftingoperatorsdescription} follows from directly from the work of \cite[Section 5.6]{dieker2018spectral}.

\begin{lemma}[Theorem 21 \cite{dieker2018spectral}]
	\label{chpt5:lem:liftingoperatorsdescription}
	Let $\lambda \vdash n$. The lifting maps as defined in Definition 
	\ref{chpt5:def:liftingmaps} are a linear combination of the shuffling and switching operators, explicitly they are:	
	\begin{eqnarray}
	\pro_{i}^{\lambda,\lambda+e_{i}} = \sum_{1\leq b_{1} <\ldots <b_{m} <b_{m+1}=i} \left( \prod_{j=1}^{m} \frac{1}{(\lambda_{i} -i) - (\lambda_{b_{j}} -b_{j})}\Theta_{b_{j},b_{j+1}} \right) \sh_{b_{1}}.	
	\end{eqnarray}
	For example, the first three lifting operators are:
	\begin{eqnarray*}
		\pro_{1}^{\lambda,\lambda+e_{1}} & = & \sh_{1} \\
		\pro_{2}^{\lambda,\lambda+e_{2}} & = &\sh_{2} + \frac{1}{(\lambda_{2}-2) - (\lambda_{1}-1)}\Theta_{1,2} \sh_{1} \\
		\pro_{3}^{\lambda,\lambda+e_{3}} & = & \sh_{3} + \frac{1}{(\lambda_{3} - 3) - (\lambda_{2}-2)} \Theta_{2,3} \circ\sh_{2} + \frac{1}{(\lambda_{3} - 3) - (\lambda_{1}-1)} \Theta_{1,3} \circ \sh_{1} \\
		& + &\frac{1}{((\lambda_{3} - 3) - (\lambda_{2}-2))((\lambda_{3} - 3) - (\lambda_{1}-1))}\Theta_{2,3} \circ \Theta_{1,2} \circ \sh_{1}
	\end{eqnarray*}
\end{lemma}

\begin{example}
	\label{chpt5:ex:difflift}
	Consider the Specht module $S^{(1)} = \langle 1 \rangle$. The 
	eigenvector of $\ALR_{1}$ belonging to this module is $1$ with 
	eigenvector $1$. Using our lifting operators $\pro_{1}^{(1),(2)}$ and $\pro_{2}^{(1,1)}$ we find:
	\begin{eqnarray*}
		\pro_{1}^{(1),(2)}(1) & = &\sh_{1} (1)  = 11 \\ 
		\pro_{2}^{(1,1)}(1) & = & \left(\sh_{2} -\frac{1}{2}\Theta_{1,2}\sh_{1}\right)(1) = 
		\frac{1}{2}\cdot(12- 21) 
	\end{eqnarray*}
	The above elements belong to the modules $S^{(2)}$ and $S^{(1,1)}$ respectively. Furthermore  $11$ is an eigenvector of $\ALR_{2}$ with eigenvalue $2$, and $12-21$ is an eigenvector of $\ALR_{2}$ with eigenvalue $1$.

\end{example}

For the analysis of the mixing time of the one-sided transposition shuffle we only use the eigenvalues computed in Theorem \ref{chpt5:thm:liftingeig}.
To end this section we give an explicit example of computing an eigenvalue using the process described in the proof Theorem \ref{chpt5:thm:liftingeig}.

\begin{example}
	\label{runningexample}
	Let $\lambda = (4,2,1)$ and choose a standard 
	Young tableau
	\[T_{1} = 
	\ytableausetup{mathmode,baseline,aligntableaux=center,boxsize=1.2em} 
	\begin{ytableau} 1 & 3 & 6 & 7\\ 2 & 4\\5 \end{ytableau}.\]
	We build up $T_{1}$ using our lifting maps, at each step Theorem 
	\ref{chpt5:thm:lift} tells us the change in our eigenvalue. We 
	start at $n=0$ with eigenvector $\omega \in S^{\emptyset}$ with eigenvalue $0$. Applying 
	the lifting map $\pro_{1}^{\emptyset, \emptyset + e_{1}}$ gives us:
	\[\pro_{1}^{\emptyset,(1)}(\omega) =1 \in S^{(1)}
	\textnormal{ with the eigenvalue } 0 + 
	\frac{2 + 0 -1}{1} = 1.\]
	On the second row, this corresponds to applying $\pro_{2}^{(1), (1)+e_{2}}$, which gives:
	\[	\pro_{2}^{(1),(1,1)}\pro_{1}^{\emptyset,(1)} (\omega) = \pro_{2}^{(1),(1,1)}(1) = \frac{1}{2}(12-21)
	\textnormal{ with eigenvalue } 1 + 
	\frac{2 + 0 -2}{2} = 1.\]
	Continuing the procedure in the proof of Lemma \ref{chpt5:thm:liftingeig} we build up $T_{1}$ keeping track of the change in 	eigenvalue at each step.
	\[
	\ytableausetup{mathmode,baseline,aligntableaux=top,boxsize=1.2em}
	\emptyset \xrightarrow{\frac{1}{1}}	
	\begin{ytableau}
	1
	\end{ytableau} \xrightarrow{\frac{0}{2}}	
	\begin{ytableau}
	1  \\2
	\end{ytableau} \xrightarrow{\frac{2}{3}}	
	\begin{ytableau}
	1 & 3\\
	2 
	\end{ytableau} \xrightarrow{\frac{1}{4}}	
	\begin{ytableau}
	1 & 3\\
	2 & 4
	\end{ytableau} \xrightarrow{\frac{-1}{5}}	
	\begin{ytableau}
	1 & 3 \\
	2 & 4\\
	5
	\end{ytableau} \xrightarrow{\frac{3}{6}}	
	\begin{ytableau}
	1 & 3 & 6 \\
	2 & 4 \\
	5
	\end{ytableau} \xrightarrow{\frac{4}{7}}	
	\begin{ytableau}
	1 & 3 & 6 & 7 \\
	2 & 4 \\
	5
	\end{ytableau} 	\]
	The arrows represent the lifting operators 
	and the value above each arrowis the change in eigenvalue at each step.
	From here computing the eigenvalue for $T_{1}$ is a matter of summing all 
	the changes then dividing by $n$ (because we scaled 
	$\RL_{n}$ in Definition \ref{RLop}). We find the eigenvalue for $T_{1}$ being
	\[\eig(T_{1}) =\frac{1}{7} \left( \frac{1}{1} + \frac{0}{2} + 
	\frac{2}{3} + 
	\frac{1}{4} -\frac{1}{5} + \frac{3}{6} + \frac{4}{7}\right)  =  
	\frac{1}{n}\sum_{\substack{\textnormal{boxes}\\ (i,j)}} 
	\frac{j-i+1}{T_{1}(i,j)} = 
	\frac{1171}{2940} .\]
\end{example}

\subsection{Lifting Eigenvectors for the Random Transposition Shuffle}

The eigenvalues of the random transposition shuffle were stated in Section \ref{chpt4:subsec:rt} and shown to correspond to Fourier transforms of the random transposition shuffle at the irreducible representations of $S_{n}$.
Given a partition $\lambda$ of $n$ we computed these to be
\[ \frac{1}{n} + 
\frac{n-1}{n}\frac{\chi_{\lambda}(\tau)}{d_{\lambda}}  = \frac{n + 2\D(\lambda)}{n^{2}} \textnormal{ with multiplicity } d_{\lambda}^{2}.\]
We show our lifting operators $\pro_{a}^{\lambda,\lambda+e_{i}}$ may be used to lift eigenvectors of the random transposition shuffle. Thus, we recover the a description of the eigenvalues of the random transposition shuffle in terms of partitions of $n$.   The analysis follows from the work of Section \ref{chpt5:subsec:algebra}, but we need to account for changes in constants needed to represent a different shuffle.
Recall that in Definition \ref{chpt4:def:RT} the random transposition shuffle was defined with driving probability $\RT$.	
We begin our analysis by turning the random transposition shuffle into an element of the group algebra $\mathfrak{S}_{n}$.

\begin{defn}
	The random transposition shuffle on $n$ cards can be viewed as the following element of $\mathfrak{S}_{n}$:
	\begin{eqnarray}
	\ART_{n} \, :=  \, n^{2} \sum_{\sigma \in S_{n}} \RT(\sigma) \, \sigma  \, =  \,n \cdot \gpid + 2 \sum_{1\leq i < j \leq n} (i \, j) \label{chpt5:eqn:algrt}.
	\end{eqnarray}
	We call this element the \emph{algebraic random transposition shuffle}. Note that here we have scaled by $n^{2}$.
\end{defn}

We remark that the eigenvalues of $\RT$ may be recovered from those of $\ART_{n}$, as we have done for the the random-to-random and the one-sided transposition shuffles (see Lemma \ref{chpt5:lem:eigensame}). The random transposition shuffle also exhibits a recursive structure shown in the equation below
\begin{eqnarray}
\ART_{n+1} - \ART_{n} = \gpid + 2 \sum_{1\leq i \leq n} ( i \hspace{0.3em} n+1). \label{chpt5:eqn:rtdifference}
\end{eqnarray}
This allows us to lift the eigenvectors of $\ART_{n}$ to those of $\ART_{n+1}$. The next result replicates Theorem \ref{chpt5:thm:master}, which relates the shuffle $\ART_{n+1}$ to a combination of $\ART_{n}$ and our adding and switching operators.

\begin{thm}
	\label{chpt5:thm:RTmaster}
	Given $n\in\mathbb{N}$, we have
	\begin{eqnarray}
	\ART_{n+1} \circ \sh_{a} - \sh_{a} \circ \ART_{n} = \sh_{a} + 2 \sum_{1 \leq b \leq n} \sh_{a}\circ \Theta_{b,a}.
	\end{eqnarray}
\end{thm}
\begin{proof}
	This follows from the same proof as Theorem \ref{chpt5:thm:master} allowing for changes in constants in equation \eqref{chpt5:eqn:algrt}
\end{proof}

From here the rest of the analysis in Section \ref{chpt5:subsec:algebra} can be followed without fundamental changes to the algebra. Keeping track of the coefficients in Theorem \ref{chpt5:thm:RTmaster} carefully we can swiftly arrive at a modified Theorem \ref{chpt5:thm:lift} for the random transposition shuffle.

\begin{thm}[Theorem 49 of \cite{dieker2018spectral}]
	\label{chpt5:thm:RTlift}

	Let $\lambda \vdash n$, and  $a \in \{1,2,\ldots,l(\lambda)+1\}$. Take $i\in [n]$ such that  $1 \leq i \leq a$ and set $\mu = \lambda +e_{i}$. Then,
	\begin{eqnarray}
	\ART_{n+1} \circ \pro_{a}^{\lambda,\mu} - 
	\pro_{a}^{\lambda,\mu} \circ \ART_{n} = (1 + 2(\lambda_{i}-i)) \pro_{a}^{\lambda,\mu} .
	\end{eqnarray}
	In particular if $v \in S^{\lambda}$ is an eigenvector of $\ART_{n}$ 
	with eigenvalue $\varepsilon$, then  
	$\pro_{a}^{\lambda,\mu}(v)$ is an eigenvector of 
	$\ART_{n+1}$ belonging to $\mu$ with eigenvalue 
	\begin{eqnarray}
	\varepsilon + 1 + 2(\lambda_{i}+1-i)
	\end{eqnarray}
	The value $\lambda_{i}+1- i$ is the diagonal index of the added box $(i,\lambda_{i}+1)$.
\end{thm}

Notice that the change in eigenvalue $(1+ 2(\lambda_{i}+1-i))$ has no dependence on $n$. Thus, the eigenvalues we recover are not dependent on the order in which we lift our eigenvectors (path we take up Young's lattice), in other words, we can label the eigenvalues  by partitions $\lambda \vdash n$, rather than standard Young tableaux of $\lambda$. Using Theorem \ref{chpt5:thm:RTlift} and the proof of Theorem \ref{mapsnew}  we can recover exactly the eigenvalues we found in Section \ref{chpt4:subsec:rt}.

\begin{lemma}
	\label{chpt5:lem:rteig}

	The eigenvalues of the random transposition shuffle $\RT$ are indexed by partitions $\lambda$. For a partition $\lambda$ the corresponding eigenvalue is:
	\[  \frac{n + 2\D(\lambda)}{n^{2}} \textnormal{ with multiplicity } d_{\lambda}^{2} \]
\end{lemma}
\begin{proof}
	This follows from following the procedure of Theorem \ref{chpt5:thm:liftingeig} keeping track of the new changes in eigenvalue from Theorem \ref{chpt5:thm:RTlift}.
\end{proof}

Comparison of the lifting for the one-sided transposition and random transposition shuffles allows us to see what parts of the method are coloured by our choice of shuffle. The equations in Theorems \ref{chpt5:thm:master} and \ref{chpt5:thm:RTmaster} are crafted to fit with our chosen shuffle, but the lifting operators $\pro_{a}^{\lambda,\mu}$ remain the same in both cases.
The key here is how our shuffles act on $n+1$ versus $n$ cards, they both show a similar structure (compare equations \eqref{chpt5:eqn:lrdifference} and \eqref{chpt5:eqn:rtdifference}) which allows the recovery of their eigenvectors using our lifting operators.
For lifting eigenvectors to work on a generic family of transposition shuffles $\{P_{n}\}$ we need the difference of $P_{n+1}$ and $P_{n}$ to only involve the movement of the new card $n+1$. This allows the adding operator $\sh_{a}$ to mimic the addition of this card, and for the recovery of the eigenvalues of $P_{n+1}$ using those known for $P_{n}$. In Section \ref{chpt5:sec:biased} we show this technique to be applicable to an entire class of transposition shuffles, called the biased one-sided transposition shuffles.

\section{Eigenvalues Analysis}
\label{chpt5:sec:eigenvalues}
In this section we establish important results about the eigenvalues of the one-sided transposition shuffle. Throughout this section we will use 
standard facts and definitions about Young diagrams and tableaux, these 
were laid out in Sections \ref{chpt4:subsec:youngdigrams}, \ref{chpt4:subsec:youngtableau}. 
We also need some specialised notation to deal with the eigenvalues of 
$\LR_{n}$ which we introduce now.

\begin{defn}
	\label{def:bestworst}
	For any $\lambda \vdash n$, define the tableau 
	$\best$ by inserting the numbers $1,\dots, n$ from left to right. Define the 
	tableau $\worst$ by inserting the numbers $1,\dots, n$ from top to bottom.
\end{defn}

Following from Section \ref{chpt5:subsec:algebra} we know the eigenvalues for $\LR_{n}$ are labelled by Young tableaux of size $n$, and Theorem \ref{chpt5:thm:liftingeig} 
gives an explicit formula for the eigenvalue associated to any given tableau. 
Before applying the classical $\ell^{2}$ bound on total variation distance, we first 
investigate relationships between the eigenvalues. 
We show that the eigenvalue corresponding to $T \in \SYT(\lambda)$ is bounded by 
the eigenvalues for $\best$ and $\worst$.  
To simplify our upper bound calculation, 
we prove that we only need to consider the partitions for which $\best$ gives a positive eigenvalue. Lastly, we prove that the eigenvalues corresponding to $\best,\worst$ decrease 
as one moves down the dominance order of partitions. 
We first illustrate the preceding definitions and discussion with an example.

\begin{example}
	Let $\lambda=(3,2)\vdash 5$.
	There are $5$ standard Young tableaux of shape $\lambda$, these tableaux together with the associated eigenvalues are given in Table~\ref{eigenvaluetable} below.  
	In this table, $\best$ is the first tableau listed and $\worst$ is the last one; we can see that the corresponding eigenvalues bound all the others.
	
	\begin{table}[H]
		\centering
		{
			\renewcommand{\arraystretch}{1.5}
			\begin{tabular}{c|c|c|c|c|c}
				$T \in \SYT((3,2))$ &
				$\ytableausetup{mathmode,baseline,aligntableaux=bottom,boxsize=1em}
				\begin{ytableau} 1 & 2 & 3\\ 4 & 5 
				\end{ytableau}$ & 
				$\ytableausetup{mathmode,baseline,aligntableaux=bottom,boxsize=1em}
				\begin{ytableau} 1 & 2 & 4\\ 3 & 5
				\end{ytableau} $ & 
				$\ytableausetup{mathmode,baseline,aligntableaux=bottom,boxsize=1em}
				\begin{ytableau} 1 & 2 & 5\\ 3 & 4
				\end{ytableau} $ & 
				$\ytableausetup{mathmode,baseline,aligntableaux=bottom,boxsize=1em}
				\begin{ytableau} 1 & 3 & 4\\ 2 & 5
				\end{ytableau} $ & 
				$\ytableausetup{mathmode,baseline,aligntableaux=bottom,boxsize=1em}
				\begin{ytableau} 1 & 3 & 5\\ 2 & 4
				\end{ytableau}$   \\\hline
				$\eig(T)$ & $0.64$	& $0.59$ & $0.57$ &  $0.52\overline{3}$	& $0.50\overline{3}$

			\end{tabular}
		}
		\caption[Eigenvalues corresponding to $T \in \SYT((3,2))$]{Eigenvalues corresponding to $T \in \SYT((3,2))$.}
		\label{eigenvaluetable}
	\end{table}
\end{example}

For a Young tableau $T$ which is not necessarily standard, define $\eig(T)$ 
to be the value given by the formula in Theorem \ref{chpt5:thm:liftingeig} (if $T$ 
is not 	standard this value has no relation to the eigenvalues of $\LR_{n}$). We  begin our analysis by showing how swapping numbers 
in a tableau affects the corresponding eigenvalue. Throughout this section 
we accompany results with explicit examples to aid understanding of the results and their proofs.

\begin{lemma}
	\label{swapping}
	Let  $T$ be a Young tableau. Suppose we form a new tableau $S$ by swapping 
	two values in $T$ which have coordinates $(i_{1},j_{1}),(i_{2},j_{2})$ in 
	$T$. WLOG assume  $T(i_{1},j_{1}) < 
	T(i_{2},j_{2})$. Then the change in corresponding eigenvalues satisfies the following inequality:
	\[\eig(S) - \eig(T) 
	\begin{cases}
	\geq 0 & \text{ if } (i_{1}-i_{2})+ (j_{2}-j_{1}) \geq 0\\
	< 0 & \text{ if }  (i_{1}-i_{2}) + (j_{2}-j_{1})<0\,.
	\end{cases}
	\]	
	Importantly, if we move the larger entry down and to the left the change 
	in eigenvalue is non-negative; if it moves up and to the right then the change is negative.
\end{lemma}

\begin{proof}
	Since $S$ and $T$ agree in all but two entries the difference in eigenvalues 
	is given by 
	\begin{align*}
	\eig(S) - \eig(T) &= \frac{1}{n}\left(\frac{j_{1}-i_{1}+1}{T(i_{2},j_{2})} + 
	\frac{j_{2}-i_{2}+1}{T(i_{1},j_{1})}\right) - 
	\frac{1}{n}\left(\frac{j_{1}-i_{1}+1}{T(i_{1},j_{1})} + 
	\frac{j_{2}-i_{2}+1}{T(i_{2},j_{2})}\right) \\
	&= \frac{(i_{1}-i_{2})+(j_{2}-j_{1})}{n} \left(\frac{1}{T(i_{1},j_{1})} - 
	\frac{1}{T(i_{2},j_{2})} 
	\right).\qedhere
	\end{align*}
\end{proof}

\begin{example}
	\label{chpt5:ex:swapingtableau}
	Let $\lambda = (3,3,3)$ and take $T$ to be the Young tableau of shape $\lambda$,
	\[T = \ytableausetup{mathmode,baseline,aligntableaux=center,boxsize=1em}
	\begin{ytableau} 1 & 2 & 3\\ 8 & 9 & 4 \\7 & 6 & 5
	\end{ytableau}  \, .\]
	To demonstrate Lemma \ref{swapping} we transpose the value $9$ with every other value in the tableau. Table \ref{table:eigenvalue3} gives the new tableaux and the corresponding change in eigenvalue.
	We see that when $9$ is moved down and to the left 
	(swapping with values $6,7,8$) the  change in eigenvalue is non-negative, and when it is moved up and to the right (swapping  with values $2,3,4$) the change in eigenvalue is negative. 
\end{example}

\begin{table}[H]
	
	\begin{adjustbox}{width=\columnwidth,center}
		\begin{tabular}{c|c|c|c|c|c|c|c|c}
			$S \in \textnormal{YT}((3,2))$ & 
			$\ytableausetup{mathmode,baseline,aligntableaux=bottom,boxsize=1em}
			\begin{ytableau} 9 & 2 & 3\\ 8 & 1 & 4 \\7 & 6 & 5
			\end{ytableau}$ & 
			$\ytableausetup{mathmode,baseline,aligntableaux=bottom,boxsize=1em}
			\begin{ytableau} 1 & 9 & 3\\ 8 & 2 & 4 \\7 & 6 & 5
			\end{ytableau}$ & 
			$\ytableausetup{mathmode,baseline,aligntableaux=bottom,boxsize=1em}
			\begin{ytableau} 1 & 2 & 9\\ 8 & 3 & 4 \\7 & 6 & 5
			\end{ytableau}$ & 
			$\ytableausetup{mathmode,baseline,aligntableaux=bottom,boxsize=1em}
			\begin{ytableau} 1 & 2 & 3\\ 8 & 4 & 9 \\7 & 6 & 5
			\end{ytableau} $
			& 
			$\ytableausetup{mathmode,baseline,aligntableaux=bottom,boxsize=1em}
			\begin{ytableau} 1 & 2 & 3\\ 8 & 5 & 4 \\7 & 6 & 9
			\end{ytableau}$
			& 
			$\ytableausetup{mathmode,baseline,aligntableaux=bottom,boxsize=1em}
			\begin{ytableau} 1 & 2 & 3\\ 8 & 6 & 4 \\7 & 9 & 5
			\end{ytableau}$
			& 
			$\ytableausetup{mathmode,baseline,aligntableaux=bottom,boxsize=1em}
			\begin{ytableau} 1 & 2 & 3\\ 8 & 7 & 4 \\9 & 6 & 5
			\end{ytableau}$
			& 
			$\ytableausetup{mathmode,baseline,aligntableaux=bottom,boxsize=1em}
			\begin{ytableau} 1 & 2 & 3\\ 9 & 8 & 4 \\7 & 6 & 5
			\end{ytableau}$
			\\\hline 
			$5670 \cdot (\eig(S)- \eig(T))$ & $0$ & $-245$ &  
			$-280$ & $-87.5$ & $0$ & $35$ & $40$ & $8.75$

		\end{tabular}
	\end{adjustbox}
	
	\caption[A list of tableaux formed from $T$ in Example \ref{chpt5:ex:swapingtableau} by transposing the value $9$ with a value from $\{1,\ldots, 8\}$]{A list of tableaux formed from $T$ in Example \ref{chpt5:ex:swapingtableau} by transposing the value $9$ with a value from $\{1,\ldots, 8\}$. To aid comparison of the eigenvalues we have found the difference with $\eig(T)$ and scaled by $5670$. }
	\label{table:eigenvalue3}
\end{table}

Lemma \ref{swapping} allows us to prove that the eigenvalue for any $T \in \SYT(\lambda)$ is bounded between those for $\worst$ and $\best$.

\begin{lemma}
	\label{seq}
	Let $\lambda\vdash n$.
	For any $T \in \SYT(\lambda)$ we have the following 
	inequality:
	\begin{eqnarray}
	\eig(\worst) \leq \eig(T) \leq 
	\eig(\best).
	\end{eqnarray}
\end{lemma}

\begin{proof}
	Reading across the rows of $T$, beginning with the first row, identify the first box in which $T$ and $\best$ have different entries; write $(i,j)$ for the coordinates of this box. Due to the way in which $\best$ is constructed, $T(i,j)>\best(i,j)$. Furthermore, the number $T(i,j)-1$ must occur strictly below and to the left of $T(i,j)$, since $T$ is a standard Young tableau. Swapping entries $T(i,j)-1$ and $T(i,j)$ in tableau $T$ produces a new element of $\SYT(\lambda)$ 
	whose corresponding eigenvalue is no smaller than $\eig(T)$, thanks to Lemma~\ref{swapping}.
	
	We iterate this procedure, swapping $T(i,j)-1$ with 
	$T(i,j)-2$ etc, until $T(i,j) = \best(i,j)$. Note that at this point the 
	entries in the first $T(i,j)$ boxes of $T$ and $\best$ must agree, moreover 
	these entries are now fixed in place. We now 
	proceed to the next box in which $T$ and $\best$ differ, and repeat: this 
	results in a sequence of swaps which make the entries of $T$ agree with 
	those in $\best$, and which can only ever cause the corresponding 
	eigenvalue to increase. This proves the second inequality in 
	Lemma~\ref{seq}, and the first one follows via an analogous argument on the 
	columns of $T$ this time with our eigenvalue decreasing after each iteration.

\end{proof}

\begin{example}
	Let $\lambda = (3,2,1,1)$, to illustrate the procedure defined in Lemma \ref{seq}, we show the first two iterations of the algorithm for turning  $\worst$ into $\best$. The first two iterations involve transforming the boxes $(1,2)$ and $(1,3)$ in $\worst$ to match the boxes of $(1,2)$ and $(1,3)$ in $\best$ (box $(1,1)$ is always matched for any pair of standard Young tableaux). 	We keep track of the eigenvalue of the tableaux at each step of the algorithm to show it is always increasing. Applying transpositions in order to fix box $(1,2)$ we find,
	
	\begin{table}[H]
		\begin{tabular}{cccccccc}
			$\ytableausetup{mathmode,baseline,aligntableaux=center, 
				boxsize=1em}
			\begin{ytableau} 1 & 5 & 7 & 8 \\ 2 & 6 \\ 3 \\ 
			4\end{ytableau}\vspace{0.1em}$ &  
			$\rightarrow$  & 
			$\ytableausetup{mathmode,baseline,aligntableaux=center, 
				boxsize=1em}
			\begin{ytableau} 1 & 4 & 7 & 8 \\ 2 & 6 \\ 3 \\ 
			5\end{ytableau}\vspace{0.1em}$ & 
			$\rightarrow$ &
			$\ytableausetup{mathmode,baseline,aligntableaux=center, 
				boxsize=1em}
			\begin{ytableau} 1 & 3 & 7 & 8 \\ 2 & 6 \\ 4 \\ 
			5\end{ytableau}\vspace{0.1em}$ & 
			$\rightarrow$ & 
			$\ytableausetup{mathmode,baseline,aligntableaux=center, 
				boxsize=1em}
			\begin{ytableau} 1 & 2 & 7 & 8 \\ 3 & 6 \\ 4 \\ 
			5\end{ytableau}\vspace{0.1em}$ 
			\\\hline
			$0.207738$ & & $0.232738$ & & $0.263988$ & & $0.305655$
		\end{tabular}
	\end{table}
	Notice that after each swap we still remain at a standard Young tableau. Once box $(1,2)$ is fixed we proceed to fix box $(1,3)$ in the same way
	\begin{table}[H]
		\begin{tabular}{cccccccccc}
			$\ytableausetup{mathmode,baseline,aligntableaux=center, boxsize=1em}
			\begin{ytableau} 1 & 2 & 7 & 8 \\ 3 & 6 \\ 4 \\ 
			5\end{ytableau}\vspace{0.1em}$ &  
			$\rightarrow$  & 
			$\ytableausetup{mathmode,baseline,aligntableaux=center, boxsize=1em}
			\begin{ytableau} 1 & 2 & 6 & 8 \\ 3 & 7 \\ 4 \\ 
			5\end{ytableau}\vspace{0.1em}$ & 
			$\rightarrow$ &
			$\ytableausetup{mathmode,baseline,aligntableaux=center, boxsize=1em}
			\begin{ytableau} 1 & 2 & 5 & 8 \\ 3 & 7 \\ 4 \\ 
			6\end{ytableau}\vspace{0.1em}$ & 
			$\rightarrow$ & 
			$\ytableausetup{mathmode,baseline,aligntableaux=center, 
				boxsize=1em}
			\begin{ytableau} 1 & 2 & 4 & 8 \\ 3 & 7 \\ 5 \\ 
			6\end{ytableau}\vspace{0.1em}$
			& 	$\rightarrow$ & 
			$\ytableausetup{mathmode,baseline,aligntableaux=center, boxsize=1em}
			\begin{ytableau} 1 & 2 & 3 & 8 \\ 4 & 7 \\ 5 \\ 
			6\end{ytableau}\vspace{0.1em}$
			\\\hline
			$0.305655$ & & $0.311607$ & & $0.33244$ & & $0.35744$ & & $0.38869$
		\end{tabular}
	\end{table}
	
	If we repeat this process until all boxes $(i,j)$ match $\best$, we only ever increase the eigenvalue and thus find an upper bound.	
	If instead we wanted to turn $\best$ into $\worst$ we would first fix box $(2,1)$ like so:
	
	\begin{table}[H]
		\begin{tabular}{cccccccc}
			$\ytableausetup{mathmode,baseline,aligntableaux=center, 
				boxsize=1em}
			\begin{ytableau} 1 & 2 & 3 & 4 \\ 5 & 6 \\ 7 \\ 
			8\end{ytableau}\vspace{0.1em}$ &  
			$\rightarrow$  & 
			$\ytableausetup{mathmode,baseline,aligntableaux=center, 
				boxsize=1em}
			\begin{ytableau} 1 & 2 & 3 & 5 \\ 4 & 6 \\ 7 \\ 
			8\end{ytableau}\vspace{0.1em}$ & 
			$\rightarrow$ &
			$\ytableausetup{mathmode,baseline,aligntableaux=center, 
				boxsize=1em}
			\begin{ytableau} 1 & 2 & 4 & 5 \\ 3 & 6 \\ 7 \\ 
			8\end{ytableau}\vspace{0.1em}$ & 
			$\rightarrow$ & 
			$\ytableausetup{mathmode,baseline,aligntableaux=center, 
				boxsize=1em}
			\begin{ytableau} 1 & 3 & 4 & 5 \\ 2 & 6 \\ 7 \\ 
			8\end{ytableau}\vspace{0.1em}$ 
			\\\hline
			$0.471726$ & & $0.446726$ & & $0.415476$ & & $0.37381$
		\end{tabular}
	\end{table}
	Now the eigenvalue is decreasing after each step, if we fix every column in turn we eventually end at $\worst$.
\end{example}

The next result and its corollary establish that when bounding eigenvalues, we only need to consider
those given by $\best$. 

\begin{lemma}
	\label{lem:sumtranspose}
	Let $\lambda\vdash n$.
	For any $T \in \SYT(\lambda)$ we have 
	\[\eig(T) + \eig(T^{\prime}) = \frac{2H_{n}}{n}.\]
\end{lemma}

\begin{proof}
	Let $T\in SYT(\lambda)$. Then
	\begin{align*} \eig(T) +\eig(T')&= \frac 1n \sum_{\substack{\textnormal{boxes}\\ (i,j)\in T}} 
	\frac{j-i+1}{T(i,j)} + \frac 1n \sum_{\substack{\textnormal{boxes}\\ (j,i)\in T'}} 
	\frac{i-j+1}{T'(j,i)} \\
	&= \frac 1n \sum_{\substack{\textnormal{boxes}\\ (i,j)\in T}} 
	\frac{j-i+1}{T(i,j)} + \frac 1n \sum_{\substack{\textnormal{boxes}\\ (i,j)\in T}} 
	\frac{-(j-i)+1}{T(i,j)} = \frac {2H_n}{n}. \qedhere
	\end{align*}

\end{proof}

\begin{corollary}
	\label{flip}
	Let $\lambda \vdash n$, and suppose we 
	have 
	$\eig(\worst) \leq 0$, then we 
	have 
	\begin{eqnarray}
	\label{case2}
	\eig(T_{\lambda^{\prime}}^{\rightarrow}) \geq 
	|\eig(\worst)| \geq  0.
	\end{eqnarray}
\end{corollary}

\begin{proof}
	It follows from Lemma~\ref{lem:sumtranspose} that $\eig(T_{\lambda'}^\rightarrow) + \eig(T_{\lambda}^\downarrow) = 2H_n/n$. Thus if $\eig(\worst) \leq 0$ then 
	\begin{eqnarray}
	\eig(T_{\lambda'}^\rightarrow) = \frac{2H_n}{n} - 
	\eig(T_{\lambda}^\downarrow) \ge - \eig(T_{\lambda}^\downarrow) = 
	|\eig(T_{\lambda}^\downarrow)| \ge 0\,. \qedhere
	\end{eqnarray}
\end{proof}

We end this section by establishing a relationship between eigenvalues and the dominance ordering on partitions.
\begin{lemma}
	\label{order}
	Let $\lambda, \mu \vdash n$. If $\lambda \trianglerighteq \mu$ then 
	\begin{eqnarray}
	\eig(\best)& \geq & 
	\eig(T_\mu^{\rightarrow}) \label{chpt5:eqn:bestdomdown}\\
	\label{down}
	\text{and }\quad \eig(\worst) & \geq &
	\eig(T_\mu^{\downarrow}). \label{chpt5:eqn:worstdomdown}
	\end{eqnarray}
\end{lemma}

\begin{proof}
	If we can show the statements hold for any partition $\mu$ which is formed from $\lambda$ by 
	moving only one box then inductively it will hold for all $\lambda \trianglerighteq \mu$.
	Suppose $\mu$ is formed from $\lambda$ by moving a box from row $a$ to row $b$, with $a<b\leq l(\lambda)+1$
	(if $b=l(\lambda)+1$ then a new row is created by placing the removed box on the very bottom of the diagram). 
	The box we move goes from coordinates $(a,\lambda_{a})$ of $\lambda$ to $(b,\lambda_{b}+1)$ in $\mu$. 
	
	We begin by proving $	\eig(\best) \geq  
	\eig(T_\mu^{\rightarrow})$. 
	Since $\best$ and  $T_\mu^\rightarrow$ are both numbered from left to right, the effect of moving a box from row $a$ to row $b$ is that $T_\mu^{\rightarrow}(i,j) = \best(i,j) -1$ for any box $(i,j) \in T_{\lambda}^{\rightarrow} \cap T_{\mu}^{\rightarrow}$  with $a<i\le b$; boxes in all other rows contain the same values in both tableaux. Using equation~\eqref{eigenvaluesum}, and remembering to include a term to account for the box being moved, we find that:
	\begin{align*}
	n(\eig(\best)  -
	\eig(T_\mu^{\rightarrow})) 
	& =  
	\left(\frac{\lambda_{a}-a+1}{\best(a,\lambda_{a})} 
	-\frac{(\lambda_{b}+1)-b+1}{T_\mu^{\rightarrow}(b,\lambda_{b}+1)}\right)  
	+ \sum_{\substack{(i,j) \in T_{\lambda}^{\rightarrow}\cap T_{\mu}^{\rightarrow} \\ \textnormal{with } a < i 
			\leq b}}\left[\frac{1}{\best(i,j)} - 
	\frac{1}{T_\mu^{\rightarrow}(i,j)}\right](j-i+1)\\
	& \geq  \left( \frac{\lambda_{a}-a+1}{\best(a,\lambda_{a})} 
	-\frac{(\lambda_{b}+1)-b+1}{T_\mu^\rightarrow(b,\lambda_b+1)}\right)  + (\lambda_a-a+1)
	\left(\frac{1}{T_\mu^{\rightarrow}(b,\lambda_{b}+1)} - 
	\frac{1}{\best(a,\lambda_{a})}\right)\\
	& =  \frac{(\lambda_{a} -\lambda_{b}) + 
		(b-a) -1}{T_\mu^{\rightarrow}(b,\lambda_{b}+1)} 	\geq 0.
	\end{align*} 
	The first inequality holds because all the square-bracketed terms in the sum are negative; 
	we upper bound $j-i+1 \le \lambda_{a}-a+1$, and the resulting sum telescopes. The final inequality holds because $(\lambda_{a}-\lambda_{b})\ge 1$ and $(b-a)\ge 1$.

	For the inequality \eqref{chpt5:eqn:worstdomdown}, recall that $\lambda \trianglerighteq \mu$ if and only if $\mu^{\prime} \trianglerighteq \lambda^{\prime}$. 
	Therefore, using inequality \eqref{chpt5:eqn:bestdomdown} we find that $\eig(T_{\mu^{\prime}}^{\rightarrow}) \geq \eig(T_{\lambda^{\prime}}^{\rightarrow})$. 
	Now Lemma~\ref{lem:sumtranspose} gives $-\eig(T_{\mu}^{\downarrow}) 
	\geq -\eig(T_{\lambda}^{\downarrow})$ and thus we recover the desired 
	inequality.

\end{proof}

\begin{example}
	The dominance ordering  on partitions of size $4$ is a linear order. Table \ref{table:eigenvalue4} lists $\best$ and $\worst$ for every partition $\lambda \vdash 4$. Reading from the table we can see the eigenvalues 
	of $\best$ and $\worst$ decrease down the dominance ordering. Also we may use the eigenvalues in the table to verify that $\eig(T) + \eig(T^{\prime}) = 2H_{4} /4 = 25/24$.
	
	\begin{table}[H]		
		\begin{tabular}{c|c|c|c|c|Sc}
			$\best$ & 
			$\ytableausetup{mathmode,baseline,aligntableaux=center} 
			\begin{ytableau} 1 & 2&3&4 \end{ytableau}$ & 
			$\ytableausetup{mathmode,baseline,aligntableaux=center} 
			\begin{ytableau} 1 & 2&3 \\4 \end{ytableau}$ & 
			$\ytableausetup{mathmode,baseline,aligntableaux=center} 
			\begin{ytableau} 1 & 2\\3&4 \end{ytableau}$ & 
			$\ytableausetup{mathmode,baseline,aligntableaux=center} 
			\begin{ytableau} 1 & 2 \\3 \\4 \end{ytableau}$ & 
			$\ytableausetup{mathmode,baseline,aligntableaux=center} 
			\begin{ytableau} 1 \\ 2\\3\\4 \end{ytableau}\vspace{0.01em}$ 
			\\[0.2cm]\hline
			$24 \cdot \eig\left(\best\right)$ & $24$ & $18$ & $13.5$  & 
			$10.5$ 
			& $1$
			\\\hline 
			$\worst$ & $\ytableausetup{mathmode,baseline,aligntableaux=center} 
			\begin{ytableau} 1 & 2&3&4 \end{ytableau}$ & 
			$\ytableausetup{mathmode,baseline,aligntableaux=center} 
			\begin{ytableau} 1 & 3&4 \\2 \end{ytableau}$ & 
			$\ytableausetup{mathmode,baseline,aligntableaux=center} 
			\begin{ytableau} 1 & 3 \\2 & 4 \end{ytableau}$ & 
			$\ytableausetup{mathmode,baseline,aligntableaux=center} 
			\begin{ytableau} 1 & 4 \\2 \\ 3 \end{ytableau}$ & 
			$\ytableausetup{mathmode,baseline,aligntableaux=center} 
			\begin{ytableau} 1 \\ 2\\3\\4 \end{ytableau}$
			\\\hline
			$24 \cdot \eig\left(\worst\right)$ & $24$ & $14.5$ & 
			$11.5$ & $7$ & 	$1$
		\end{tabular}
		
		\caption[The eigenvalues $\best$, $\worst$ for every partition 
		$\lambda \vdash 4$]{The eigenvalues $\best$, $\worst$ for every partition 
			$\lambda \vdash 4$. The dominance order follows from left to right. 
			We have scaled the values by $24$ for ease of comparison}
		\label{table:eigenvalue4}
	\end{table}
	
\end{example}

\section{Upper Bound for the One-sided Transposition Shuffle}
\label{chpt5:sec:upperbound}
In this section we complete the proof of the upper bound present in Theorem 
\ref{chpt5:thm:cutoff}, making use of the results of 
Section~\ref{chpt5:sec:eigenvalues}. 
The analysis splits into two parts, dealing separately with those partitions $\lambda$ having either \emph{large} or \emph{small} first row.

Theorem \ref{chpt2:thm:classicL2} allows us to upper bound the total 
variation distance in terms of the non-trivial eigenvalues of the 
transition matrix. Using Theorem~		\ref{chpt5:thm:liftingeig} we see that 
the trivial 
eigenvalue corresponds to the one-dimensional partition $\lambda=(n)$, and 
so Theorem~\ref{chpt5:thm:liftingeig} implies that 
\begin{eqnarray}
4\lVert \LR^{t}_{n} -\pi_{n} \rVert_{\textnormal{TV}}^{2}\,\le\, 
\sum_{\substack{\lambda\vdash n\\\lambda \neq (n)}} 
d_{\lambda} \sum_{T \in 
	SYT(\lambda)} 
\eig(T)^{2t} \,. \label{eqn:UB_start}
\end{eqnarray}
Recall from Lemma~\ref{seq} that for any $T \in \SYT(\lambda)$ the 
eigenvalue corresponding to $T$ may be bounded by those corresponding to 
$\worst$ and $\best$. With this in mind, we let $\Lambda_n^\rightarrow = 
\{\lambda\vdash n\,:\, |\eig(T_\lambda^\downarrow)|\,\le\, 
|\eig(T_\lambda^\rightarrow)|\}$ and $\Lambda_n^\downarrow = 
\{\lambda\vdash n\,:\, |\eig(T_\lambda^\downarrow)|\,>\, 
|\eig(T_\lambda^\rightarrow)|\}$; note that these are disjoint sets, with 
$\Lambda_n^\rightarrow \subseteq \{\lambda\vdash n\,:\, \eig(\best) 
\,\ge\,0 \}$ and $\Lambda_n^\downarrow \subseteq \{\lambda\vdash n\,:\, 
\eig(\worst) \,<\,0 \}$. Using Lemma~\ref{seq} and then 
Corollary~\ref{flip} we relax the upper bound as follows:
\begin{eqnarray}
4\lVert \LR^{t}_{n} -\pi_{n} \rVert_{\textnormal{TV}}^{2} & \leq & 
\eig\left(T_{(1^n)}\right)^{2t} \,+ \,
\sum_{\substack{\lambda\in \Lambda_n^\rightarrow\\\lambda \neq(n)}} 
d_{\lambda} \sum_{T \in SYT(\lambda)} \eig(T)^{2t}
+ \sum_{\substack{\lambda\in \Lambda_n^\downarrow\\\lambda \neq(1^n)}} 
d_{\lambda} \sum_{T \in SYT(\lambda)} \eig(T)^{2t}\nonumber \\
& \leq &   \eig\left(T_{(1^n)}\right)^{2t} \,+ \,
\sum_{\substack{\lambda\in \Lambda_n^\rightarrow\\\lambda \neq(n)}} 
d_{\lambda}^2 \,\eig(\best)^{2t}
+ \sum_{\substack{\lambda\in \Lambda_n^\downarrow\\\lambda \neq(1^n)}} 
d_{\lambda}^2 \,\eig(\worst)^{2t} \nonumber\\	
& \leq & \eig\left(T_{(1^n)}\right)^{2t} \,+ \,
\sum_{\substack{\lambda\,:\,\eig(\SEQ) \geq 0\\ \lambda \neq (n)}} 
d_{\lambda}^{2}  \,	\eig(\best)^{2t} 
+ \sum_{\substack{\lambda\,:\,\eig(\OSEQ) < 0 \\ \lambda \neq(1^n)}} 
d_{\lambda}^{2} \,\eig(\worst)^{2t} \nonumber\\
& \leq &   \eig\left(T_{(1^n)}\right)^{2t} \,+ \,
\sum_{\substack{\lambda\,:\,\eig(\SEQ) \geq 0\\ \lambda \neq (n)}} 
d_{\lambda}^2 \,\eig(\best)^{2t}
+ \sum_{\substack{\lambda\,:\,\eig(\OSEQ) < 0 \\ \lambda'\neq(1^{n})}} 
d_{\lambda'}^2 \,\eig(T^\rightarrow_{\lambda'})^{2t} \nonumber \\
& \leq & \eig\left(T_{(1^n)}\right)^{2t} \,+ \,
2\,\sum_{\substack{\lambda\,:\,\eig(\SEQ) \geq 0\\ \lambda \neq (n)}} 
d_{\lambda}^{2}  \,	\eig(\best)^{2t} \,. \label{eqn:UB_eigs}
\end{eqnarray}
In the penultimate line we have used Corollary~\ref{flip} and the fact that 
$d_{\lambda'}=d_\lambda$. The final inequality follows by a second 
application of Corollary~\ref{flip}: if $\lambda$ satisfies $\eig(\OSEQ)<0$ 
then $\eig(T^\rightarrow_{\lambda'})$ must be non-negative.

The first term in \eqref{eqn:UB_eigs} is simple to deal with at time $t=n\log(n)+cn$. We have already observed that $\eig\left(T_{(n)}\right) = 1$, and so Lemma~\ref{flip} implies that $\eig\left(T_{(1^{n})}\right) = 2 H_{n}/n-1$. This means that
\begin{equation}\label{eqn:first_partition}
\eig\left(T_{(1^n)}\right)^{2t} = 
\left(1 - 
\frac{2H_{n}}{n} \right)^{2(n\log n + cn)} \leq e^{-4 H_{n} (\log + c)} \textnormal{ for $n\geq 5$ }
\end{equation}
here we have used the bound $1-x \leq e^{-x}$, we see that this tends to 
zero 
for any positive $c$ as $n\to\infty$.

It therefore remains to bound the sum in \eqref{eqn:UB_eigs}. The partitions with 
the biggest eigenvalues are be those with large first rows $\lambda_{1}$, and so we split the analysis into two parts according to this value; by \emph{large} partitions we mean those with 
$\lambda_{1} \geq 3n/4$, and \emph{small} partitions are those with 
$\lambda_{1} <3n/4$. Large partitions give the biggest eigenvalues for 
$\LR_{n}$ and must be dealt with carefully; it is these which will 
determine the mixing time of the shuffle. Small partitions have 
correspondingly larger dimensions, but eigenvalues which are small enough 
to give control around time of order $n\log(n)$. 
We begin by identifying the partition at the top of the dominance ordering 
for any fixed value of $\lambda_1$, which allows us to employ 
Lemma~\ref{order}.

\begin{defn}
	Let $k\in[n]$, define the partition $(n-k,\star)$ to have as many rows of $n-k$ as possible, with the last row being formed of  $n - (n-k) \lfloor n/(n-k) \rfloor$ boxes. For example,
	\begin{eqnarray*}
		(n-k,\star) = \begin{cases}
			(n-k,k) & \textnormal{ if }k \in[1,n/2] \\
			(n-k,n-k,2k-n) & \textnormal{ if } k \in(n/2,n/3]
		\end{cases}
	\end{eqnarray*}
\end{defn}
\begin{lemma}
	Suppose $\lambda\vdash n$ has first row equal to 
	$\lambda_{1}=n-k$, then by moving boxes up and to the right it follows 
	trivially that  $\lambda \trianglelefteq (n-k,\star)$.
\end{lemma}

The notation of $(n-k,\star)$ is driven by our analysis as it will transpire that only the size of the first two 
rows are important for our bounds. 	For each $k$ we also need a bound 
on sum of the squared dimensions of all partitions with 
$\lambda_{1}=n-k$, and for this recall that Lemma \ref{chpt4:lem:youngdimbound} tells us
\[\sum_{\substack{\lambda \vdash n \\ \lambda_{1} =n-k} }
d_{\lambda}^{2} \leq {n \choose k}^{2} k!\leq   \frac{n^{2k}}{k!} \,.\]

\subsubsection{The Eigenvalues of Partition $(n-1,1)$}

Before we proceed with the analysis of large and small partitions let us look at the partition $(n-1,1)$ in more detail. This partition gives the first and largest term in the sum in	\eqref{eqn:UB_eigs}, and controls the mixing time of the one-sided transposition shuffle.

Label the different standard Young tableau of shape $(n-1,1)$ by $T_{i}$, where $i$ denotes the value in the second row of $(n-1,1)$.
The eigenvalue of $T_{i}$ is given by,
\[\eig(T_{i}) = 1 -\frac{1}{n}\left(1 + \sum_{j=i+1}^{n} \frac{1}{j}\right)\]
We have $\eig(T_{n}) = \eig(T_{(n-1,1)}^{\rightarrow})$, and a clear linear order on the eigenvalues given by 
\[\eig(T_{n}) > \eig(T_{n-1}) > \ldots > \eig(T_{3}) > \eig(T_{2}) .\]

The tableau $T_{n}$ gives us the second biggest (in absolute value) eigenvalue of the one-sided transposition shuffle. The contribution of this tableau to the sum \eqref{eqn:UB_start} is,
\begin{eqnarray}
(n-1) \left(1 - \frac{1}{n}\right)^{2t} \label{chpt5:eqn:lr2bigeig}
\end{eqnarray}
and we have seen previously that this eigenvalue is killed at time $t =(n/2)\log n$. We have said before that the time needed to kill the second biggest eigenvalue often tells us the time expected for a random walk to converge to its stationary distribution. For the one-sided transposition shuffle the eigenvalue $T_{n}$ is tightly grouped around the other eigenvalues for partition $(n-1,1)$, with the biggest difference being $T_{n} -T_{2} =(H_{n}-1)/n$. This means that in order to estimate the mixing time it is not good enough to consider the single eigenvalue \eqref{chpt5:eqn:lr2bigeig}; rather, we must consider the sum,
\begin{eqnarray}
\sum_{i=2}^{n} (n-1) \eig(T_{i})^{2t} = (n-1) \sum_{i=2}^{n} \left(1 -\frac{1}{n}\left(1 + \sum_{j=i+1}^{n} \frac{1}{j} \right)\right)^{2t}. \label{chpt5:eqn:lrbigeigsum}
\end{eqnarray}
This sum is bounded at time close to $n\log n$ as $n\to \infty$. Therefore, we expect the mixing time of the one-sided transposition shuffle to be around $n\log n$, as stated in Theorem \ref{chpt5:thm:cutoff}. Note that in sum \eqref{eqn:UB_eigs} we reduce the sum \eqref{chpt5:eqn:lrbigeigsum} of eigenvalues for $(n-1,1)$, to
\[(n-1)^{2} \left(1 - \frac{1}{n}\right)^{2t} \]
which is bounded in $n$ exactly at time $t=n\log n$.

\subsection{Large Partitions}
In this subsection we prove that the sum of large partitions is bounded 
with a decay of $e^{-2c}$ at time $n \log n +cn$. 
Let $\lambda$ be a partition satisfying $\eig(\best)\ge 0$, and for which $\lambda_1 = n-k$ for some  $k\leq n/4$. We have observed that $\lambda \trianglelefteq(n-k,k)$, and so Lemma~\ref{order} suggests that we look at the eigenvalue of $T_{(n-k,k)}^{\rightarrow}$. Using our eigenvalue formula from Theorem \ref{chpt5:thm:liftingeig} we calculate this as follows, with the first and second sum corresponding to the first and second row of $T_{(n-k,k)}^{\rightarrow}$ respectively:
\begin{align}
\eig\left(T_{(n-k,k)}^{\rightarrow}\right) & = \frac 1n \sum_{j=1}^{n-k} 
\frac{j}{T_{(n-k,k)}^{\rightarrow}(1,j)} 
\,+\, \frac 1n\sum_{j=1}^{k} 
\frac{j-1}{T_{(n-k,k)}^{\rightarrow}(2,j)}\nonumber \\
&=
\frac{n-k}{n} \,+\, \frac 1n \sum_{j=1}^{k} 
\frac{j-1}{n-k+j} \label{eqn:large} \\
&=  1 -\frac{(n-k+1)}{n}(H_{n}-H_{n-k+1}) - \frac 1n \,.
\end{align}
We now use this, along with the inequality $1-x \leq e^{-x}$, to bound the contribution of large partitions to the sum in \eqref{eqn:UB_eigs}:
\begin{align}
\sum_{k=1}^{n/4}
\sum_{\substack{\lambda\,:\,\eig(\best)\ge 0\\ \lambda_{1}=n-k}} 
d_{\lambda}^{2} \text{eig}(T_{\lambda}^{\rightarrow})^{2t}
& \, \leq\, \sum_{k=1}^{n/4} \,
\text{eig}\left(T_{(n-k,k)}^{\rightarrow}\right)^{2t} 
\sum_{\substack{\lambda\,:\,\eig(\best)\ge 0\\ \lambda_{1}=n-k}} 
d_{\lambda}^{2}\nonumber \\
& \, \leq\,   	\sum_{k=1}^{n/4} \frac{n^{2k}}{k!}\,
\text{eig}\left(T_{(n-k,k)}^{\rightarrow}\right)^{2t}  \qquad\text{(by Lemma \ref{chpt4:lem:youngdimbound})} \nonumber \\
& \,\leq\, 	\sum_{k=1}^{n/4} \frac{n^{2k}}{k!}\, \left( 1 
-\frac{(n-k+1)}{n}(H_{n}-H_{n-k+1}) - \frac{1}{n}\right)^{2t} \nonumber \\
& \,\leq\, 	\sum_{k=1}^{n/4}  \frac{n^{2k}}{k!}\, e^{-2t\left(\frac{(n-k+1)}{n}(H_{n}-H_{n-k+1})+ 
	\frac{1}{n}\right)} \nonumber \\
& \,\leq\,  e^{-2c}	\sum_{k=1}^{n/4}
\frac{n^{2k-2(n-k+1)(H_{n}-H_{n-k+1})-2}}{k!}\,, \label{sum_large}
\end{align}
in the last step we have substituted $t=n\log n + cn$. When $k=1$ we 
get the following term
\[\left(n^{2-2(n)(H_{n}-H_{n})-2}\right)/1! = 1.\]
Thus if the ratio between consecutive terms is less than $1$ for $n$ 
suitably 
large we may bound the 
sum via a geometric series.
The ratio of the $(k+1)^{\textnormal{th}}$ term to the $k^{\textnormal{th}}$ term in 
\eqref{sum_large} is given 
by 
\begin{eqnarray}
\label{ratio1}
\frac{n^{2(H_n-H_{n-k})}}{k+1}\,.
\end{eqnarray}
For large $n$ this ratio is approximated by $(n^{2\log (n/n-k)})/(k+1)$. To 
see that this ratio is less than $1$ for all $k\in \{1,\dots,n/4\}$ as 
$n\rightarrow \infty$ we consider the two cases, $k =O(1)$ and $k=O(n)$.  In both cases $\lim_{n\to\infty} n^{2(H_{n} - H_{n-k})}  = \lim_{n\to\infty} n^{2\log n/(n-k)}$. Now if $k$ is 
constant then we have $\lim_{n\to\infty} n^{2\log n/ (n-k)} = 1$ and so 
\[\lim_{n\to\infty} (n^{2\log (n/n-k)})/(k+1) = \frac{1}{k+1}.\] If $k = an$ with $a\in(0,1/4]$ then $n^{2\log (n/(n-an))} = n^{2\log 
	1/(1-a)} 
\leq n^{2\log (4/3)} = n^{0.58}$. Taking the limit in $n$ we have,
\[\lim_{n\to\infty} (n^{2\log (n/n-k)})/(k+1) \leq \lim_{n\to\infty} 
n^{0.58-1}/a = 0.\]
Therefore, for large enough $n$ the ratio \eqref{ratio1} is less than $1$ 
for all $k \in \{1,\dots,n/4\}$. Indeed for $n\to \infty$ the largest value 
of the ratio over this range of $k$ 
is achieved when $k=1$, at which point it equals $n^{2/n}/2$. For 
sufficiently large $n$ this ratio is thus bounded above by $3/4$, say, 
which permits us to bound the sum in \eqref{sum_large} by a geometric 
series with initial	term 1:
\begin{equation}
\label{eqn:sum_large_1}
e^{-2c}\,	\sum_{k=1}^{n/4}
\frac{n^{2k-2(n-k+1)(H_{n}-H_{n-k+1})-2}}{k!} \, \le \,
e^{-2c}	\sum_{k=1}^{n/4} (3/4)^{k-1} \,\le\, 4e^{-2c}\,.
\end{equation}

\subsection{Small Partitions}

Now consider a partition $\lambda$ satisfying $\eig(\best)\ge 0$ and for which 
$\lambda_1 = n-k$ with  $n/4<k\le n-2$. Suppose first of all that $n/4<k\le n/2$; as in the large partition case, any such partition is dominated by $(n-k,k)$, and the same calculation as in equation \eqref{eqn:large} shows that 
\begin{equation}\label{eqn:medium}
\eig\left(T_{(n-k,k)}^{\rightarrow}\right)  
\,=\, \frac{n-k}{n} \,+\, \frac 1n \sum_{j=1}^{k} 
\frac{j-1}{n-k+j} \,. 
\end{equation}

Now consider the case when $k>n/2$. We have already identified that 
$\lambda\trianglelefteq(n-k,\star)$, and so we proceed by calculating the 
eigenvalue of $T^\rightarrow_{(n-k,\star)}$. Note first that for any box on 
row three or below -- that is $(i,j)$ with $i\geq 3$ -- its contribution to equation \eqref{eigenvaluesum} may be bounded by: 
\begin{align*}
\frac{j-i+1}{T_{(n-k,\star)}^{\rightarrow}(i,j)} &= \frac{j-i+1}{(i-1)(n-k) + 
	j} \le \frac{(n-k)}{(i-1)(n-k) + (n-k)} \le \frac 13\,.
\end{align*}
Using this inequality in conjunction with Theorem \ref{chpt5:thm:liftingeig} we bound 
$\eig\left(T_{(n-k,\star)}^{\rightarrow}\right)$ as follows:
\begin{align}
\eig(T_{(n-k,\star)}^{\rightarrow}) & = \frac 1n \sum_{j=1}^{n-k} 
\frac{j}{T_{(n-k,\star)}^{\rightarrow}(1,j)} 
\,+\, \frac 1n\sum_{j=1}^{n-k} 
\frac{j-1}{T_{(n-k,\star)}^{\rightarrow}(2,j)} \, + \, \frac 1n 
\sum_{\substack{(i,j) \nonumber \\ 
		i\geq 3}} \frac{j-i+1}{T_{(n-k,\star)}^{\rightarrow}(i,j)} \\
& \leq \frac{n-k}{n} \,+\, \frac 1n \sum_{j=1}^{n-k} 
\frac{j-1}{n-k+j} \, + \,  \frac{n-2(n-k)}{3n}\,. \label{eqn:small1}
\end{align}
Observe that if we substitute $k \in(n/4,n/2]$ in \eqref{eqn:small1} 
it provides an upper bound for the expression in \eqref{eqn:medium}. 
Indeed, 
for $n/4<k\le n/2$ we may write
\begin{align}
\frac{n-k}{n} \,+\, \frac 1n \sum_{j=1}^{n-k} 
\frac{j-1}{n-k+j} \, + \,  \frac{n-2(n-k)}{3n} - 
\eig(T_{(n-k,k)}^{\rightarrow}) &\,=\, \frac 1n \sum_{j=k+1}^{n-k} 
\left(\frac{j-1}{n-k+j} -\frac 13 \right) \nonumber
\\
&\,= \,  \frac{2(n-2k)}{3n} - \frac{(n-k+1)}{n}(H_{2(n-k)}-H_n) \nonumber 
\,. 
\end{align}

Substituting $k=\gamma n$, the final expression is bounded below by
\begin{eqnarray}
\frac{2(1-2\gamma)}{3} - \left(1-\gamma 
+\frac{1}{n}\right)\log(2(1-\gamma)) 
\label{eq:smalldiff}.
\end{eqnarray}
For $n\geq 15$ the expression \eqref{eq:smalldiff} is non-negative for all 
$\gamma \in [1/4,1/2]$ thus completing our claim. We have just shown that \eqref{eqn:small1} provides a bound on 
$\eig\left(T_{(n-k,\star)}^{\rightarrow}\right)$ for all $k\in(n/4,n-2]$, 
therefore, it provides a bound for all $\eig\left(\best \right)$ with 
$\lambda_{1} =n-k$.
Working with our new bound we rearrange it to:
\begin{align}
\eig\left(T_{\lambda}^{\rightarrow}\right) \,&\le\, \frac{n-k}{n} \,+\, 
\frac 1n \sum_{j=1}^{n-k} 
\frac{j-1}{n-k+j} \, + \,  \frac{n-2(n-k)}{3n} \nonumber\\
&= \frac{n-k}{n} + \frac{n-k-1 - (n-k+1)(H_{2(n-k)} -H_{n-k +1})}{n} +  
\frac{2k-n}{3n} \nonumber\\ 
& = 1 -\frac{(4k-2n+3)}{3n} - \frac{(n-k+1)}{n}(H_{2(n-k)} -H_{n-k +1})\, 
\label{eqn:small2}.
\end{align}
Using the inequalities $1-x \leq e^{-x}$ for all $x$, and, our eigenvalue 
bound 
\eqref{eqn:small2},
we are able to bound the 
contributions of small partitions in the sum \eqref{eqn:UB_eigs} at time $t=n\log n +cn$ as follows:
\begin{align}
\sum_{k=n/4}^{n-2}
\sum_{\substack{\lambda\,:\,\eig(\best)\ge 0\\ \lambda_{1}=n-k}} 
d_{\lambda}^{2} \eig\left(T_{\lambda}^{\rightarrow}\right)^{2t}
& \, \leq\,   	\sum_{k=n/4}^{n-2} \frac{n^{2k}}{k!}\,
e^{-\frac{2t}{n}\left(\frac{4k-2n+3}{3} + (n-k+1)(H_{2(n-k)} -H_{n-k +1})\right)}\nonumber \\
& \leq   
e^{-2c}\sum_{k =n/4}^{n-2} 
\frac{n^{\frac{4n-2k-6}{3}-2(n-k+1)(H_{2(n-k)} -H_{n-k +1})}}{k!}\,. 
\label{sum_small2}
\end{align}

To analyse this bound further we require the following inequality.
\begin{lemma}
	\label{lem:smallbound}
	$(n+1)\left(H_{2n} -H_{n+1}\right) > (n-1)\log 2$ for all integers $n\geq 1$ 
\end{lemma}
\begin{proof}
	Rearranging our inequality we have to show that
	$(n+1)\left(H_{2n} -H_{n+1} - \log 2\right) > -2 \log 2$.
	The $n^{\textnormal{th}}$ harmonic number may be bounded by $\log n + 
	\gamma + \frac{1}{2n+1} \leq H_{n} \leq \log n 
	+\gamma + \frac{1}{2n-1}$ where $\gamma$ is the Euler–Mascheroni constant.
	
	Using the lower bound for $H_{2n}$ and upper 
	bound for $H_{n+1}$, we have:

	\begin{eqnarray*}
		(n+1)\left(H_{2n} -H_{n+1} - \log 2\right) & \geq &
		(n+1) \left( \frac{1}{4n+1} - \frac{1}{2n+1} 
		+\log\left(\frac{n}{n+1}\right)\right)\\
		& = &  -\frac{2n(n+1)}{(4n+1)(2n+1)}  -(n+1)\log\left(\frac{n+1}{n}\right) \\
		& \geq & -\frac{2n(n+1)}{(4n+1)(2n+1)} -\frac{(n+1)}{n} \\
		& \geq & -2 \log 2 \textnormal{ for all } n\geq 8.
	\end{eqnarray*}
	The last inequality comes from $\frac{2n(n+1)}{(4n+1)(2n+1)} +\frac{(n+1)}{n}$ being a decreasing function, and $n=8$ is the first time it passes $2\log 2$.	The original inequality may be verified for the remaining integer values $1 \leq n \leq 7 $ 
\end{proof}

Using the new bound provided by Lemma \ref{lem:smallbound} we may bound 
\eqref{sum_small2} via the following:
\begin{eqnarray}
e^{-2c}\sum_{k =n/4}^{n-2} 
\frac{n^{\frac{4n-2k-6}{3}-2(n-k-1)\log 2}}{k!}\,. 
\label{sum_small}
\end{eqnarray}
Once again writing $k=\gamma n$, now for $\gamma\in[1/4,1]$, the terms in 
the summand of \eqref{sum_small} may be rewritten as
\begin{eqnarray} n^{\frac{2n(2-\gamma)}{3}-2n(1-\gamma)\log 2 - 0.5 +2\log 
	2} /(\gamma n!) \label{eqn:smallterm}
\end{eqnarray}
An application of Stirling's formula to $\gamma n!$ tells us that
$\gamma n! > \sqrt{2\pi} (\gamma n)^{\gamma n +0.5} e^{-\gamma n}$.
Combining this with equation \eqref{eqn:smallterm} we get
\begin{eqnarray}
\frac{e^{\gamma n} n^{2\log 2 -1}}{\gamma^{\gamma 
		n}\sqrt{2\pi}}n^{\frac{n}{3}\left( 4 -5\gamma -6(1-\gamma)\log 2 \right)} 
\label{eqn:smallterm2}
\end{eqnarray}

Thus the dominant term 
of \eqref{eqn:smallterm2} takes the form $n^{\frac{n}{3}g(\gamma)}$,
where $g(\gamma) = 4-5\gamma - 6(1-\gamma)\log 2 < 0$ for all
$\gamma\in[1/4,1]$. It follows that, for any positive $c$, 
\begin{eqnarray}
\lim_{n\rightarrow \infty} e^{-2c}\sum_{k =n/4}^{n-2} 
\frac{n^{\frac{4n-2k-6}{3}-2(n-k-1)\log 2}}{k!} =0 \label{eqn:smallfinal}
\end{eqnarray}	
This completes the analysis of small partitions.

\subsubsection*{Proof of the Upper Limit in Theorem \ref{chpt5:thm:cutoff}}
Combining the results and bounds of \eqref{eqn:UB_eigs}, 
\eqref{eqn:first_partition}, \eqref{eqn:sum_large_1} and 
\eqref{eqn:smallfinal} we find at time $t=n\log n +cn$, 
\begin{eqnarray*}
	4 \,\limsup_{n\to \infty}\lVert \LR^{t}_{n} -\pi_{n} \rVert_{\tiny \textnormal{TV}}^{2} 
	& \leq & \limsup_{n\to \infty}
	\eig\left(T_{(1^n)}\right)^{2t} \,+ \,
	2\,\limsup_{n\to\infty} \sum_{\substack{\lambda\,:\,\eig(\SEQ) \geq 0\\ 
			\lambda \neq (n)}} 
	d_{\lambda}^{2}  \,	\eig(\best)^{2t} \,. \\
	&\leq & 2\,\limsup_{n\to\infty} \sum_{k=1}^{n/4} 
	d_{\lambda}^{2} \eig\left(T_{(n-k,k)}^{\rightarrow}\right)^{2t} +
	+ 2\,\limsup_{n\to\infty} \sum_{k=n/4}^{n-2} 
	d_{\lambda}^{2} \eig\left(T_{(n-k,\star)}^{\rightarrow}\right)^{2t}\\
	& \leq & 8e^{-2c}
\end{eqnarray*}
This completes the proof that   $\limsup_{n\to \infty} \lVert \LR^{n\log n +cn}_{n} -\pi_{n} 
\rVert_{\tiny \textnormal{TV}} \leq \sqrt{2}e^{-c}$. Thus, we have an upper bound on the mixing time of the one-sided transposition shuffle of time $n\log n$.

\section{Lower Bound for the One-sided Transposition Shuffle}
\label{chpt5:sec:lowerbound}

To complete Theorem \ref{chpt5:thm:cutoff} we need to prove our lower limit 
on total variation distance. To do this we employ the usual trick of finding a set of 
permutations $F_n\subseteq S_n$  which has significantly different 
probability under the uniform distribution $\pi_n$ and the one-sided 
transposition measure $\RL_n^t$ before a certain time $t$. 
The definition of total variation distance 
then immediately yields a simple lower bound:
\[ \lVert \RL_n^t -\pi_n \rVert_{\textnormal{TV}}\, \ge\, \RL_n^t(F_n) - 
\pi_n(F_n) \,. \]
In particular, we follow in the steps of the random transposition shuffle from Section \ref{chpt4:subsec:rt} and 
find a suitable set $F_n$ by considering the number of fixed points within 
- a certain part of - the deck. Estimation of $\RL_n^t(F_n)$ then reduces 
to a novel variant of the classical coupon collector's problem.

\paragraph{}
Recall that one step of the one-sided transposition shuffle may be modelled 
by 
firstly choosing a position $\righthand \sim_{d} U\{1,\dots n\}$ with our right 
hand, and 
then choosing a position $\lefthand\sim_{d} U\{1,\dots ,\righthand\}$ with our 
left 
hand and 
transposing the cards in the chosen positions.
Since the left hand always chooses a position below that of the right hand, 
it is intuitively clear that our shuffle is relatively unlikely to 
transpose two cards near to the top of the deck. For example, taking $n>3$,  $\LR_{n}((1\,2)) = 
1/2n$, whereas $\LR_{n}((n-1 \, n)) = 1/n^{2}$, therefore we are $(n/2)$ times as likely	to apply transposition $(1\,2)$, than $(n-1 \, n)$. 

This leads us to focus  the
attention of our analysis on a set of positions at the top of the deck:  write $V_n$ for 
the top part of the deck, 
\[ V_n = \{n-n/m+1,\dots,n-1,n\} \,, \]
where $m=m(n)$ is growing in $n$ and to be chosen later. We want to keep track of 
fixed points within this part of the deck, let $ F_n = \{\sigma\in S_n \, | \, \text{$\sigma$ has at least 1 fixed point in 	$V_n$}\}.$
Note that $V_n$ contains $n/m$ positions, and so we find an easy upper 
bound on the uniform distribution $\pi_{n}(F_n) \leq 1/m \to 0$ as $n\to 
\infty$.

To bound the value of $\RL_n^t(F_n)$ we  
reduce the problem to studying a simpler Markov chain linked to coupon 
collecting. When either of our hands $(\righthand, \lefthand)$ picks a new 
(previously 
untouched) card we say that this card gets \emph{collected}. The  
uncollected cards in $V_n$ at time $t$ are those which have not yet been 
picked 
by either hand, and thus the size of this set gives us a lower bound on the 
number of fixed points in $V_n$. Writing $U_n^t$ for the set of uncollected 
cards in $V_n$ after $t$ steps of the one-sided transposition shuffle, it follows that 
\begin{eqnarray}
\label{Ineq1}
\RL_{n}^{t}(F_n) \,\geq\, \mathbb{P}(|U^{t}_{n}| \ge 1 )\,.
\end{eqnarray}

We wish to show that at time $t=n\log n - n\log\log n$ the probability on 
the 
right hand side of \eqref{Ineq1} is large, thus recovering a lower bound on 
$\LR_{n}^{t}$. For the one-sided transposition shuffle the value of $\lefthand$ is clearly not independent of the value $\righthand$. Importantly if $\righthand \notin V_{n}$ then $\mathbb{P}(\lefthand \in V_{n}) =0$.   This means that a standard 
coupon-collecting argument for the time taken to collect all of the 
cards/positions in $V_n$ cannot be applied to our shuffle, and a little 
more work is therefore required.

Note that at each step there are four possibilities: both hands collect new 
cards, only one hand does (left or right) or neither does. This permits us 
to bound the change in the number of \emph{collected} cards as follows:
\begin{align}
|V_n\setminus U_n^{t+1}| & \,=\, |V_n\setminus U_n^{t}| + 
|\{L^{t+1},R^{t+1}\}\cap U_n^t| \nonumber \\
& \,\le\, |V_n\setminus U_n^{t}| + 2 \cdot \ind [L^{t+1}\in U_n^t] + 
\ind[L^{t+1}\notin U_n^t,R^{t+1}\in 
U_n^t]\,,\label{eqn:indicator_fns} 
\end{align}
where $\ind[\cdot]$ is an 
indicator function. We may find a useable 
upper bound for $\mathbb{P}(L^{t+1}\in U_n^t)$ by conditioning only on the size of $U_{n}^{t}$. Furthermore, we shall show that $\mathbb{P}(L^{t+1}\in U_n^t)$ is relativity small compared with $\mathbb{P}(L^{t+1}\notin U_n^t,R^{t+1}\in 
U_n^t)$, thus our approximation does not stop us finding the correct mixing time.

The probability $ \mathbb{P}(L^{t+1}\in U_n^t)$, naturally depends on what 
positions the 
uncollected cards are in at time $t$. However, our left hand is 
more likely to choose positions towards the bottom of the pack so, letting 
$\hat U_n^t = \{n-n/m+1,\dots,n-n/m+ |U_n^t|\}$, i.e. the $|U_n^t|$ 
lowest numbered positions in $V_n$, we may form an upper bound on our 
probabilities as follows:
\begin{align*}
\mathbb{P}(L^{t+1}\in U_n^t) &\,\le\, \mathbb{P}(L^{t+1}\in \hat U_n^t) \,.
\end{align*}
Given the number of uncollected cards,the 
the set $\hat U_{n}^{t}$ is completely determined, so we may compute a bound in terms of $|U_{n}^{t}|$, 
\begin{align}
\mathbb{P}(L^{t+1}\in \hat U_n^t) &\,=\, \frac 1n\sum_{k\in \hat 
	U_n^t}\mathbb{P}\left(L^{t+1}\in \hat U_n^t\,|\, R^{t+1}=k\right)+  
\frac 
1n\sum_{k\in V_n\setminus \hat U_n^t}\mathbb{P}\left(L^{t+1}\in \hat 
U_n^t\,|\, 
R^{t+1}=k\right) \nonumber \\
&\,=\, \frac 1n\sum_{k\in \hat U_n^t} \frac{k-(n-n/m)}{k} +  \frac 
1n\sum_{k\in 
	V_n\setminus \hat U_n^t}\frac{|U_n^t|}{k} \nonumber \\
&\,\le\, \frac{\sum_{k=1}^{|U_n^t|} k + (n/m -|U_n^t|)|U_n^t|}{n(n-n/m)} 
\,\le\,\frac{|U_n^t|}{(m-1)n}\,. \label{eqn:bound_minus2}
\end{align}
The probability of the event $\{L^{t+1}\notin U_n^t,R^{t+1}\in 
U_n^t\}$ in \eqref{eqn:indicator_fns} is simple 
to bound:
\begin{equation}\label{eqn:bound_minus1}
\mathbb{P}(L^{t+1}\notin U_n^t,R^{t+1}\in U_n^t) \,\le\, 
\mathbb{P}(R^{t+1}\in U_n^t) \,\le\, \frac{|U_n^t|}{n}\,.
\end{equation}
Using \eqref{eqn:indicator_fns}, \eqref{eqn:bound_minus2} and 
\eqref{eqn:bound_minus1} together, we now define a counting process $M_n^t$ 
which stochastically dominates the number of collected cards $|V_n\setminus 
U_n^{t}|$ at all times: 
\begin{align}
M_n^0 \,&=\, 0 \,; \nonumber \\
\mathbb{P}({M}^{t+1}_{n} = {M}^{t}_{n}+k) \,&=\, 
\begin{cases} 
\frac{1}{(m-1)n}\left(\frac{n}{m}-{M}^{t}_{n}\right) & \text{ 
	if } k=2\\
\frac{1}{n}\left(\frac{n}{m}-{M}^{t}_{n}\right) & \text{ if } k =1\\
1- \frac{m}{(m-1)n}\left(\frac{n}{m}-{M}^{t}_{n}\right) & 
\text{ if } k =0\,. \label{eqn:M_change}
\end{cases}
\end{align}

Combining this with \eqref{Ineq1} we obtain the following bound on 
$\LR^{t}_{n}(F_n)$:
\begin{equation}
\label{Ineq3}
\RL^{t}_{n}(F_n)  \geq 
\mathbb{P}\left(M^{t}_{n} 
< n/m\right) \,.
\end{equation}

The idea behind this counting process is that it increases by $1$ 
whenever $R^{i}$ collects a card (i.e. the event $\{R^{t+1} \in U^{t+1}_{n}\}$), and increases by $2$ whenever 
$L^{i}$ collects a card (i.e. the event $\{L^{t+1} \in U^{t+1}_{n}\}$). We know that for $R^{i}$ to collect 
every card in $V_{n}$ would take $O(n\log n)$ time by standard 
coupon collecting. In our case $L^{i}$ is also helping to collect cards 
speeding up this time but we mitigate this additional help by focusing on collecting the cards in $V_{n}$.

From this point on we are interested in the time at which the process 
${M}^{t}_{n}$ first reaches level $n/m$, we now take $m=m(n) =\log n$ for the remainder of this section. 
\begin{lemma}\label{lem:lower_bound}
	Let $\mathcal{T} = \min\{t\,:\, M_n^t \ge n/\log n\}$. Then for any 
	$c>2$,
	\[\mathbb{P}(\mathcal{T} \leq n\log n - n \log 
	\log n 
	-cn)\, \leq \, \frac{\pi^{2}}{6(c-2)^2} \,. \] 
\end{lemma}

\begin{proof}
	Let $\mathcal{T}_{i}$ be the time spent by the process 
	${M}^{t}_{n}$ in each state $i\ge0$. We have $\mathcal{T} = 
	\mathcal{T}_{0} + \mathcal{T}_{1} + \dots + 
	\mathcal{T}_{(n/m)-1}$. Define $p_{i}$ to be the probability that we leave state $i$ after one step of $M_{n}^{t}$, from \eqref{eqn:M_change} we see that 
	\begin{equation}\label{eqn:M_up}
	p_i\,:=\,\mathbb{P}({M}^{t+1}_{n} > {M}^{t}_{n}\,|\, M_n^t = i) \,=\, 
	\frac{m}{(m-1)n}\left(\frac{n}{m}-i\right)\,.
	\end{equation}
	In the standard coupon collector's problem each of the random variables 
	$\mathcal T_i$ has a geometric distribution with success probability 
	$p_i$ (see Section \ref{chpt4:subsec:t2r}). Here, however, we have to take into account the chance that our 
	counting process $M_n$ increments by two, leading it to spend zero time 
	at some state. Note first that 
	\[\mathbb{P}({M}^{t+1}_{n} ={M}^{t}_{n} + 2 \,|\, {M}^{t+1}_{n} > 
	{M}^{t}_{n}) \,=\, \frac{1}{m}\,,\]
	independently of the value of $M_n^t$, taking $m=\log n$ the probability of us skipping a state is tending to $0$ as $n\to \infty$.

	Prior to spending any time in state $i$, the process ${M}^{t}_{n}$ must 
	visit (at least) one of the states $i-1$ or $i-2$. A simple argument 
	shows that 
	\[ \mathbb{P}(\mathcal{T}_{i}>0\, | \,\mathcal{T}_{i-1}>0) \,=\, 1- 
	\frac{1}{m}\,,\quad\text{and}\quad
	\mathbb{P}(\mathcal{T}_{i}>0\, |\, \mathcal{T}_{i-2}>0) \,=\, 1- 
	\frac{1}{m}\left(1-\frac{1}{m}\right) \ge 1-\frac 1m\,. \]
	Therefore $\mathbb{P}(\mathcal T_i>0) \ge 1-\frac{1}{m}$ for all states 
	$i$, and so $\mathcal{T}_i$ stochastically dominates the random 
	variable $\mathcal T'_i$ with mass function
	\begin{equation}\label{eqn:Tprime}
	\mathbb P(\mathcal T'_i = k) = \begin{cases}
	1/m &\quad k=0 \\
	(1-1/m)p_i(1-p_i)^{k-1} &\quad k\ge 1\,.
	\end{cases}
	\end{equation}
	Define  $\mathcal T' = \mathcal{T}'_{0} + \mathcal{T}'_{1} + 
	\dots + \mathcal{T}'_{(n/m)-1}$, it  follows that $\mathbb P(\mathcal T<t) \le \mathbb P(\mathcal T'<t)$ 
	for any $t$.
	
	Substituting $m = m(n)=\log n$ we may bound the expectation and variance of 
	$\mathcal T'$ using the variables $\mathcal T'_{i}$:
	\begin{eqnarray*}
		\E[\mathcal{T}^{\prime}] & = & 
		\sum_{i=0}^{n/m-1} 
		\frac{m-1}{mp_i} =  \left(\frac{m-1}{m}\right)^{2} 
		n\log(n/m)\,\geq\, n\log n  
		-n\log \log n -2n \,;\\
		\Var[\mathcal{T}^{\prime}] & \le & \sum_{i=0}^{n/m-1}  
		\frac{1}{p_i^2} \,\le\, 
		\sum_{i=1}^{n/m}\frac{n^{2}}{i^{2}} 
		\,\leq\, 	\frac{\pi^{2}}{6}n^{2}\,.
	\end{eqnarray*}
	Finally, applying Chebyshev's inequality yields the following for any 
	$c>2$:
	\[	\mathbb{P}\left(\mathcal{T}^{\prime} \leq n\log(n) - n \log \log n 
	-cn\right) 
	\,\leq\, 
	\mathbb{P}\left(|\mathcal{T}^{\prime} -\E[\mathcal{T}^{\prime}]\,|\, \geq  
	(c-2)n \right)
	\,\leq\,  \frac{\pi^2}{6(c-2)^2}\,.  \qedhere		\]

\end{proof}

\subsubsection{Proof of the Lower Limit in Theorem \ref{chpt5:thm:cutoff}}	
Lemma~\ref{lem:lower_bound} 
quickly leads to a proof of the lower bound in Theorem~\ref{chpt5:thm:cutoff}. Setting $t = n\log n - n \log \log n $, and $c>2$ we obtain
\begin{align*}
\lVert \RL_n^{t-cn} -\pi_n \rVert_{\tiny \textnormal{TV}}\, &\ge\, 
\RL_n^{t-cn}(F_n) - \pi_n(F_n) 
\,\ge\, \mathbb{P}\left(M^{t-cn}_{n} < n/m(n)\right)- 1/m(n)\\
&= \, \mathbb{P}(\mathcal{T} >t-cn) - \frac{1}{\log n } 
\,\ge\,  1 -\frac{\pi^{2}}{6(c-2)^2} - \frac{1}{\log n }.
\end{align*}
Therefore taking a limit in $n$, we recover $\liminf_{n\to\infty} \lVert \RL_n^{t-cn} -\pi_n \rVert_{\textnormal{TV}} \geq 1 - \frac{\pi^{2}}{6(c-2)^{2}}$ as required.

\paragraph{}
This finishes our proof of Theorem \ref{chpt5:thm:cutoff}, showing the unbiased one-sided transposition shuffle exhibits a total variation cutoff at time $n\log n$. The mixing time for the unbiased one-sided transposition shuffle is therefore half as slow as the mixing time of the random transposition shuffle. This fact fits our intuition that restricting the choice of our left hand should slow down the randomisation of our deck. We now explore the behaviour of a new class of one-sided transposition shuffles where our left hand retains its behaviour but the right hand may be biased on $[n]$.

\section{Biased One-sided Transposition shuffles}
\label{chpt5:sec:biased}

The unbiased one-sided transposition shuffle had our right hand choosing uniformly from the support $[n]$. In this section we allow the right hand to choose a position according to a weighted distribution on $[n]$, and we call these \emph{biased one-sided transposition shuffles}. We generalise the work of Section \ref{chpt5:sec:algebra} to compute the eigenvalues for our biased shuffles and show that within this class certain natural choices for the weights lead to shuffles which exhibit a cutoff in total variation distance.

\begin{defn}
	\label{chpt5:def:biasedLR}
	Given a \emph{weight function} $w:\mathbb{N} \rightarrow (0,\infty)$, let $\NW = \sum_{i=1}^{n} w(k)$ denote the cumulative weight up to $n$. Then the \emph{biased one-sided transposition shuffle} $\BLR_{n,w}$ is the random walk on $S_n$ generated by the following distribution on transpositions:
	\begin{eqnarray}
	\BLR_{n,w}(\tau) = 
	\begin{cases}
	\frac{w(j)}{\NW}\cdot\frac{1}{j} &  \text{if } \tau = (i\,j) 
	\text{ for some } 
	1\leq i 
	\leq j 
	\leq n\\
	0 & \text{otherwise.}
	\end{cases}
	\end{eqnarray}
	Note that if $w(j)=1$ for all $j$, we recover the unbiased one-sided transposition shuffle.
\end{defn}

According to a biased one-sided transposition shuffle with weight function $w$ the random variables $R^{t}, L^{t}$ defined by the choices of our right and left hands follow the distributions,
\begin{eqnarray}
\mathbb{P}(R^{t}=j) & = &\frac{w(j)}{\NW}, \text{ for } 1 \leq j 
\leq n\ \\
\mathbb{P}(L^{t}=i \, | \, R^{t}=j) & = &\begin{cases}
\frac{1}{j} & \textnormal{ if } i \leq j \\
0 & \textnormal{ otherwise }
\end{cases}.
\end{eqnarray}

Notice the left hand $L^{t}$ still chooses a position uniformly on the set $\{1,\ldots,R^{t}\}$. Importantly in Definition \ref{chpt5:def:biasedLR} the weight of each position $w(j)$ may only 
depend on $j$ and not the size of the deck $n$. 
This setup preserves the recursive algebraic structure identified in Section \ref{chpt5:sec:algebra}, and allows us to recover the eigenvalues of the biased shuffles by minor modifications to our analysis. 

\subsubsection{Lifting Eigenvectors for the Biased One-sided Transposition Shuffle}
Our first step is to turn the shuffle $\BLR_{n,w}$ into an element of our group algebra.

\begin{defn}
	\label{chpt5:def:BLRop}
	Let $n \in \mathbb{N}$.	The biased one-sided transposition shuffle on $n$ cards, with bias $w(j)$, may be viewed as the following element of the group algebra $\mathfrak{S}_n$:
	\begin{eqnarray}
	\ABLR_{n,w}:= \sum_{1 \leq i \leq j \leq n}\frac{w(j)}{j} (i\,j) .
	\label{chpt5:eqn:biasedalgebrashuffle}
	\end{eqnarray}
	Note that to form $\ABLR_{n,w}$ we have scaled $\BLR_{n,w}$ by a factor of $\NW$.
\end{defn}

The eigenvalues of the biased one-sided transposition shuffle $\BLR_{n,w}$ may be recovered from the eigenvalues of the element $\ABLR_{n,w}$ acting on the simple Specht modules $S^{\lambda}$. Thus, as before we may focus our attention on the action of our shuffle as an element of the group algebra. The biased one-sided transposition shuffle inherits the recursive structure present in the unbiased case. Compare 
\begin{eqnarray}
\ABLR_{n+1,w} - \ABLR_{n,w} = \frac{w(n+1)}{n+1} \sum_{1\leq i \leq n+1} (i \hspace{0.2cm} n+1) \label{chpt5:eqn:biasedcompare}
\end{eqnarray}
to equation \eqref{chpt5:eqn:lrdifference}, we can see that the only change is a new factor of $w(n+1)$. If our weight function is allowed to depend on things other than the position $j$ then we run into issues with our lifting because $w(n+1)$ may no longer solely be a function of $(n+1)$.  Now we establish an analogue of Theorem \ref{chpt5:thm:master} for our biased shuffles.

\begin{thm}
	\label{chpt5:thm:biasedmaster}
	Given $n \in \mathbb{N}$, we have
	\begin{equation}
	\ABLR_{n+1,w} \circ \sh_{a} - \sh_{a} \circ \ABLR_{n,w} = 
	\frac{w(n+1)}{n+1}\sh_{a} + 
	\frac{w(n+1)}{n+1}\sum_{1 \leq b \leq n} \sh_{b}\circ\Theta_{b,a} \,. 
	\label{chpt5:eqn:biasedmastereq}
	\end{equation}
\end{thm}
\begin{proof}
	This follows the same proof as Theorem \ref{chpt5:thm:master} up to changes in constants.
\end{proof}

If our bias is uniform, i.e. $w(j)=1$, then we recover exactly the equation used for the unbiased one-sided transposition shuffle. Overall we can see that our bias affects Theorem \ref{chpt5:thm:master} minimally. From here we may use the previously established lifting operators $\pro_{a}^{\lambda, \mu}$ to skip to the conclusion of lifting eigenvectors. We restate Theorem \ref{chpt5:thm:lift} for our biased one-sided transposition shuffle.

\begin{thm}[Theorem 49 of \cite{dieker2018spectral}]
	\label{chpt5:thm:biasedlifting}

	Let $\lambda \vdash n$, and  $a \in \{1,2,\ldots,l(\lambda)+1\}$. Take $i\in [n]$ such that  $1 \leq i \leq a$ and set $\mu = \lambda +e_{i}$. Then,
	\begin{eqnarray}
	\ABLR_{n+1,w} \circ \pro_{a}^{\lambda,\mu} - 
	\pro_{a}^{\lambda,\mu} \circ \ABLR_{n,w} = \frac{w(n+1)(2 + 
		\lambda_{i} -i)}{n+1}\pro_{a}^{\lambda,\mu} .
	\end{eqnarray}
	In particular if $v \in S^{\lambda}$ is an eigenvector of $\ABLR_{n,w}$ 
	with 
	eigenvalue $\varepsilon$, then  
	$\pro_{a}^{\lambda,\mu}(v)$ is an eigenvector of 
	$\ABLR_{n+1,w}$ with eigenvalue 
	\begin{eqnarray}
	\varepsilon + \frac{w(n+1)(2 +\lambda_{i} -i )}{n+1}.
	\end{eqnarray}
\end{thm}
\begin{proof}
	This follows from Theorem \ref{chpt5:thm:lift} with changes in constants for our biased shuffles.
\end{proof}

\begin{lemma}
	\label{biasedCompute}
	The eigenvalues for the biased one-sided transposition shuffle $P_{n,w}$ on $n$ cards 
	are indexed by standard Young tableaux of size $n$. Moreover, the 
	eigenvalue corresponding to a tableau $T$ is given by
	\begin{eqnarray*}
		\label{eigenvaluesum2}
		\eig(T) = \frac{1}{\NW}\sum_{\substack{\textnormal{boxes}\\ 
				(i,j) \in T}} 
		\frac{j-i+1}{T(i,j)} \cdot w(T(i,j))
	\end{eqnarray*}
\end{lemma}

\subsubsection{Mixing Times of Biased One-sided Transposition Shuffles}

We have recovered the eigenvalues for the biased one-sided transposition shuffles. To find the mixing time of the shuffle $\BLR_{n,w}$, we need some extra conditions on the weight function $w$.
We focus on a class of weight functions with form $w(j) = 
j^{\alpha}$ for $\alpha \in \mathbb{R}$; we denote the resulting shuffle as $\BLR_{n,\alpha}$, and 
write $\NA$ in place of $\NW$. For $\alpha=0$ we recover our original one-sided 
transposition shuffle $\LR_{n}$, while if $\alpha>0$ ($\alpha<0$) the right 
hand is biased towards the top (respectively, bottom) of the deck. Restating Lemma \ref{biasedCompute} for this class of weight functions we find the eigenvalue for $\BLR_{n,\alpha}$ given by tableau $T$ to be,
\begin{eqnarray}
\eig(T) = \frac{1}{\NA}\sum_{\substack{\textnormal{boxes}\\ 
		(i,j) \in  T}} (j-i+1)\cdot T(i,j)^{\alpha-1} = \frac{1}{\NA}\sum_{m=1}^{n} 
T(m) \, m^{\alpha-1},			\label{chpt5:eqn:eigenvaluesumalpha}
\end{eqnarray}
where $T(m)$ is the index $j-i+1$ of the box containing value $m$ in $T$.

\paragraph{}
In the following sections we analyse the the mixing time of biased one-sided transposition shuffles generalising work from Sections \ref{chpt5:sec:upperbound} and \ref{chpt5:sec:lowerbound}. The conlusion of this section is the biased one-sided transposition shuffles exhibit a cutoff for all real choices of $\alpha$. 

\begin{thm}
	\label{chpt5:thm:biasedbounds}
	Define the time $t_{n,\alpha}$ as, 
	\[t_{n,\alpha}= \begin{cases}
	\NA/n^{\alpha} & \textnormal{ if } \alpha \leq 1 \\
	\NA/N_{\alpha-1}(n) & \textnormal{ if } \alpha \geq 1 
	\end{cases}. \]
	The biased one-sided transposition shuffle $\BLR_{n,\alpha}$ satisfies the following bounds for any $c_{1} >5/2, c_{2} >\max(2, 3-\alpha)$:
	
	\begin{eqnarray*}
		\limsup_{n\rightarrow \infty}\, \lVert 
		\BLR_{n,\alpha}^{t_{n,\alpha}\left(\log n + c_{1}\right)} -\pi_{n} 
		\rVert & 
		\leq & 	Ae^{-2c_{1}} 
		\textnormal{ 	for a universal constant } A, \textnormal{ for all } 
		\alpha \\
		\text{ and } \quad 
		\liminf_{n\to\infty}\, \lVert 
		\BLR_{n,\alpha}^{t_{n,\alpha}\left(\log n - \log \log n - 
			c_{2}\right)} 
		-\pi_{n} \rVert & 
		\geq 
		& \begin{cases}
			1 - \frac{\pi^{2}}{6(c_{2}-3+\alpha)^{2}} & \textnormal{ if } \alpha 
			\leq 1 \\
			1 - \frac{\pi^{2}}{6(c_{2}-2)^{2}} & \textnormal{ if } \alpha \geq 1 
		\end{cases} 
	\end{eqnarray*}
	Thus, the biased one-sided transposition shuffle  exhibits a total variation cutoff at time $t_{n,\alpha}\log n$ for all $\alpha$ with a window of size $t_{n,\alpha}\log \log n$. 
\end{thm}

The asymptotics of the times $t_{n,\alpha}$ as $n\to\infty$ for the one-sided transposition shuffle are summarised in Table \ref{chpt5:table:biasedtime}

\begin{table}[H]
	\begin{tabular}{c|c|c|c|c}
		& $\alpha \in (-\infty,-1)$ & $\alpha =-1$ & $\alpha \in (-1,1]$ & $\alpha \in (1,\infty)$
		\\\hline
		$t_{n,\alpha}\log n$ & $\zeta(-\alpha) n^{-\alpha} \log n$ & $n(\log n)^{2}$ & $ \frac{1}{1+\alpha} n \log n$ & $\frac{\alpha}{1+\alpha} n \log n$ 
	\end{tabular}
	
	\caption[Asymptotics of $t_{n,\alpha} \log n$ as $n\to\infty$]{Asymptotics of $t_{n,\alpha} \log n$ as $n\to\infty$.}
	\label{chpt5:table:biasedtime}
\end{table}

The fastest mixing time of a biased one-sided transposition shuffle is obtained when $\alpha = 1$; using this weight function the shuffle is constant on the conjugacy class 
of transpositions, with probability close to that of the random 
transposition shuffle, $P_{n,1}((i \, j)) = 2/(n(n+1))$.
We obtain a mixing time of $t_{n,1} \sim (n/2)\log n$ which agrees with that of 
the random transposition shuffle. The mixing time increases as $\alpha$ moves 
away from $1$ in either direction but as $\alpha \to \infty$ the time
stays bounded above by $n\log n$ whereas if $\alpha \to -\infty$ then the mixing time is unbounded. As $\alpha \to -\infty$ our right hand is choosing positions near to the bottom of the deck frequently, and since the left hand is restricted to the range $\{1,\ldots, \righthand\}$, this leads to more and more mass being placed on the identity element as $\alpha \to\infty$, thus slowing the mixing time.

\subsection{Cutoff for Biased One-sided Transposition Shuffles with $\alpha \leq 1$}
\label{chpt5:subsec:smallalpha}

The proof of a total variation cutoff for the biased one-sided transposition shuffle with $\alpha \leq1$ follows from generalisations of Sections~\ref{chpt5:sec:upperbound} and \ref{chpt5:sec:lowerbound}. For the upper bound we establish bounds on the eigenvalues of $\BLR_{n,\alpha}$ for large and small partitions, and reduce the analysis to a previously bounded summation. The lower bound follows from the same argument as Section \ref{chpt5:sec:lowerbound} with careful attention paid to how the weighting affects the probabilities of $\lefthand, \righthand$.

\subsubsection{Upper Bound}
\label{GENUPPERBOUNDSMALLALPHA}
First of all we once again use Lemma~\ref{chpt2:thm:classicL2} to form an upper bound on the total variation distance of $\BLR_{n,\alpha}$ from $\pi_{n}$.
Furthermore, for $\alpha\leq 1$ every result of Section~\ref{chpt5:sec:eigenvalues} holds after suitable adjustments to account for the new eigenvalues. In particular, Lemma \ref{lem:sumtranspose} now becomes: for any $T \in \SYT(\lambda)$ we have $\eig(T) + \eig(T^{\prime}) = 2 N_{\alpha-1}(n) /\NA $. Applying the same analysis as equation \eqref{eqn:UB_eigs} we reduce our bound to:

\begin{eqnarray}
4\lVert \BLR_{n,\alpha}^{t} - \pi_{n} \rVert^{2}_{\textnormal{TV}} & \leq & \sum_{\substack{ \lambda \vdash n \\ \lambda \neq (n)}} \sum_{T \in \SYT(\lambda)} d_{\lambda} \left(\eig(T) \right)^{2t} \nonumber\\ 
& \leq & \left(\eig(T_{(1^{n})})\right)^{2t} + 2\sum_{\substack{\lambda : \eig(T_{\lambda}^{\rightarrow}) \geq 0 \\ \lambda \neq (n)}} d_{\lambda}^{2} \eig(T_{\lambda}^{\rightarrow})^{2t}.\label{chpt5:eqn:biasedbigalpha}
\end{eqnarray}

We bound the eigenvalues of large and small partitions separately.

\begin{lemma}
	\label{lem:genbound1}
	Let $\lambda \vdash n$ with $\lambda_{1} = n-k$. Then the eigenvalue $\eig(\best)$ for the shuffle $\BLR_{n,\alpha}$ with $\alpha \leq1$ may be bounded as follows:
	\[
	\eig(\best) \,\le\, 
	\begin{cases}
	1 - \frac{(n-k+1)kn^{\alpha}}{n\NA} & \quad \text{if $k \le n/4$} \\
	1 - \frac{k n^{\alpha}}{2\NA} & \quad \text{if $k >n/4$.} \\
	\end{cases}
	\]
\end{lemma}

\begin{proof}
	For $k\le n/2$, the maximum partition in the dominance order of 
	partitions with $\lambda_{1}=n-k$ is $(n-k,k)$, and so $\eig(\best) \le \eig\left(T_{(n-k,k)}^{\rightarrow}\right)$. The eigenvalue of $T_{(n-k,k)}^{\rightarrow}$ may be calculated by summing over the rows of the partition $(n-k,k)$ and 
	using equation \eqref{eigenvaluesum2}, as follows:
	\begin{eqnarray*}
		\NA \, 
		\textnormal{eig}\left(T_{(n-k,k)}^{\rightarrow}\right) & = &
		\sum_{m=1}^{n-k} m^{\alpha} + \sum_{m=n-k+1}^{n} (m-n+k-1) m^{\alpha-1} \\
		& = & \sum_{m=1}^{n} m^{\alpha} -(n-k+1) \sum_{m=n-k+1}^{n} m^{\alpha-1}
		\\
		& \leq & \NA -\frac{(n-k+1)kn^{\alpha}}{n}.
	\end{eqnarray*}
	This immediately proves the desired inequality for $k\le n/4$. The above bound also holds for $k \in (n/4, n/2]$, and in this case
	\[\NA -\frac{(n-k+1)kn^{\alpha}}{n} \leq\NA - \frac{kn^{\alpha}}{2} \]
	giving us our required bound.
	
	For $k>n/2$ we once again need to bound 
	$\eig\left(T_{(n-k,\star)}^{\rightarrow}\right)$. Letting $\nu 
	=(n-k,\star)$ for ease of notation we calculate as follows:
	\begin{eqnarray}
	\NA \, 
	\eig(T_{\nu}^{\rightarrow}) & = & \sum_{j=1}^{n-k} 
	(j-1+1)^{\alpha}  + \sum_{i=2}^{l(\nu)} 
	\sum_{j=1}^{\nu_{i}} (j-i+1)  \frac{((i-1) (n-k)
		+j)^{\alpha}}{(i-1) (n-k) +j} \nonumber\\
		& = & \NA - \sum_{i=2}^{l(\nu)} 
	\sum_{j=1}^{\nu_{i}} (i-1) (n-k+1)  
	\frac{((i-1) (n-k)+j)^{\alpha}}{(i-1) (n-k)+j}\nonumber \\
	& \leq & \NA - 
	\frac{n^\alpha}{n}\sum_{i=2}^{l(\nu)} 
	\sum_{j=1}^{\nu_{i}} \left((i-1)\cdot (n-k+1)\right) \nonumber\\
	& \leq & \NA -\frac{(n-k+1)n^{\alpha}}{n} 
	\sum_{i=2}^{l(\nu)} (i-1)  \nu_{i} \, \label{chpt5:eqn:smallalphastep1}. 
	\end{eqnarray}
	By definition of the partition $\nu$, each row but the last has size $n-k$, and the final row has size $\nu_{l(\nu)} = n-(l(\nu)-1)(n-k)$. In addition, since $l(\nu)$ is equal to the ceiling of $n/(n-k)$ we may write $l(\nu) = n/(n-k)+x$ for some $0\le x<1$.	Substituting these values into equation \eqref{chpt5:eqn:smallalphastep1} we obtain:
	\begin{eqnarray*}
		\NA \, 
		\eig(T_{\nu}^{\rightarrow}) 
		& \leq & \NA
		-\frac{(n-k)n^{\alpha}}{n} 
		\frac{(l(\nu)-1)( 2n - (n-k)l(\nu))}{2}\\
		& = & \NA
		-\frac{n^{\alpha}}{2n} (n - (1-x)(n-k)) (n-x(n-k))\\
		& = & \NA
		-\frac{n^{\alpha}}{2n} (nk+x(1-x)(n-k)^2)\\
		& \leq & 
		\NA -\frac{kn^{\alpha}}{2}\,. 
	\end{eqnarray*}	
\end{proof}

Using these eigenvalue bounds we complete the proof of the upper bound for the biased one-sided transposition shuffles with $\alpha \leq 1$. Splitting the sum in \eqref{chpt5:eqn:biasedbigalpha} around the big and small partitions we find,

\begin{eqnarray}
4\lVert \BLR_{n,\alpha}^{t} - \pi_{n} \rVert^{2}_{\tiny \textnormal{TV}} 
&\leq &  \left(\eig(T_{1^{n}})\right)^{2t} + 2\sum_{k=1}^{n/4} 
{n\choose k}^{2} k!	\left(1 - 
\frac{2(n-k+1)kn^{\alpha}}{2n\NA}\right)^{2t} \nonumber\\ &+& 
2\sum_{k>n/4}^{n-1} {n\choose k}^{2} k!
\left( 1 - 	\frac{kn^{\alpha}}{2\NA}\right)^{2t}. \nonumber
\end{eqnarray}
The singular term $\eig(T_{1^{n}})$ may be seen to be tending to $0$ as $n\to$ $\infty$ at time $t=t_{n,\alpha}(\log n  +  c)$,
\[\lim_{n\to\infty} \eig(T_{1^{n}})^{2t} = \lim_{n\to\infty}\left(\frac{1}{\NA} \sum_{i=1}^{n} (2-i)^{\alpha-1} \right)^{2t} \leq\lim_{n\to\infty} \left( \frac{n}{\NA}\right)^{2t_{n,\alpha}(\log n  +c)} = 0\]
We are left with two sums to control, taking $t = t_{n,\alpha} (\log n + c)$ these sums may be reduced too:
\begin{eqnarray}
e^{-2c}\sum_{k=1}^{n/4}  \frac{n^{\frac{2(k(k-1))}{n}}}{k!} + \sum_{k>n/4}^{n-1} {n\choose k}^{2} k! n^{-k} e^{-k c}
\label{chpt5:eqn:biasedsmallalphasum}
\end{eqnarray}
The first summation in \eqref{chpt5:eqn:biasedsmallalphasum} was shown to be bounded by a universal constant by Diaconis \cite[Chapter 3D Theorem 5]{Diaconis1988}. The ratio between consecutive terms in the  second summation is decreasing in $k$ and less than $1$ at its start ($k=n/4$) if $c>2.5$. Thus we may bound the summation by bounding it by $n$ times its first term; this reduces to bounding

\begin{eqnarray}
\sum_{k>n/4}^{n-1} {n\choose k}^{2} k! n^{-k} e^{-k c} \leq  n {n\choose n/4}^{2} ( n/4)!  n^{-n/4} e^{-(n/4) c} .\label{chpt5:eqn:stirling}
\end{eqnarray}
We now compute the asymptotics of the binomial and factorial terms using Stirling's approximation, these results are given in Table \ref{chpt4:table:stirlings}. 
\begin{table}[H]

	\begin{adjustbox}{width=\columnwidth,center}
		\renewcommand{\arraystretch}{1.3}
		\begin{tabular}{c|c|c}
			Term    & ${n\choose \delta n}^{2}$ & $(\delta n)!$  \\\hline
			Stirling's approximation & $(\delta^{-2\delta n} (1-\delta)^{-2(1-\delta)n} ) / (2\pi\delta(1-\delta)n)$ & $\sqrt{2\pi \delta n}\left(\delta n / e\right)^{\delta n}$ \\\hline
			$\log (\textnormal{Stirling's approximation})$ & $-\log(2\pi \delta) -2\delta n\log (\delta) - 2(1-\delta)n\log (1-\delta)$ & $\frac{1}{2}\log(2\pi\delta n) + \delta n \log(\delta n) - \delta n $
		\end{tabular}
		
	\end{adjustbox}
	\caption[Asymptotics of binomial and factorial terms using Stirling's approximation]{Asymptotics of binomial and factorial terms computed using Stirling's approximation.}
	\label{chpt4:table:stirlings}
\end{table}
Using the results in Table \ref{chpt4:table:stirlings} we find the asymptotics of the logarithm of \eqref{chpt5:eqn:stirling} to be:
\begin{eqnarray}
-(n/4)\left( c+ 1 +  2 \log (0.25) +  6\log (0.75) \right) + o(n).\label{chpt5:eqn:stirling2}
\end{eqnarray}
The coefficient of the leading order term
\[c+ 1 +  2 \log (0.25) +  6\log (0.75) \]
is positive if $c>5/2$.	Therefore for $c>5/2$, the term \eqref{chpt5:eqn:stirling2} is tending to $-\infty$ as $n\to\infty$. Hence, the second summation in \eqref{chpt5:eqn:biasedsmallalphasum} is tending to $0$ as $n\to\infty$. Putting together all of the bounds above we find,
\[\limsup_{n\to \infty} \lVert \BLR_{n,\alpha}^{t_{n,\alpha}(\log n+ + c)} - \pi_{n} \rVert_{\tiny \textnormal{TV}} \leq Ae^{-c} \textnormal{ for some universal constant $A$} .\]
This completes the upper bound for $\alpha \leq 1$ present in Theorem \ref{chpt5:thm:biasedbounds}.

\subsubsection{Lower Bound}
We use a coupon-collecting argument as in Section~\ref{chpt5:sec:lowerbound}, once again letting $V_n = 
\{n-n/m+1,\dots,n-1,n\}$ with $m =\log n$, and considering the set $F_n = \{\sigma\in S_n \, | \, \sigma$ has at least 1 fixed point in $V_n\}$. We have seen previously that under the uniform distribution we have $\pi_{n}(F_{n}) \leq 1/m = 1/\log n$. For the biased one-sided transposition shuffle we modify the bounds \eqref{eqn:bound_minus2} and \eqref{eqn:bound_minus1} as follows, using the inequality 
$k^\alpha \le (\frac{m}{m-1})^{1-\alpha}n^\alpha$ for all $k\in V_n$ (which holds for all $\alpha\le 1$):
\begin{align}
\mathbb{P}(L^{t+1}\in \hat U_n^t) &\,=\, \sum_{k\in \hat 
	U_n^t} \frac{w(k)}{\NA}\frac{k-(n-n/m)}{k} +  \sum_{k\in 
	V_n\setminus \hat U_n^t}\frac{w(k)}{\NA}\frac{|U_n^t|}{k} \nonumber \\
&\,\le\, \frac{ (\frac{m}{m-1})^{1-\alpha}n^\alpha}{\NA}\frac{\sum_{k=1}^{|U_n^t|} k + (n/m -|U_n^t|)|U_n^t|}{(n-n/m)} 
\,\le\,\frac{(\frac{m}{m-1})^{1-\alpha}n^\alpha|U_n^t|}{\NA(m-1)}\,; \label{eqn:bound_minus2biased} \\
\mathbb{P}(R^{t+1}\in U_n^t) &\,\le\, \frac{(\frac{m}{m-1})^{1-\alpha}n^\alpha |U_n^t|}{\NA}\,. \label{eqn:bound_minus1biased}
\end{align}
Using these as before we construct a counting process $M_{n,\alpha}^t$
which stochastically dominates the number of collected cards $|V_n\setminus 
U_n^{t}|$ at all times: 
\begin{align}
M_{n,\alpha}^0 \,&=\, 0 \,; \nonumber \\
\mathbb{P}(M_{n,\alpha}^{t+1}= M_{n,\alpha}^t+k) \,&=\, 
\begin{cases} 
\frac{1}{(m-1)} \left(\frac{\left(\frac{m}{m-1}\right)^{1-\alpha} n^{\alpha} |U^{t}_{n}|}{\NA} \right) & \text{ 
	if } k=2\\
\left(\frac{\left(\frac{m}{m-1}\right)^{1-\alpha} n^{\alpha} |U^{t}_{n}|}{\NA} \right) & \text{ if } k =1\\
1- \frac{m}{(m-1)} \left(\frac{\left(\frac{m}{m-1}\right)^{1-\alpha} n^{\alpha} |U^{t}_{n}|}{\NA} \right) & 
\text{ if } k =0\,. \label{chpt5:eqn:biassmallM_change}
\end{cases}
\end{align}
Our counting process $M_{n,\alpha}^{t}$ gives a lower bound on the probability of being in $F_{n}$ as follows,
\[\BLR_{n,\alpha}^{t}(F_{n}) \geq \mathbb{P}(M_{n,\alpha}^{t} < n/m).\]
We are now interested in the time it takes $M_{n,\alpha}^{t}$ to pass value $n/m$.

\begin{lemma}
	\label{chpt5:lem:smallalphalowerbound}
	Let $\mathcal{T} = \min \{t \, : \, M_{n,\alpha}^{t} \geq n / \log n\}$. Then for any $c>3-\alpha$,
	\[\mathbb{P}(\mathcal{T} \leq t_{n,\alpha}(\log n - \log \log n - c)) \geq 1 - \frac{\pi^{2}}{6(c -(3-\alpha))^{2}}.\]
\end{lemma}

\begin{proof}
	Construct the time $\mathcal{T} = \mathcal{T}_{0} + \mathcal{T}_{1} + \ldots + \mathcal{T}_{n/m -1}$ with $T_{i}$ being the time spent in state $i$. Denote the escape probability of state $i$ as $p_{i}$, the expression for $p_i$ modified from \eqref{eqn:M_up} becomes 
	\[ p_i = \left(\frac{m}{m-1}\right)^{2-\alpha}\frac{n^\alpha}{\NA} \left( \frac nm-i\right) \,, \]
	and this is easily checked to be strictly less than one for all values of $\alpha \le 1$ if  $n$ is sufficiently large.
	The remainder of the analysis mirrors the unbiased case: using the new expression for $p_i$ the distribution of the random variable $\mathcal T'_i$ is exactly as given in \eqref{eqn:Tprime}, that is
	\begin{equation*}
	\mathbb P(\mathcal T'_i = k) = \begin{cases}
	1/m &\quad k=0 \\
	(1-1/m)p_i(1-p_i)^{k-1} &\quad k\ge 1\,.
	\end{cases}
	\end{equation*}
	Setting $\mathcal{T}^{'} = \mathcal T'_0 + \ldots + \mathcal T'_{n/m -1}$ then $\mathcal{T}$ stochastically dominates $\mathcal T'$. Analysing the expectation and variance of $\mathcal T'$ we arrive at: 
	\begin{eqnarray*}
		\E[\mathcal{T}^{\prime}]  & = &
		\sum_{i=0}^{n/m-1} \,
		\frac{m-1}{mp_i} = \frac{\NA}{n^{\alpha}}\left(1 - \frac{1}{m}\right)^{3-\alpha}  \sum_{i=1}^{n/m} \frac{1}{i} \ge  \frac{\NA}{n^\alpha}(\log n  -\log \log n -(3-\alpha)) \,;\\
		\Var[\mathcal{T}^{\prime}] & \le &  \sum_{i=0}^{n/m-1} \, \frac{1}{p_i^2}  = \frac{\NA^2}{n^{2\alpha}} \left(1 - \frac{1}{m}\right)^{2(2-\alpha)}\sum_{i=1}^{n/m} \frac{1}{i^{2}}\le 
		\frac{\NA^2}{n^{2\alpha}} \frac{\pi^{2}}{6} \,.
	\end{eqnarray*}
	An application of Chebyshev's inequality give us:
	\begin{eqnarray*}
		\mathbb{P}(\mathcal T' \leq t_{n,\alpha}(\log n -\log \log n -(3-\alpha))) &\leq &\mathbb{P}(|\mathcal T' - \E[\mathcal T']| \geq (c- (3-\alpha))n)\\
		& \leq & \frac{\pi^{2}}{6(c - (3-\alpha))^{2}} .
	\end{eqnarray*}
\end{proof}
Lemma \ref{chpt5:lem:smallalphalowerbound} quickly leads to a lower bound, taking $c> 3-\alpha$ we have
\begin{eqnarray*}
	\lVert \BLR_{n,\alpha}^{t_{n,\alpha}(\log n - \log \log n -c)} -\pi_{n} \rVert_{\tiny \textnormal{TV}} &\geq &\mathbb{P}(\mathcal{T} > t_{n,\alpha}(\log n - \log \log n -c)) - \pi_{n} (F_{n})\\
	& \geq &1 - \frac{\pi^{2}}{6(c-(3-\alpha))^{2}} - \frac{1}{\log n}.
\end{eqnarray*}
Taking limits in $n$ recovers the limit present in Theorem \ref{chpt5:thm:biasedbounds}.

\subsection{Cutoff for Biased One-sided Transposition Shuffles with $\alpha \geq 1$}
\label{chpt5:subsec:bigalpha}

The proof of a total variation cutoff for biased one-sided transposition shuffle with $\alpha \geq 1$  requires more work than the $\alpha \leq 1$ case. To prove the upper limit we construct  a new special tableau, $T_{\lambda}^{\searrow}$, which bounds all the eigenvalues associated to  partition $\lambda$. This new tableau is then used to reduce our analysis to summations which we have previous seen to be bounded at the correct mixing time.
The lower bound follows from the same argument as Section \ref{chpt5:sec:lowerbound} up to the replacement of a normalising factor.

\subsubsection{Upper Bound}
Unlike the  $\alpha \leq 1$ previous case not all results of Section \ref{chpt5:sec:eigenvalues} generalise to the bias $\alpha \geq 1$. The results which do generalise include Lemma \ref{seq} and \ref{flip}, we restate these below with all eigenvalues calculated for a biased one-sided transposition shuffle with $\alpha \geq 1$.

\begin{lemma}
	Let $\lambda\vdash n$. For any $T \in \SYT(\lambda)$ we have the following 
	inequality :
	\begin{eqnarray}
	\eig(\best) \leq \eig(T) \leq 
	\eig(\worst).
	\end{eqnarray}
\end{lemma}
\begin{proof}
	This follows from the proof of Lemma \ref{seq} with the roles of $\best$ and $\worst$ reversed.
\end{proof}

\begin{corollary}
	Let $\lambda \vdash n$, and suppose we 
	have 
	$\eig(\best) \leq 0$, then we 
	have 
	\begin{eqnarray}
	\eig(T_{\lambda^{\prime}}^{\downarrow}) \geq 
	|\eig(\best)| \geq  0.
	\end{eqnarray}
\end{corollary}
\begin{proof}
	This follows from the proof of Lemma \ref{flip} with the roles of $\best$ and $\worst$ reversed.
\end{proof}

Using the above lemmas allows us to reduce our analysis to looking at $T_{\lambda}^{\downarrow}$ but we may not solely focus on partitions at the top of the dominance order like before. Applying the same analysis as equation \eqref{eqn:UB_eigs} we reduce our bound to:

\begin{eqnarray}
4\lVert \BLR_{n,\alpha}^{t} - \pi_{n} \rVert^{2}_{\textnormal{TV}} & \leq & \sum_{\substack{ \lambda \vdash n \\ \lambda \neq (n)}} \sum_{T \in \SYT(\lambda)} d_{\lambda} \left(\eig(T) \right)^{2t} \nonumber\\ 
& \leq & \left(\eig(T_{1^{n}})\right)^{2t} + 2\sum_{\substack{\lambda : \eig(T_{\lambda}^{\downarrow}) \geq 0 \\ \lambda \neq (n)}} d_{\lambda}^{2} \eig(T_{\lambda}^{\downarrow})^{2t}.\label{chpt5:eqn:biasedsmallalpha}
\end{eqnarray}

To bound the eigenvalue of $\worst$ when $\lambda_{1}=n-k$ we introduce a new  special tableau for a Young diagram $\lambda$.

\begin{defn} 
	Let $\lambda \vdash n$. Define the $\lambda$-tableau $T_{\lambda}^{\searrow}$  by 
	filling in the diagonals of $\lambda$ from left to right, with each 
	diagonal filled in top to bottom. For example, the tableaux $T_{(4^{4})}^{\searrow}, \, T_{(4,3,1^{2})}^{\searrow}, \, T_{(6,4)}^{\searrow}$ are given below,
	\[	\ytableausetup{mathmode,baseline,aligntableaux=center,boxsize=1.2em} 
	\begin{ytableau} 
	7 & 11 & 14 & 16 \\
	4 &  8 & 12 & 15 \\
	2 & 5  & 9 & 13 \\ 
	1 & 3 & 6 & 10\end{ytableau} \hspace{2cm} 
	\begin{ytableau} 
	4 & 6 & 8 & 9 \\
	3 &  5 & 7  \\
	2  \\ 
	1 \end{ytableau}  \hspace{2cm} 
	\begin{ytableau} 
	2& 4& 6& 8& 9& 10\\
	1& 3& 5& 7
	\end{ytableau} \]
\end{defn}

\begin{lemma}
	\label{lem:eigbound}
	For $\alpha \geq1$, and  $\lambda \vdash n$ with $\lambda_{1} =n-k$ we 
	have, 
	\[\eig(T_{\lambda}^{\downarrow}) \leq 
	\eig\left(T_{(n-k,\star)}^{\searrow}\right) .\]
\end{lemma}
\begin{proof}

	The eigenvalue associated to the tableau $\worst$ is given by,
	\[\NA \cdot \eig(T) = 	\sum_{m=1}^{n} T(m) m^{\alpha-1}.\]
	Now consider all indexes  $\worst(m)$ for $m \in [n]$ including repeats and order them from smallest to largest as $c_{i}$ for $i \in [n]$.  The fact that $\alpha \geq 1$ allows us to upper bound our eigenvalue in the following 
	way: \[\NA \,\eig(\worst)  = \sum_{m=1}^{n} \worst(m) \, 
	m^{\alpha-1} \leq \sum_{m=1}^{n} c_{m} \, 
	m^{\alpha-1} =\sum_{m=1}^{n} T_{\lambda}^{\searrow}(m) \, 
	m^{\alpha-1} .\]
	The first inequality holds because we have matched up the pairs of values $c_{i}$ and $m^{\alpha-1}$ in such a way as to maximise the summations value.
	The last equality comes from the fact that every diagonal in $\lambda$ has the same index $(j-i+1)$ and this increases from left to right matching the ordering of	$T_{\lambda}^{\searrow}$.
	
	To complete our proof we need to show that moving from $T_{\lambda}^{\searrow}$ to 	$T_{(n-k,\star)}^{\searrow}$ increases the eigenvalue. To get from $\lambda$ to $(n-k,\star)$ we move boxes up and to the right. Every box that moves increases its index $(j-i+1)$ and thus its value in our eigenvalue summation. Hence when we put the values of the boxes of 
	$(n-k,\star)$  in via $\searrow$, we find $T_{\lambda}^{\searrow}(m) \leq 
	T_{(n-k,\star)}^{\searrow}(m)$ for all $m$, 
	establishing
	\[\sum_{m=1}^{n} T_{\lambda}^{\searrow}(m) \cdot 
	m^{\alpha-1} \leq \sum_{m=1}^{n} T_{(n-k,\star)}^{\searrow}(m) \cdot 
	m^{\alpha-1} = 	\NA \cdot \eig\left(T_{(n-k,\star)}^{\searrow}\right). 
	\qedhere\] 
\end{proof}

\begin{example}
	As a demonstration of Lemma \ref{lem:eigbound} consider the partition $\lambda = (4,2)$. The tableaux $T_{(4,2)}^{\downarrow}$, and $T_{(4,3)}^{\searrow}$ are given below:
	\[T_{(4,2)}^{\downarrow} = 
	\ytableausetup{mathmode,baseline,aligntableaux=center,boxsize=1.2em} \begin{ytableau} 1 & 3 & 5 & 6 \\ 2 & 4 
	\end{ytableau}  \hspace{2cm}  T_{(4,2)}^{\searrow} = \begin{ytableau} 
	2& 4& 5& 6\\
	1& 3  	\end{ytableau}.\]
	For any $\alpha$ we may compute the eigenvalues of both tableaux using equation \eqref{chpt5:eqn:eigenvaluesumalpha}:
	\begin{eqnarray*}
		N_{\alpha}(8) \,\eig\left(T_{(5,3)}^{\downarrow}\right) & = & 1 (1^{\alpha-1}) + 0 (2^{\alpha-1}) +2 (3^{\alpha-1}) + 1 (4^{\alpha-1})+3(5^{\alpha-1}) + 4(6^{\alpha-1}) \\  
		N_{\alpha}(8) \, \eig\left(T_{(5,3)}^{\searrow}\right) & = & 0 (1^{\alpha-1}) + 1(2^{\alpha-1}) + 1(3^{\alpha-1}) + 2(4^{\alpha-1}) + 3(5^{\alpha-1}) + 4(6^{\alpha-1}) .
	\end{eqnarray*}
	For $\alpha \geq 1$ the values $i^{\alpha-1}$ are increasing in $i$. Therefore rearranging the coefficients given by the diagonal value of each box  (in our case the values $\{0,1,1,2,3,4\}$) to be in increasing order only increases our eigenvalue of the tableau.
\end{example}

\begin{lemma}
	\label{lem:boxindex}
	Let $n \in \mathbb{N}$, $k \in \{1,\dots, n-2\}$. Suppose $\lambda = (n-k,\star)$, then for every box $(i,j) \in \lambda$ we have
	\begin{align}
	j-i+1 \leq \frac{n-k}{n} T_{\lambda}^{\searrow}(i,j) \label{eqn:nkbound}
	\end{align}
\end{lemma}

\begin{proof}	
	Let us write $l(n-k,\star) = l^{+} = \lceil \frac{n}{n-k} \rceil$ and 
	$l^{-} = \lfloor \frac{n}{n-k} \rfloor$.  	
	If $(i_{1},j_{1}), (i_{2},j_{2})$ belong to 
	the same diagonal of $\lambda$ then $j_{1} - i_{1} +1 = j_{2} - i_{2} +1$, and hence if \eqref{eqn:nkbound} holds for the smallest value $T_{\lambda}^{\searrow}(i,j)$ on a diagonal it holds for every entry of that diagonal.  Furthermore, the bound trivially holds for any box $(i,j)$ which satisfies $j-i+1 \leq0$ (for which the left hand side of \eqref{eqn:nkbound} is non-positive). 	
	Combining these two observations, we see that it suffices to prove the bound for boxes $(1,j)$, which appear on the first row of $(n-k,\star)$ as they must contain the smallest value on their given diagonal.
	
	Note that no diagonal can contain more than $l^{+}$ boxes: call diagonals with $l^{-}$ or fewer boxes \emph{short} diagonals, and all others \emph{long} diagonals. Note that long diagonals can only exist when $l^+=l^-+1$. Any long diagonals clearly occur strictly before the short ones, when working left to right along the first row. If the box $(1,j)$ lies on a long diagonal, then the numbering pattern for $T^{\searrow}$ implies that this box will contain the integer $T_{(n-k,\star)}^{\searrow}(i,j)=\binom{l^+}{2} + 1 + (j-1)l^+$. For this value of $m$, the left hand side of \eqref{eqn:nkbound} becomes
	\begin{align*}
	\frac{n-k}{n} \left( \frac{l^+(l^+-1)}{2} + 1 + (j-1)l^+ \right)\, & = \,
	\frac{(n-k)l^+}{n} \left( \frac{l^+-1}{2} + \frac{1}{l^{+}}-1 + j \right)\, \ge \, j 
	\end{align*}
	thanks to the definition of $l^+$ and the fact that $(x-1)/2 + 1/x  \geq 1$ if  $x\geq2$. 
	
	It remains to deal with the short diagonals which contain a box on the first row. For these diagonals we now work from right $(j=n-k)$ to left  $(j=1)$. For rightmost box on the first row $(1,n-k)$ we know that $T_{(n-k,\star)}^{\searrow}((1,n-k)) = n$, and it is clear that \eqref{eqn:nkbound} holds for box $(1,n-k)$. Suppose \eqref{eqn:nkbound} holds for a box $(1,j)$, then for box $(1,j-1)$ on a small diagonal it is straightforward to see that
	\begin{eqnarray}
	j-1 \leq  \frac{n-k}{n}T_{\lambda}^{\searrow}(1,j) -1  \leq \frac{n-k}{n} T_{\lambda}^{\searrow}(1,j-1)
	\end{eqnarray}
	because $T_{\lambda}^{\searrow}(1,j) - n/(n-k)  \leq T_{\lambda}^{\searrow}(1,j-1)$.	The result for all short diagonals now follows quickly by induction.	
\end{proof}

\begin{example}
	\label{chpt5:ex:bigalphadiag}
	To illustrate Lemma \ref{lem:boxindex} consider the partition $\lambda = (10,10,9)$, which may be represented as $\lambda = (n-k ,\star)$ for $n=29, k =19$. For $\lambda = (n-k,\star)$ we have $(n-k/n) = 10/29$. Define 
	$T$ a $\lambda$-tableau, such that $T(i,j) = \lfloor (10/29) T_{\lambda}^{\searrow}(i,j)\rfloor$ for all $(i,j) \in \lambda$. The two tableaux $T_{\lambda}^{\searrow}, T$ are given below,
	
	\[	\ytableausetup{mathmode,baseline,aligntableaux=center,boxsize=1.2em}  
	T_{\lambda}^{\searrow} =
	\begin{ytableau} 
	4 & 7 & 10 & 13 & 16 & 19 & 22 & 25 & 27 & 29 \\
	2 & 5 & 8 & 11 & 14 & 17 & 20 & 23 & 26 & 28 \\
	1 & 3 & 6 & 9 & 12 &15 & 18 & 21 & 24 
	\end{ytableau}  
	\hspace{2cm}  
	T =
	\begin{ytableau} 
	1 & 2 & 3 & 4 & 5 & 6 & 7 & 8 & 9 & 10 \\
	0 & 1 & 2 & 3 & 4 & 5 & 6 & 7 & 8 & 9 \\
	0 & 1 & 2 & 3 & 4 & 5 & 6 & 7 & 8 
	\end{ytableau}  \,.\]
	We may clearly see that for any box $(i,j) \in \lambda$ we have $j-i+1 \leq T (i,j)$.

\end{example}

\begin{lemma}
	\label{lem:genbound2}
	Let $\lambda \vdash n$ with $\lambda_{1} = n-k$. Then the eigenvalue 
	$\eig(\worst)$ for the shuffle $\BLR_{n,\alpha}$ with $\alpha \geq1$ may be 
	bounded as follows:
	\[
	\eig(\worst) \,\le \, 
	\begin{cases}
	1 -	\frac{k(n-k) N_{\alpha-1}(n)}{n \NA} & \quad \text{if $k \le n/4$} \\
	1 - \frac{k}{n} & \quad \text{for all $k$}  \\
	\end{cases}
	.\]
\end{lemma}

\begin{proof}
	We begin by quickly proving the second bound for all $k$ using Lemmas \ref{lem:eigbound} and \ref{lem:boxindex}:
	\begin{eqnarray*}
		\eig(T_{\lambda}^{\downarrow}) \leq  
		\eig\left(T_{(n-k,\star)}^{\searrow}\right) & = & \frac{1}{N_\alpha(n)} \sum_{(i,j) \in T_{(n-k,\star)}^{\searrow}} (j-i +1)  \, T_{(n-k,\star)}^{\searrow}(i,j)^{\alpha-1} \\
		& \leq &	\frac{1}{N_\alpha(n)} \frac{n-k}{n} \sum_{(i,j) \in T_{(n-k,\star)}^{\searrow}} T_{(n-k,\star)}^{\searrow}(i,j)^{\alpha} =1 -\frac{k}{n}\,.
	\end{eqnarray*}
	Now we prove the bound for $k\leq n/4$, by finding a tighter bound on 
	$T_{(n-k,\star)}^{\searrow}$ and once again using Lemma \ref{lem:eigbound}. For $k\le n/4$, we know 
	$(n-k,\star) = (n-k,k)$, and thanks to the order in which the boxes are filled, we see that $T_{(n-k,\star)}^{\searrow}(m) \leq m/2$ for $1\leq m\leq 2k$, and $T_{(n-k,\star)}^{\searrow}(m) = m-k$ for $2k+1\leq m\leq n$. This gives us the following simple bound:
	\begin{eqnarray}
	\NA \cdot \eig(T_{(n-k,\star)}^{\searrow}) & \leq & \sum_{m=1}^{2k} 
	\frac{m}{2} m^{\alpha-1} + \sum_{m=2k+1}^{n} (m-k) m^{\alpha-1} \nonumber \\ 
	& = & \NA - k N_{\alpha-1}(n)+ \left( k N_{\alpha-1}(2k)-  \frac12 N_\alpha(2k)\right) \,. \label{eqn:pre-int}
	\end{eqnarray}
	Now, 
	\begin{eqnarray}
	k N_{\alpha-1}(2k)-  \frac12 N_\alpha(2k) & = & k \sum_{m=1}^{2k} m 
	^{\alpha-1}\left(1 - \frac{m}{2k}\right) 
	\leq k \int_{0}^{2k} x^{\alpha-1}\left(1 - 
	\frac{x}{2k}\right)\mathrm{d} x  \nonumber\\
	&=& \frac{k(2k)^\alpha}{\alpha(1+\alpha)} \nonumber \\
	&\leq & \frac{k^2 n^{\alpha-1}}{\alpha} \frac{2^{2-\alpha}}{(1+\alpha)} \qquad\text{(since $k\leq n/4$)} \nonumber \\
	&\leq & \frac{k^2 N_{\alpha-1}(n)}{n}  \label{chpt5:eqn:bigalphafinallemma}
	\end{eqnarray}
	in the last step we have used the inequality $n^{\alpha}/\alpha = \int^{n}_{0} j^{\alpha-1} \leq \sum^{n}_{0} j^{\alpha-1}  =  N_{\alpha-1}(n) $ which holds for $\alpha \geq 1$. Combining \eqref{chpt5:eqn:bigalphafinallemma} and \eqref{eqn:pre-int} yields the desired result.	
\end{proof}

Note that the  bound for $k\leq n/4$ in Lemma \ref{lem:genbound2}  is not the same as our previous bound (Lemma \ref{lem:genbound1}) when $\alpha=1$. The accuracy 
lost here helps us establish one bound for all 
$\alpha \geq 1$. Using Lemma \ref{lem:genbound2} above we complete our upper 
bound argument for $P_{n,\alpha}$ with $\alpha\geq1$. The total variation 
distance is bounded as follows,
\begin{eqnarray*}
	4\lVert \BLR_{n,\alpha}^{t} - \pi_{n} \rVert^{2}_{\textnormal{TV}} 
	& \leq & 
	\left(\eig(1^{n})\right)^{2t} + 2\sum_{k=1}^{(n/4)} 
	{n\choose k}^{2} k!
	\left(\eig(T_{(n-k,k)}^{\searrow})\right)^{2t}\\
	& + &
	2\sum_{k>(n/4)}^{n-2} {n\choose k}^{2} k!
	\left(\eig(T_{(n-k,\star)}^{\searrow})\right)^{2t}
\end{eqnarray*}
Substituting $t=t_{n,\alpha}(\log n +  c)$, the first term $\eig(T_{1^{n}})$ disappears as $n\to\infty$,
\[\lim_{n\to\infty} \eig(T_{1^{n}})^{2t} = \lim_{n\to\infty} \left(\frac{1}{\NA} \sum_{i=1}^{n} (2-i)^{\alpha-1} \right)^{2t} \leq \lim_{n\to\infty} \left( \frac{n^{\alpha}}{\NA}\right)^{2t_{n,\alpha}(\log n + c)} = 0.\]
We are left with the following two sums to control:
\begin{eqnarray*}
	2\sum_{k=1}^{(n/4)} 
	{n\choose k}^{2} k!
	\left(1 - 
	\frac{(n-k)kN_{\alpha-1}(n)}{ n\NA}\right)^{2t} + 
	2\sum_{k>(n/4)}^{n-2} {n\choose k}^{2} k!
	\left( 1 - 
	\frac{k}{n}\right)^{2t}.
\end{eqnarray*}
The first sum at time $t = t_{n,\alpha}(\log n +c)$ may be reduced to:
\[\sum_{k=1}^{(n/4)} 
{n\choose k}^{2} k!
\left(1 - 
\frac{(n-k)kN_{\alpha-1}(n)}{ n\NA}\right)^{2t_{n,\alpha}(\log n + c)} \leq e^{-2c}\sum_{k=1}^{n/4}   \frac{n^{2k^{2}/n}}{k!}.\]
Following the work of Diaconis presented in Section \ref{chpt4:subsec:rt} this may be	shown to be bounded by $Ae^{-2c}$ for sufficiently large $n$. 
For the second sum note 
that $\NA/N_{\alpha-1}(n)$ is increasing in $\alpha$, and at $\alpha=1$ we have $\NA/N_{\alpha-1}(n) = (n+1)/2$. The second summation has been shown to be tending to $0$ at time $(n/2)\log n + cn$ as $n\to \infty$ by Diaconis \cite[Chapter 3D Theorem 5]{Diaconis1988}.

Therefore, the sum tending to $0$ at the larger time $t = (\NA/N_{\alpha-1}(n))(\log n  +  c)$ as $n\to\infty$. Putting these bounds together we get find,
\[\limsup_{n\to \infty} \lVert \BLR_{n,\alpha}^{t_{n,\alpha}(\log n  + c)} - \pi_{n} \rVert_{\tiny \textnormal{TV}} \leq Ae^{-c} \textnormal{ for some universal constant $A$} .\]
This completes the upper bound for $\alpha \geq 1$ present in Theorem \ref{chpt5:thm:biasedbounds}.

\subsubsection{Lower Bound}
For the case of $\alpha \geq 1$ a straightforward generalisation of Section \ref{chpt5:sec:lowerbound} does not work. Instead, we are going 
to present a symmetric argument, swapping the roles of the left and right hands. 
For $\alpha \geq 1$, let $V_{n} = \{1,\dots,n/m\}$ (i.e. the bottom $n/m$ cards) and
\[F_n = \{\sigma\in S_n \, | \, \text{$\sigma$ has at least 1 fixed point in 
	$V_n$}\}.\]
We have already seen that $\pi_{n}(F_n) \leq 1/m$. We again want to bound the probability of $\BLR_{n,\alpha}(F_{n})$ via a coupon collector's argument. Writing $U_n^t$ for the set of 
uncollected cards in $V_n$ after $t$  steps of the biased one-sided transposition shuffle, it 
follows that 
\begin{eqnarray}
\label{biasedIneq1}
\RL_{n,\alpha}^{t}(F_n) \,\geq\, \mathbb{P}(|U^{t}_{n}| \ge 1 )\,.
\end{eqnarray}
Given $\alpha \geq 1$ our 
right hand is now choosing the top of the deck with higher probability meaning 
it is not likely to touch any card in $U_{n}^{t}$. Whereas, our left hand is comparatively more likely to touch cards in $U_{n}^{t}$. In one step we still only have four 
choices: both hands touch uncollected cards, only one does (left or right), or neither does. This permits us to bound the change 
in the number of \emph{collected} cards as follows:
\begin{align}
|V_n\setminus U_n^{t+1}| & \,=\, |V_n\setminus U_n^{t}| + 
|\{L^{t+1},R^{t+1}\}\cap U_n^t| \nonumber \\
& \,\le\, |V_n\setminus U_n^{t}| + 2 \cdot \ind [R^{t+1}\in U_n^t] +  	\ind [R^{t+1}\notin U_n^t,L^{t+1}\in U_n^t]\,,\label{chpt5:eqn:bigalphaindicator} 	\end{align}
This is the `reverse' of equation \eqref{eqn:indicator_fns}, this time we increase by two if our right hand collects a card. For sufficiently large $n$ we may bound the 
probability of the events in \eqref{chpt5:eqn:bigalphaindicator} as follows:
\begin{align}
\mathbb{P}(R^{t+1}\in U_n^t) & \leq \frac{1}{\NA} \sum_{i=(n/m)- 
	|U_n^t|}^{n/m} i^{\alpha} \leq \frac{|U_n^t|}{\NA} 
\left(n/m\right)^{\alpha}  \leq \frac{|U_n^t|}{\NA} 
\frac{n^{\alpha}}{\alpha (m-1)} \nonumber \\
& \leq \frac{|U_n^t| N_{\alpha-1}(n)}{\NA (m-1)} 
\label{eqn:bound_minusbiga1} \\[0.5cm]
\mathbb{P}(R^{t+1}\notin U_n^t,L^{t+1}\in U_n^t) & \le\, \frac{|U_n^t|}{\NA} 
\sum_{i=|U_n^t|+1}^{n} \frac{i^{\alpha}}{i} \leq 
\frac{|U_n^t| N_{\alpha-1}(n)}{\NA}.
\label{eqn:bound_minusbiga2}
\end{align}
In the last inequality of \eqref{eqn:bound_minusbiga1} we have used the fact that $n^{\alpha}/\alpha \leq N_{\alpha-1}(n)$ for $\alpha\geq1$.
Using \eqref{chpt5:eqn:bigalphaindicator}, \eqref{eqn:bound_minusbiga2} and 
\eqref{eqn:bound_minusbiga1} together, we now define a counting process $M_{n,\alpha}^t$ 
which stochastically dominates the number of collected cards $|V_n\setminus 
U_n^{t}|$ at all times: 
\begin{align}
M_{n,\alpha}^0 \,&=\, 0 \,; \nonumber \\
\mathbb{P}({M}^{t+1}_{n, \alpha} = {M}^{t}_{n,\alpha}+k) \,&=\, 
\begin{cases} 
\frac{N_{\alpha-1}(n) 
}{\NA(m-1)}\left(\frac{n}{m}-{M}^{t}_{n,\alpha}\right) 
& 
\text{ 
	if } k=2\\
\frac{N_{\alpha-1}(n) }{\NA}\left(\frac{n}{m}-{M}^{t}_{n,\alpha}\right) 
& 
\text{ if } 
k 
=1\\
1-\frac{m}{m-1}\frac{N_{\alpha-1}(n) 
}{\NA}\left(\frac{n}{m}-{M}^{t}_{n,\alpha}\right)  & 
\text{ if } k =0\,. \label{eqn:biasedbigM_change}
\end{cases}
\end{align}
Setting $m= \log n$, this forms a valid 
probability distribution for $n$ sufficiently large.
The counting process $M_{n,\alpha}$ is analogous to the process 
$M_{n}$ defined in \eqref{eqn:M_change} with $1/n$ replaced with $N_{\alpha-1}(n) /\NA$. We recover the bound,
\[\BLR_{n,\alpha}^{t}(F_{n}) \geq \mathbb{P}(M_{n,\alpha}^{t} < n/m).\]

\begin{lemma}
	\label{chpt5:lem:bigalphalowerboundtime}
	Let $\mathcal{T} = \min\{t\,:\, M_n^t \ge n/\log n\}$. Then for any 
	$c>2$,
	\[\mathbb{P}(\mathcal{T}\leq t_{n,\alpha}(\log n - n \log 
	\log n 	-cn)\, \leq \, \frac{\pi^{2}}{6(c-2)^2} \,. \] 
\end{lemma}
\begin{proof}
	This follows from the proof of Lemma \ref{lem:lower_bound} accounting for the change of $1/n$ to $N_{\alpha-1}(n)/\NA$.
\end{proof}
Following from Lemma \ref{chpt5:lem:bigalphalowerboundtime} and taking $c>2$ we find,
\begin{eqnarray*}
	\lVert \BLR_{n,\alpha}^{t_{n,\alpha}(\log n - \log \log n -c)} -\pi_{n} \rVert_{\tiny \textnormal{TV}} &\geq &\mathbb{P}(\mathcal{T} > t_{n,\alpha}(\log n - \log \log n -c)) - \pi_{n} (F_{n})\\
	& \geq &1 - \frac{\pi^{2}}{6(c-2)^{2}} - \frac{1}{\log n}.
\end{eqnarray*}
Taking limits in $n$ recovers the limit present in Theorem \ref{chpt5:thm:biasedbounds} for biased shuffles with $\alpha \geq 1$.

\subsection{Cutoff for Biased One-sided Transposition Shuffles with General Weight Functions}

We have established a total variation cutoff for the biased one-sided transposition shuffles where the weight function takes the form $w(j) = j^{\alpha}$ for some $\alpha$.  In this section we look to extend our previous results to a more general class of weight functions.
Notice that $j^{\alpha} / j$ is monotonically increasing  if and only if $\alpha \geq 1$ and monotonically decreasing if and only if $\alpha \leq 1$.
We now consider the class of weight functions $w: \mathbb{N} \to (0,\infty)$ with the property $w(j)/j$ is monotonic. We first consider the case where $w(j)/j$ is monotonically decreasing (corresponding to $\alpha \leq 1$). In this situation we can prove the following extension of Lemma \ref{lem:genbound1}.

\begin{lemma}
	\label{chpt5:lem:monodecreasingbound}
	Let $\lambda \vdash n$ with $\lambda_{1} = n-k$. Then the eigenvalue $\eig(\best)$ for the shuffle $\BLR_{n,w}$ with 	$w(j)/j$ monotonically decreasing may be bounded as follows:
	\[
	\eig(\best) \,\le\, 
	\begin{cases}
	1 - \frac{(n-k+1)k w(n)}{n\NW} & \quad \text{if $k \le n/4$} \\
	1 - \frac{k w(n)}{2\NW} & \quad \text{if $k >n/4$.} \\
	\end{cases}
	\]
\end{lemma}

\begin{proof}
	This follows the same proof as Lemma \ref{lem:genbound1} with $m^{\alpha}/m$ (or $m^{\alpha-1}$) replaced with $w(m)/m$.

\end{proof}

\begin{corollary}
	The mixing time of the biased one-sided transposition shuffle $\BLR_{n,w}$ with $w(j)/j$ monotonically decreasing is at most $\left(\NW /w(n) \right) \log n$. 
\end{corollary}
\begin{proof}
	Applying the same argument as Section \ref{chpt5:subsec:smallalpha} Upper Bound with the new bounds from Lemma \ref{chpt5:lem:monodecreasingbound} immediately yields the desired result.
\end{proof}

To complete the existence of a total variation cutoff we need to find a matching lower bound of the same time. However, the argument presented in Section \ref{chpt5:subsec:smallalpha} is not easily generalised for our new choices of weight function. The main issue being the relationship between $w(k)$ and $w(n)$ for $k \in V_{n} = \{n - n/m,\ldots, n\}$. In the previous case with $w(j) = j^{\alpha}$ we can easily see that $w(k) = w(k/n) w(n)$. In the current situation we only know $w(n)/n \leq w(k) /k$; this is not enough information to establish a relationship between $w(k)$ and $w(n)$ for use in equations \eqref{eqn:bound_minus2biased}, and \eqref{eqn:bound_minus1biased}. A new approach is likely to be needed to prove a lower bound on the mixing time of the biased one-sided transposition shuffle for this class of weight functions.

\begin{conj}
	The biased one-sided transposition shuffle $\BLR_{n,w}$ with $w(j)/j$ monotonically decreasing exhibits a total variation cutoff at time $\left(\NW /w(n) \right) \log n$. 
\end{conj}

\paragraph{}
We now consider the class of weight functions where $w(j)/j$ is monotonically increasing (corresponding to $\alpha \geq 1$). Previously we proved the mixing time for the biased one-sided transposition shuffle with $\alpha \geq 1$ to be $\left( \NA / N_{\alpha-1}(n) \right) \log n$. For a general function $w$ we replace the factor of $N_{\alpha-1}(n)$ by the summation $ N_{w}^{\prime}(n):= \sum_{i=1}^{n} w(i)/i$. Therefore, we now look to prove a cutoff at time $\left( \NW / N_{w}^{\prime}(n) \right) \log n$. The lower bound argument presented in Section \ref{chpt5:subsec:bigalpha} Lower Bound, can easily be modified for this new class of weight functions. 

\begin{lemma}
	The mixing time of the biased one-sided transposition shuffle $\BLR_{n,w}$ with $w(j)/j$ monotonically decreasing is at least $\left(\NW / N_{w}^{\prime}(n) \right) \log n$. 
\end{lemma}
\begin{proof}
	The result follows from the work of Section \ref{chpt5:subsec:bigalpha} Lower Bound with $N_{\alpha-1}(n) /\NA$ replaced by $N_{w}^{\prime}(n) / N_{w}(n)$, and the inequalities \eqref{eqn:bound_minusbiga1} and \eqref{eqn:bound_minusbiga2} replaced with:
	\begin{align}
	\mathbb{P}(R^{t+1}\in U_n^t) & \leq \frac{1}{\NW} \sum_{i=(n/m)- 
		|U_n^t|}^{n/m} w(i)\leq  \frac{|U_{n}^{t}|}{(m-1)\NW} (m-1)w(n/m)  \nonumber \\ 
	& \leq  \frac{|U_n^t|}{(m-1)\NW} \sum_{i=1}^{n} \frac{w(i)}{i} = \frac{|U_n^t| N_{w}^{\prime}(n)}{(m-1)\NW}
	\\[0.5cm]
	\mathbb{P}(R^{t+1}\notin U_n^t,L^{t+1}\in U_n^t) & \le\, \frac{|U_n^t|}{\NW} 
	\sum_{i=|U_n^t|+1}^{n} \frac{w(i)}{i} \leq \frac{|U_n^t| N_{w}^{\prime}(n)}{\NW}.
	\end{align}
\end{proof}
To find a matching upper bound on the mixing time we look to generalise the results of Section \ref{chpt5:subsec:bigalpha} Upper Bound.
For a general weight function with $w(j)/j$ monotonically increasing  Lemma \ref{lem:eigbound} still holds and so we can bound the eigenvalues associated to a partition $\lambda$ using the value $\eig\left(T_{\lambda}^{\searrow} \right)$. For small partitions we are able to use the bound constructed in Lemma \ref{lem:genbound2} without any modifications. For large partitions we encounter an issue because the proof of the bound in Lemma \ref{lem:genbound2} relies on exact knowledge of the weight function $w$. In particular, we use the integral of $w$ in order to obtain an upper bound on the value of $\eig\left(T_{\lambda}^{\searrow} \right)$,  this argument can not be generalised to a generic weight function $w$. However, we conjecture that an upper bound of the correct mixing time is possible.

\begin{conj}
	The biased one-sided transposition shuffle $\BLR_{n,w}$ with $w(j)/j$  monotonically decreasing  exhibits a total variation cutoff at time $\left(\NW / N_{w}^{\prime}(n) \right) n\log n$. 
\end{conj}

The class of functions $\{w:\mathbb{N} \to (0,\infty) \, | \,  w(j)/j \textnormal{ is monotonic}\}$ is still quite restrictive. We suspect that cutoffs in total variation distance may be proved for biased one-sided transposition shuffles with more general weight function $w$ (for example, slowly and regularly varying functions), but any proof of this would certainly need different techniques to those developed in this section.

\section{Separation Distance for the One-sided Transposition Shuffle}
\label{chpt5:sec:separation}
\subsection{Cutoff for the Unbiased One-sided Transposition Shuffle}

The unbiased one-sided transposition 
shuffle exhibits a cutoff in total variation distance at time $n\log n$. To prove that we also have a cutoff in separation distance at time $n\log n$ we identify a strong stationary time for the unbiased one-sided transposition shuffle. The aim of this section is to prove the following result.
\begin{thm}
	\label{thm:SSTtime}
	There exists a strong stationary time $T$ for the one-sided 
	transposition shuffle, with		
	$\mathbb{P}(T> n\log n + cn) \leq e^{-c}$.
\end{thm}

From Theorem \ref{thm:SSTtime} we may quickly establish that the one-sided transposition shuffle exhibits a cutoff in separation distance at time $n\log n$.

\begin{thm}
	\label{thm:sepcutoff}
	The one-sided transposition shuffle exhibits a cutoff in separation distance at time $n\log n$.

\end{thm}

\begin{proof}[Proof of Theorem \ref{thm:sepcutoff}]
	The lower limit on separation distance mixing time follows from Lemma \ref{chpt2:lem:sepupperbound} and Theorem \ref{chpt5:thm:cutoff}. 
	The upper limit on separation distance mixing time follows from  Theorem 
	\ref{thm:SSTtime}.

\end{proof}

\subsubsection{A Strong Stationary Time Argument}
\label{chpt5:subsec:ostsst}

Recall that we may view the elements of $S_{n}$ as the possible permutations of deck of cards which is made up of cards and positions both indexed by $[n]$. Any permutation $\sigma \in S_{n}$ is a bijection from cards to positions and $\sigma^{-1}$ is a bijection from positions to cards, i.e. $\sigma(i)$ tells us the position of card $i$ whereas $\sigma^{-1}(i)$ tells us what card is in position $i$.  We begin all our random walks at the identity permutation with positions and labels fully matched. Throughout the rest of this section let $(X^{t})_{t\in 	\mathbb{N}}$ denote a Markov chain driven by the unbiased one-sided transposition shuffle, and let $(Y^{t})_{t \in \mathbb{N}}$ be a Markov chain defined by setting $Y^{t} = (X^{t})^{-1}$ for all times $t$. For any $\sigma \in S_{n}$ and time $t$ we have $\mathbb{P}(X^{t} = \sigma) = \mathbb{P}(Y^{t} = \sigma^{-1})$.  The Markov chains $(X^{t})$ and $(Y^{t})$ represent two different ways to view the one-sided transposition shuffle.

Let $\tau^{t}$ be the transposition chosen at step $t$ of the unbiased one-sided transposition shuffle. To construct our strong stationary time we need to condition on the exact permutation of cards in positions above position $j$ at time $t$, that is the random variables $Y^{t}(i)$ for $j< i  \leq n$.
Given this information we also know which cards can be in positions $1$ to $j$ at time $t$, define this set as, 
\[A_{j}^{t} = [n] \setminus \{ Y^{t}(i)  \, | \, j < i \leq n\}.\]

\begin{defn}
	We say the random walk $(Y^{t})_{t\in \mathbb{N}}$ satisfies \emph{property $\mathcal{P}_{j}$ at time $t$} if we have:
	\begin{eqnarray}
	\mathbb{P}\left(Y^{t} (j) = l \, | \, Y^{t}(i) \textnormal{ for all } j < i \leq n \right) = 
	\begin{cases}
	1/j & \textnormal{ if } l \in A_{j}^{t} \\
	0 & \textnormal{ otherwise } 
	\end{cases} \label{chpt5:eqn:propertyj}
	\end{eqnarray} 
	This property tells us that given total information about the deck strictly above position $j$, the card in position $j$ is equally likely to be any of the remaining cards.
\end{defn}

\begin{lemma}
	\label{lem:SUTpropj}
	Let $T_{j}$ be the first time our right hand chooses the position
	$j$ when performing the unbiased one-sided transposition shuffle. 
	If $T_{j} \leq t$ then the Markov chain $(Y^{t})_{t\in \mathbb{N}}$ satisfies property $\mathcal{P}_{j}$ at time $t$.
\end{lemma}
\begin{proof}
	We prove this by induction: once property $\mathcal{P}_{j}$ holds for some $t$, it holds for all times after $t$.

	Consider the time $T_{j}$, at this step of our Markov chain we must have applied  a transposition 
	$(i\,j)$ with $i\leq j$. The probability of picking any one of the transpositions $(i \,j)$ at time $T_{j}$ is $\mathbb{P}(\tau^{T_{j}} = (i\,j)) = 1/j$ for all $i\leq j$. Therefore, the card in position $j$ at time $T_{j}$ has a uniform chance of being any of the cards in $A_{j}^{T_{j-1}} = A_{j}^{T_{j}}$. Thus, we may clearly see that,
	\begin{eqnarray}
	\mathbb{P}\left(Y^{T_{j}} (j) = l \, | \, Y^{T_{j}}(i) \textnormal{ for all } j < i \leq n \right)= 
	\begin{cases}
	1/j & \textnormal{ if } l \in A_{j}^{T_{j}}\\
	0 & \textnormal{ otherwise } 
	\end{cases}
	\end{eqnarray} 
	so $\mathcal{P}_{j}$ holds at time $T_{j}$.

	Now suppose property $\mathcal{P}_{j}$ holds at time $t$. We study the time $t+1$ and split the analysis 
	into cases based on which transposition $(a \,b)$ (with $a\leq b$) was applied at time 
	$t+1$,
	\begin{eqnarray}
	& & \mathbb{P}\left(Y^{t+1}(j) =  l \, | \, Y^{t+1}(i) \textnormal{ for all } j < i \leq n\right) \nonumber\\
	& = &
	\sum_{\substack{1\leq a \leq b \leq n}} \mathbb{P}(\tau^{t+1} = 
	(a\,b)) \, \mathbb{P}\left(Y^{t+1} (j) 
	= l \, | \, Y^{t+1}(i) \textnormal{ for all } j < i \leq n, \, \tau^{t+1} = (a\,b) 
	\right).\hspace{0.5cm}\label{eqn:SSTinduction} 
	\end{eqnarray}
	Using knowledge of the transposition $\tau^{t+1}$ we evolve our deck backwards in time ($X^{t} = \tau^{t+1} X^{t+1}$) to recover the random variables $Y^{t}(i)$ from $Y^{t+1}(i)$ for $j < i \leq n$, and relate $Y^{t+1}(j)$ to $Y^{t}(j)$. This allows us to use our inductive hypothesis.

	If $b=j$ then our random walk satisfies property $\mathcal{P}_{j}$ at time $t+1$ for the same reasoning as time $T_{j}$. Suppose $a \leq b < j$ then  we know that $Y^{t}(i) = Y^{t+1}(i)$ for all $j<i \leq n$. Suppose instead that $a,b >j$, then we have $Y^{t}(b) = Y^{t+1}(a)$ and $Y^{t}(a) = Y^{t+1}(b)$, with $Y^{t}(i) = Y^{t+1}(i)$ for all other $j<i\leq n$. In either case we know $A_{j}^{t+1} = A_{j}^{t}$, and the card in position $j$ has not moved from time $t$ to $t+1$. Therefore, we have
	
	\begin{eqnarray}
	& & \mathbb{P}\left(Y^{t+1} (j) = l \, | \, Y^{t+1}(i) \textnormal{ for all } j < i \leq n, \, \tau^{t+1}  =(a \, b )\right) \\	
	& = & \mathbb{P}\left(Y^{t} (j) = l \, | \, Y^{t}(i) \textnormal{ for all } j < i \leq n \right)  = 
	\begin{cases}
	1/j & \textnormal{ if } l \in A_{j}^{t} = A_{j}^{t+1}\\
	0 & \textnormal{ otherwise } 
	\end{cases} \nonumber.
	\end{eqnarray}
	
	We now study in detail the effects of the remaining transpositions $(a \, b)$ with $a\leq j<b$.
	In this case  we can not fully recover the random variables $Y^{t}(i)$ with $i < j \leq n$ without extra assumptions. To this end fix $b>j$,  a card $C \in A_{j}^{t+1}$, and suppose that $\tau^{t+1} = (X^{t+1}(C) \hspace{0.15cm} b)$, i.e., card $C$ is moved from position $b$ into a position below $j$ by $\tau^{t+1}$. Letting $C$ range over all choices in $A_{j}^{t+1}$ will recover every transposition $(a \, b)$ with $b>j$ fixed and  $a\leq j$. In the case that  $\tau^{t+1} = (X^{t+1}(C) \hspace{0.15cm} b)$, we know that $Y^{t}(b) = C$ and the other positions above $j$ have $Y^{t}(i) = Y^{t+1}(i)$ for $j < i \leq n$ and $i \neq b$.  Therefore, for this choice of $\tau^{t+1}$ we know that $A_{j}^{t} = \left(A_{j}^{t+1} \sqcup \{Y^{t+1}(b)\} \right) \setminus \{C\}$. Now consider the probability:
	\begin{eqnarray}
	\mathbb{P}\left(Y^{t+1}(j) = l \, | \, Y^{t+1}(i) \textnormal{ for all } j < i \leq n , \, \tau^{t+1} = (X^{t+1}(C) \hspace{0.2cm} b)\right) \nonumber.
	\end{eqnarray}
	
	If $l = C$ then the event in question can only occur if  the card currently in position $b$, i.e. $Y^{t+1}(b)$, was in position $j$ at time $t$. Noting that $Y^{t+1}(b) \in A_{j}^{t}$, and using our inductive hypothesis we find,
	\begin{eqnarray}
	& & \mathbb{P}\left(Y^{t+1}(j) = C \, | \, Y^{t+1}(i) \textnormal{ for all } j < i \leq n , \, \tau^{t+1} = (X^{t+1}(C) \hspace{0.2cm} b)\right) \nonumber \\
	& = & \mathbb{P}\left(Y^{t}(j) = Y^{t+1}(b) \, | \, Y^{t}(i) \textnormal{ for all } j < i \leq n  \right) = 1/j \label{chpt5:eqn:SSThardind1}
	\end{eqnarray}
	Alternatively, suppose $l \in A_{j}^{t+1} \setminus \{C\}$, we know the card $l$ does not move from its position at time $t$ to time $t+1$, and we know $l \in A_{j}^{t}$. Therefore, we find
	\begin{eqnarray}
	& & \mathbb{P}\left(Y^{t+1}(j) = l \, | \, Y^{t+1}(i) \textnormal{ for all } j < i \leq n , \, \tau^{t+1} = (X^{t+1}(C) \hspace{0.2cm} b)\right) \nonumber \\
	& = & \mathbb{P}\left(Y^{t}(j) = l \, | \, Y^{t}(i) \textnormal{ for all } j < i \leq n  \right) = 1/j \label{chpt5:eqn:SSThardind2}
	\end{eqnarray}
	Putting the equations \eqref{chpt5:eqn:SSThardind1} and \eqref{chpt5:eqn:SSThardind2} together gives us,
	\begin{eqnarray*}
		\mathbb{P}\left(Y^{t+1}(j) = l \, | \,  Y^{t+1}(i) \textnormal{ for all } j < i \leq n, \, \tau^{t+1} = (X^{t+1}(C) \hspace{0.2cm} b)\right) & = & 
		\begin{cases}
			1/j & \textnormal{ if }  l \in A_{j}^{t+1} \setminus\{ C \}\\
			1/j & \textnormal{ if }  l = C \\
			0 & \textnormal{ otherwise }
		\end{cases}.
	\end{eqnarray*}
	Letting $C$ range over all possible choices of card in $A_{j}^{t+1}$ while keeping $b> j$ fixed, we cover the desired probability for every transposition $(a\, b)$ with $a\leq j < b$. 
	Finally applying every separate case to \eqref{eqn:SSTinduction} we have established that 
	\[\mathbb{P}\left(Y^{t+1}(j) = l \, | \, Y^{t+1}(i) \textnormal{ for all } j < i \leq n
	\right) = \begin{cases}
	1/j & \textnormal{ if } l \in A_{j}^{t+1}\\
	0 & \textnormal{ otherwise } 
	\end{cases} \]
	as required, thus by induction our hypothesis holds for all 
	$t \geq T_{j}$.
	
\end{proof}

\begin{lemma}
	\label{lem:SUTresutl}
	Let $T$ be the first time our right hand has chosen every position. Then $T$ is a strong uniform time for $(X^{t})_{t\in\mathbb{N}}$ and $(Y^{t})_{t\in\mathbb{N}}$.
	
\end{lemma}

\begin{proof}
	Note that $X^{t}$ is uniformly distributed if and only if $Y^{t}$ is uniformly distributed. Lemma \ref{lem:SUTpropj} implies that by time $T \geq T_{j}$ the Markov chain $(Y^{t})_{t \in \mathbb{N}}$ satisfies all properties $\mathcal{P}_{j}$. Hence, for any $\sigma \in  S_{n}$,  

	\begin{eqnarray*}
		\mathbb{P}(Y^{t} = \sigma^{-1} \, | \, T\leq t) & = &\mathbb{P}\left( 
		\cap_{j=1}^{n} \{Y^{t}(j) = \sigma^{-1}(j) \} \, | \, T\leq t\right)\\  
		& = & \prod_{j=1}^{n}\mathbb{P}\left(  Y^{t}(j) = \sigma^{-1}(j) \, | \, 
		\cap_{i=j+1}^{n} \{
		Y^{t}(i) = \sigma^{-1}(i) \},\, T\leq t\right) \\
		& = &	\prod_{j=1}^{n} 
		\frac{1}{j} = \pi_{n}(\sigma).
	\end{eqnarray*}
\end{proof}

We have found a strong stationary time for the unbiased one-sided transposition shuffle. Following quickly from this we may prove Theorem \ref{thm:SSTtime}, and thus establish a cutoff in separation distance for the one-sided transposition shuffle.

\begin{proof}[Proof of Theorem \ref{thm:SSTtime}]
	Let $T$ be the first time our right hand has chosen every position 
	$j$. 
	Our right hand is choosing positions via a uniform probability on 
	$[n]$. Thus $T$ is modelled by the standard coupon collector's problem with $n$ 	coupons. To complete our argument recall (Section \ref{chpt4:subsec:t2r}, equation \eqref{chpt4:eqn:classiccoupon})  that for the standard coupon collector's problem on $n$ cards we have $\mathbb{P}(T_{n} >n\log n +cn) \leq e^{-c}$. 
\end{proof}

\paragraph{}
The
strong stationary time $T$ is novel in a few ways. 
Often strong stationary 
times rely on a sequence of stopping times $(T_{i})_{i\in \mathbb{N}}$ where at 
time $T_{i}$ we know there exists a subset of $i$ cards which are uniformly distributed with respect to $S_{i}$. This is exactly the technique employed to find a strong stationary time for the top-to-random shuffle which we covered in Section \ref{chpt4:subsec:t2r}.
The strong stationary time $T$ can not be split up in this way. Suppose we require that we 
touch positions in order from $1$ to $n$, then after we 
choose position $i$ the cards in the subset of positions 
$[i]\subseteq [n]$ are uniformly distributed on $S_{i}$, similarly to the top-to-random shuffle. 	
However, if instead we require that positions are touched in reverse order, i.e. from $n$ to $1$, then we do not find any uniform subsets of $S_{n}$ before becoming completely uniform when we touch card $1$. Therefore, because we do not impose a condition on the order in which we must touch positions to get to the time $T$,  at any particular time 
before $T$ there is no certainty that we have a uniformly distributed subset of the cards $[n]$. 
This sudden uniformity is the kind of behaviour which drives the existence of a cutoff in separation distance.

\subsection{Generalising for Biased One-sided Transposition Shuffles}
\label{chpt5:subsec:generalsst}

In Section \ref{chpt5:sec:biased} we were able to prove cutoff in total variation for all biased one-sided transposition shuffles. Given that the unbiased one-sided transposition shuffle $\BLR_{n,0}$ also shows a cutoff 
in separation at the same time as for total variation we may hope the same can 
be proven for all biased shuffles $\BLR_{n,\alpha}$. However, in this section 
we show that simply extending the method used for $\BLR_{n,0}$ is not enough to prove the existence of a cutoff in separation distance for all $\alpha$. Recall that for a biased 
one-sided transposition shuffle $\BLR_{n,\alpha}$ our right hand chooses a position $j$ with probability $j^{\alpha} / \NA$. 

\begin{lemma}
	\label{lem:SSTtimebias}
	Let $T_{\alpha}$ be the time for the right hand of the biased shuffle 
	$\BLR_{n,\alpha}$ to choose every position in $[n]$ at least once. Then 
	$T_{\alpha}$ is a strong stationary time for the shuffle 
	$\BLR_{n,\alpha}$, for all $\alpha$.
\end{lemma}
\begin{proof}
	For the biased shuffle $\BLR_{n,\alpha}$ our left hand is still choosing uniformly below our right hand, therefore this result follows from the proofs of Lemmas \ref{lem:SUTpropj} and \ref{lem:SUTresutl}.
\end{proof}

Given $T_{\alpha}$ is still a strong stationary time all we have left to do is analyse the probability $\mathbb{P}(T_{\alpha}> t)$. This reduces to a biased coupon collector's problem. The next lemma establishes a bound on the time $T_{\alpha}$ which will allow us to prove a cutoff in separation distance for biased one-sided transposition shuffles with $\alpha\leq 0$.

\begin{lemma}
	\label{chpt5:lem:SSTnegalpha}
	Let $t_{n,\alpha}$ be the time defined in Theorem \ref{chpt5:thm:biasedbounds}, that is
	\[t_{n,\alpha}= \begin{cases}
	\NA/n^{\alpha} & \textnormal{ if } \alpha \leq 1 \\
	\NA/N_{\alpha-1}(n) & \textnormal{ if } \alpha \geq 1 
	\end{cases}. \]
	For the biased one-sided transposition shuffle with $\alpha\leq0$ we have 
	that $\mathbb{P}(T_{\alpha} > t_{n,\alpha}(\log n + c)) \leq e^{-c}$.
\end{lemma}
\begin{proof}
	We say position $j$ is \emph{collected} the first time it is chosen by our right hand. Its collection probability at each step of our random walk is $j^{\alpha}/\NA$. We bound the time $T_{\alpha}$ by a coupon collector's argument with inspiration taken from the top to random shuffle (Lemma \ref{chpt4:lem:ttrtimebound}). Let $C_{j}^{t}$ be the event that we have not collected position $j$ by time $t$. We may find a simple upper bound as follows:
	\begin{eqnarray}
	\mathbb{P}(T_{\alpha} > t)  = \mathbb{P}(\cup_{j=1}^{n} C_{j}^{t}) \leq \sum_{j=1}^{n} \mathbb{P}(C_{j}^{t}) = \sum_{j=1}^{n}\left( 1 - \frac{j^{\alpha}}{\NA} \right)^{t}. \label{chpt5:eqn:biasedsep}
	\end{eqnarray}
	Now substituting $t = t_{n,\alpha}(\log n +c) = (\NA /n^{\alpha})(\log n +c)$ for $\alpha\leq 0$, and applying the bound $(1-x)^{t} \leq e^{-tx}$, we find,
	
	\begin{eqnarray*}
		\mathbb{P}(T_{\alpha} > t)  \leq  \sum_{j=1}^{n}\left( 1 - \frac{j^{\alpha}}{\NA} \right)^{ t_{n,\alpha}(\log n +c)} \leq e^{-c}\sum_{j=1}^{n} n^{-(j/n)^{\alpha}} \leq e^{-c} n^{1 -(1)^\alpha} \leq e^{-c},
	\end{eqnarray*} 
	with the second to last inequality following from $(j/n)^{\alpha}  = (n/j)^{-\alpha} \geq 1^{-\alpha}$ for $\alpha \leq 0$.
	
\end{proof}

\begin{thm}
	\label{chpt5:thm:biasedsepcutoff}
	The biased one-sided transposition shuffle $\BLR_{n,\alpha}$ with $\alpha\leq 0$ exhibits a cutoff in separation distance at time $t_{n,\alpha} \log n$.	
\end{thm}

\begin{proof}
	The lower bound on separation distance mixing time follows from Lemma \ref{chpt2:lem:sepupperbound} and Theorem \ref{chpt5:thm:biasedbounds}. The upper bound on separation distance mixing time follows from  Theorem 
	\ref{thm:SSTtime}.
\end{proof}

The problem stopping us from copying Lemma 
\ref{chpt5:lem:SSTnegalpha} for shuffles with $\alpha>0$ is the approximation in equation \eqref{chpt5:eqn:biasedsep}.
For $\alpha \in (0,1]$ we still find that
\[ 	\mathbb{P}(T_{\alpha} > t) \leq \sum_{j=1}^{n}\left( 1 - \frac{j^{\alpha}}{\NA} \right)^{t_{n,\alpha}(\log n +c)} \leq e^{-c} \sum_{j=1}^{n} n^{-(j/n)^{\alpha}}. \]
The summation present above is unbounded as $n\to\infty$ for $\alpha \in (0,1]$. Similarly for $\alpha \in [1,\infty)$ we may find the bound \eqref{chpt5:eqn:biasedsep} to be,
\[	\mathbb{P}(T_{\alpha} > t) \leq \sum_{j=1}^{n}\left( 1 - \frac{j^{\alpha}}{\NA} \right)^{t_{n,\alpha}(\log n +c)} \leq \sum_{j=1}^{n}n^{-j^{\alpha}/N_{\alpha-1}(n)} e^{-c \,(j^{\alpha}/N_{\alpha-1}(n))}\]
and this summation is still unbounded as $n\to\infty$. Therefore, to use bound \eqref{chpt5:eqn:biasedsep} for $\alpha >0$ we would need to increase the time we consider from $t_{n,\alpha}\log n$ until the respective summation is bounded in $n$.  For example if $\alpha =1$ the summation  \eqref{chpt5:eqn:biasedsep} becomes
\[e^{-c}\sum_{j=1}^{n}\left( 1 - \frac{2j}{n(n+1)} \right)^{t} \]
which is bounded at time $t = O(n^{2})$, which is of a greater order than our lower bound of $(n+1)/2\log n$ from Theorem \ref{chpt5:thm:biasedbounds}.

Overall the bound presented in the proof of Lemma \ref{chpt5:lem:SSTnegalpha} is not good enough to establish an upper bound on our separation distance mixing time that matches our lower bound of $t_{n,\alpha} \log n$ for the biased one-sided transposition shuffles with $\alpha>0$.

\paragraph{}
The shuffle $\BLR_{n,1}$ is closely related to the random transposition shuffle.  Strong stationary times for the random transposition shuffle have historically been difficult to find. In \cite{broderthesis} Broder found a strong stationary time for random transpositions which gave an upper bound on the mixing time of $2n\log n$. Matthews claimed to have improved this to $n\log n$ and then to $(n/2)\log n$, thus establishing a cutoff in separation distance \cite{matthews1988strong}. However, the argument Matthew presented had a subtle mistake which was only recently discovered by White and explored in his thesis \cite{whitethesis}. White has constructed a strong stationary time for random transpositions which gives the correct upper bound of $(n/2)\log n$ \cite{white2019strong}. The proof of this involves the creation of a sophisticated stationary time involving keeping track of the possible cycle structure the deck could be in at certain times. In light of this it is not a 
surprise that a simple generalisation does not yield a satisfactory upper bound on separation distance for the shuffle $\BLR_{n,1}$. Despite not finding a strong stationary time for the biased one-sided transposition shuffles with $\alpha >0$, we conjecture that they show a cutoff in separation distance at the same time as their respective total variation cutoff. 

\begin{conj}
	The biased one-sided transposition shuffle $\BLR_{n,\alpha}$ with $\alpha \in (0,\infty)$ exhibits a cutoff in separation distance at $t_{n,\alpha}\,\log n$.
\end{conj}

\chapter{The Hyperoctahedral Group and Random Walks}
\label{chpt6:chpt}

The hyperoctahedral group $B_{n}$ is a natural extension of the symmetric group. We may view the hyperoctahedral group as the arrangement of a deck of $n$ cards with the extra information of knowing whether a given card is face up or face down. The group has a lot of structural similarly with the symmetric group; its elements may still be decomposed into a product of different cycles, and the conjugacy classes are determined by cycle type. This allows us to extend shuffles on $S_{n}$ into those on $B_{n}$, in particular we focus on extensions of the random transposition shuffle and one-sided transposition shuffle.
The module structure of the hyperoctahedral group also resembles that of $S_{n}$ with its simple modules being generalised Specht modules $S^{\bl}$ labelled now by bi-partitions of $n$. This enables us to modify the technique of lifting eigenvectors for the hyperoctahedral group and extended shuffles. 

In the first few sections we cover the basics of the hyperoctahedral group and its module structure, good references for the material we cover are \emph{The Representations of the Weyl Groups of Type $B_{n}$} by Aamily, Morris and Peel \cite{al1981representations}, and \emph{Representations of the hyperoctahedral Groups} by Geissinger and Kinch \cite{geissinger1978representations}. The last two sections detail lifting eigenvectors for the hyperoctahedral group, we recover the full spectrum of the random transposition shuffle and one-sided transposition shuffle.

\section{The Hyperoctahedral Group}
\label{chpt6:subsec:hyper}
Define the set $[\pm n ]:= \{1,\ldots, n\} \cup \{-1,\ldots, -n\}$
The hyperoctahedral group on $n$ elements, denoted $B_{n}$, may be defined as the group of all bijections $\sigma: [\pm n] \to [\pm n]$ such that $\sigma(-i) = - \sigma(i)$, there are $2^{n} n!$ such bijections. Any map $\sigma \in B_{n}$ is completely determined by its values on the positive elements $\{1,\ldots, n\}$. Another common way to think about the group $B_{n}$ is the arrangements of a deck of $n$ cards where we may now distinguish between cards that are face up or face down. The bijection $\sigma$ tell us what position card $i$ is in and whether it is face up ($\sigma(i)$ positive), or face down ($\sigma(i)$ negative). The last way to view the hyperoctahedral group is as the wreath product $\mathbb{Z}_{2} \wr S_{n}$. We may identify important subgroups $S_{n} = \{\sigma \, | \, \sigma(i) >0 \textnormal{ if $i > 0$}\}$ and $\mathbb{Z}_{2}^{n} = \{\sigma \, | \, \sigma(i) = \pm i \textnormal{ for all $i$}\}$ inside the hyperoctahedral group.

Define multiplication performed in the group $B_{n}$ to be the composition of functions from right to left. We may decompose elements of $B_{n}$ into a product of negative transpositions, and a permutation from $S_{n}$. This allows us to use cycle notation from $S_{n}$  (Section \ref{chpt4:sec:symgp}) to represent the elements of $B_{n}$.
\begin{defn}
	The conjugacy class of \emph{negative transpositions} in $B_{n}$ is composed of the elements $\xi_{i}$ for $1 \leq i \leq n$, with:
	\[\xi_{i}(j) = \begin{cases}
	-i & \textnormal{ if } j= i\\
	+i & \textnormal{ if } j=-i\\
	j & \textnormal{ otherwise }
	\end{cases}
	.\]
	The elements $\xi_{i}$ can be thought of as flipping card $i$ over.
\end{defn}
\begin{lemma}
	\label{chpt6:lem:decom}
	Let $\sigma \in B_{n}$. We may decompose $\sigma$ into a permutation from $S_{n}$ and a product of negative transpositions. I.e., there exists $x_{i} \in \{0,1\}$ and $\eta \in S_{n} \subset B_{n}$ such that,
	\[\sigma= \left(\prod_{i=1}^{n} \xi_{i}^{x_{i}}\right) \eta .\]
\end{lemma}
\begin{proof}
	By applying negative transpositions to $\sigma$ we may recover an element of $S_{n}$ inside of $B_{n}$, i.e., there exists a collection of $y_{i} \in \{0,1\}$ such that, $\left(\prod_{i=1}^{n}\xi_{i}^{y_{i}}\right) \sigma = \eta$ for some $\eta \in S_{n}$. From here taking inverses of the elements $\xi_{i}$ gives the desired form.
\end{proof}

The conjugacy class of \emph{positive transpositions} for $B_{n}$ is formed from the elements $(i \, j)$ and $\xi_{i} \xi_{j}( i \, j)$ with $i<j$. The set of \emph{transpositions} for the hyperoctahedral group is the union of the sets of positive and negative transpositions. Lemma \ref{chpt6:lem:decom} may by restated as; every element of $B_{n}$ may be decomposed into a product of transpositions.
Previously we described the conjugacy classes of $S_{n}$ in terms of their cycle type. Any element of $B_{n}$ may be split into positive and negative cycles which define the cycle type of that element. To label the possible cycles types for $B_{n}$ we introduce the concept of bi-partitions. 

\begin{defn}
	Let $n\in\mathbb{N}^{0}$, a \emph{bi-partition} of $n$, denoted $\bl$, is a pair of partitions $\bl = (\lambda^{1},\lambda^{2})$ such that $\lambda^{1} \vdash m$ and $\lambda^{2} \vdash n-m$ for some $0\leq m \leq n$. If $\bl$ is a bi-partition of $n$ we write $\bl \vdash n$.
\end{defn}

\begin{defn}
	Let $\sigma \in B_{n}$ with decomposition $\sigma = \left(\prod_{i=1}^{n} \xi_{i}^{x_{i}}\right) \eta$. We call a cycle in $\eta$ a \emph{positive cycle} if $\sum_{j}x_{j} \equiv 0 \mod 2$ for elements $j$ in the cycle. We call a cycle in $\eta$ a \emph{negative cycle} if $\sum_{j}x_{j} \equiv 1\mod 2$ for elements $j$ in the cycle.
	Form two non-increasing tuples $\lambda^{1}$, $\lambda^{2}$ of the lengths of positive and negative cycles respectively. The \emph{cycle type} of $\sigma$ is defined as the bi-partition $(\lambda^{1},\lambda^{2})$.
\end{defn}

\begin{example}
	Let $n =6$ and consider $\sigma \in B_{6}$ as defined by,
	\[\sigma =	\left(\xi_{1}\xi_{2} \xi_{4} \xi_{5} \right) (1)(2 \, 3 \, 4) (5 \, 6).\]
	The element $\sigma$ only has one positive cycle of length $3$, being $\xi_{2}\xi_{4} (2 \, 3 \, 4)$. We also have two negative cycles, the element $\xi_{1}$ forms a cycle of length $1$, and the elements $\xi_{5}(5 \,6)$ giving a cycle of length $2$. Therefore, $\sigma$ has cycle type $((3),(2,1))$.
\end{example}

\begin{lemma}
	Elements of $B_{n}$ belong to the same conjugacy class if and only if they share the same cycle type. Hence, bi-partitions label the conjugacy classes of $B_{n}$.
\end{lemma}

The conjugacy classes of negative and positive transpositions have cycle types $((1^{n-1}),(1))$ and  $((2,1^{n-2}),(0))$ respectively.
Following from the decomposition of Lemma \ref{chpt6:lem:decom}, we may introduce the generalised sign function for the hyperoctahedral group, which counts the number of transpositions (positive and negative) needed to form an element.
\begin{defn}
	Let $\sigma \in B_{n}$ with decomposition $ \sigma = \left(\prod_{i=1}^{n} \xi_{i}^{x_{i}} \right) \eta$. Define the  \emph{sign function} for the hyperoctahedral group, denoted $\sign:B_{n} \to \{-1,1\}$, as follows,
	\[\sign(\sigma) = (-1)^{\sum_{i=1}^{n}x_{i}} \, \sign(\eta).\]
	An element in $B_{n}$ is called \emph{even} if it has positive sign, and is called \emph{odd} if it has negative sign. The sign function on $B_{n}$ is multiplicative, that is $\sign(\sigma \tau) = \sign(\sigma) \sign(\tau)$. Thus, the set of positive elements of $B_{n}$ forms a normal subgroup of index $2$.
\end{defn}

\section{The Structure of Modules for the Hyperoctahedral Group.}
\label{chpt6:subsec:modules}
Define the group algebra of the hyperoctahedral group as $\mathfrak{B}_{n} := \mathbb{C}[B_{n}]$.
To construct the modules of the hyperoctahedral group we extend the notions of partitions and tableaux to bi-partitions and bi-tableaux. This new notation allows for the formulation of permutation and Specht modules for $B_{n}$ following a combinatoric method similar to that of $S_{n}$.

\subsection{Bi-partitions}

Every bi-partition has an associated \emph{Young diagram}, which is formed from the Young diagrams for partitions $\po, \pt$ separately. For example, the bi-partition $((3,1), (2^{2},1))$ has Young diagram
\[\ytableausetup{mathmode,baseline,aligntableaux=center,boxsize=1em}\left( \ydiagram{3,1}\, , \, \ydiagram{2,2,1} \right) .\]
In a Young diagram we denote the empty partition $(0)$ as $\emptyset$. 
We refer to the boxes of $\bl$ by triples $(i,j,k)$, which mean box $(i,j)$ in $\lambda^{k}$. For any bi-partition $\bl$ we may form the \emph{transpose} of $\bl$, denoted $\bl^{\prime}$, by swapping the partitions $\lambda^{1},\lambda^{2}$ and taking their transposes separately, i.e., $\bl^{\prime} = (\lambda^{2 \,\prime}, \lambda^{1 \, \prime})$. For example the transpose of bi-partition $((3,1),(2^{2},1))$ is $((3,2),(2,1,1))$, which has corresponding Young diagram,
\[\ytableausetup{mathmode,baseline,aligntableaux=center,boxsize=1em}\left( \ydiagram{3,2}\, , \, \ydiagram{2,1,1} \right) .\]
We may extend the dominance ordering on partitions to a partial ordering on bi-partitions in the following way.

\begin{defn}
	\label{chpt6:def:bidom}
	Let $\bl, \bm$ be bi-partitions.
	Define the \emph{dominance ordering} on bi-partitions as follows: 
	\[\bl \trianglerighteq \bm
	\Leftrightarrow \begin{cases} \text{ we have } & |\lambda^{1}| > 
	|\mu^{1}|,\\ 
	\text{ or } & |\lambda^{1}| = 
	|\mu^{1}|  \text{ and } \lambda^{1}
	\trianglerighteq \mu^{1}, \lambda^{2} \trianglerighteq 
	\mu^{2}.
	\end{cases} 
	\]
	In the case $\bl \trianglerighteq \bm$ we say $\bl$ \emph{dominates} $\bm$. 
\end{defn}
Consider the bi-partitions $\bl = ((2,1),(2,2))$, $\bm = ((3),(2,1,1))$, we do not have $\bl$ dominating $\bm$ nor vice versa, because $\mu^{1} \trianglerighteq \lambda^{1}$ and $\lambda^{2} \trianglerighteq \mu^{2}$. Therefore, the dominance order of bi-partitions is not necessarily a total ordering (however, it may be for small $n$). The dominance ordering does have maximal and minimal elements, which are given by $((n),(0))$ and $((0),(1^{n}))$ respectively.   
Note that if $\bl \trianglerighteq \bm$ we may form $\bm$ from $\bl$ by moving boxes from $\po$ to $\pt$, or moving boxes down and to the left within each separate partition. 

\begin{lemma}
	Let $\bl, \mu$ be bi-partitions of $n$. Then $\bl \trianglerighteq \bm$ if and only if $\bm^{\prime} \trianglerighteq \bl^{\prime}$. 
\end{lemma}
\begin{proof}
	Suppose $\bl \trianglerighteq \bm$, if $|\lambda^{1}| > |\mu^{1}|$ then the lemma follows immediately from the definitions of transpose and dominance ordering. If $|\lambda^{1}| = |\mu^{1}|$ then  also $|\lambda^{2}| = |\mu^{2}|$ and we know $\lambda^{k} \trianglerighteq \mu^{k}$ if and only if $\mu^{k \, \prime} \trianglerighteq \lambda^{k \, \prime}$. Hence, $\bm^{\prime} \trianglerighteq \bl^{\prime}$. 
\end{proof}

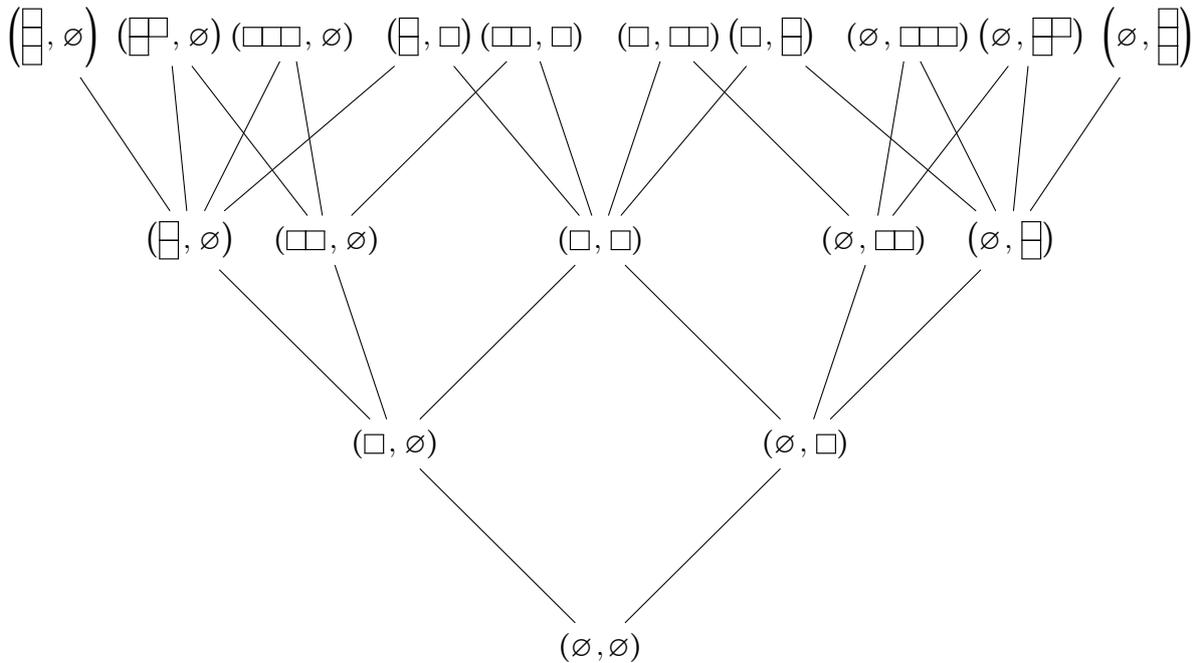
\begin{figure}[H]
	\makebox[\linewidth]{
		\begin{tikzpicture}[scale= 0.9]
		\node (start) at (0,0) {$\left(\emptyset\,,\emptyset \right)$};
		
		\node (11) at (3,3) {$\ytableausetup{mathmode,baseline,aligntableaux=center,boxsize=0.6em}\left( \emptyset\, , \, \ydiagram{1} \right)$};
		\node (2) at (-3,3) {$\ytableausetup{mathmode,baseline,aligntableaux=center,boxsize=0.6em}\left( \ydiagram{1}\, , \, \emptyset \right)$};
		
		\node (left) at (-6,6) {$\ytableausetup{mathmode,baseline,aligntableaux=center,boxsize=0.6em}\left( \ydiagram{1,1}\, , \, \emptyset \right)$};
		\node (3) at (-4,6) {$\ytableausetup{mathmode,baseline,aligntableaux=center,boxsize=0.6em}\left( \ydiagram{2}\, , \, \emptyset \right)$};
		\node (21) at (0,6) {$\ytableausetup{mathmode,baseline,aligntableaux=center,boxsize=0.6em}\left( \ydiagram{1}\, , \, \ydiagram{1} \right)$};
		\node (111)  at (4,6)   {$\ytableausetup{mathmode,baseline,aligntableaux=center,boxsize=0.6em}\left( \emptyset\, , \, \ydiagram{2} \right)$};;
		\node (right) at (6,6)  {$\ytableausetup{mathmode,baseline,aligntableaux=center,boxsize=0.6em}\left( \emptyset\, , \, \ydiagram{1,1} \right)$};
		
		\node (left111) at (-8,9){$\ytableausetup{mathmode,baseline,aligntableaux=center,boxsize=0.6em}\left( \ydiagram{1,1,1}\, , \, \emptyset \right)$};
		\node (left21) at (-6.3,9) {$\ytableausetup{mathmode,baseline,aligntableaux=center,boxsize=0.6em}\left( \ydiagram{2,1}\, , \, \emptyset \right)$};
		\node (left3) at (-4.5,9) {$\ytableausetup{mathmode,baseline,aligntableaux=center,boxsize=0.6em}\left( \ydiagram{3}\, , \, \emptyset \right)$};
		
		\node (r111) at (8,9) {$\ytableausetup{mathmode,baseline,aligntableaux=center,boxsize=0.6em}\left( \emptyset\, , \, \ydiagram{1,1,1} \right)$};
		\node (r21) at (6.3,9) {$\ytableausetup{mathmode,baseline,aligntableaux=center,boxsize=0.6em}\left( \emptyset\, , \, \ydiagram{2,1} \right)$};
		\node (r3) at (4.5,9) {$\ytableausetup{mathmode,baseline,aligntableaux=center,boxsize=0.6em}\left( \emptyset\, , \, \ydiagram{3} \right)$};
		
		\node(m21) at (-1,9) {$\ytableausetup{mathmode,baseline,aligntableaux=center,boxsize=0.6em}\left( \ydiagram{2}\, , \, \ydiagram{1} \right)$};
		\node(m111) at (-2.5,9) {$\ytableausetup{mathmode,baseline,aligntableaux=center,boxsize=0.6em}\left( \ydiagram{1,1}\, , \, \ydiagram{1} \right)$};
		\node(m12) at (1,9) {$\ytableausetup{mathmode,baseline,aligntableaux=center,boxsize=0.6em}\left( \ydiagram{1}\, , \, \ydiagram{2} \right)$};
		\node(m1112) at (2.5,9) {$\ytableausetup{mathmode,baseline,aligntableaux=center,boxsize=0.6em}\left( \ydiagram{1}\, , \, \ydiagram{1,1} \right)$};	
		
		\draw[-] (start) -- (2);
		\draw[-] (start) -- (11);
		
		\draw[-] (2) -- (21);
		\draw[-] (11) -- (111);
		\draw[-] (11) -- (21);
		\draw[-] (2) -- (3);
		\draw[-] (2) -- (left);
		\draw[-] (11) -- (right);
		
		\draw[-] (left) -- (left111);
		\draw[-] (left) -- (left21);
		\draw[-] (left) -- (left3);
		\draw[-] (3) -- (left21);
		\draw[-] (3) -- (left3);
		
		\draw[-] (right) -- (r111);
		\draw[-] (right) -- (r21);
		\draw[-] (right) -- (r3);
		\draw[-] (111) -- (r21);
		\draw[-] (111) -- (r3);
		
		\draw[-] (21) -- (m21);
		\draw[-] (21) -- (m111);
		\draw[-] (21) -- (m12);
		\draw[-] (21) -- (m1112);
		
		\draw[-] (left) -- (m111);
		\draw[-] (3) -- (m21);
		\draw[-] (right) -- (m1112);
		\draw[-] (111) -- (m12);
		\end{tikzpicture}
		
		\caption[Young's lattice for bi-partitions of size $n \in \{0,1,2, 3\}$]{Young's lattice for bi-partitions of size $n \in \{0,1,2, 3\}$.}
		
		\label{chpt6:fig:lattice}
	}
\end{figure}

Given two bi-partitions $\bl, \bm$ of possibly different sizes, we say $ \bl \subseteq \bm$ if $\lambda^{1} \subseteq \mu^{1}$ and $\lambda^{2} \subseteq \mu^{2}$. 
Given a bi-partition of $n$ we may turn in into a bi-partition of $n+1$ by adding a single box to one of the partitions $\po,\pt$. Let the element $e^{1}_{i}$ denote the tuple $(e_{i},(0))$ and $e_{i}^{2} = ((0),e_{i})$, with $e_{i}$ as defined in Section \ref{chpt4:subsec:youngdigrams}. When we add a box to $\bl$ on the $i^{th}$ row of partition $k$ we form a new bi-tuple $\bl + e_{i}^{k}$. We may extend Young's lattice (Figure \ref{chpt4:fig:lattice}) to bi-partitions by restriction our attention to choices of $e_{i}^{k}$ that result in a new bi-partition. This structure may be seen in Figure \ref{chpt6:fig:lattice}.

Each path upwards (respectively downwards) in Figure \ref{chpt6:fig:lattice} represents the placement (respectively removal) of a box to form a Young diagram. This is an important structure as it connects the bi-partitions of $n$ to bi-partitions of $n+1$, this in turn allows us to link the modules of $\mathfrak{B}_{n}$ and $\mathfrak{B}_{n+1}$.

\subsection{Bi-tableaux}

Given a bi-partition $\bl$ we may form a \emph{Young bi-tableau} (alternatively a \emph{$\bl$-tableau}), denoted $\bt = (T^{1},T^{2})$, by putting the numbers $\pm1, \ldots, \pm n$ into the boxes of (the Young diagram of) $\bl$ such that the value $a$ (positive or negative) only appears once. We say the value $a$ \emph{occurs} in a bi-tableau $\bt$ if $\pm a$ appears in a box of $\bt$. In this chapter we refer to bi-tableaux as just tableaux with the understanding that for the hyperoctahedral group we are working with bi-partitions, and bi-tableaux. 
The set of Young tableaux of shape $\bl$ is denoted $\YT(\bl)$. A \emph{standard Young tableau} $\bt$ is a Young tableau where the values in the boxes of $\bt$ 
are all positive, and increasing across rows and down columns within $T^{1}$ and $T^{2}$. 
The set of standard Young tableaux of shape $\bl$ is denoted $\SYT(\bl)$. The size of the set $\SYT(\bl)$ is denoted $d_{\bl}$ and called the dimension of $\bl$.

\begin{lemma}
	The dimension $d_{\bl}$ is equal to the product ${n\choose |\lambda^{1}|} d_{\lambda^{1}}d_{\lambda^{2}}$ where $d_{\lambda^{k}}$ is the dimension of $\lambda^{k}$ as a partition of $|\lambda^{k}|$.
\end{lemma} 
\begin{proof}
	To form a standard Young tableau of shape $\bl$ first split the values of $[n]$ across the two Young diagrams, there are ${n \choose |\lambda^{1}|}$ ways to do this. Then form standard Young tableaux of shape $\lambda^{1}, \lambda^{2}$ with their assigned values, there are $d_{\lambda^{1}}d_{\lambda^{2}}$ way to do this.
\end{proof}

The transpose of a tableau $\bt$, denoted $\bt^{\prime}$, is defined by swapping the tableaux $T^{1},T^{2}$ and taking each transpose separately, i.e., $\bt^{\prime}= (T^{2 \, \prime}, T^{1 \, \prime})$. A tableau $\bt$  is a standard Young tableau if and only if its transpose $\bt^{\prime}$ is a standard Young tableau. If a tableau $\bt$ has a box in position $(i,j,k)$, we let $\bt(i,j,k)$ denote the value in that box; otherwise $\bt(i,j,k)$ is undefined. 
\begin{example}
	The bi-partition $\bl = ((3,1), (2,1))$ has ${7 \choose 4} d_{(3,1)}d_{(2,1)} = 35*6 = 210$ standard Young tableaux, $6$ of these are given below:
	\[
	\ytableausetup{mathmode,baseline,aligntableaux=center,boxsize=1.1em}\left(
	\begin{ytableau} 1 & 2 & 3\\ 4
	\end{ytableau}, \hspace{0.1cm}\begin{ytableau} 5 & 6  \\
	7	\end{ytableau}\right) \hspace{1cm} \left(
	\begin{ytableau} 1 & 2 & 4\\ 3
	\end{ytableau}, \hspace{0.1cm}\begin{ytableau} 5 & 6  \\
	7	\end{ytableau}\right) \hspace{1cm} \left(
	\begin{ytableau} 1 & 3 & 4\\ 2
	\end{ytableau}, \hspace{0.1cm}\begin{ytableau} 5 & 6  \\
	7	\end{ytableau}\right) \]
	\[
	\ytableausetup{mathmode,baseline,aligntableaux=center,boxsize=1.1em}\left(
	\begin{ytableau} 1 & 2 & 3\\ 4
	\end{ytableau}, \hspace{0.1cm}\begin{ytableau} 5 & 7  \\
	6	\end{ytableau}\right) \hspace{1cm} \left(
	\begin{ytableau} 1 & 2 & 4\\ 3
	\end{ytableau}, \hspace{0.1cm}\begin{ytableau} 5 & 7  \\
	6	\end{ytableau}\right) \hspace{1cm} \left(
	\begin{ytableau} 1 & 3 & 4\\ 2
	\end{ytableau}, \hspace{0.1cm}\begin{ytableau} 5 & 7  \\
	6	\end{ytableau}\right).\]
	The remaining standard Young tableaux may be found by choosing different divisions of the numbers $\{1,2,3,4,5,6,7\}$ between the partitions $(3,1)$ and $(2,1)$.
\end{example}

Each standard Young tableau of shape $\bl$ corresponds to one path up Young's lattice (Figure \ref{chpt6:fig:lattice}) starting at $\emptyset$ and ending at $\bl$. To form this correspondence take any $\bt \in \SYT(\bl)$ and create a path up Young's lattice by adding boxes to $\emptyset$ in the order given by the entries in $\bt$. 

The following result is an extension of Lemma \ref{chpt4:lem:tabdom} which establishes a link between Young tableaux and their respective shapes. In particular by comparing two young tableaux of possibly different shapes we may learn whether one bi-partition dominates the other.

\begin{defn}
	Let $\tl$ be a $\bl$-tableau and $\tm$ be a $\bm$-tableau. We say that $\tl$ 
	\emph{agrees} with $\tm$, if for all $a$ that occur in $\tl^{k}$ we have $\pm a$ occurring in $\tm^{k}$ for $k\in \{1,2\}$. If the tableaux $\tl$, $\tm$ agree then we must have $|\lambda^{k}| = |\mu^{k}|$ for $k\in\{1,2\}$.
\end{defn}

\begin{lemma}
	\label{chpt6:lem:tabdom}
	Let $\tl$ be a $\bl$-tableau and 
	$\tm$ be a $\bm$-tableau that agree. Suppose for each index $i$ the elements in row $i$ of $\tm^{k}$ appear in different columns of $\tl^{k}$. Then 
	$\bl \trianglerighteq \bm$. 
\end{lemma}

\begin{proof}
	We know $|\lambda^{k}| = |\mu^{k}|$ for $k \in \{1,2\}$ and thus using the 
	same argument as Lemma \ref{chpt4:lem:tabdom}  on the first and second tableaux separately we have $\lambda^{1} \trianglerighteq \mu^{1}$, $\lambda^{2} \trianglerighteq \mu^{2}$.
\end{proof}

\subsection{The Module Structure of the Hyperoctahedral Group}

We are now in a position to define the permutation and Specht modules for the hyperoctahedral group. We find one permutation module for every bi-partition $\bl$ which contains the simple Specht modules as special submodules. The following proofs are adapted from several sources, including Geissinger and Kinch \cite{geissinger1978representations}, Aamily, Morris, and Peel \cite{al1981representations}, and Can \cite{can1996representations}. We have adapted the constructions presented in the literature to allow us to re-establish Young's rule (Lemma \ref{chpt4:lem:youngrule}) for the permutation modules of the hyperoctahedral group.

\subsubsection{Permutation Modules}

To define the permutation modules of $B_{n}$ we reintroduce the idea of a row permutation for a tableau $\bt$, this defines an equivalence relation on the Young tableaux of a fixed bi-partition.

\begin{defn}
	\label{chpt6:def:tabaction}
	Let $\bt$ be a Young tableau of size $n$, and $\sigma \in B_{n}$.  Define an action of $\sigma$ on $\bt$ by applying $\sigma$ to the values of $\bt$ box wise. 
\end{defn}

\begin{defn}
	Let $\bt$ be a bi-tableau of size $n$. Define the subgroup of \emph{row permutations} of $\bt$, denoted $R_{\bt}$, as the set of elements in $B_{n}$ which fix rows of $T^{1}$ up to changes in sign, and fix rows in $T^{2}$ completely. Define the subgroup of \emph{column permutations} of $\bt$, denoted $C_{\bt}$, as the set of elements in $B_{n}$ which fix columns in $T^{2}$ up to changes in sign, and fix columns in $T^{1}$ completely.
\end{defn}

\begin{defn}
	Define an equivalence relation on tableaux of shape $\bl$ called \emph{row equivalence}, denoted $\sim_{R}$, in the following way:
	\[\overline{T_{1}} \sim_{R} \overline{T_{2}} 	\Leftrightarrow \exists\sigma \in R_{\overline{T_{1}}}\textnormal{ such that } \sigma \,  \overline{T_{1}} = \overline{T_{2}}.\]
	The row equivalence class of a tableau $\overline{T}$ is denoted by $\{\overline{T}\}$ and called a  \emph{(bi-)tabloid}.

\end{defn}

\begin{example}
	Let $\bl = ((2,1), (2))$. Consider the following three $\bl$-tableaux,
	\[\ytableausetup{mathmode,baseline,aligntableaux=center,boxsize=1.4em} \overline{T}_{1} = \left(	\begin{ytableau} 1 & 5 \\
	3\end{ytableau}\, , 	\begin{ytableau} 4 & 2\end{ytableau} \right) \hspace{1cm} \overline{T}_{2} = \left(	\begin{ytableau} 5 & -1 \\ 	3\end{ytableau}\, , 	\begin{ytableau} -4 & 2 \end{ytableau} \right) \hspace{1cm} \overline{T}_{3} = \left(	\begin{ytableau} -1 & -5 \\		-3\end{ytableau}\, , 	\begin{ytableau} 2 & 4\end{ytableau} \right) .\]
	The tableaux $\overline{T}_{1}$, $\overline{T}_{2}$ are not row equivalent because we have changed the sign of $4$ in the second tableau. The tableaux $\overline{T}_{1}$, $\overline{T}_{3}$ are row equivalent for the permutation $\left(\xi_{1} \xi_{3} \xi_{5} (2 \, 4) \right) \in R_{T_{1}}$, therefore, $\{\overline{T}_{1}\} = \{\overline{T}_{3}\}$.
\end{example}

The action of $B_{n}$ on tableaux, naturally extends to an action on tabloids, given by $\sigma \, \{\bt\} = \{\sigma \, \bt\}$. We now construct the permutation modules for $\mathfrak{B}_{n}$.
\begin{defn}
	\label{chpt6:def:permmod}
	Let $\bl \vdash n$. The \emph{permutation module} for $\mathfrak{B}_{n}$, denoted $M^{\bl}$, is the following vector space
	\[M^{\bl} = \langle \{\overline{T} \} \, | \, T\textnormal{ is a tableau of shape } \bl \rangle \]
	with action of $\mathfrak{B}_{n}$ given by extending the action of $B_{n}$ on tabloids (Definition \ref{chpt6:def:tabaction}) linearly.
\end{defn}

Note that just simply extending the previous definition of row permutations (Definition \ref{chpt4:def:rowequiv}) for single tableaux $T$ to both tableaux $T^{1}, T^{2}$ is not strong enough to distinguish all the permutation modules of $B_{n}$ from one another. For example, if we study row permutations of bi-tableaux which fix rows completely in both separate tableaux, then for any partition $\lambda \vdash n$ the modules $M^{(\lambda,(0))}$ and $M^{((0),\lambda)}$ would be isomorphic to one another.

\begin{lemma}
	The permutation module $M^{\bl}$ is cyclic and has dimension
	\[	2^{|\lambda^{2}|}\frac{n!}{\prod_{i=1}^{2}\left(\prod_{j=1}^{l(\lambda^{i})} \lambda_{j}^{i}! \right)} .\]
\end{lemma}
\begin{proof}
	The main fraction above is a multinomial coefficient formed from $n$ and the size of every row in $\bl$. The extra factor $2^{|\lambda^{2}|}$ represents the fact that all entries of $\lambda^{2}$ have fixed signs. The permutation module being cyclic follows from the row equivalence relation.
\end{proof}

\subsubsection{Specht Modules}
We construct the Specht modules $S^{\bl}$ for the hyperoctahedral group as stable subspaces of the corresponding permutation modules $M^{\bl}$.
We form a basis for the Specht module $S^{\bl}$ by taking linear combinations of tabloids, we call these new elements \emph{polytabloids}.

\begin{defn}
	Given a Young tableau $\bt$, form the element $k_{\bt} \in \mathfrak{B}_{n}$ as the following sum:
	\begin{eqnarray}
	k_{\bt} = \sum_{\sigma \in C_{\bt}} \sign(\sigma) \, \sigma.
	\end{eqnarray} The \emph{polytabloid} associated to tabloid $\bt$, denoted $e_{\bt}$, is the element $e_{\bt} = k_{\bt} \, \{\overline{T}\}$. We say that the tabloid $\{\bt^{\prime}\}$ is contained in the polytabloid $e_{\bt}$ if it appears with a non-zero coefficient.
\end{defn}

\begin{lemma}
	\label{chpt6:lem:spechtcyclic}
	Let $\bt$ be a Young tableau of size $n$, and $\eta \in B_{n}$. Then $\eta \, e_{\bt} = e_{\eta \, \bt}$.  
\end{lemma}
\begin{proof}
	\begin{eqnarray}
	\eta \, e_{\bt} =  \sum_{\sigma \in C_{\bt}}  \sign(\sigma) \, \eta \,  \sigma \,\eta^{-1} \{\eta \bt\}  =  \sum_{\sigma \in C_{\eta \, \bt}}  \sign(\sigma)  \,  \sigma \, \{\eta \bt\} = e_{\eta \,\bt}. \nonumber \qedhere 
	\end{eqnarray}
\end{proof}

\begin{defn}
	Let $\bl \vdash n$. The \emph{Specht module} for $\mathfrak{B}_{n}$, denoted $S^{\bl}$, is defined as the following vector space,
	\[S^{\bl} = \langle e_{\bt} \, | \, \bt \in \YT(\bl) \rangle \]
	with the natural action of $\mathfrak{B}_{n}$ on tabloids. 	In fact we may restrict the spanning set $\{ e_{\bt} \, | \, \bt \in \YT(\bl) \}$ to a basis of polytabloids formed from standard Young tableaux of $\bl$ (see \cite[Sections 4 and 5]{can1996representations}), that is
	\begin{eqnarray}
	S^{\bl} = \langle e_{\bt} \, | \, \bt \in \SYT(\bl) \rangle.
	\end{eqnarray}
	Thus, the Specht module $S^{\bl}$ has dimension $d_{\bl}$, and is cyclic by applying Lemma \ref{chpt6:lem:spechtcyclic}.
\end{defn}

\begin{example}
	\label{chpt6:ex:specht}
	Let $\bl = ((2,1), (2))$, the Specht module $S^{\bl}$ has dimension $d_{\bl} = {5 \choose 3} d_{(2,1)}d_{2} = 20$. Take $\overline{T}$ to be the following standard Young tableau of shape $\bl$,
	\[\ytableausetup{mathmode,baseline,aligntableaux=center,boxsize=1.1em} \overline{T} = \left(	\begin{ytableau} 1 & 2 \\3\end{ytableau}\, , 	\begin{ytableau} 4 & 5\end{ytableau} \right).\]
	The group of column permutations for any tableau of shape $\bl$ has size $8$, as we may swap elements in the first column of $\overline{T}^{1}$, and swap the sign of every value in $\overline{T}^{2}$. Starting from the tabloid $\bt$ we can compute its associated polytabloid:
	\begin{eqnarray*}
		e_{\overline{T}} = k_{\overline{T}} \, \{\overline{T}\} & = & \ytableausetup{mathmode,baseline,aligntableaux=center,boxsize=1.3em} \left \lbrace	\begin{ytableau} 1 & 2 \\3\end{ytableau}\, , 	\begin{ytableau} 4 & 5\end{ytableau} \right\rbrace \, - \,  \left \lbrace	\begin{ytableau} 3 & 2 \\1\end{ytableau}\, , 	\begin{ytableau} 4 & 5\end{ytableau} \right\rbrace \, - \,  \left \lbrace	\begin{ytableau} 1 & 2 \\3\end{ytableau}\, , 	\begin{ytableau} -4 & 5\end{ytableau} \right\rbrace \, + \,  \left \lbrace	\begin{ytableau} 3 & 2 \\1\end{ytableau}\, , 	\begin{ytableau} -4 & 5\end{ytableau} \right\rbrace \\ 
		& - & \ytableausetup{mathmode,baseline,aligntableaux=center,boxsize=1.3em} \left \lbrace	\begin{ytableau} 1 & 2 \\3\end{ytableau}\, , 	\begin{ytableau} 4 & -5\end{ytableau} \right\rbrace \, + \,  \left \lbrace	\begin{ytableau} 3 & 2 \\1\end{ytableau}\, , 	\begin{ytableau} 4 & -5\end{ytableau} \right\rbrace \, + \,  \left \lbrace	\begin{ytableau} 1 & 2 \\3\end{ytableau}\, , 	\begin{ytableau} -4 & -5\end{ytableau} \right\rbrace \, - \,  \left \lbrace	\begin{ytableau} 1 & 2 \\3\end{ytableau}\, , 	\begin{ytableau} -4 & -5\end{ytableau} \right\rbrace . 
	\end{eqnarray*}
\end{example}
\paragraph{}
We have seen that the conjugacy classes of $B_{n}$ are labelled by bi-partitions, and we have found exactly one Specht module for each bi-partition $\bl$. We now prove that the Specht modules form a complete set of simple modules for the hyperoctahedral group. To do this we study the properties of the groups action on polytabloids. We first show every Specht module to be simple and then prove that $S^{\bl} \cong S^{\bm}$ if and only if $\bl = \bm$. The following arguments take inspiration from the construction of Specht modules for the symmetric group (Section \ref{chpt4:subsec:spechtmod}), and work by Can \cite{can1996representations}.

\begin{lemma}
	\label{chpt6:lem:polyzero}
	Let $\bl,\bm \vdash n$. Let $\tl$ be a $\bl$-tableau and $\tm$ be a $\bm$-tableau. Suppose that at least one of the following conditions holds: 1) there exists $ a\in \tm^{1}$ such that $\pm a \in \tl^{2}$, 2) there exists $a,b$ in the same row of $\tm^{k}$, with  $c =\pm a ,d = \pm b$ in the same column of $\tl^{k}$. Then $k_{\tl} \{\tm\} = 0$.
\end{lemma}

\begin{proof}
	First we establish that $C_{\tl} \cap R_{\tm} \neq \{\gpid\}$. If the first condition is met then clearly $\xi_{a} \in C_{\tl} \cap R_{\tm}$. If the second condition is met with $k=1$ then $(a \, b) \in C_{\tl} \cap R_{\tm}$, if $k=2$ then $(c \, d) \in C_{\tl} \cap R_{\tm}$.  Let $\tau$ be a transposition in $C_{\tl} \cap R_{\tm} $, then $(\gpid - \tau) \{\tm \}=0$. 
	All choices of $\tau$ generate a subgroup of order $2$ in $C_{\tl}$, taking signed coset representatives $\sigma_{1}, \ldots \sigma_{m}$ for this subgroup. Then,
	\[k_{\tl} \{\tm\} = \left(\sum_{\sigma \in C_{\tl}} 
	\sign(\sigma) \,\sigma \right) \{\tm\} =  \left(\sum_{i=1}^{m} 
	\sigma_{i} (\gpid - \tau)\right) \{\tm\} = 0. \qedhere \]
\end{proof}

\begin{lemma}
	\label{chpt6:lem:notdom}
	Let $\bl,\bm \vdash n$. Let $\tl$ be a $\bl$-tableau and $\tm$ be a $\bm$-tableau. Suppose $\bl \ntrianglerighteq \bm$, then $k_{\tl} \{\tm\} = 0$.
\end{lemma}
\begin{proof}
	Suppose $\tl$ does not agree with $\tm$ 
	then by our hypothesis there must exist $a \in \tm^{1}$ such 
	that $\pm a\in \tl^{2}$ thus the conclusion holds by 
	Lemma \ref{chpt6:lem:polyzero}. However, if $T_{\bl}$ agrees with $T_{\bm}$ the conclusion holds by  Lemma \ref{chpt6:lem:polyzero} and Lemma \ref{chpt6:lem:tabdom}. 
\end{proof}

\begin{corollary}
	\label{chpt6:cor:dom}
	Let $\bl,\bm \vdash n$. Let $\tl$ be a $\bl$-tableau and $\tm$ be a $\bm$-tableau. Suppose 
	$k_{\tl} \{\tm\} \neq 0$, then $\bl \trianglerighteq \bm$.  Furthermore, if $\bl = \bm$ then $k_{\tl} \{\tm\} = \pm e_{\tl}$.
\end{corollary}

\begin{proof}
	The first statement is the contrapositive of Lemma \ref{chpt6:lem:notdom}.
	Let $\bl = \bm$ and $k_{\tl} \{\tm\} \neq 0$. The tableaux $\tl, \tm$ must agree and furthermore, if $a,b$ belong to the same row of $\tl^{k}$ then $c= \pm a, \,d =\pm b$ belong to different columns of $\tm$. This implies that 	
	there exists a column permutation $\eta \in C_{\tl}$ such that $\{\tm \}= \eta \{\tl\}$ for more details see \cite[Lemma 3.6]{can1996representations}. Therefore,
	\[k_{\tl} \{\tm\} = \left( \sum_{\sigma \in C_{\tl}} \sign(\sigma) \sigma \eta \right) \{\tl\}  = \left( \sum_{\sigma \in C_{\tl}} \sign(\eta^{-1}) \sign(\sigma)\sigma \right) \{\tl\}  = \sign(\eta) e_{\tl} . \qedhere \]
\end{proof}

\begin{corollary}
	\label{chpt6:cor:agreeconstant}
	Let $\tl$ be a $\bl$-tableau and $v \in 
	M^{\bl}$. Then $k_{\tl} \,v$ is a multiple (possibly zero) of $e_{\tl}$.
\end{corollary}
\begin{proof}
	Write $v = \sum_{i} c_{i} \{\overline{T}_{i} \}$, where $\overline{T}_{i}$ is a $\tl$-tableau. Each summand $k_{\tl} \{\overline{T}_{i}\}$ is either zero or a multiple of $e_{\tl}$ by Corollary \ref{chpt6:cor:dom}.
\end{proof}

We are now in a position to prove that every Specht module is simple, and that they are pairwise non-isomorphic to each other. 	The proofs of the next theorems follow closely from their symmetric group counterparts (Theorems \ref{chpt4:thm:submodthm} and \ref{chpt4:thm:spechtmaps}).

\begin{thm}[Submodule Theorem]
	Let $V \subseteq M^{\bl}$ be a submodule. Then $V \supseteq S^{\bl}$ or $V \subseteq (S^{\bl})^{0}$. Therefore, the Specht module $S^{\bl}$ is simple.
\end{thm}
\begin{proof}
	Take $v \in V$, and $\tl$ a $\bl$-tableau. By Corollary \ref{chpt6:cor:agreeconstant} we have $k_{\tl} v = c e_{\tl}$ for some $c\in\mathbb{C}$. Suppose we have $c\neq 0$ for some choice of  $\tl$, then $c^{-1} k_{\tl} \,v = e_{\tl} \in V$ hence we may generate $S^{\bl}$ inside $V$.
	
	Now suppose that $c = 0$ for all choices of tableau $\tl$, then consider the inner product on $M^{\bl}$ defined on its basis by $\langle \{\overline{T}_{1} \} \,, \, \{\overline{T}_{2} \} \rangle = \delta_{\{\overline{T}_{1} \} \,, \, \{\overline{T}_{2} \}}$. This inner product is invariant under the action of $B_{n}$, so we find
	\[\langle v, e_{\tl} \rangle = \sum_{\sigma \in C_{\tl}} \langle v, \sign(\sigma) \, \sigma \, \{\tl\} \rangle = \sum_{\sigma \in C_{\tl}} \langle \sign(\sigma) \, \sigma^{-1} \, v, \{\tl\} \rangle =  \langle k_{T_{\bl}} \,v , \{\tl\} \rangle= \langle 0, \{\tl\} \rangle= 0.\]
	Therefore no polytabloid appears in $M^{\bl}$, and because any single polytabloid would span $S^{\bl}$, we have $v \notin S^{\bl} \Rightarrow v \in (S^{\bl})^{0}$. 
	
\end{proof}

\begin{thm}
	\label{chpt6:thm:noniso}
	Let $\bl,\bm \vdash n$. Suppose there is a non-zero homomorphism $\psi: S^{\bl} \to M^{\bm}$. Then $\bl \trianglerighteq \bm$, and if $\bl = \bm$ then $\psi$ is multiplication by a scalar.
\end{thm}

\begin{proof}
	Take a basis vector $e_{\tl} \in S^{\bl}$ such that $\psi(e_{\tl}) \neq 0$. Extend the homomorphism $\psi$ to the whole module $M^{\bl}$ by setting it to be $0$ on the complement of $S^{\bl}$. Then
	\[0 \neq \psi(e_{\tl}) = k_{\tl} \psi( \{\tl\}) = k_{\tl} \left( \sum_{i} c_{i} \{\overline{T}_{i}\} \right) \]
	where the $\overline{T}_{i}$ are $\bm$-tableau. We must have at least one $c_{i}$ and product $k_{\tl} c_{i} \{\overline{T}_{i}\}$ being non-zero. Therefore, by Corollary \ref{chpt6:cor:dom} we have $\bl \trianglerighteq \bm$. 
	
	If $\bl = \bm$ we have that $\psi(e_{\tl}) = c \cdot e_{\tl}$ for some constant $c \in \mathbb{C}$. Hence, for any $\sigma \in B_{n}$, using Lemma \ref{chpt6:lem:spechtcyclic} we have
	\[\psi(e_{\sigma \tl}) = \psi( \sigma \, e_{\tl}) = \sigma \psi(e_{\tl}) = c\cdot \sigma e_{\tl} = c\cdot e_{\sigma \, \tl}. \qedhere\]
\end{proof}

\begin{corollary}
	\label{chpt6:cor:simp}
	The Specht modules $S^{\bl}$ for $\bl \vdash n$ form a complete set of pairwise non-isomorphic simple modules for $\mathfrak{B}_{n}$.
\end{corollary}
\begin{proof}
	The Submodule Theorem tells us the Specht modules are simple. Now if $\bl = \bm$ we may easily see that $S^{\bl} \cong S^{\bm}$. Conversely if $S^{\bl} \cong S^{\bm}$ then we have non-zero homomorphisms $\psi:S^{\bl} \to M^{\bm}$ and $\phi:S^{\bm} \to M^{\bl}$. Therefore, by Theorem \ref{chpt6:thm:noniso} we must have $\bl \trianglerighteq \bm$ and $\bm \trianglerighteq \bl$ which implies $\bl = \bm$.
\end{proof}

Given a permutation module $M^{\bl}$ we would like to decompose it into its simple submodules. An application of Theorem \ref{chpt6:thm:noniso} allows us to re-establish Young's rule for the the permutation modules of the hyperoctahedral group.

\begin{lemma}[Young's Rule]
	\label{chpt6:lem:youngsrule}
	For $\bm \vdash n$ we have,
	\begin{eqnarray}
	M^{\bm} \cong \bigoplus_{\substack{\bl \vdash n \\ \bl \trianglerighteq \bm }}K_{\bl,\bm} S^{\bl},
	\end{eqnarray}
	where $K_{\bl,\bm} S^{\bl}$ denotes a direct sum of $K_{\bl,\bm}$ copies of $S^{\bm}$. We call the coefficients $K_{\bl,\bm} \in \mathbb{N}^{0}$ \emph{generalised Kostka numbers}. We know that $K_{\bl,\bl} =1$ for all $\bl \vdash n$. Furthermore, if $|\lambda^{1}| = |\mu^{1}|$ then $K_{\bl,\bm} = K_{\lambda^{1},\mu^{1}} \, K_{\lambda^{2},\mu^{2}}$.  
\end{lemma}
\begin{proof}
	Theorem \ref{chpt6:thm:noniso} tells us that if $S^{\bl}$ appears as a summand of $M^{\bm}$ we must have $\bl \trianglerighteq \bm$. If $\bl=\bm$ we know that any homomorphism $S^{\bl} \to M^{\bl}$ is multiplication by a scalar hence there must only be one copy of $S^{\bl}$ in $M^{\bl}$. The last assertion follows from work by Geissinger and Kinch \cite[Corollary II.3]{geissinger1978representations}.
\end{proof}

Similarly to the symmetric group case the regular module $\mathfrak{B}_{n}$ is isomorphic to the permutation module associated to the bi-partition at the bottom of the dominance ordering, that is the partition $((0),(1^{n}))$.

\begin{lemma}
	\label{chpt6:lem:regdecomp}
	The permutation module $M^{((0),(1^{n}))} \cong \mathfrak{B}_{n}$ as $\mathfrak{B}_{n}$-modules. Thus, it has canonical decomposition,
	\[M^{((0),(1^{n}))} \cong \bigoplus_{\bl \vdash n} d_{\bl} S^{\bl}\textnormal{ as $\mathfrak{B}_{n}$-modules} .\]
\end{lemma}
\begin{proof}
	Define a isomorphism of $\mathfrak{B}_{n}$-modules $\psi:\mathfrak{B}_{n} \to M^{((0),(1^{n}))}$ by its action on a single permutation $\sigma \in B_{n}$:
	\[\psi(\sigma) = 
	\ytableausetup{mathmode,baseline,aligntableaux=center,boxsize=2.2em}
	\left\lbrace \emptyset, \, \begin{ytableau} \sigma(1) \\ \dots  \\ 
	\sigma(n)\end{ytableau} \,
	\right\rbrace. \qedhere\]
\end{proof}

The version of Young's rule we stated for the hyperoctahedral group is weaker than the result for the symmetric group (Lemma \ref{chpt4:lem:youngrule}). For the purpose of lifting eigenvectors we do not require that every Specht module $S^{\bl}$ with  $ \bl \trianglerighteq \bm$ appears in the decomposition of $M^{\bm}$. 

The hyperoctahedral group has a natural inclusion structure $B_{m} \hookrightarrow B_{n}$ for $m<n$, by extending every permutation of $B_{m}$ to a permutation of $B_{n}$ by choosing it to fix all elements of $[\pm m] \setminus [\pm n]$. The group algebras of the hyperoctahedral group and subsequently $\mathfrak{B}_{n}$-modules inherit this recursive structure.	To end this section we state the branching rules for the Specht modules of the hyperoctahedral group, a proof and further details of this result may be found in \cite[Section III]{geissinger1978representations}.

\begin{thm}[Branching rules for $B_{n}$]
	\label{chpt6:thm:branching}
	Let $n \geq 2$, and $\bl \vdash n$. The \emph{branching rules} for the simple modules of the hyperoctahedral group are as follows:
	\begin{eqnarray}
	\res_{B_{n-1}}^{B_{n}} S^{\bl} & \cong & \bigoplus_{\substack{\bm \vdash n-1 \\ \bm \subseteq \bl}} S^{\bm} \textnormal{ as $\mathfrak{B}_{n-1}$-modules}\label{chpt6:eqn:restrict}\\
	\induce_{B_{n}}^{B_{n+1}} S^{\bl} & \cong &\bigoplus_{\substack{\bm \vdash n+1 \\ \bl \subseteq \bm }} S^{\bm} \textnormal{ as $\mathfrak{B}_{n+1}$-modules} \label{chpt:eqn:induce}
	\end{eqnarray}
\end{thm}

The branching rules for Specht modules are closely related to the edges of Young's lattice. When we induce a Specht module $S^{\bl}$ one step we create a direct sum of the Specht modules associated to bi-partitions directly above $\bl$ in Young's lattice, conversely when we restrict one step we find a direct sum of Specht modules associated to bi-partitions directly below $\bl$ in Young's lattice.

\subsection{Switching to Words}

Let $\bl \vdash n$, and take a tableau $\bt $ of shape $\bl$. If we add a box to $\bt$ containing the value $\pm(n+1)$ we transform $\bt$ to a tableau of size $n+1$. Similarly we can transform a tabloid $\{\bt\}$ if size $n$ into a tabloid of size $n+1$ by adding a single box. This is exactly how we lift vectors from permutation module $M^{\bl}$ to permutation modules of $\mathfrak{B}_{n+1}$. Notationally this process is cumbersome, therefore we instead form a correspondence between the study of tabloids to the study of words over a finite alphabet. This allows the algebra of lifting eigenvectors in Section \ref{chpt6:sec:randomwalks} to be presented in a concise way. The construction of the correspondence follows closely to that of the symmetric group, however the alphabet and words now contain signed and unsigned letters.

\paragraph{}
Define the set $[n^{\pm}] = \{1^{+}, \ldots, n^{+}\} \cup \{1^{-}, \ldots, n^{-}\}$, viewed as $2n$ distinct symbols.  We call the letters in $[n^{\pm}]$ \emph{signed}, and letters in $[n]$  \emph{unsigned}.
Given $n\in \mathbb{N}^{0}$ we denote by $W^{\overline{n}}$ the set of words of length $n$ 
with letters from the set $[\overline{n}] = [n] \cup [n^{\pm}]$. For $n=0$ we allow the empty word $\omega$ as the unique word of length zero. The size of $W^{\overline{n}}$ is $(3 n)^{n}$.

The hyperoctahedral group has a natural action on $W^{\overline{n}}$. Take $\sigma \in B_{n}$ with decomposition $\sigma = \left( \prod_{i=1}^{n} \xi_{i}^{x_{i}} \right) \eta $, then $\sigma$ acts on the word $w \in W^{\overline{n}}$ in the following way: permute the letters in $w$ by $\eta$ (as defined in Section \ref{chpt4:subsec:words}), then for all $x_{i}=1$ flip the sign of $w_{i}$ if possible.
For example, take $ w = 1^{+}2 4^{-}3^{+}3$ and $\sigma = \left(\xi_{1}\xi_{4}\xi_{5}\right) (1 \, 2 \, 3)$, then $\sigma w = 4^{+}1^{+}23^{-}3$, note that $\xi_{5}$ has no affect on $w_{5} = 3$ because it is an unsigned letter. Let $M^{\overline{n}}$ be the vector space over $\mathbb{C}$ spanned by words in $W^{\overline{n}}$. The action of $B_{n}$ on $W^{\overline{n}}$ extends linearly to an action of $\mathfrak{B}_{n}$ on $M^{\overline{n}}$, thus the vector space $M^{\overline{n}}$ is a $\mathfrak{B}_{n}$-module.

\paragraph{}
To each word $w \in W^{\overline{n}}$ we may associate a bi-tuple of size $n$ of non-negative integers called its \emph{evaluation}, denoted $\overline{eval}$, as follows. Define $\eval_i(w)$ to count the number of occurrences of unsigned symbol $i$ in $w$, and let $\eval(w) := (\eval_{1}, \ldots, \eval_{n})$. Define $\eval_{i}^{\pm}(w)$ to count the number of occurrences of signed symbol $i^{\pm}$ in $w$ and let $\eval^{\pm}(w) := (\eval_{1}^{\pm}, \ldots, \eval_{n}^{\pm})$. Then for a word $w \in W^{\overline{n}}$ define its evaluation as the tuple $\overline{\eval}(w) := (\eval(w), \eval^{\pm}(w))$.  Note that the evaluation of any word in $W^{\overline{n}}$ sums to $n$. For example, the word $1^{+}24^{-}3^{+}3$ has evaluation $((0,1,1,0,0),(1,0,1,1,0))$. If the bi-tuple $\overline{eval}(w)$ is non-decreasing we associate it with the corresponding bi-partition, e.g. $1^{+}112^{-}2$ has evaluation $((2,1,0^{2}),(1,1,0^{3}))$, and we associate this bi-tuple with the bi-partition $((2,1),(1,1))$.
The evaluation of any word is stable under the action of $B_{n}$ because we only permute symbols and change the sign of signed symbols, both of which do not affect $\eval^{\pm}(w)$. Therefore, we may form a submodule of $M^{\overline{n}}$ by restricting to words with a given evaluation.

\begin{defn}
	Let $\overline{\nu}$ be a bi-tuple of non-negative integers which sum to $n$. Define the $\mathfrak{B}_{n}$-module $M^{\overline{\nu}}$ as the following stable vector space
	\[M^{\overline{\nu}} = \langle w \in W^{\overline{n}} \, | \, \overline{\eval}(w) = \overline{\nu} \, \rangle \subseteq M^{\overline{n}}.\]
\end{defn}

Let $\overline{\nu} \vdash n$, then the module $M^{\overline{\nu}}$ is isomorphic to the permutation module defined in Definition \ref{chpt6:def:permmod}. To establish the equivalence of these two definitions we construct a bijection from tabloids of shape $\bl$ to words with evaluation $\bl$. 

\begin{defn}
	Let $\bl \vdash n$. Define a map $w:\YT(\bl) \to W^{\overline{n}}$ as follows: for each tableau $\bt$ of shape $\bl$, let $w(\bt) = w_{1} \, \ldots \, w_{n}$ be the word with $w_{\bt(i,j,1)} = i$ for each box $(i,j,1) \in \bt$, and $W_{\bt(i,j,2)} = i^{\pm}$ for each box $(i,j,2) \in \bt$ with the sign of $i$ given by the sign present in box $(i,j,2)$ of $\bt$.
	The word $w(\bt)$ formed by this process has evaluation $\bl$. 
	The map $w$ respects the action of $B_{n}$ on tableaux and words, thus forming a homomorphism between the $\mathfrak{B}_{n}$-modules of $\bl$-tableaux and words in $M^{\bl}$. 		
	
\end{defn}

\begin{example}
	Let $\bl = ((2,1),(1,1))$, and take $\bt$ a $\bl$-tableau, 
	\[\overline{T} =\ytableausetup{mathmode,baseline,aligntableaux=center,boxsize=1.3em} \left(\begin{ytableau} 2 & 1 \\-4\end{ytableau}\, , \, \begin{ytableau} -3 & 5\end{ytableau}\right) \textnormal{ which has corresponding word } w(\overline{T}) =  111^{-}21^{+}.\]
	
	Now take $\sigma \in B_{n}$, the action of $\sigma$ commutes with the linear map $w$. For example take  $\sigma = \left(\xi_{1}\xi_{2}\right)(1 \, 3) (2\, 4\,5)$, then
	\[\sigma \overline{T} =  \ytableausetup{mathmode,baseline,aligntableaux=center,boxsize=1.3em} \left(\begin{ytableau} 4 & 3 \\-5\end{ytableau}\, , \, \begin{ytableau} 1 & -2\end{ytableau}\right) \textnormal{ which has corresponding word } w(\sigma \overline{T}) = 1^{+}1^{-}112 = \sigma \, w(\overline{T}).\]
\end{example}

The homomorphism $w:\YT(\bl) \to M^{\bl}$ is surjective but not injective. To form an isomorphism between the two definitions of the permutation modules $M^{\bl}$ we restrict the domain of $w$ from tableaux of shape $\bl$, to tabloids of shape $\bl$ which are in one-to-one-correspondence with words of evaluation $\bl$.

\begin{lemma}
	Let $\bl \vdash n$ and $\overline{T}_{1}, \overline{T}_{2}$ be two $\bl$-tableaux. Then
	\[\{\overline{T}_{1}\} = \{\overline{T}_{2}\}  \Leftrightarrow \overline{T}_{1} \sim_{R} \overline{T}_{2} \Leftrightarrow w(\overline{T}_{1}) = w(\overline{T}_{2}).\]
\end{lemma}
\begin{proof}
	This follows from the definition of $w$. Any row permutation fixes the occurrences of letters in our word, but may change the sign of those in the first tableau. However, the letters represented by the first tableau are unsigned, therefore changing the sign does not create a new word.
\end{proof}

\begin{lemma}
	Let $\bl \vdash n$, the permutation module $M^{\bl}$  may be seen as a vector 
	space over the following bases:
	\begin{eqnarray*}
		M^{\bl} & = & \langle \{\overline{T}\} \,|\, \overline{T} \textnormal{ is a tableau of 
			shape 
		} 	\bl \, \rangle \\
		& \cong & \langle w(\overline{T}) \, | \,\overline{T} \textnormal{ is a tableau of shape } 	
		\bl \,\rangle \\ 
		& \cong & \langle w \in M^{\overline{n}} \, | \,\overline{\eval}(w) = \bl \, \rangle
	\end{eqnarray*}  
\end{lemma}
\begin{proof}
	The map $w: M^{\bl} \to M^{\bl}$ forms a bijection between tabloids of shape $\bl$ and 
	words of evaluation $\bl$, which respects the action of $\mathfrak{B}_{n}$.	
\end{proof}

Using the one-to-one correspondence between tabloids and words we may also use words to describe the Specht modules for $\mathfrak{B}_{n}$. In example \ref{chpt6:ex:specht} we saw an explicit description of a polytabloid $e_{\overline{T}}$ belonging to the module $S^{((2,1),(2))}$. We can restate this element replacing every tabloid by its corresponding word: 
\begin{eqnarray}
e_{\overline{T}}  & = & \ytableausetup{mathmode,baseline,aligntableaux=center,boxsize=1.3em} \left \lbrace	\begin{ytableau} 1 & 2 \\3\end{ytableau}\, , 	\begin{ytableau} 4 & 5\end{ytableau} \right\rbrace \, - \,  \left \lbrace	\begin{ytableau} 3 & 2 \\1\end{ytableau}\, , 	\begin{ytableau} 4 & 5\end{ytableau} \right\rbrace \, - \,  \left \lbrace	\begin{ytableau} 1 & 2 \\3\end{ytableau}\, , 	\begin{ytableau} -4 & 5\end{ytableau} \right\rbrace \, + \,  \left \lbrace	\begin{ytableau} 3 & 2 \\1\end{ytableau}\, , 	\begin{ytableau} -4 & 5\end{ytableau} \right\rbrace \nonumber \\ 
& - & \ytableausetup{mathmode,baseline,aligntableaux=center,boxsize=1.3em} \left \lbrace	\begin{ytableau} 1 & 2 \\3\end{ytableau}\, , 	\begin{ytableau} 4 & -5\end{ytableau} \right\rbrace \, + \,  \left \lbrace	\begin{ytableau} 3 & 2 \\1\end{ytableau}\, , 	\begin{ytableau} 4 & -5\end{ytableau} \right\rbrace \, + \,  \left \lbrace	\begin{ytableau} 1 & 2 \\3\end{ytableau}\, , 	\begin{ytableau} -4 & -5\end{ytableau} \right\rbrace \, - \,  \left \lbrace	\begin{ytableau} 1 & 2 \\3\end{ytableau}\, , 	\begin{ytableau} -4 & -5\end{ytableau} \right\rbrace \nonumber\\ 
& = &1121^{+}1^{+} - 2111^{+}1^{+} - 1121^{-}1^{+} + 2111^{-}1^{+} \nonumber \\
& - &  1121^{+}1^{-} + 2111^{+}1^{-} + 1121^{-}1^{-} - 2111^{-}1^{-}. \label{chpt6:eqn:polytabexample}
\end{eqnarray}

The structure of words allow us to describe the lifting from modules of $\mathfrak{B}_{n}$ to $\mathfrak{B}_{n+1}$ by adding letters to words instead of adding boxes to tabloids. The lifting operators for the Specht modules of $\mathfrak{S}_{n}$ were described by appending letters to the end of a word. We perform a similar analysis for the Specht modules of $\mathfrak{B}_{n}$ but we now we have a choice of adding a signed or unsigned letter. Overall the correspondence with words helps to simplify the presentation of the lifting operators for $\mathfrak{B}_{n}$ and the results in Section \ref{chpt6:subsec:rtlifting}.

\section*{Random Walks on The Hyperoctahedral Group}
\label{chpt6:sec:randomwalks}

We may extend any shuffle on the symmetric group to a shuffle on $B_{n}$ by adding an additional action by negative transpositions to the end of the shuffle. In the following sections we study extended versions of the random transposition shuffle and one-sided transposition shuffle. We first extend the random transposition shuffle to $B_{n}$ and show how to recover its eigenvalues by lifting eigenvectors. Then we extend the one-sided transposition shuffle to $B_{n}$ and recover its eigenvalues using the same technique.

\section{The Random Transposition Shuffle on The Hyperoctahedral Group}
\label{chpt6:subsec:rt}

The random transposition shuffle for the hyperoctahedral group is described by the following procedure: pick two positions in $[n]$ uniformly at random and switch the cards in these positions, then flip a fair coin; if heads do nothing, if tails flip the moved cards over to their opposite sides. There are other ways we could choose to extend the random transposition shuffle to the hyperoctahedral group -- Schoolfield \cite{schoolfield2002random} studied an alternative description where each card's flip are done independently at the end of the shuffle, e.g., we could flip one card and not the other. Schoolfield went on to analyse the mixing time of this shuffle using Fourier transforms and found it to be tightly bounded at time $(n/2)\log n$, i.e., the same time as the random transposition shuffle on $S_{n}$. Instead of using Fourier transformations we analyse the random transposition shuffle using the technique of lifting eigenvectors.

\begin{defn}
	The \emph{random transposition shuffle for $B_{n}$}, denoted $\brt$,  is driven by the following probability distribution:
	\begin{eqnarray}
	\brt( \sigma ) = \begin{cases}
	1/2n & \textnormal{ if } \sigma = \gpid \\
	1/2n^{2} & \textnormal{ if } \sigma = \xi_{i} \textnormal{ for } i \in [n] \\
	1/n^{2} & \textnormal{ if } \sigma = (i \, j)  \textnormal{ for } i,j \in [n] \textnormal{ with } i<j\\
	1/n^{2} & \textnormal{ if } \sigma = \xi_{i} \xi_{j} (i \, j)  \textnormal{ for } i,j \in [n] \textnormal{ with } i<j \\
	0 & \textnormal{ otherwise }
	\end{cases}.
	\end{eqnarray}
\end{defn}

This shuffle is defined on the conjugacy classes of positive and negative transpositions.
Lifting eigenvectors for the random transposition shuffle on $B_{n}$ requires modification from the symmetric group case in order to account for the new Specht modules associated to bi-partitions. In particular we require two different lifting operators to take eigenvectors of $\mathfrak{B}_{n}$ to those of $\mathfrak{B}_{n+1}$. Throughout this section we will be working towards the following result.

\begin{thm}
	\label{chpt6:thm:rteig}
	The eigenvalues for the random transposition shuffle $\brt$ are indexed by bi-partitions
	$\bl \vdash n$. The eigenvalue corresponding to partition $\bl$ occurs with multiplicity $d_{\bl}^{2}$, and is given by
	\begin{eqnarray}
	\eig(\bl) = \frac{1}{2n^{2}}\left(2|\lambda^{1}| + 
	4 \, \D(\lambda^{1}) + 
	4 \, \D(\lambda^{2})\right). \label{chpt6:eqn:rteig}
	\end{eqnarray}

\end{thm}

The random transposition shuffle for $B_{n}$ is closely related to that for $S_{n}$, and in fact by considering a group homomorphism from $B_{n}$ to $S_{n}$ we recover an easy lower bound on the total variation distance between $\brt^{t}$ and $\pi_{n} := 1/ (2^{n} \, n!)$. 

\begin{lemma}
	\label{chpt6:lem:rtlower}
	The random transposition shuffle $\brt$ satisfies the following bound for any $c>0$:
	\begin{eqnarray}
	\liminf_{n\to\infty}\lVert \brt^{(n/2)\log n -c n} - \pi_{n} \rVert_{\tiny\textnormal{TV}} & \geq & \frac{1}{e} -e^{-e^{2c}}  .\label{chpt6:eqn:brtlowerbound}
	\end{eqnarray}
\end{lemma}
\begin{proof}
	The proof of this lemma follows from the fact that total variation distance can only decrease under projections, see \cite[Lemma 7.10]{Levin2017}.
	Form a surjective homomorphism $\psi:B_{n} \to S_{n}$ by ignoring the signs on permutations in $B_{n}$, i.e., for $\sigma\in B_{n}$ with decomposition $\sigma = \left( \sum_{i=1}^{n} \xi_{i}^{x_{i}}\right) \eta$ we have $\psi(\sigma) = \eta$. Consider the preimage of $\psi$ at a permutation in $\eta \in S_{n}$, we find that
	\[\sum_{\sigma \in \psi^{-1}(\eta)} \brt(\sigma) = \RT(\eta) .\]
	Extending this equality forward in time we have, 
	\[\sum_{\sigma \in \psi^{-1}(\eta)} \brt^{t}(\sigma) = \RT^{t}(\eta) \textnormal{ for all $t\geq 0$ and $\eta \in S_{n}$} .\]
	Thus, we may reduce the total variation distance of random transposition shuffle on $B_{n}$ to random transposition shuffle on $S_{n}$:
	\begin{eqnarray*}
		\lVert \brt^{t} - \pi_{n} \rVert_{\tiny\textnormal{TV}} & = & \frac{1}{2}\sum_{\sigma \in B_{n}}|\brt^{t}(\sigma) - \pi_{n}(\sigma)| \\
		& \geq &\frac{1}{2}\sum_{\eta \in S_{n}} \left\lvert\sum_{ \substack{\sigma \in B_{n} \\ \sigma \in \psi^{-1}(\eta)}}\brt^{t}(\sigma) - \pi_{n}(\sigma)\right\rvert
		= 
		\lVert \RT^{t} - \pi_{S_{n}}\rVert_{\tiny\textnormal{TV}}.
	\end{eqnarray*}
	Using the above inequality the lower bound is then established from Theorem \ref{chpt4:thm:RTcutoff}.
\end{proof}

\subsection{Upper Bound for the Random Transposition Shuffle}
\label{chpt6:subsec:upperboundrt}

To complete a cutoff argument for the random transposition shuffle on $B_{n}$ we need to show that its mixing time is at most $(n/2)\log n$.
The random transposition shuffle is reversible so we can use the eigenvalues given by Theorem \ref{chpt6:thm:rteig} to upper bound the total variation distance between $\brt^{t}$ and $\pi_{n}$ (see Theorem \ref{chpt2:thm:classicL2}). 
\begin{eqnarray}
4 \lVert \brt^{t} - \pi_{n} \rVert_{\tiny\textnormal{TV}}^{2} & \leq & \sum_{\substack{\bl \vdash n \\ \bl \neq ((n),(0))}} d_{\bl}^{2} \,  \eig(\bl)^{2t} \\
& = & \sum_{k=0}^{n} \sum_{\substack{\bl \vdash n \\ \bl \neq ((n),(0)) \\ |\lambda^{1} | = k}}  {n \choose k}^{2} \,  d_{\lambda^{1}}^{2} d_{\lambda^{2}}^{2} \left( \frac{k +2\D(\lambda^{1})+2\D(\lambda^{2})}{n^{2}}\right)^{2t}. \label{chpt6:eqn:rtbound1}
\end{eqnarray}

\paragraph{}

The first step in bounding the summation \eqref{chpt6:eqn:rtbound1} is to decrease the number of bi-partitions $\bl$ we have to consider.
For any bi-partition $\bl = (\lambda^{1}, \lambda^{2})$, we can create at most eight other bi-partitions by taking transposes of $\lambda^{1}$, $\lambda^{2}$, as well as swapping the positions of $\lambda^{1}$ and $\lambda^{2}$. For example two of the bi-partitions we can form following these rules are $(\lambda^{2 \, \prime}, \,\lambda^{1 \, \prime})$, and $(\lambda^{1 \, \prime}, \lambda^{2} )$. Furthermore, all the bi-partitions we can form have the same dimension because $d_{\lambda} = d_{\lambda^{\prime}}$ and ${n \choose |\lambda^{1}|} = {n \choose |\lambda^{2}|}$. If we consider the eigenvalues associated to these eight possible bi-partitions, the largest eigenvalue in magnitude must come from the bi-partition $\bl$ with $|\lambda^{1}| \geq |\lambda^{2}|$ and $\D(\lambda^{1}) \geq 0$, $\D(\lambda^{2}) \geq 0$. Therefore, the contribution of any bi-partition to \eqref{chpt6:eqn:rtbound1} may be upper bounded by a bi-partition $\bl$ with $|\lambda^{1}| \geq |\lambda^{2}|$, $\D(\lambda^{1}) \geq 0$, and $\D(\lambda^{2}) \geq 0$. 	
This allows us to reduce the number of bi-partitions we have to consider by roughly a factor of $8$. However, we need to make a special exception for the bi-partitions $((1^{n}),(0)), ((0),(n)), ((0),(1^{n}))$, because $((n),(0))$ does not belong to our sum. The eigenvalues associated to these bi-partitions are: $\eig((1^{n}),(0)) = -1 +(2/n)$ and $\eig((0),(n)) = -\eig ((0),(1^{n}))= 1 - (1/n)$. These eigenvalues all have dimension $1$ and so can easily be seen to be bounded at time $(n/2)\log n$,
\begin{eqnarray}
\limsup_{n\to \infty} \; \left( 1 - \frac{2}{n}\right)^{n \log n} + 2\left( 1 - \frac{1}{n}\right)^{n \log n} =0 . \nonumber
\end{eqnarray}
Altogether applying our reduction we may see that the summation \eqref{chpt6:eqn:rtbound1} is upper bounded by
\begin{eqnarray}
3 \left(1 - \frac{1}{n}\right)^{2t}\, +  8 \sum_{k\geq n/2}^{n}   \sum_{\substack{\bl \vdash n \\ \bl \neq ((n),(0)) \\ |\lambda^{1} | = k \\ \D(\lambda^{1}) \geq 0 \\ \D(\lambda^{2}) \geq 0 }} {n \choose k}^{2} \, d_{\lambda^{1}}^{2} d_{\lambda^{2}}^{2} \left( \frac{k +2\D(\lambda^{1})+2\D(\lambda^{2})}{n^{2}}\right)^{2t}. \label{chpt6:eqn:rtbound3}
\end{eqnarray}
The values $\D(\lambda^{1})$ and $\D(\lambda^{2})$ are closely related to the eigenvalues $\eig(\lambda^{1})$ and $\eig(\lambda^{2})$ for the random transposition shuffle on $S_{k}$ and $S_{n-k}$ respectively. In fact we can upper bound $\eig(\bl)$ by
\[\eig(\bl ) \leq  \frac{k +2\D(\lambda^{1})+ (n-k) + 2\D(\lambda^{2})}{n^{2}} =  \frac{k^{2}}{n^{2}}\eig(\lambda^{1}) + \frac{(n-k)^{2}}{n^{2}}\eig(\lambda^{2}) .\]
We have previously bounded the value of $\D(\lambda)$ for $\lambda$ a partition of $n$ using information about $\lambda_{1}$. The following bound was given in Lemma \ref{chpt4:lem:diagbounds} for $\lambda \vdash n$:
\begin{eqnarray}
2\D(\lambda) \leq \begin{cases}
(n-1)n - 2(n-\lambda_{1})(\lambda_{1}+1) & 
\textnormal{ 
	if } \lambda_{1} \geq \frac{n}{2}\\
(\lambda_{1}-1)n & \textnormal{ for all } \lambda  
\end{cases}.
\end{eqnarray} 
Using this result we may establish bounds on the value $\eig(\bl)$ for particular choices of $\lambda^{1}$ and $\lambda^{2}$. We would like to analyse equation \eqref{chpt6:eqn:rtbound3} by reducing it to previously studied bounds for the random transposition shuffle on the symmetric group.
This technique was used by Schoolfield to analyse the mixing time of the random transposition shuffle on $B_{n}$ with independent card flips \cite{schoolfield2002random}, and recently by Ghosh to analyse the mixing time of the flip-transpose top with random shuffle \cite{ghosh2019total}.
However, for the shuffle $\brt$ this technique leads to problems which stem from how our eigenvalues $\eig(\bl)$ are composed from the eigenvalues $\eig(\lambda^{1})$ and $\eig(\lambda^{2})$. 	

To see where this approach fails let us focus on an example. 
Suppose that $\lambda^{1} =k-i$ with $i \leq k/4$, and $\lambda^{2} = (n-k)-j$ with $j< (n-k)/4$, in other words both $\lambda^{1}$, $\lambda^{2}$ are large partitions of $k$ and $n-k$ respectively (see Section \ref{chpt5:sec:upperbound}). In this case the eigenvalue $\eig(\bl)$ may be bounded as follows:
\begin{eqnarray}
\frac{k +2\D(\lambda^{1})+2\D(\lambda^{2})}{n^{2}} \leq 1 -\frac{2k(n-k)}{n^{2}} - \frac{2i(k-i+1)}{n^{2}} - \frac{2j(n-k-j+1)}{n^{2}}.
\end{eqnarray} 
In the above bound there are three different negative terms, we want to use the first to bound the binomial coefficient ${n \choose k }^{2}$, the second to bound  $d_{\lambda^{1}}^{2}$, and the third to bound $d_{\lambda^{2}}^{2}$. Making use of Lemma \ref{chpt4:lem:youngdimbound} and setting $t = (n/2)\log n +cn$ we can reduce the contributions of large partitions in \eqref{chpt6:eqn:rtbound3} to:
\begin{eqnarray}
& & \sum_{k\geq n/2}^{n}  \,\sum_{i=0}^{k/4} \sum_{j=0}^{(n-k)/4} \sum_{\substack{\bl \vdash n \\ \bl \neq ((n),(0)) \\ |\lambda^{1} | = k \\ \lambda_{1}^{1} = k-i \\ \lambda_{1}^{2} = n-k -j }} {n \choose k}^{2}  \, d_{\lambda^{1}}^{2} d_{\lambda^{2}}^{2} \left(1 -\frac{2k(n-k)}{n^{2}} - \frac{2i(k-i+1)}{n^{2}} - \frac{2j(n-k-j+1)}{n^{2}}\right)^{2t} \nonumber\\
& = & \sum_{k\geq n/2}^{n}  \, \sum_{i=0}^{k/4} \sum_{j=0}^{(n-k)/4} {n \choose k}^{2} \, \left(1 -\frac{2k(n-k)}{n^{2}} - \frac{2i(k-i+1)}{n^{2}} - \frac{2j(n-k-j+1)}{n^{2}}\right)^{2t} \sum_{\substack{\lambda^{1} \vdash k \\ \lambda_{1}^{1} =k-i}}  d_{\lambda^{1}}^{2} \sum_{\substack{\lambda^{2} \vdash n-k \\ \lambda_{1}^{2} =n-k-j} } d_{\lambda^{2}}^{2} \nonumber\\
& \leq & \sum_{k\geq n/2}^{n}  \,\sum_{i=0}^{k/4} \, \sum_{j=0}^{(n-k)/4}  {n \choose k}^{2} \, {k \choose i}^{2} i! \, {n-k \choose j}^{2} j! \, \left(1 -\frac{2k(n-k)}{n^{2}} - \frac{2i(k-i+1)}{n^{2}} - \frac{2j(n-k-j+1)}{n^{2}}\right)^{2t} \nonumber\\ 
& \leq & \sum_{k\geq n/2}^{n} \, \sum_{i=0}^{k/4} \, \sum_{j=0}^{(n-k)/4} {n \choose k}^{2} \, {k \choose i}^{2} i! \, {n-k \choose j}^{2} j! \, e^{- (n\log n  +2cn )\left( \frac{2k(n-k)}{n^{2}} + \frac{2i(k-i+1)}{n^{2}} + \frac{2j(n-k-j+1)}{n^{2}} \right) }\nonumber\\
& \leq & e^{-2c} \,\sum_{k\geq n/2}^{n} \, \sum_{i=0}^{k/4} \, \sum_{j=0}^{(n-k)/4} {n \choose k}^{2} \, {k \choose i}^{2} i! \, {n-k \choose j}^{2} j! \, n^{-\frac{2k(n-k)}{n} - \frac{2i(k-i+1)}{n} - \frac{2j(n-k-j+1)}{n} } \label{chpt6:eqn:rtbound6}
\end{eqnarray}
From here we split equation \eqref{chpt6:eqn:rtbound6} into three separate summations, one for the binomial coefficients and one each for the dimensions of $\lambda^{1}$ and $\lambda^{2}$,
\begin{eqnarray}
e^{-2c} \sum_{k\geq n/2}^{n}  {n \choose k}^{2}  n^{-\frac{2k(n-k)}{n}}\, \sum_{i=0}^{k/4} {k \choose i}^{2} i! n^{- \frac{2i(k-i+1)}{n}}\, \sum_{j=0}^{(n-k)/4}  \, {n-k \choose j}^{2} j! \, n^{ - \frac{2j(n-k-j+1)}{n}}. \nonumber
\end{eqnarray}
We can now try to show that equation \eqref{chpt6:eqn:rtbound6} is bounded in $n$ by analysing each of these summations separately. The sums corresponding to partitions $\lambda^{1}$ and $\lambda^{2}$ given by 
\begin{eqnarray}
\sum_{i=0}^{k/4} {k \choose i}^{2} i! n^{- \frac{2i(k-i+1)}{n}} \hspace{0.5cm} \textnormal{ and }  \hspace{0.5cm} \sum_{j=0}^{(n-k)/4}  \, {n-k \choose j}^{2} j! \, n^{ - \frac{2j(n-k-j+1)}{n}} \label{chpt6:eqn:rtbound8}
\end{eqnarray}
are closely related to the bounds analysed by Diaconis \cite[Chapter 3D]{Diaconis1988} for the random transposition shuffle on $S_{k}$ and $S_{n-k}$, these are respectively given by,
\begin{eqnarray}
\sum_{i=0}^{k/4} {k \choose i}^{2} i! k^{- \frac{2i(k-i+1)}{k}} \hspace{0.5cm} \textnormal{ and }  \hspace{0.5cm} \sum_{j=0}^{(n-k)/4}  \, {n-k \choose j}^{2} j! \, (n-k)^{ - \frac{2j(n-k-j+1)}{(n-k)}}. \label{chpt6:eqn:rtbound7}
\end{eqnarray}
Diaconis showed the summations in \eqref{chpt6:eqn:rtbound7} are bounded by a universal constant as $k \to \infty$ and $n-k \to \infty$ respectively. The difference between equations \eqref{chpt6:eqn:rtbound8} and \eqref{chpt6:eqn:rtbound7} comes from the eigenvalue $\eig(\bl)$ being of a different dimension compared with those of $\eig(\lambda^{1})$, $\eig(\lambda^{2})$. 
These differences make the sums in \eqref{chpt6:eqn:rtbound8} challenging to bound. For example, if $n-k = O(n)$ we find that the summation given by $\lambda^{2}$ is unbounded as $n\to \infty$. 
A similar issue is found when we consider $\lambda^{1}$ and $\lambda^{2}$ to be any other combination of large or small partitions. Hence, trying to reduce the analysis of \eqref{chpt6:eqn:rtbound3} to separate bounds related to the  symmetric group case would appear to fail.

There are several ways we could try to circumvent this issue. The first way is to consider a more holistic approach to bounding equation \eqref{chpt6:eqn:rtbound6}. The binomial term given by 
${n \choose k}^{2}  n^{-\frac{2k(n-k)}{n}}$
is not only bounded in $n$ but is decaying as $k \to n/2$, whereas the terms ${k \choose i}^{2} i! n^{- \frac{2i(k-i+1)}{n}}$ and ${n-k \choose j}^{2} j! \, n^{ - \frac{2j(n-k-j+1)}{n}}$ corresponding to $\lambda^{1}$ and $\lambda^{2}$ respectively are growing as $k \to n/2$. By considering all these terms together  instead of as three separate sums we may be able to use the decay of the binomial term to counteract the growth of the  $\lambda^{1}$ and $\lambda^{2}$ terms.
Alongside this it may be useful to consider a more fine-grain analysis of the eigenvalues, either by further restricting the choices of $\lambda^{1}$ and $\lambda^{2}$, or by splitting the analysis into cases depending on whether $n-k = O(1)$ or $n-k = O(n)$.

To date we have been unable to successfully use the eigenvalues to establish an upper bound on the mixing time of the random transposition shuffle on $B_{n}$, but we conjecture that an upper bound of time $(n/2)\log n$ may be found. Thus, we expect the random transposition shuffle on $B_{n}$ to exhibit a total variation cutoff at the same time as the random transposition shuffle on $S_{n}$. 
\begin{conj}
	\label{chpt6:conj:upper}
	The random transposition shuffle $\brt$ satisfies the following bound:
	\begin{eqnarray}
	\lim_{c\to\infty} \liminf_{n\to\infty}\lVert \brt^{(n/2)\log n +c n} - \pi_{n} \rVert_{\tiny\textnormal{TV}} & = & 0  \label{chpt6:eqn:brtupperbound} 
	\end{eqnarray}
	Thus, the random transposition shuffle on the hyperoctahedral group exhibits a cutoff in total variation distance at time $(n/2)\log n$.
\end{conj}

\subsection{Lifting Eigenvectors for the Random Transposition Shuffle}
\label{chpt6:subsec:rtlifting}
We lift the eigenvectors of $\brt$ by reducing the lifting of bi-partitions $\bl$ to the lifting of the individual partitions $\lambda^{1}$ and $\lambda^{2}$. This allows us to extend the technique of lifting eigenvectors that we developed for the random transposition shuffle on $S_{n}$ to the random transposition shuffle on $B_{n}$. The first step as before is to turn the probability $\brt$ into an element of the group algebra $\mathfrak{B}_{n}$.

\begin{defn}
	The random transposition shuffle on $B_{n}$ may be viewed as the following element of the group algebra $\mathfrak{B}_{n}$.
	\begin{eqnarray}
	\abrt = n \cdot \gpid + \sum_{1\leq i \leq n} \xi_{i} + 2 
	\sum_{1\leq i < j \leq n} (i \, j)  + 2\sum_{1\leq i <j\leq n} \xi_{i}\xi_{j}(i \, j ).
	\end{eqnarray}
	Note that here we have scaled by $2n^{2}$.
\end{defn}

The eigenvectors and eigenvalues of $\brt$ are in one-to-on correspondence with those given by $\abrt$ acting on the regular module $\mathfrak{B}_{n}$. Using the canonical decomposition of the regular module we now focus our attention on finding the eigenvectors of $\abrt$ belonging to the Specht modules $S^{\bl}$. If we compare the elements $\abrtp$ and $\abrt$ we can see they have a close relationship to one another,
\begin{eqnarray}
\abrtp - \abrt = \gpid + \xi_{n+1} + 2\sum_{1\leq i \leq n} (i \hspace{0.2cm} n+1) + 2\sum_{1\leq i \leq n} \xi_{i}\xi_{n+1}(i \hspace{0.2cm} n+1) .\label{chpt6:eqn:nn1diff}
\end{eqnarray}
The equation \eqref{chpt6:eqn:nn1diff} only depends on transpositions involving $n+1$, this is similar to equation \eqref{chpt6:eqn:nn1diff}, and gives us an insight into why lifting eigenvectors works for the random transposition shuffle on $B_{n}$. To turn eigenvectors for the module $S^{\bl}$ into those for $S^{\bl + e_{i}^{k}}$ we introduce new adding and switching operators.

\begin{defn}
	Let $w$ be a word in $M^{\bl}$, we may add a letter to 
	$w$ to form a word of size $n+1$.	Define the linear \emph{adding 
		operators}, denoted $\ad_{a}$, for letter $a\in [n]$ as follows,
	\[\ad_{a}(w) = w  a, \hspace{0.8cm} \ad_{a}^{+}( w ) = w
	a^{+}, \hspace{0.8cm} \ad_{a}^{-}( w ) = w
	a^{-}.\]
	Note that the first operator takes us from $M^{\bl}$ to 
	space $M^{\bl +e_{a}^{1} }$ whereas the second and third operators takes us to $M^{\bl 
		+e_{a}^{2} }$. Using the above basic adding operators define two new adding operators:
	\[\sh_{a}^{1} = \Phi_{a},  \hspace{0.8cm} \sh_{a}^{2} = \Phi_{a}^{+} - 
	\Phi_{a}^{-}.\]
	These new adding operators are critical to our analysis of the random transposition shuffle.
\end{defn}

\begin{defn}
	\label{chpt6:def:swtich}
	Let $w$ be a word in $M^{\bl}$ and $a,b \in [\overline{n}]$. We define the linear 
	\emph{switching operators} as follows:
	\begin{eqnarray}
	\Theta_{b,a}(w) & = &\sum_{\substack{1 \leq i \leq n \\ w_{i} =b}} w_{1} \ldots  w_{i-1}  \,a \, w_{i+1} \ldots  
	w_{n} \label{eqn:switch1}
	\end{eqnarray}
	Using this single operator we define three other switching operators for $a,b \in [n]$:
	\begin{eqnarray}
	\Theta_{b^{\pm},a}(w) & = & \Theta_{b^{-},a} + \Theta_{b^{+},a}	\label{eqn:switch2}\\
	\Theta_{b^{\pm},a^{\pm}}^{+} (w)& = & \Theta_{b^{+},a^{+}} + \Theta_{b^{-},a^{-}}	\label{eqn:switch3} \\
	\Theta_{b^{\pm},a^{\pm}}^{-} (w)& = & \Theta_{b^{+},a^{-}} + \Theta_{b^{-},a^{+}}	\label{eqn:switch4}
	\end{eqnarray}
	The operator $\Theta_{b^{\pm},a}$ takes all signed occurrences of $b$ and replaces them by the unsigned $a$. The operators $\Theta_{b^{ \pm},a^{\pm}}^{\pm}$ swap signed occurrences of $b$ for signed occurrences of $a$, with $\Theta^{+}$ fixing the signs, and $\Theta^{-}$ swapping them. Define two new operators $\Theta^{1}_{b,a},\Theta^{2}_{b,a}$  with domain $M^{\bl}$ for the unsigned symbols $a,b \in [n]$ as follows:
	\begin{eqnarray}
	\Theta^{1}_{b,a} & = & \Theta_{b,a} \\ 
	\Theta^{2}_{b,a} & = & \Theta_{b^{\pm},a^{\pm}}^{+} - \Theta_{b^{\pm},a^{\pm}}^{-}.
	\end{eqnarray}
	The operator $\Theta^{1}_{b,a}$ acts only on unsigned letters, and $\Theta^{2}_{b,a}$ acts only on signed letters.
	The image of $\Theta_{b,a}^{1}$ is $\bl - e_{b}^{1} + e_{a}^{1}$, and the image of $\Theta^{2}_{b,a}$ is $\bl -e_{b}^{2} + e_{a}^{2}$. 
\end{defn}

\begin{lemma}
	The switching operators $\Theta_{b^{\pm},a}, \Theta_{b,a}^{1}$ and $\Theta_{b,a}^{2}$ are $\mathfrak{B}_{n}$-module morphisms.
\end{lemma}
\begin{proof}
	The operators $\Theta_{b,a}^{1}$ and $\Theta_{b^{\pm},a}$, swap a letter $b$ (signed or unsigned) for an unsigned letter $a$, this is unaffected by movement of the letter and changes in sign (as it always ends up unsigned), therefore it commutes with elements $B_{n}$. The operator $\Theta^{+}_{b^\pm,a^\pm}$ swaps a signed letter $b^{\pm}$ for $a^{\pm}$ while preserving the current sign of $b^{\pm}$ so it does not matter if the letter is moved or the sign of the letter is changed before or after applying the operator. Similar logic holds for the operator $\Theta^{-}_{b^\pm,a^\pm}$.  Hence, the signed lifting operator $\Theta^{2}_{b,a}$ is a  $\mathfrak{B}_{n}$-module morphism.
\end{proof}

Notice that we have not defined a switching operator which takes unsigned letters and turns them into a signed letters, the reason behind this is there is no way to give a letter a sign and respect the action of $B_{n}$. The adding operators are not module morphisms, however, they can be seen to commute with the switching operators with an extra adjustment term.

\begin{lemma}
	\label{chpt6:lem:addswitchrelation}
	The adding and switching operators satisfy the following equalities:
	\begin{eqnarray}
	\Phi_{b}^{1} \circ \Theta_{b,a}^{1} & = &\Theta_{b,a}^{1} \circ \Phi_{b}^{1} - \Phi_{a}^{1} \label{chpt6:eqn:addshufflerelation1}\\
	\Phi_{b}^{2} \circ \Theta_{a, b}^{2} & = & \Theta_{a,b}^{2} \circ \Phi_{b}^{2} - 2\cdot \Phi_{a}^{2} \label{chpt6:eqn:addshufflerelation2}.
	\end{eqnarray}
\end{lemma}

We can already see that lifting eigenvectors for $B_{n}$ requires more careful tools than for $S_{n}$; we have gone from requiring one adding and switching operator to two of them each. We now establish a relationship between the modules of $\mathfrak{B}_{n}$ and $\mathfrak{B}_{n+1}$ using the newly defined operators. This is similar to Theorem \ref{chpt5:thm:master} however we now consider two cases, based on if we lift the partition $\lambda^{1}$ or $\lambda^{2}$.

\begin{thm}
	\label{chpt6:thm:master1}
	Let $n\in \mathbb{N}$ and $\bl \vdash n$. For words in $M^{\bl}$ we have the following equalities: 
	\begin{eqnarray}
	\abrtp \circ \sh_{a}^{1} - \sh_{a}^{1} \circ 
	\abrt & = & 2\sh_{a}^{1} +  
	2\sum_{1 \leq b \leq n} 
	2 \cdot \sh_{b}^{1}\circ\Theta_{b,a}^{1} + \left(\sh_{b}^{+} + \sh_{b}^{-}\right) 
	\circ\Theta_{b^{\pm},a}.\label{chpt6:eqn:master1} \\
	\abrtp \circ \sh_{a}^{2} - \sh_{a}^{2} \circ 
	\abrt & = & 2 \sum_{1 \leq b \leq n} 
	\sh_{b}^{2}\circ\Theta_{b,a}^{2}. \label{chpt6:eqn:master2}
	\end{eqnarray}

\end{thm}

\begin{proof}
	Take $w \in M^{\bl}$ a generic element to verify these 
	equations with. 
	In both equations all non-identity elements of $\abrtp$ which do not 
	involve $n+1$ commute with the adding operators $\sh^{1}_{a}, \sh^{2}_{a}$, and therefore cancel with $\abrt$. 
	This just leaves us the elements present in equation \eqref{chpt6:eqn:nn1diff} to work with. 
	First we 
	prove 
	equation \eqref{chpt6:eqn:master1}, using the above observation we reduce it to,
	\begin{eqnarray*}
		\abrtp \circ \sh_{a}^{1} - \sh_{a}^{1} \circ 
		\abrt & = & \left(\gpid + \xi_{n+1} + 2\sum_{1\leq i \leq n} (i \hspace{0.2cm} n+1) + 2\sum_{1\leq i \leq n} \xi_{i}\xi_{n+1}(i \hspace{0.2cm} n+1) \right) \sh_{a}^{1}(w).
	\end{eqnarray*}
	Clearly $(\gpid + \xi_{n+1})\sh_{a}^{1}(w) = 2 \sh_{a}^{1}(w)$. 
	The letter we add to $w$ is unsigned, define the word $v_{i} = w_{1} \ldots w_{i-1} \, a \, w_{i+1} \ldots w_{n}$. Then the positive transpositions act on $w \, a$ in the following way:
	\[\left( (i \hspace{0.25cm} n+1) + \xi_{i}\xi_{n+1}(i \hspace{0.25cm} 
	n+1)\right)\sh_{a}^{1}(w)= \begin{cases}
	v_{i}b + v_{i}b = 2v_{i}b & \text{ if } w_{i} = b\\
	v_{i}b^{+} + v_{i}b^{-} & 
	\text{ if } w_{i} = b^{+}\\
	v_{i}b^{-} + v_{i}b^{+}& 
	\text{ if } w_{i} = b^{-}
	\end{cases}
	.\]
	Thus, summing over all $i$ or equivalently all signed and unsigned letters 
	we recover  \eqref{chpt6:eqn:master1}. 	Now we prove equation 
	\eqref{chpt6:eqn:master2}, again it may be reduced it to
	\begin{eqnarray*}
		\abrtp \circ \sh_{a}^{2} - \sh_{a}^{2} \circ 
		\abrt & = & \left(\gpid + \xi_{n+1} + 2\sum_{1\leq i \leq n} (i \hspace{0.2cm} n+1) + 2\sum_{1\leq i \leq n} \xi_{i}\xi_{n+1}(i \hspace{0.2cm} n+1) \right) \sh_{a}^{2}(w) . 
	\end{eqnarray*}
	Looking at the negative transposition we find $(\gpid + \xi_{n+1})\sh_{a}^{2}(w) = (\gpid + \xi_{n+1}) (wa^{+}-wa^{-}) = 0$. For the positive transpositions we must be careful because the added letter $a$ is now signed. Define  $v_{i}^{+} = w_{1}\dots w_{i-1} \, a^{+} \, w_{i+1}\dots w_{n}$, and $v_{i}^{-} = 
	w_{1}\dots w_{i-1} \, a^{-} \, w_{i+1}\dots w_{n}$. The positive transpositions act on $(wa^{+} - wa^{-})$ in the following way:
	\[\left( (i \hspace{0.2cm} n+1) + \xi_{i}\xi_{n+1}(i \hspace{0.2cm} 
	n+1)\right) \sh_{a}^{2}(w) = \begin{cases}
	(v_{i}^{+}b - v_{i}^{-}b) 
	+ (v_{i}^{-}b - v_{i}^{+}b) = 0 & \text{ if } w_{i} = b\\
	(v_{i}^{+}b^{+} - v_{i}^{-}b^{+}) 
	+ (v_{i}^{-}b^{-} - v_{i}^{+}b^{-}) = \Phi_{b}^{2}\left(v_{i}^{+} - v_{i}^{-}\right)  & 
	\text{ if } w_{i} = b^{+}\\
	(v_{i}^{+}b^{-} - v_{i}^{-}b^{-}) 
	+ (v_{i}^{-}b^{+} - v_{i}^{+}b^{+}) = \Phi_{b}^{2}\left(v_{i}^{-} - v_{i}^{+}\right) & 
	\text{ if } w_{i} = b^{-}
	\end{cases}
	.\]
	Thus summing over all $1\leq i \leq n$ or equivalently all signed letters $b^{\pm}$ in $w$ completes equation \eqref{chpt6:eqn:master2}. 
\end{proof}

Once refined, equations \eqref{chpt6:eqn:master1}, \eqref{chpt6:eqn:master2} will enable us to turn eigenvectors of $\abrt$ into those for $\abrtp$. The next step is to restrict their domains to a Specht module $S^{\bl}$.

\begin{lemma}
	\label{chpt6:lem:thetadom}
	Let $\bl \vdash n$ be such that $\bm + e_{b}^{k} = \bl + e_{a}^{l}$ for some $k,l \in \{1,2\}$ and $a,b \in [n]$.  
	If $k>l$ then $\Theta_{b^{\pm},a}(S^{\bl}) =0$. 
	If $k=l$ then  $\Theta_{b,a}^{k}$ is non-zero on $S^{\bl}$ if and only if $\lambda^{k}$ dominates the non-increasing rearrangement of $\mu^{k}$. In particular, if $a<b$, then $\Theta_{b,a}^{k}(S^{\bl}) = 0$.

\end{lemma}
\begin{proof}
	We know that for both cases the respective switching operators $\Theta \, (S^{\bl})$ must belong in $M^{\bm}$, and so for the map to be non-zero we must have $\bl$ dominating the non-increasing rearrangement of $\bm$  by Lemma \ref{chpt6:lem:youngsrule}. If $k>l$ then $\bm \trianglerighteq \bl$ therefore $\Theta_{b^{\pm},a}(S^{\bl}) = 0$. If $k=l$ the assertion holds by the same reasoning as Lemma \ref{chpt5:lem:restrict}.
\end{proof}

For an illustration of Lemma \ref{chpt6:lem:thetadom} take the polytabloid $e_{\overline{T}}$ belonging to module $S^{((2,1),(2))}$ from Example \ref{chpt6:ex:specht} (also equation \eqref{chpt6:eqn:polytabexample}), the switching operators $\Theta_{1^{\pm},1}$, $\Theta_{2,1}^{1}$, and $\Theta_{1,3}^{1}$ applied to $e_{\overline{T}}$ gives us:
\begin{eqnarray}
\Theta_{1^{\pm},1} (e_{\overline{T}})& = & (11211^{+} + 1121^{+}1) - (21111^{+} + 2111^{+}1) - (11211^{+} + 1121^{-}1) + (21111^{+} + 2111^{-}1) \nonumber \\
& - &  (11211^{-} + 1121^{+}1) +  (21111^{-} + 2111^{+}1)  +  (11211^{-} + 1121^{-}1) -  (21111^{-} + 2111^{-}1) = 0 \nonumber \\ 
\Theta_{2,1}^{1} (e_{\overline{T}}) & = & 1111^{+}1^{+} - 1111^{+}1^{+} - 1111^{-}1^{+} + 1111^{-}1^{+}
\nonumber\\
& - & 1111^{+}1^{-} + 1111^{+}1^{-} + 1111^{-}1^{-} - 1111^{-}1^{-} = 0 \nonumber\\
\Theta_{1,3}^{1} (e_{\overline{T}}) & = & (1321^{+}1^{+} + 3121^{+}1^{+}) - (2131^{+}1^{+} + 2311^{+}1^{+})  \nonumber \\ 
& - & ( 1321^{-}1^{+} + 3121^{-}1^{+})+ (2131^{-}1^{+} + 2311^{-}1^{+})  \nonumber \\
& - &  ( 1321^{+}1^{-} + 3121^{+}1^{-}) + (2131^{+}1^{-} + 2311^{+}1^{-})  \nonumber \\
& + & ( 1321^{-}1^{-} + 3121^{-}1^{-}) - (2131^{-}1^{-} + 2311^{-}1^{-}) \neq 0 \nonumber.
\end{eqnarray}
Applying Lemma \ref{chpt6:lem:thetadom} we now restrict equations \eqref{chpt6:eqn:master1} and \eqref{chpt6:eqn:master2} to the domain of a Specht module $S^{\bl}$.
		
\begin{corollary}
	\label{chpt6:cor:restrict}
	Let $\bl \vdash n$, restricting equations 
	\eqref{chpt6:eqn:master1} and \eqref{chpt6:eqn:master2} to have domain 
	$S^{\bl}$ we find the following: 
	\begin{eqnarray}
	\abrtp \circ \sh_{a}^{1} - \sh_{a}^{1} \circ 
	\abrt |_{S^{\bl}} & = &2\sh_{a}^{1}|_{S^{\bl}}   +  
	4  \sum_{1 \leq b \leq a} 
	\sh_{b}^{1}\circ\Theta_{b,a}^{1} |_{S^{\bl}}  
	\label{chpt6:eqn:masterres1}\\
	\abrtp \circ \sh_{a}^{2} - \sh_{a}^{2} \circ 
	\abrt |_{S^{\bl}}  & = & 2 \sum_{1 \leq b \leq a}
	\sh_{b}^{2}\circ  \Theta_{b,a}^{2} |_{S^{\bl}}  .
	\label{chpt6:eqn:masterres2}
	\end{eqnarray}

\end{corollary}

Notice that equations \eqref{chpt6:eqn:masterres1} and \eqref{chpt6:eqn:masterres2} now only depend on unsigned or signed letters, respectively,  being added or moved. 	
This allows us to stop focusing on the bi-partition $\bl$ and instead focus on the two single partitions $\lambda^{1}, \lambda^{2}$. Equations \eqref{chpt6:eqn:masterres1}, \eqref{chpt6:eqn:masterres2} are used to lift eigenvectors by adding a box to $\lambda^{1}$ or $\lambda^{2}$ respectively. Before we project onto the Specht modules $S^{\bl+e_{a}^{1}}$ or $S^{\bl+e_{a}^{2}}$ we must be sure they are contained as a submodule of the image. The next lemma is an analogue of Lemma \ref{chpt5:lem:lives} for the hyperoctahedral group.

\begin{lemma}
	\label{chpt6:lem:phimapscorrect}

	The subspace $\sh_{a}^{k}(S^{\bl})$ is contained in a 
	$\mathfrak{B}_{n+1}$ submodule of $M^{\bl + e_{a}^{k}}$, 
	isomorphic to $\bigoplus S^{\bm}$ where $\bm$ ranges 
	over the partitions obtained from $\bl$ by adding a box $ e_{i}^{k}$ with $i \leq a$. 	
\end{lemma}

\begin{proof}
	We begin by proving our hypothesis for $k=1$, that is adding a box to $\lambda^{1}$.
	Let $w$ be a word in $M^{\bl}$, let $a \in [n+1]$, 
	if letter $b$ does not appear in $w$ then
	\[\sh_{a}^{1}(w) = \Theta_{b,a}^{1} (\sh_{b}^{1}(w)) = \Theta_{b,a}^{1}(w \,b).\]
	Let $b = l(\lambda^{1})+1$, and consider the $\mathfrak{B}_{n+1}$-module 
	\[N_{1} = \langle x \, b \, : \, x \in S^{\bl} 
	\rangle 
	\cong 
	\textnormal{Ind}_{\mathfrak{B}_{n} \times 
		\mathfrak{B}_{1}}^{\mathfrak{B}_{n+1}} ( S^{\bl} \otimes 
	S^{((1),\emptyset)}) \cong \bigoplus_{\substack{\bm \vdash n+1 \\ \lambda^{1} \subset \mu^{1} \\\lambda^{2} = \mu^{2}  }} S^{\bm}\]
	where the last isomorphism follows from the branching rules for $B_{n}$ (see \cite[Theorem III.2]{geissinger1978representations}).	
	Using the observation at the start of the proof, we obtain
	\[\Phi_{a}^{1}(S^{\bl}) = 
	\Theta_{b,a}^{1}(\Phi_{b}^{1}(S^{\bl})) 
	\subseteq \Theta_{b,a}^{1}(\langle\Phi_{b}^{1}(S^{\bl})\rangle ) 
	= 
	\Theta_{b,a}^{1}(N_{1}) \cong \bigoplus_{\substack{\bm \vdash n+1 \\ \lambda^{1} \subset \mu^{1} \\\lambda^{2} = \mu^{2}  }} \Theta_{b,a}^{1} (S^{\bm}).\]
	Now note that $\Theta_{b,a}^{1}$ sends any word with evaluation 
	$\bl+e_b^{1}$ to a word with evaluation $\bl+e_a^{1}$, and hence
	$\Theta_{b,a}^{1}(N_{1}) \subseteq\Theta_{b,a}^{1}(M^{\bl +e_{b}^{1}}) \subseteq M^{\bl +e_{a}^{1}}$. 
	It follows that all nonzero summands $S^{\bm}$ appearing on the right 
	hand side occur for $\bm\vdash n+1$ dominating the non-decreasing rearrangement of $\bl+e_a^{1}$, and then
	by Lemma \ref{chpt6:lem:thetadom} we can conclude that $\bm$ is obtained from 
	$\bl$ by adding a box $e_{i}^{1}$ with $i \leq a$, as required.

	Suppose $k=2$, let $w$ be a word in $\bl$, let 
	$a \in 
	[n+1]$, 
	if the letter $b^{\pm}$ does not appear in $w$ then
	\[\Phi_{a}^{2}(w) = \frac{1}{2}\Theta_{b,a}^{2} \circ \Phi_{b}^{2}(w) 
	=\frac{1}{2}\Theta_{b,a}^{2}\left( w \, b^{+} - w \, b^{-}\right) .\]
	Let $b = l(\lambda^{2})+1$, consider the set 
	\[N_{2} = \langle x \, b^{+} - x \, b^{-} \, : \, x \in 
	S^{\bl} 
	\rangle \cong 
	\textnormal{Ind}_{\mathfrak{B}_{n} \times 
		\mathfrak{B}_{1}}^{\mathfrak{B}_{n+1}} ( S^{\bl} \otimes 
	S^{(\emptyset,(1))}) \cong \bigoplus_{\substack{\bm \vdash n+1 \\ \lambda^{1} = \mu^{1} \\\lambda^{2} \subseteq \mu^{2}   }} S^{\bm}\]	
	where the last isomorphism follows from the branching rules for $B_{n}$ (see \cite[Theorem III.2]{geissinger1978representations})	
	Using the observation at the start of the proof,  we obtain as modules over 
	$\mathfrak{B}_{n+1}$: 
	\[\Phi_{a}^{2}(S^{\bl}) =
	\Theta_{b,a}^{2}(\Phi_{b}^{2}(S^{\bl})) 
	\subseteq 
	\Theta_{b,a}^{2}(\langle\Phi_{b}^{2}(S^{\bl})\rangle) =
	\Theta_{b,a}^{2}(N_{2}) \cong \bigoplus_{\substack{ \bm \vdash n+1 \\ \lambda^{1} = \mu^{1} \\ \lambda^{2} \subseteq \mu^{2}   }} \Theta_{b,a}^{2} (S^{\bm}).\]
	By the same reasoning as before the final direct sum must be over all $\bm$ with $\mu^{2}$  formed by adding a box $e_{i}^{2}$ to $\bl$ with $i\leq a$. 
\end{proof}

\begin{defn}
	\label{chpt6:def:proj}
	Let $\iso^{\bm}$ be the isotypic projector onto the simple 
	module $S^{\bm}$. Define the \emph{lifting operator} for 
	$k \in \{1,2\}$ as:
	\begin{align}
	\pro_{a}^{\bl,\bm,k} = (\pi^{\bm} \circ \sh_{a}^{k}) 
	\big|_{S^{\bl}}:S^{\bl} 
	\rightarrow  M^{\bl + e_{a}^{k}} .
	\end{align}
	Lemma \ref{chpt6:lem:phimapscorrect} tells us that the image of $\sh_{a}^{k}$ may contain at most one copy of any Specht module $S^{\bm}$. Therefore, if there exists a copy in the image the element  $\pro_{a}^{\bl,\bm,k}(w)$ belongs exactly to the Specht module $S^{\bm}$ within $M^{\bl + e_{a}^{k}}$.
\end{defn}	

\begin{corollary}
	\label{chpt6:cor:prononzero}
	For any $\bl \vdash n$, $a \in \{1,2,\ldots, l(\lambda^{k})+1\}$ and $1 \leq i\leq a$, there exists $v \in 
	S^{\bl}$ such that
	\[\pro_{a}^{\bl,\bl+e_{i}^{k},k}(v) \neq 	0.\]
\end{corollary}

\begin{proof}
	Suppose 
	$\pro_{a}^{\bl,\bl+e_{i}^{k},k}(v) =0$, 
	then the image $\Phi_{a}^{k}(S^{\bl})$ lies in the kernel of 
	the projection $\pi^{\bl+ e_{i}^{k}}: M^{\bl 
		+e_{i}^{k}} \to S^{\bl + e_{i}^{k}}$ which is an 
	$\mathfrak{B}_{n+1}$-submodule with no component equal to 
	$S^{\bl+e_{i}^{k}}$. Hence, the submodule generated by $\sh_{a}^{k}(S^{\bl})$ has no competent equal to $S^{\bl+e_{i}^{k}}$. But we previously observed that (with 
	notation from Lemma \ref{chpt6:lem:phimapscorrect} and the two cases presented 
	together)
	\begin{eqnarray}
	\langle \Phi_{a}^{k}(S^{\bl}) \rangle =\langle 
	\Theta_{b,a}^{k} (\Phi_{b}^{k}(S^{\bl})) \rangle 
	= \Theta_{b,a}^{k} (\langle\Phi_{b}^{k}(S^{\bl})\rangle ) 
	= \Theta_{b,a}^{k}(N_{k}) \cong \bigoplus_{1 \leq i \leq a} S^{\bl+ e_{i}^{k}}
	\end{eqnarray}	
	Since the corresponding right hand side contains $S^{\bl+e_{i}^{k}}$ as a summand, we have a contradiction.
\end{proof}

\begin{lemma}
	\label{chpt6:lem:pro2inject}
	The linear operators 
	$\pro_{a}^{\bl,\bl+e^{k}_{i},k}$ for $1 \leq i \leq a$ are 
	$\mathfrak{B}_{n}$-module morphisms with trivial kernels. Therefore, these lifting operators are injective.
\end{lemma}
\begin{proof}
	Our key observation  is that $\sigma (\Phi_{a}^{k}(w)) = 
	\Phi_{a}^{k}(\sigma w)$ for any $k\in \{1,2\}$ and $\sigma \in 
	\mathfrak{B}_{n} \subset \mathfrak{B}_{n+1}$, which fixes the coordinate $n+1$. Thus 
	our lifting operators $\pro_{a}^{\bl,\bl+e^{k}_{i},k}$ are a composition of two $\mathfrak{B}_{n}$ 
	module morphisms. Finally we know from Corollary \ref{chpt6:cor:prononzero} that 
	the map is non-zero and therefore by Schur's Lemma must be injective.
\end{proof}

We are now ready to lift eigenvectors for the random transposition shuffle. We state two theorems, one for lifting $\lambda^{1}$, and another for lifting $\lambda^{2}$.
Afterwards we combine the theorems to recover all the eigenvalues for the random transposition shuffle on $B_{n}$. The results below follow from work presented in \cite{dieker2018spectral} and Section \ref{chpt5:subsec:algebra} (Lemma \ref{chpt5:lem:lifting} and Theorem \ref{chpt5:thm:lift}) with suitable modifications to the new setting of the hyperoctahedral group.

\begin{lemma}
	\label{chpt6:lem:lifting1}
	Let $\bl \vdash n$, and  $a \in \{1,2,\ldots,l(\lambda^{1})+1\}$. Take $i\in [n]$ such that  $1 \leq i \leq a$ and set $\bm = \bl +e_{i}^{1}$. Then,
	\begin{eqnarray*}
		\abrtp\circ \pro^{\bl,\bm,1}_{a} - 
		\pro^{\bl,\bm,1}_{a} \circ \abrt & = & (2 + 
		4(\lambda^{1}_{a}+1-a)) \pro^{\bl,\bm,1}_{a} +  4 \sum_{i\leq b < a} \Theta_{b,a}^{1} \circ 
		\pro^{\bl,\bm,1}_{b}.
	\end{eqnarray*}
\end{lemma}
\begin{proof}
	Continuing from Corollary \ref{chpt6:cor:restrict} we know
	\[			\abrtp \circ \sh_{a}^{1} - \sh_{a}^{1} \circ 
	\abrt |_{S^{\bl}}  = 2\sh_{a}^{1}|_{S^{\bl}}   +  
	4  \sum_{1 \leq b \leq a } 	\sh_{b}^{1}\circ\Theta_{b,a}^{1} |_{S^{\bl}} .\]
	Apply the isotypic projection $\pi^{\bm}$ to both sides of the equation. Since $\abrtp$ is given by the action of an element of the group algebra $\mathfrak{B}_{n}$ and $\pi^{\bm}$ is an $\mathfrak{B}_{n+1}$-module morphism, these operators commute and so we have			
	\begin{eqnarray}
	\abrtp \circ \pro^{\bl,\bm,1}_{a} - \pro^{\bl,\bm,1}_{a} \circ 
	\abrt & = &2 \pro^{\bl,\bm,1}_{a}   +  
	4  \sum_{1 \leq b \leq a} (\iso^{\bm} \circ	\sh_{b}^{1}\circ\Theta_{b,a}^{1} )|_{S^{\bl}}. \label{chpt6:eqn:importlem1}
	\end{eqnarray}	
	Applying equation \eqref{chpt6:eqn:addshufflerelation1} we may see that
	\[(\iso^{\bm} \circ	\sh_{b}^{1}\circ\Theta_{b,a}^{1} ) =  \Theta_{b,a}^{1} \circ  \pro^{\bl,\bm,1}_{b} -  \pro^{\bl,\bm,1}_{a}.\]
	The right hand side of equation \eqref{chpt6:eqn:importlem1} now becomes,
	\begin{eqnarray}
	2\pro^{\bl,\bm,1}_{a} +  
	4  \sum_{1 \leq b \leq a} \left(\Theta_{b,a}^{1} \circ \pro^{\bl,\bm,1}_{b} -  \pro^{\bl,\bm,1}_{a} \right) \label{chpt6:eqn:importantinglemmalast}
	\end{eqnarray}
	Notice that if $b=a$ then $\Theta_{a,a}^{1} (w)$ acts as a scalar by the number of occurrences of the symbol $a$ in $w$. All words in $\pro^{\bl,\bm,1}_{a}(S^{\bl})$ contain $\lambda_{a}^{1} +1$ occurrences of $a$.
	Finally if $b<i$ we know that $\sh_{b}^{1}(S^{\bl})$ does not contain $S^{\bl+e_{i}^{1}}$ as a submodule (by Lemma \ref{chpt6:lem:phimapscorrect}), so $\pro^{\bl,\bm,1}_{b} =0$. Thus, equation \eqref{chpt6:eqn:importantinglemmalast} is equal to
	\[( 2 +4(\lambda_{a}^{1} + 1 -a)) \pro^{\bl,\bm,1}_{a} +  
	4  \sum_{i \leq b < a} \Theta_{b,a}^{1} \circ \pro^{\bl,\bm,1}_{b}.\]
\end{proof}

\begin{lemma}
	\label{chpt6:lem:lifting2}
	Let $\bl \vdash n$, and  $a \in \{1,2,\ldots,l(\lambda^{2})+1\}$. Take $i\in [n]$ such that  $1 \leq i \leq a$ and set $\bm = \bl +e_{i}^{2}$. Then,
	\begin{eqnarray*}
		\abrtp \circ \pro^{\bl,\bm,2}_{a} - 
		\pro^{\bl,\bm,2}_{a}  \circ  \abrt& = & 
		4(\lambda^{2}_{a}+1-a) \pro^{\bl,\bm,2}_{a} +   2
		\sum_{i\leq b < a} 
		\Theta_{b,a}^{2} \circ 
		\pro^{\bl,\bm,2}_{a}.
	\end{eqnarray*}
\end{lemma}
\begin{proof}
	This follows from the same proof as Lemma \ref{chpt6:lem:lifting1} replacing the use of the first equality in Lemma \ref{chpt6:lem:thetadom} with the second equality.
\end{proof}

\begin{thm}[Lifting for $\lambda^{1}$]
	\label{chpt6:thm:lifting1}
	Let $\bl \vdash n$, and  $a \in \{1,2,\ldots,l(\lambda^{1})+1\}$. Take $i\in [n]$ such that  $1 \leq i \leq a$ and set $\bm = \bl +e_{i}^{1}$. Then,
	\begin{eqnarray*}
		\abrtp\circ \pro^{\bl,\bm,1}_{a} - 
		\pro^{\bl,\bm,1}_{a} \circ \abrt & = & (2 
		+ 
		4(\lambda_{i}^{1}+1-i))\pro^{\bl,\bm,1}_{a}
	\end{eqnarray*}
	Thus if we have an eigenvector $v \in S^{\bl}$ with 
	eigenvalue $\epsilon$, we find 
	$\pro^{\bl,\bm,1}_{a}(v)$ to be an eigenvector of
	$S^{\bm}$ with eigenvalue $\epsilon + (2 + 
	4(\lambda^{1}_{i}+1-i))$.
\end{thm}

\begin{proof}
	For $i=a$ the result follows from Lemma \ref{chpt6:lem:lifting1}. Let $\mu = \bl + e_{i}^{1}$, again from Lemma \ref{chpt6:lem:lifting1} we know that 
	\[\abrtp \circ \pro_{i}^{\bl,\bm,1} - \pro_{i}^{\bl,\bm,1}\circ \abrt = (2 + 4(\lambda^{1}_{i} + 1 -i))\pro_{i}^{\bl,\bm,1} .\]
	Applying the linear operator $\Theta_{i,a}^{1}$ to the above equation,
	\begin{eqnarray}
	\abrtp \circ \Theta_{i,a}^{1}\circ \pro_{i}^{\bl,\bm,1} - \Theta_{i,a}^{1}\circ \pro_{i}^{\bl,\bm,1}\circ \abrt = (2 + 4(\lambda^{1}_{i} + 1 -i)) \Theta_{i,a}^{1} \circ \pro_{i}^{\bl,\bm,1} . \label{chpt6:eqn:thm1}
	\end{eqnarray}
	Consider the left hand side of \eqref{chpt6:eqn:thm1}, break up the lifting operator into $\pro_{i}^{\bl,\bm,1} = \pi^{\bm} \circ \sh_{i}^{1}|_{S^{\bl}}$. The projection $\pi^{\bm}$ commutes with the $\mathfrak{B}_{n+1}$-module morphism $\Theta_{i,a}^{1}$. Performing this and using equation \eqref{chpt6:eqn:addshufflerelation1} we obtain the equation below, restricted to $S^{\bl}$;
	\begin{eqnarray}
	&  &\abrtp \circ \pi^{\bm} \circ \Theta_{i,a}^{1}\circ \sh_{i}^{1} - \pi^{\bm} \circ \Theta_{i,a}^{1}\circ \sh_{i}^{1} \circ \abrt  \nonumber\\
	& = &  \abrtp \circ \pi^{\bm} \circ (\sh_{i}^{1}\circ \Theta_{i,a}^{1} + \sh_{i}^{1}) - \pi^{\bm} \circ(\sh_{i}^{1}\circ \Theta_{i,a}^{1} + \sh_{i}^{1}) \circ \abrt  \nonumber\\
	& = & (\abrtp \circ \pro_{a}^{\bl,\bm,1} -\pro_{a}^{\bl,\bm,1} \circ \abrt) + (\abrtp \circ \pro_{i}^{\bl,\bm,1} - \pro_{i}^{\bl,\bm,1} \circ \abrt) \Theta_{i,a}^{1} \nonumber \\ 
	& = &(\abrtp \circ \pro_{a}^{\bl,\bm,1} -\pro_{a}^{\bl,\bm,1} \circ \abrt) + (2 + 4(\lambda_{i}^{1} +1 -i)) \pro_{i}^{\bl,\bm,1} \circ \Theta_{i,a}^{1} \label{chpt6:eqn:importthm2},
	\end{eqnarray}
	Manipulating the right hand side of \eqref{chpt6:eqn:importlem1} by again splitting our lifting operator $\pro_{i}^{\bl,\bm,1}$ and applying \eqref{chpt6:eqn:addshufflerelation1} we recover,
	\begin{eqnarray}
	\eqref{chpt6:eqn:importlem1}  & = &(2 + 4(\lambda^{1}_{i} + 1 -i)) \pi^{\bm}\circ \Theta_{i,a}^{1} \circ \sh_{i}^{1} \nonumber \\
	& = & (2 + 4(\lambda^{1}_{i} + 1 -i)) \pi^{\bm} \circ (\sh_{i}^{1}\circ \Theta_{i,a}^{1} + \sh_{i}^{1}) \nonumber \\
	& = &  (2 + 4(\lambda_{i}^{1} +1 -i)) \pro_{i}^{\bl,\bm,1} \circ \Theta_{i,a}^{1} + (2 + 4(\lambda_{i}^{1} +1 -i)) \pro_{a}^{\bl,\bm,1}. \label{chpt6:eqn:importthm3}
	\end{eqnarray}
	Combining equations \eqref{chpt6:eqn:importthm2} and \eqref{chpt6:eqn:importthm3} completes the proof.
\end{proof}

\begin{thm}[Lifting for $\lambda^{2}$]
	\label{chpt6:thm:lifting2}
	Let $\bl \vdash n$, and  $a \in \{1,2,\ldots,l(\lambda^{2})+1\}$. Take $i\in [n]$ such that  $1 \leq i \leq a$ and set $\bm = \bl +e_{i}^{2}$. Then,
	\begin{eqnarray*}
		\abrtp\circ \pro^{\bl,\bm,2}_{a}  - 
		\pro^{\bl,\bm,2}_{a}  \circ \abrt & = & 
		4(\lambda^{2}_{i}+1-i) \pro^{\bl,\bm,2}_{a} 
	\end{eqnarray*}
	Thus if we have an eigenvector $v \in S^{\bl}$ with 
	eigenvector $\epsilon$, we find 
	$\pro^{\bl,\bm,2}_{a} (v)$ to be an eigenvector of
	$S^{\bm}$ with eigenvalue $\epsilon +	4(\lambda^{2}_{i}+1-i)$.
\end{thm}
\begin{proof}
	This follows the same proof as Theorem \ref{chpt6:thm:lifting1}.
\end{proof}

Notice that there is a quantifiable difference in the change in eigenvalue depending on whether we lift $\lambda^{1}$ or $\lambda^{2}$. The value $\lambda_{i}^{k} +1 -i$ refers to the diagonal index of the new box added to $\lambda^{k}$, this value does not depend on current size of $\bl$. 
Combining Theorems \ref{chpt6:thm:lifting1} and \ref{chpt6:thm:lifting2} we recover all the eigenvectors for the random transposition shuffle belonging to each Specht module, this allows us to give a swift proof of Theorem \ref{chpt6:thm:rteig}

\begin{lemma}
	\label{chpt6:lem:liftingunique}
	For any $\bm \vdash n+1$ we may find a basis of eigenvectors of $\abrtp$ 
	for the module $S^{\bm}$ by lifting the eigenvectors of $\abrt$ 
	belonging in the modules $S^{\bl}$ with $\bl \vdash n$ and $\bl \subset \bm$.
\end{lemma}
\begin{proof}
	We proceed by induction. For $n=1$ we know the simple modules $\langle a \rangle = S^{((1),(0))}, 
	\langle a^{+} - a^{-} \rangle = S^{((0),(1))}$ of $\mathfrak{B}_{1}$ are both one dimensional, and their elements are eigenvectors for $\abrtn_{1}$.

	Consider the simple module $S^{\bm}$ with $\bm\vdash 
	n+1$. We know classically 
	from the branching rules of the hyperoctahedral group (Theorem \ref{chpt6:thm:branching}) that the restriction of this module to 
	$\mathfrak{B}_{n}$ is given by
	\[\textnormal{Res}_{\mathfrak{B}_{n}}^{\mathfrak{B}_{n+1}} \left( 
	S^{\bm}\right) \cong \bigoplus_{\substack{\bl 
			\vdash n \\ \bl \subset \bm}} 
	S^{\bl} .\]
	Now suppose we have a basis of eigenvectors for every 
	$S^{\bl}$. By Lemma \ref{chpt6:lem:pro2inject} the map 
	$\pro^{\bl,\bm}(S^{\bl})$ gives a 
	basis for the submodule $S^{\bl}$ inside of the vector space of 
	$\textnormal{Res}_{\mathfrak{B}_{n}}^{\mathfrak{B}_{n+1}} \left( 
	S^{\bm}\right)$ which vector space is isomorphic to the vector space 	$S^{\bm}$. Hence, considering all of the lifted eigenvectors 
	from every $S^{\bl}$ together we find a basis for 
	$S^{\bm}$. By Theorems \ref{chpt6:thm:lifting1} and \ref{chpt6:thm:lifting2} 
	the lifted eigenvectors form a basis of eigenvectors for 
	$S^{\bm}$.
\end{proof}

\begin{proof}[Proof of Theorem \ref{chpt6:thm:rteig}]
	We have found a basis of eigenvectors of the random transposition shuffle for every Specht module $S^{\bl}$. Each eigenvector in this basis may be constructed by repeated applications of Theorems \ref{chpt6:thm:lifting1} and \ref{chpt6:thm:lifting2}, starting at the empty word $\omega$ belonging to $S^{((0),(0))}$, and ending at the eigenvector belonging to $S^{\bl}$. When we add box $e_{i}^{1}$ in the process of forming $\bl$ the change in eigenvalue is $2 + 4(\lambda^{1}_{i}+1-i)$, the second coefficient here is the diagonal index of the box $(i,\lambda_{i}^{1} +1)$ added to $\lambda^{1}$. Similarly when we add box $ e_{i}^{2}$ the change in eigenvalue is the change in eigenvalue is $4(\lambda^{2}_{i}+1-i)$, which is the diagonal index of the box  $(i,\lambda_{i}^{2} +1)$ added to $\lambda^{2}$. Both of these changes are independent of when in the lifting we choose to add each box. Therefore, building up $\bl$ in any order produces an eigenvalue $2|\lambda^{1}| + 4\D(\lambda^{1}) + 4\D(\lambda^{2})$ with multiplicity $d_{\bl}^{2}$. Normalising the eigenvalue by $2n^{2}$ gives the expression \eqref{chpt6:eqn:rteig}.

\end{proof}

If we choose to ignore the second partition $\lambda^{2}$ when we lift eigenvalues we recover exactly the eigenvalues for the random transposition shuffle on $S_{n}$. That is, for a partition $\lambda \vdash n$, the bi-partition $(\lambda,(0))$ has corresponding eigenvalue,
\[\eig\left((\lambda,(0)) \right) = \frac{1}{n^{2}} \left( n + 2 \D(\lambda)\right) .\]
as seen in Lemma \ref{chpt5:lem:rteig}.
We end the lifting of the random transposition shuffle by  providing an explicit description of the lifting operators $\pro_{i}^{\bl,\bl+e_{i}^{k},k}$, this result follows from work of Dieker and Saliola \cite{dieker2018spectral}.

\begin{lemma}[Theorem 21 \cite{dieker2018spectral}]
	\label{chpt6:lem:liftingoperatorsdescription}
	Let $\lambda \vdash n$. The lifting maps as defined in Definition 
	\ref{chpt6:def:proj} are a linear combination of the shuffling and switching operators, explicitly they are:	
	\begin{eqnarray}
	\pro_{i}^{\bl,\bl+e_{i}^{k},k} & = & \sum_{1\leq b_{1} <\ldots <b_{m} <b_{m+1}=i} \left( \prod_{j=1}^{m} \frac{1}{k\left((\lambda_{i}^{k} -i) - (\lambda_{b_{j}}^{k} -b_{j}) \right)}\Theta_{b_{j},b_{j+1}}^{k} \right) \sh_{b_{1}}^{k}.	
	\end{eqnarray}
\end{lemma}

\begin{example}
	Let $\bl = ((1),(1))$, its associated Specht module is 2-dimensional and spanned by $S^{\bl} =   \langle 1 1^{+} - 11^{-}, 1^{+} 1 - 1^{-} 1\rangle$. Both these basis vectors of $S^{\bl}$ are eigenvectors for $\abrtn_{2}$. Take $w =  1 \,1^{+} - 1\,1^{-}$ using the formulae given in Lemma \ref{chpt6:lem:liftingoperatorsdescription} we lift $w$ to an eigenvector of $\abrtn_{3}$ belonging to the Specht modules $S^{((1,1),(1))}$ and $S^{((1),(1,1))}$ by adding the box $e_{2}^{1}$, $e_{2}^{2}$ respectively.
	\begin{eqnarray*}
		\pro_{2}^{\bl,\bl+e_{2}^{1},1}(w) & = & \left( \sh_{2}^{1} - \frac{1}{2}\Theta_{1,2}^{1} \sh_{1}^{1}\right) (w) = \frac{1}{2}\left[ (1 \, 1^{+} \, 2 - 2 \, 1^{+} \, 1 ) - (1 \, 1^{-} \, 2 - 2 \, 1^{-} \, 1) \right]  \\
		\pro_{2}^{\bl,\bl+e_{2}^{2},2}(w) & = & \left( \sh_{2}^{2} - \frac{1}{4}\Theta_{1,2}^{2} \sh_{1}^{2}\right) (w) \\ 
		&= &\frac{1}{2}\left[(11^{+}2^{+}  - 12^{+}1^{+}) + (11^{-}2^{-} - 12^{-}1^{-}) -  (11^{+}2^{-} - 12^{-}1^{+}) - (11^{-}2^{+} - 12^{+}1^{-} )\right].
	\end{eqnarray*}
\end{example}

\section{One-sided Transposition Shuffles on The Hyperoctahedral Group}

The one-sided transposition shuffle for the hyperoctahedral group is described by the following procedure: apply a transposition chosen according to the one-sided transposition shuffle for $S_{n}$, i.e. $\BLR_{n,w}$, then flip a fair coin; if heads do nothing, if tails flips the cards that were moved this step to their opposite sides.

\begin{defn}
	The \emph{biased one-sided transposition shuffle for $B_{n}$} with bias $w(j)$, denoted $\bblr_{n,w}$,  is driven by the following probability distribution:
	\begin{eqnarray}
	\bblr_{n,w}( \sigma ) = \begin{cases}
	\frac{\sum_{j=1}^{n}w(j)/j}{2N_{w}} & \textnormal{ if } \sigma = \gpid \\
	\frac{w(j)}{N_{w}} \frac{1}{2j}& \textnormal{ if } \sigma = \xi_{j} \textnormal{ for } j \in [n] \\
	\frac{w(j)}{N_{w}} \frac{1}{2j}& \textnormal{ if } \sigma = (i \, j)  \textnormal{ for } i,j \in [n] \textnormal{ with } i <j\\
	\frac{w(j)}{N_{w}} \frac{1}{2j} & \textnormal{ if } \sigma = \xi_{i} \xi_{j} (i \, j)  \textnormal{ for } i,j \in [n] \textnormal{ with } i< j \\
	0 & \textnormal{ otherwise }
	\end{cases} .
	\end{eqnarray}
\end{defn}
The one-sided transposition shuffle is a transitive, aperiodic, and reversible random walk on  $B_{n}$, and is not constant on the conjugacy classes of $B_{n}$.  The method of lifting eigenvectors for the random transposition shuffle on $B_{n}$ may be modified to be applicable for the biased one-sided transposition shuffle on $B_{n}$. 
In Section \ref{chpt6:subsec:ostlifting} we prove the following result.

\begin{thm}
	\label{chpt6:thm:lreig}
	The eigenvalues for the biased one-sided transposition shuffle $\bblr_{n,w}$ are indexed by standard Young tableaux of shape $\bl \vdash n$, and the eigenvalue represented by a tableau of shape $\bl$ has multiplicity $d_{\bl}$. For a standard Young tableau $\bt$ of shape $\bl$ the eigenvalue 
	corresponding to $\bt$ is given by	
	\[\eig(\bt) = \frac{1}{N_{w}} \left(\sum_{(i,j,1) \in \overline{T}} \frac{j-i+1}{T^{1}(i,j)} w(T^{1}(i,j)) +\sum_{(i,j,2) \in \overline{T}} \frac{j-i}{T^{2}(i,j)} w(T^{2}(i,j)) \right).\]
\end{thm}

If we focus on bi-tableau $\overline{T} = (T^{1},T^{2})$ with $T^{2} = \emptyset$ (those we can form by lifting only the first partition) then we recover the eigenvalues for the biased one-sided transposition shuffle on $S_{n}$ (Lemma \ref{biasedCompute}). We once again restrict our attention to weight functions of the form $w(j) = j^{\alpha}$, and define the biased one-sided transposition shuffle with this weight function as $\bblr_{n,\alpha}$. Recall that the time $t_{n,\alpha}$ was defined as, 
\[t_{n,\alpha}= \begin{cases}
\NA/n^{\alpha} & \textnormal{ if } \alpha \leq 1 \\
\NA/N_{\alpha-1}(n) & \textnormal{ if } \alpha \geq 1 
\end{cases}. \]
Applying the same projection argument used in Lemma \ref{chpt6:lem:rtlower} we may establish a lower bound of $t_{n,\alpha} (\log n)$ on the mixing time of the biased one-sided transposition shuffle $\bblr_{n,\alpha}$. Using the eigenvalues established in Theorem \ref{chpt6:thm:lreig} we can compute an upper bound on the total variation distance between $\bblr_{n,\alpha}$ and $\pi_{n}$. We suspect that any analysis of this bound is likely to encounter similar problems to those described in Section \ref{chpt6:subsec:upperboundrt} for the random transposition shuffle. We conjecture that it is possible to find a matching upper bound on the mixing time of $\bblr_{n,\alpha}$ for all $\alpha$.

\begin{lemma}
	\label{chpt6:lem:ostlower}
	The biased one-sided transposition shuffle on the hyperoctahedral group $\bblr_{n,\alpha}$ satisfies the following bound for any $c >\max(2, 3-\alpha)$:
	\begin{eqnarray*}
		\liminf_{n\to\infty}\, \lVert 
		\bblr_{n,\alpha}^{t_{n,\alpha}\left(\log n - \log \log n - 
			c\right)} 
		-\pi_{n} \rVert_{\tiny \textnormal{TV}} & 
		\geq 
		& \begin{cases}
			1 - \frac{\pi^{2}}{6(c-3+\alpha)^{2}} & \textnormal{ if } \alpha 
			\leq 1 \\
			1 - \frac{\pi^{2}}{6(c-2)^{2}} & \textnormal{ if } \alpha \geq 1 
		\end{cases} 
	\end{eqnarray*}
\end{lemma}

\begin{conj}
	\label{chpt6:conj:ostupper}
	The biased one-sided transposition shuffle on the hyperoctahedral group  $\bblr_{n,\alpha}$ satisfies the following bound:
	\begin{eqnarray*}
		\lim_{c\to \infty} \limsup_{n\rightarrow \infty}\, \lVert 
		\bblr_{n,\alpha}^{t_{n,\alpha}\left(\log n + c \right)} -\pi_{n} 
		\rVert_{\tiny \textnormal{TV}} & 
		= & 0 
	\end{eqnarray*}
	Thus, the biased one-sided transposition shuffle  exhibits a cutoff in total variation distance at time $t_{n,\alpha}\log n$ for all $\alpha$.
\end{conj}

\subsection{Lifting Eigenvectors for One-sided Transposition Shuffles}
\label{chpt6:subsec:ostlifting}

The eigenvectors of the one-sided transposition shuffle $\bblr_{n,w}$ may be computed using the same lifting operators $\pro_{a}^{\bl,\bm,k}$ as the random transposition shuffle. The difference in the lifting is reflected by changes to Theorems \ref{chpt6:thm:master1} and \ref{chpt6:thm:lifting1}. We begin by transforming $\bblr_{n,w}$ into an element of the group algebra $\mathfrak{B}_{n}$.

\begin{defn}
	The one-sided transposition shuffle on $B_{n}$ may be viewed as the following element of the group algebra
	$\mathfrak{B}_{n}$,
	\begin{eqnarray}
	\abblr = \sum_{1 \leq j \leq n} \frac{w(j)}{2j} \, \gpid + \sum_{1\leq j \leq n} \frac{w(j)}{2j} \xi_{i} +  
	\sum_{1\leq i < j \leq n} \frac{w(j)}{2j} (i \, j)  +  \sum_{1\leq i <j\leq n} \frac{w(j)}{2j} \, \xi_{i}\xi_{j}(i \, j ). \label{chpt6:eqn:blr}
	\end{eqnarray}
	Note that above we have scaled our probability by $\NW$.
\end{defn}
The one-sided transposition shuffle on $B_{n}$ has a similar recursive structure to the shuffle on $S_{n}$. Taking the difference of $\abblrp$ and $\abblr$ we find 
\[\abblrp - \abblr = \frac{w(n+1)}{2(n+1)} \, \left(\gpid + \xi_{n+1} \right)+  
\sum_{1\leq i < \leq n} \frac{w(n+1)}{2(n+1)} \left( (i \hspace{0.2em} n+1)  +  \xi_{i}\xi_{n+1}(i \hspace{0.2em} n+1) \right)\]
which only depends the movement of the new card $n+1$ (compare this with equation \eqref{chpt5:eqn:biasedcompare}).
Using this relationship we recover a new version of Theorem \ref{chpt6:thm:master1} for the one-sided transposition shuffle.

\begin{thm}
	\label{chpt6:thm:master3}
	Let $\bl \vdash n$ for $n\in \mathbb{N}$, for words in $M^{\bl}$ we have the following equalities: 
	\begin{eqnarray}
	\abblrp \circ \sh_{a}^{1} - \sh_{a}^{1} \circ 
	\abblr & = & \frac{w(n+1)}{n+1} \left(\sh_{a}^{1} +  
	\sum_{1 \leq b \leq n} 
	\sh_{b}^{1}\circ\Theta_{b,a}^{1} + \frac{1}{2}\sum_{1\leq b\leq n} 
	\left(\sh_{b}^{+} + \sh_{b}^{-}\right)  
	\circ\Theta_{b^{\pm},a} \right). \nonumber
	\\
	\abblrp \circ \sh_{a}^{2} - \sh_{a}^{2} \circ 
	\abblr & = &\frac{w(n+1)}{2(n+1)} \sum_{1 \leq b \leq n} 
	\sh_{b}^{2}\circ \Theta_{b,a}^{2}. \nonumber
	\end{eqnarray}
	
\end{thm}

\begin{proof}
	This follows the same proof as Theorem \ref{chpt6:thm:master1} with changes in constants to reflect equation \eqref{chpt6:eqn:blr}.
\end{proof}

Following from Theorem \ref{chpt6:thm:master3} the work of Section \ref{chpt6:subsec:rtlifting} may be replicated for the one-sided transposition shuffle. We summarise the results of the lifting in the following theorems whose proofs follow from those of the random transposition shuffle with changes in coefficients from Theorem \ref{chpt6:thm:master3}.

\begin{thm}[Lifting for $\lambda^{1}$]
	\label{chpt6:thm:lifting3}
	Let $\bl \vdash n$, and  $a \in \{1,2,\ldots,l(\lambda^{1})+1\}$. Take $i\in [n]$ such that  $1 \leq i \leq a$ and set $\bm = \bl +e_{i}^{1}$. Then,
	\begin{eqnarray*}
		\abblrp \circ \pro^{\bl,\bm,1}_{a} - 
		\pro^{\bl,\bm,1}_{a} \circ \abblr & = & \frac{w(n+1)(1 + (\lambda^{1}_{i}+1-i))}{n+1}\pro^{\bl,\bm,1}_{a}
	\end{eqnarray*}
	Thus if we have an eigenvector $v \in S^{\bl}$ with 
	eigenvalue $\epsilon$, we find 
	$\pro^{\bl,\bm,1}_{a}(v)$ to be an eigenvector of
	$S^{\bm}$ with eigenvalue $\epsilon + \frac{w(n+1)(1 + (\lambda^{1}_{i}+1-i))}{n+1}$.
\end{thm}

\begin{thm}[Lifting for $\lambda^{2}$]
	\label{chpt6:thm:lifting4}
	Let $\bl \vdash n$, and  $a \in \{1,2,\ldots,l(\lambda^{2})+1\}$. Take $i\in [n]$ such that  $1 \leq i \leq a$ and set $\bm = \bl +e_{i}^{2}$. Then,
	\begin{eqnarray*}
		\abblrp \circ \pro^{\bl,\bm,2}_{a}  - 
		\pro^{\bl,\bm,2}_{a}  \circ \abblr & = & 
		\frac{w(n+1)(\lambda^{2}_{i}+1-i)}{n+1} \pro^{\bl,\bm,2}_{a} 
	\end{eqnarray*}
	Thus if we have an eigenvector $v \in S^{\bl}$ with 
	eigenvalue $\epsilon$, we find 
	$\pro^{\bl,\bm,2}_{a} (v)$ to be an eigenvector of
	$S^{\bm}$ with eigenvalue $\epsilon +\frac{w(n+1)(\lambda^{2}_{i}+1-i)}{n+1}$.
\end{thm}

Note that Lemma \ref{chpt6:lem:liftingunique} still holds for the one-sided transposition shuffle because the lifting operators have not changed. This leads us to a proof of Theorem \ref{chpt6:thm:lreig}

\begin{proof}[Proof of Theorem \ref{chpt6:thm:lreig}]
	For any bi-tableau $\bt$ of shape $\bl$ we build up a distinct eigenvector belonging to $S^{\bl}$ by lifting from partition $((0),(0))$ to $\bl$ in the order specified by the tableau $\bt$. Summing the changes in eigenvalue given by Theorems \ref{chpt6:thm:lifting3} and \ref{chpt6:thm:lifting4} and normalising by $\NW$ we recover the expression for the eigenvalue corresponding to $\bt$ given in the theorem.
\end{proof}

\subsection{A Strong Stationary Time for the Unbiased One-sided Transposition Shuffle on The Hyperoctahedral Group.}
\label{chpt6:subsec:ostsst}
In this section we prove that the unbiased one-sided transposition shuffle on $B_{n}$ exhibits a cutoff in separation distance at time $n\log n$.	To establish a cutoff in separation distance we prove that the strong stationary time introduced in Section \ref{chpt5:sec:separation} is also a strong stationary time for the unbiased one-sided transposition shuffle on $B_{n}$. We record our main results below before constructing the strong stationary time in detail.

\begin{thm}
	\label{chpt6:thm:ostSSTtime}
	There exists a strong stationary time $T$ for the unbiased one-sided 
	transposition shuffle for $B_{n}$. Furthermore, for $c>0$ we have	$\mathbb{P}(T> n\log n + cn) \leq e^{-c}$.
\end{thm}

From Theorem \ref{chpt6:thm:ostSSTtime} we may quickly establish that the one-sided transposition shuffle exhibits a cutoff in total variation distance and separation distance at time $n\log n$. This proves Conjecture \ref{chpt6:conj:ostupper} for the case $\alpha=0$.

\begin{thm}
	\label{chpt6:thm:ostsepcutoff}
	The unbiased one-sided transposition shuffle $\bblr_{n,0}$ satisfies the following bounds, for $c_{1}>0, c_{2}>2$:
	\begin{eqnarray} 
	\limsup_{n\to\infty} \lVert \bblr_{n,0}^{n\log n +c_{1}n} 
	-\pi_{n} 	\rVert_{\textnormal{sep}} & \leq & e^{-c_{1}} \label{chpt5:eqn:sepupper}\\ 
	\liminf_{n\to\infty} \lVert \bblr_{n,0}^{n\log n -n \log \log n - c_{2}n} 
	-\pi_{n} 	\rVert_{\tiny \textnormal{TV}} & \geq & 1 -\frac{\pi^{2}}{6(c_{2}-2)^{2}} \label{chpt5:eqn:seplower}
	\end{eqnarray}
	Thus, unbiased one-sided transposition shuffle exhibits a cutoff in total variation distance and separation distance at 	time $n\log n$.
\end{thm}

\begin{proof}[Proof of Theorem \ref{chpt6:thm:ostsepcutoff}]
	
	The lower limit was given in Lemma \ref{chpt6:lem:ostlower}. 
	The upper limit follows from  Theorem 
	\ref{chpt6:thm:ostSSTtime}, using the strong stationary time $T$, we see that,
	\[\limsup_{n\to\infty} \lVert \bblr_{n,0}^{n\log n +c_{1}n} 
	-\pi_{n} 
	\rVert_{\textnormal{sep}} \leq \limsup_{n\to\infty} 
	\mathbb{P}(T> n\log n + c_{1}n) \leq e^{-c_{1}} .\]
	Separation distance gives an upper bound on total variation distance (Lemma \ref{chpt2:lem:sepupperbound}), therefore combining the two bounds establishes a cutoff in both total variation distance and separation distance.
\end{proof}

\subsubsection{A Strong Stationary Time Argument}
Recall that we may view the elements of $B_{n}$ as a subset of permutations of a deck of cards made up of cards and positions both indexed by $[\pm n]$.
Any permutation $\sigma \in B_{n}$ is a bijection from cards to positions and $\sigma^{-1}$ is a bijection from positions to cards, i.e., $\sigma(i)$ tells us the position of card $i$ whereas $\sigma^{-1}(i)$ tells us what card is in position $i$.
We begin all our random walks at the identity permutation with positions and labels fully matched. 
Throughout the rest of this section let $(X^{t})_{t\in 	\mathbb{N}}$ denote a Markov chain on $B_{n}$ driven by the unbiased one-sided transposition shuffle, and let $(Y^{t})_{t \in \mathbb{N}}$ be a Markov chain on $B_{n}$ defined by setting $Y^{t} = (X^{t})^{-1}$ for all $t$. The Markov chains $(X^{t})$ and $(Y^{t})$ represent two different ways to view the unbiased one-sided transposition shuffle.

Let $\tau^{t}$ be the transposition chosen at step $t$ of the unbiased one-sided transposition shuffle. 
To construct our strong stationary time we need to condition on the exact permutation of cards in positions (both positive and negative) above position $j$  at time $t$, that is the random variables $Y^{t}(\pm i)$ for $j< i  \leq n$. Note that if we know the random variable $Y^{t}(i)$ we may compute $Y^{t}(-i) = -Y^{t}(i)$.
Given this information we also know which cards can be in positions $[\pm j]$ at time $t$, define this set as, 
\[A_{j}^{t} = [\pm  n] \setminus \{ Y^{t}(\pm i)  \, | \, j < i \leq n\}.\]

\begin{defn}
	We say the random walk $(Y^{t})_{t\in \mathbb{N}}$ satisfies \emph{property $\mathcal{P}_{j}$ at time $t$} if we have:
	\begin{eqnarray}
	\mathbb{P}\left(Y^{t} (j) = l \, | \, Y^{t}(i) \textnormal{ for all } j < i \leq n \right) = 
	\begin{cases}
	1/2j & \textnormal{ if } l \in A_{j}^{t} \\
	0 & \textnormal{ otherwise } 
	\end{cases} \label{chpt6:eqn:propertyj}
	\end{eqnarray} 
	This property tells us that given total information about the deck strictly above position $j$, the card in position $j$ is equally likely to be any of the remaining cards.
\end{defn}

\begin{lemma}
	\label{chpt6:lem:SUTpropj}
	Let $T_{j}$ be the first time our right hand chooses the position
	$j$ when performing the unbiased one-sided transposition shuffle. 
	If $T_{j} \leq t$ then the Markov chain $(Y^{t})_{t\in \mathbb{N}}$ satisfies property $\mathcal{P}_{j}$ at time $t$.
\end{lemma}
\begin{proof}
	We prove this by induction: once property $\mathcal{P}_{j}$ holds for some $t$, it holds for all times after $t$.

	Consider the time $T_{j}$, at this step of our Markov chain we must have applied  a transposition 
	$(i\,j)$ or $\xi_{i}\xi_{j}(i\, j)$ with $i\leq j$. The probability of picking any one of the transpositions $(i \,j)$ at time $T_{j}$ is $\mathbb{P}(\tau^{T_{j}} = (i\,j)) = \mathbb{P}(\tau^{T_{j}} = \xi_{i} \xi_{j}(i\,j))= 1/2j$ for all $i\leq j$. Therefore, the card in position $j$ at time $T_{j}$ has a uniform chance of being any of the cards in $A_{j}^{T_{j-1}} = A_{j}^{T_{j}}$ (which contains both positive and negative cards). Thus, we may clearly see that,
	\begin{eqnarray}
	\mathbb{P}\left(Y^{T_{j}} (j) = l \, | \, Y^{T_{j}}(i) \textnormal{ for all } j < i \leq n \right)= 
	\begin{cases}
	1/2j & \textnormal{ if } l \in A_{j}^{T_{j}}\\
	0 & \textnormal{ otherwise } 
	\end{cases}
	\end{eqnarray} 
	so $\mathcal{P}_{j}$ holds at time $T_{j}$.

	Now suppose property $\mathcal{P}_{j}$ holds at time $t$. We study the time $t+1$ and split the analysis into cases based on which permutation was applied at time $t+1$,
	\begin{eqnarray}
	& & \mathbb{P}\left(Y^{t+1}(j) =  l \, | \, Y^{t+1}(i) \textnormal{ for all } j < i \leq n\right) \nonumber\\
	& = &
	\sum_{\sigma \in B_{n}} \mathbb{P}(\tau^{t+1} = 
	\sigma) \, \mathbb{P}\left(Y^{t+1} (j) 
	= l \, | \, Y^{t+1}(i) \textnormal{ for all } j < i \leq n, \, \tau^{t+1} = \sigma 
	\right).\hspace{0.5cm}\label{chpt6:eqn:SSTinduction} 
	\end{eqnarray}
	Using knowledge of the transposition $\tau^{t+1}$ we evolve our deck backwards in time ($X^{t} = \tau^{t+1} X^{t+1}$) to recover the random variables $Y^{t}(i)$ from $Y^{t+1}(i)$ for $j < i \leq n$, and relate $Y^{t+1}(j)$ to $Y^{t}(j)$. This allows us to use our inductive hypothesis.

	If $\tau^{t+1} \in \{(a\,j), \, \xi_{a}\xi_{j} (a \,j) \, | \, a\leq j \}$ then our random walk satisfies property $\mathcal{P}_{j}$ at time $t+1$ for the same reasoning as time $T_{j}$. 
	Suppose that $\tau^{t+1} \in \{(a\,b), \, \xi_{a}\xi_{b} (a \,b ) \, | \,  a,b< j \}$. Then we know that $Y^{t}(i) = Y^{t+1}(i)$ for all $j<i \leq n$. Suppose instead that $\tau^{t+1} \in \{(a\,b), \, \xi_{a}\xi_{b} (a \,b ) \, | \,  j <a,b\}$, then we have $ Y^{t}(\pm b) = Y^{t+1}( \pm a)$ and $ Y^{t}( \pm a) = Y^{t+1}(\pm b)$ (with the signs swapping in the second case), with $Y^{t}(i) = Y^{t+1}(i)$ for all other $j<i\leq n$. In either case we know $A_{j}^{t+1} = A_{j}^{t}$, and the card in position $j$ has not moved from time $t$ to $t+1$. Therefore, we have

	\begin{eqnarray}
	& & \mathbb{P}\left(Y^{t+1} (j) = l \, | \, Y^{t+1}(i) \textnormal{ for all } j < i \leq n, \, \tau^{t+1}  \in \{(a\,b), \, \xi_{a}\xi_{b} (a \,b ) \, | \,  a,b< j \textnormal{ or both } b,a>j\} \right) \nonumber \\	
	& = & \mathbb{P}\left(Y^{t} (j) = l \, | \, Y^{t}(i) \textnormal{ for all } j < i \leq n \right)  = 
	\begin{cases}
	1/2j & \textnormal{ if } l \in A_{j}^{t} = A_{j}^{t+1}\\
	0 & \textnormal{ otherwise } 
	\end{cases} \nonumber.
	\end{eqnarray}
	
	We now study in detail the effects of the remaining transpositions $(a \, b)$ and $\xi_{a},\xi_{b}(a \,b)$ with $a\leq j<b$.
	In this case  we can not fully recover the random variables $Y^{t}(i)$ with $i < j \leq n$ without extra assumptions. To this end fix $b>j$,  a card $C \in A_{j}^{t+1}$, and suppose that $\tau^{t+1} = (X^{t+1}(C) \hspace{0.15cm} b)$, i.e., card $C$ is moved from position $b$ into a position below $j$ by $\tau^{t+1}$. Letting $C$ range over all choices in $A_{j}^{t+1}$ will recover every transposition $(a \, b)$ and $\xi_{a}\xi_{b} (a \, b)$ with $b>j$ fixed and  $a\leq j$. In the case that  $\tau^{t+1} = (X^{t+1}(C) \hspace{0.15cm} b)$, we know that $Y^{t}(b) = C$ and the other positions above $j$ have $Y^{t}(i) = Y^{t+1}(i)$ for $j < i \leq n$ and $i \neq b$.  Therefore, for this choice of $\tau^{t+1}$ we know that $A_{j}^{t} = \left(A_{j}^{t+1} \sqcup \{Y^{t+1}(\pm b)\} \right) \setminus \{\pm C\}$. Now consider the probability:
	\begin{eqnarray}
	\mathbb{P}\left(Y^{t+1}(j) = l \, | \, Y^{t+1}(i) \textnormal{ for all } j < i \leq n , \, \tau^{t+1} = (X^{t+1}(C) \hspace{0.2cm} b)\right) \nonumber.
	\end{eqnarray}
	If $l = C$ then the event in question can only occur if  the card currently in position $b$, i.e. $Y^{t+1}(b)$, was in position $j$ at time $t$. Noting that $Y^{t+1}(b) \in A_{j}^{t}$, and using our inductive hypothesis we find,
	\begin{eqnarray}
	& & \mathbb{P}\left(Y^{t+1}(j) = C \, | \, Y^{t+1}(i) \textnormal{ for all } j < i \leq n , \, \tau^{t+1} = (X^{t+1}(C) \hspace{0.2cm} b)\right) \nonumber \\
	& = & \mathbb{P}\left(Y^{t}(j) = Y^{t+1}(b) \, | \, Y^{t}(i) \textnormal{ for all } j < i \leq n  \right) = 1/2j \label{chpt5:eqn:SSThardind1}
	\end{eqnarray}
	If $l = -C$ then the event in question can only occur if the card currently in position $-b$, i.e. $Y^{t+1}(-b)$, was in position $j$ at time $t$. Noting that $Y^{t+1}(-b) \in A_{j}^{t}$, and using our inductive hypothesis we find,
	\begin{eqnarray}
	& & \mathbb{P}\left(Y^{t+1}(j) = -C \, | \, Y^{t+1}(i) \textnormal{ for all } j < i \leq n , \, \tau^{t+1} = (X^{t+1}(C) \hspace{0.2cm} b)\right) \nonumber \\
	& = & \mathbb{P}\left(Y^{t}(j) = Y^{t+1}(-b) \, | \, Y^{t}(i) \textnormal{ for all } j < i \leq n  \right) = 1/2j \label{chpt5:eqn:SSThardind1.5}
	\end{eqnarray}
	Alternatively, suppose $l \in A_{j}^{t+1} \setminus \{\pm C\}$, we know the card $l$ does not move from its position at time $t$ to time $t+1$, and we know $l \in A_{j}^{t}$. Therefore, we find
	\begin{eqnarray}
	& & \mathbb{P}\left(Y^{t+1}(j) = l \, | \, Y^{t+1}(i) \textnormal{ for all } j < i \leq n , \, \tau^{t+1} = (X^{t+1}(C) \hspace{0.2cm} b)\right) \nonumber \\
	& = & \mathbb{P}\left(Y^{t}(j) = l \, | \, Y^{t}(i) \textnormal{ for all } j < i \leq n  \right) = 1/2j \label{chpt5:eqn:SSThardind2}
	\end{eqnarray}
	Putting the equations above together gives us,
	\begin{eqnarray*}
		\mathbb{P}\left(Y^{t+1}(j) = l \, | \,  Y^{t+1}(i) \textnormal{ for all } j < i \leq n, \, \tau^{t+1} = (X^{t+1}(C) \hspace{0.2cm} b)\right) & = & 
		\begin{cases}
			1/2j & \textnormal{ if }  l \in A_{j}^{t+1} \setminus\{ \pm C \}\\
			1/2j & \textnormal{ if }  l \in \{ \pm C\} \\
			0 & \textnormal{ otherwise }
		\end{cases}.
	\end{eqnarray*}
	Letting $C$ range over all possible choices of card in $A_{j}^{t+1}$ while keeping $b> j$ fixed, we cover the desired probability for all transpositions $(a\, b)$, $\xi_{a}\xi_{b}(a\, b)$ with $a\leq j < b$. 
	Finally applying every separate case to \eqref{chpt6:eqn:SSTinduction} we have established that 
	\[\mathbb{P}\left(Y^{t+1}(j) = l \, | \, Y^{t+1}(i) \textnormal{ for all } j < i \leq n
	\right) = \begin{cases}
	1/2j & \textnormal{ if } l \in A_{j}^{t+1}\\
	0 & \textnormal{ otherwise } 
	\end{cases} \]
	as required, thus by induction our hypothesis holds for all 
	$t \geq T_{j}$.
	
\end{proof}

\begin{lemma}
	\label{chpt6:lem:SUTresutl}
	Let $T =\min\{t \geq 0  \, | \, t \geq T_{1},\ldots T_{n}\}$ be the first time our right hand has chosen every position $j$. Then $T$ is a strong uniform time for $(X^{t})_{t \in \mathbb{N}}$ and $(Y^{t})_{t \in \mathbb{N}}$.
\end{lemma}

\begin{proof}
	Note that $X^{t}$ is uniformly distributed if and only if $Y^{t}$ is uniformly distributed.
	Lemma \ref{chpt6:lem:SUTpropj} implies that by time $T$ our Markov chain $(Y^{t})_{t \in \mathbb{N}}$ satisfies all properties $\mathcal{P}_{j}$. Hence, we have
	\begin{eqnarray*}
		\mathbb{P}(Y^{t} = \sigma^{-1} \, | \, T\leq t) & = &\mathbb{P}\left(
		\cap_{j=1}^{n} \{Y^{t}(j) = \sigma^{-1}(j) \} | T\leq t\right)\\  
		& = & \prod_{j=1}^{n}\mathbb{P}\left(  Y^{t}(j) = \sigma^{-1}(j) | 
		\cap_{i=j+1}^{n} \{
		Y^{t}(i) = \sigma^{-1}(i) \}, T\leq t\right)\\
		& = & 	\prod_{j=1}^{n} \frac{1}{2j} = \frac{1}{2^{n}n!} = \pi_{n}(\sigma).
	\end{eqnarray*}
\end{proof}

We have found a strong stationary time for the unbiased one-sided transposition shuffle. Following quickly from this we may prove Theorem \ref{chpt6:thm:ostSSTtime}, and thus establish a cutoff in separation distance for the one-sided transposition shuffle. 

\begin{proof}[Proof of Theorem \ref{chpt6:thm:ostSSTtime}]
	Let $T$ be the first time our right hand has chosen every position 
	$j$. 
	Our right hand is choosing positions via a uniform probability on 
	$[n]$. Thus $T$ is modelled by the uniform coupon collectors problem with $n$ 	coupons. To complete our argument recall (Section \ref{chpt4:subsec:t2r}, equation \eqref{chpt4:eqn:classiccoupon})  that for the uniform coupon collectors problem on $n$ cards we have $\mathbb{P}(T >n\log n +cn) \leq e^{-c}$. 
\end{proof} 

The strong stationary time we have constructed in this section is actually a strong stationary time for any biased one-sided transposition shuffle. Let $T_{\alpha}$ be the first time our right hand has chosen every position $j$ following the biased one-sided transposition shuffle $\bblr_{n,\alpha}$. In Lemma \ref{chpt5:lem:SSTnegalpha} we showed that for $\alpha \leq 0$ we have the bound $\mathbb{P}(T_{\alpha} > t_{n,\alpha}(\log n + c)) \leq e^{-c}$. Therefore, extending the proof of Theorem \ref{chpt6:thm:ostsepcutoff} we are able to prove the following.

\begin{corollary}
	\label{chpt6:thm:biasedsepcutoff}
	The biased one-sided transposition shuffle $\bblr_{n,\alpha}$ with $\alpha\leq 0$ exhibits a cutoff in separation distance and total variation distance at time $t_{n,\alpha} \log n$.	
	Thus, conjecture \ref{chpt6:conj:ostupper} holds for all $\alpha \leq 0$.
\end{corollary}

\section{Further Work}

In this chapter we have explored the hyperoctahedral group as an extension of the symmetric group. The hyperoctahedral group is just one instance of the generalised symmetric group.

\begin{defn}
	Let $r,n \in \mathbb{N}$, and let $\xi$ be a $r^\textnormal{th}$ root of unity. The \emph{generalised symmetric group}, denoted $G_{r,n}$, is defined as the group of all bijections $\sigma$ on elements $\{ \xi^{k} i \, | \, i \in [n],k\in[k]\}$ such that $\sigma(\xi^{k} i ) = \xi^{k} \sigma(i)$. The generalised symmetric group $G_{r,n}$ is isomorphic to the wreath product $\mathbb{Z}_{r} \wr S_{n}$.
\end{defn}

The symmetric and hyperoctahedral groups are isomorphic to the generalised symmetric groups $G_{1,n}$ and $G_{2,n}$ respectively.
We were able to describe the module structure of the symmetric group and the hyperoctahedral group using a similar set of techniques in both cases.  The arguments presented in Sections  \ref{chpt4:sec:modules} and \ref{chpt6:subsec:modules} can be modified to work for the generalised symmetric group $G_{r,n}$ for any $r \in\mathbb{N}$. Below we describe the module structure of $G_{r,n}$ over the field $\mathbb{C}$, a full construction of the permutation and Specht modules for any choice of $r \in \mathbb{N}$ may be found in \cite{can1996representations}.

\begin{defn}
	Let $r,n\in\mathbb{N}$. An \emph{$r$-partition of $n$}, denoted $\bl$, is a tuple of partitions $\bl = (\lambda^{1},\ldots, \lambda^{r})$ such that $\sum_{k=1}^{r}|\lambda^{k}| = n$.	Let $\bl,\bm$ be $r$-partitions. Define the \emph{dominance ordering} on $r$-partitions as follows:
	\[\bl \trianglerighteq \bm \Leftrightarrow \begin{cases}
	& |\lambda^{1}| >  |\mu^{1}| \\
	\textnormal{ or } & |\lambda^{k}| = |\mu^{k}| \textnormal{ for  $1\leq k \leq i$, and } |\lambda^{i+1}| > |\mu^{i+1}| \textnormal{ for some }  i\in \{1,\ldots,r-2\} \\
	\textnormal{ or } & |\lambda^{k}| = |\mu^{k}| \textnormal{ and } \lambda^{k} \trianglerighteq \mu^{k} \textnormal{ for all $k \in [r]$} 
	\end{cases}
	.\]
	
\end{defn}

\begin{lemma}
	The permutation and simple modules for the generalised symmetric group $G_{r,n}$ are index by $r$-partitions of $n$. The permutation module and simple module corresponding to a $r$-partition $\bl$ are denoted $M^{\bl}$ and $S^{\bl}$ respectively. 
	Furthermore, the permutation and simple modules respect Young's rule, that is 
	\[M^{\bm} \cong \bigoplus_{\bl \trianglerighteq \bm} K_{\bl,\bm} S^{\bl},\]
	for constants $K_{\bl,\bm} \in \mathbb{N}^{0}$.
\end{lemma}

\begin{conj}
	Let $n \geq 2$, and $\bl \vdash n$. The \emph{branching rules} for the simple modules of $G_{r,n}$ are as follows:
	\begin{eqnarray}
	\res_{G_{r,n-1}}^{G_{r,n}} S^{\bl} & \cong & \bigoplus_{\substack{\bm \vdash n-1 \\ \bm \subseteq \bl}} S^{\bm} \textnormal{ as $G_{r,n-1}$-modules} \\
	\induce_{G_{r,n}}^{G_{r,n+1}} S^{\bl} & \cong &\bigoplus_{\substack{\bm \vdash n+1 \\ \bl \subseteq \bm }} S^{\bm} \textnormal{ as $G_{r,n+1}$-modules}.
	\end{eqnarray}
\end{conj}

\paragraph{}
The structure of $G_{r,n}$ allows us to generalise the random transposition shuffle and the one-sided transposition shuffle for any choice of $r$. To do this we first need to define the negative transpositions $\xi_{i}$ for the group $G_{r,n}$. Let $\xi$ be an $r^{\textnormal{th}}$-root of unity, and define the negative transpositions $\xi_{i} \in G_{r,n}$ for $i\in [n]$ as follows:
\[\xi_{i}(j) = \begin{cases}
\xi \,i & \textnormal{ if } j= i\\
j & \textnormal{ otherwise }
\end{cases}
.\]
The random transposition shuffle on the group $G_{r,n}$ is defined by the following procedure: pick a transposition $(i\,j)$ of $S_{n} \subseteq G_{r,n}$ via the random transposition shuffle and apply it, then uniformly at random choose a element $k \in [r]$ and apply the negative transpositions $\xi_{i}^{k}, \xi_{j}^{k}$ to the positions $i,j$. Overall we have applied an element of $G_{n,r}$ of the form $\xi_{i}^{k}\xi_{i}^{k} (i \,j)$ or $\xi_{j}^{k} (j \, j)$. If $r=1,2$ we recover the random transposition shuffle on the symmetric group and hyperoctahedral group respectively which were described in Sections \ref{chpt4:subsec:rt} and \ref{chpt6:subsec:rt}. Similarly we may form the biased one-sided transposition shuffles on the group $G_{r,n}$ by the following procedure:
pick a transposition $(i\,j)$ of $S_{n}$ via the biased one-sided transposition shuffle and apply it, then uniformly at random choose a element $k \in [r]$ and apply the negative transpositions $\xi_{i}^{k}, \xi_{j}^{k}$ to the positions $i,j$.
It is natural to think how we might extend the technique of lifting eigenvectors to the random transposition and the one-sided transposition shuffles on the group $G_{r,n}$ for $r>2$. In line with the $r=2$ case we expect there to be $r$ different lifting operators, one for adding a box to each separate partition $\lambda^{k}$ which together constitute an $r$-partition. Furthermore, we expect the eigenvectors for each shuffle to correspond to standard Young tableaux of $r$-partitions.

\begin{conj}
	The eigenvalues for the random transposition shuffle on the generalised symmetric group $G_{r,n}$ are labelled by $r$-partitions of $n$, and may be described by the technique of lifting eigenvectors. 
\end{conj}

\begin{conj}
	The eigenvalues for the biased one-sided transposition shuffle on the  generalised symmetric group $G_{r,n}$ are labelled by standard Young tableaux of $r$-partitions of $n$, and may be described by the technique of lifting eigenvectors.
\end{conj}

Once the eigenvalues for each shuffle have been found they could be used to analyse the mixing time of each shuffle on the generalised symmetric group $G_{r,n}$. We conjecture that the random transposition shuffle and the one-sided transposition shuffles on the group $G_{r,n}$ exhibit a cutoff in total variation distance at the same time as their symmetric group counterparts.

\begin{conj}
	The random transposition shuffle on the generalised symmetric group $G_{r,n}$ exhibits a cutoff in total variation distance at time $(n/2) \log n$.
\end{conj}

\begin{conj}
	The biased one-sided transposition shuffle on the generalised symmetric group $G_{r,n}$ with weight function $w$ such that $w(j) = j^\alpha$ exhibits a cutoff in total variation distance at time $t_{n,\alpha} \log n$.
\end{conj}

In fact by using the strong stationary time introduced in Section \ref{chpt5:sec:separation} we may prove the above conjecture in the case of the unbiased one-sided transposition shuffle ($\alpha=0$) on the group $G_{r,n}$.
\begin{thm}
	The unbiased one-sided transposition shuffle on the  generalised symmetric group $G_{r,n}$ exhibits a cutoff in total variation distance and separation distance at time $ n \log n$.
\end{thm}

\begin{proof}[Sketch Proof]
	The proof of this Theorem follows from extending the techniques of Sections 	\ref{chpt5:subsec:ostsst} and \ref{chpt6:subsec:ostsst}. Let $\xi$ be an $r^{\textnormal{th}}$ root of unity. Let $(X^{t})_{t\in 	\mathbb{N}}$ denote a Markov chain on $G_{r,n}$ driven by the unbiased one-sided transposition shuffle, and let $(Y^{t})_{t \in \mathbb{N}}$ be a Markov chain on $G_{r,n}$ defined by setting $Y^{t} = (X^{t})^{-1}$ for all $t$.  To construct a strong stationary time we need to condition on exact knowledge of the positions above position $j$  at time $t$, that is the random variables $Y^{t}(\xi^{k} i)$ for $j< i  \leq n$ and $k\in [r]$. Note that if we know the random variable $Y^{t}(i)$ we may compute $Y^{t}(\xi^{k} i) = \xi^{k} Y^{t}(i)$ for any $k \in [r]$.
	Given this information we also know the value of positions $\xi^{k} 1$ to $\xi^{k}j$ at time $t$, define this set as, 
	\[A_{j}^{t} = \{\xi^{k} i  \,| \, i\in [n], k \in [r] \} \setminus \{ Y^{t}(\xi^{k} i)  \, | \, j < i \leq n, \, k \in [r]\}.\]
	
	We say our random walk satisfies \emph{property $\mathcal{P}_{j}$ at time $t$} if we have,
	\begin{eqnarray}
	\mathbb{P}\left(Y^{t} (j) = l \, | \, Y^{t}(i) \textnormal{ for all } j<i \leq n \right) = 
	\begin{cases}
	1/rj & \textnormal{ if } l \in A_{j}^{t} \\
	0 & \textnormal{ otherwise } 
	\end{cases} .\label{chpt6:eqn:propertyjgeneral}
	\end{eqnarray}
	Let $T_{j}$ be the first time we apply a transposition $\xi_{i}^{k}\xi_{j}^{k}(i \, j)$ with $i \leq j$ and  $k \in [r]$ following the unbiased one-sided transposition shuffle on $G_{r,n}$. Extending the arguments of Lemmas \ref{lem:SUTpropj} and \ref{chpt6:lem:SUTpropj}, we may see that at any time $t \geq T_{j}$ the Markov chain $(Y^{t})_{t\in\mathbb{N}}$ satisfies property $\mathcal{P}_{j}$.
	
	Let $T= \min \{t\geq 0 \,: \, t \geq T_{j} \textnormal{ for all $j$}\}$, this is a stopping time for our random walk. At any time after $T$ the random walk $(Y^{t})_{t\in\mathbb{N}}$, satisfies all properties $\mathcal{P}_{j}$, therefore for any $\sigma \in G_{r,n}$ we have
	\begin{eqnarray*}
		\mathbb{P}(Y^{t} = \sigma^{-1} \, | \, T\leq t) & = &\mathbb{P}\left( 
		\cap_{j=1}^{n} \{Y^{t}(j) = \sigma^{-1}(j) \} \, | \, T\leq t\right)\\  
		& = & \prod_{j=1}^{n}\mathbb{P}\left(  Y^{t}(j) = \sigma^{-1}(j) \, | \, 
		\cap_{i=j+1}^{n} \{
		Y^{t}(i) = \sigma^{-1}(i) \}, T\leq t\right)\\
		& = &\prod_{j=1}^{n} 
		\frac{1}{rj} = \frac{1}{r^{n}n!} = \pi_{n}(\sigma).
	\end{eqnarray*}
	Hence, the time $T$ is a strong stationary time for $(Y^{t})_{t \in \mathbb{N}}$ and so $(X^{t})_{t \in\mathbb{N}}$. Using the same coupon collector's argument as before we know that $\mathbb{P}(T>n\log n + cn) \leq e^{-c}$, hence we have an upper bound of $n\log n$ on the on separation distance mixing time of the unbiased one-sided transposition shuffle.
	To find a matching lower bound we reduce the analysis of the unbiased one-sided transposition shuffle on $G_{r,n}$ to the $r=1$ (symmetric group case) using the same argument as Lemma \ref{chpt6:lem:rtlower}. This gives a lower bound of $n\log n$ on the total variation distance mixing time of the unbiased one-sided transposition shuffle.  Combining the two bounds above proves the existence of a cutoff in total variation distance and separation distance at time $n\log n$.
\end{proof}

\newpage
\addcontentsline{toc}{chapter}{Bibliography}
\bibliography{ThesisBibClean}
\bibliographystyle{plain}

\end{document}